\pdfoutput=1
\documentclass[a4paper]{article}
\usepackage[margin=30mm]{geometry}
\usepackage[svgnames]{xcolor}
\usepackage{amsmath, amssymb, amsthm, mathtools, comment, here}
\usepackage{authblk}
\usepackage[normalem]{ulem}
\usepackage[setpagesize=false]{hyperref}
\hypersetup{
 colorlinks=true,
 linkcolor=FireBrick,
 citecolor=DarkBlue,
}
\usepackage{tikz}
\usetikzlibrary{cd}

\title{Frobenius algebras associated with the $\alpha$-induction for equivariantly braided tensor categories}
\author{Mizuki Oikawa \thanks{The author is supported by JSPS KAKENHI Grant Number JP23KJ0540.}}
\affil{Graduate School of Mathematical Sciences \\ The University of Tokyo, Komaba, Tokyo 153-8914, Japan \\ E-mail: \texttt{moikawa@ms.u-tokyo.ac.jp}}
\date{}

\theoremstyle{definition}
\newtheorem{defi}{Definition}[section]
\newtheorem{thm}[defi]{Theorem}
\newtheorem{lem}[defi]{Lemma}
\newtheorem{prop}[defi]{Proposition}
\newtheorem{cor}[defi]{Corollary}
\newtheorem{rem}[defi]{Remark}
\newtheorem{eg}[defi]{Example}

\makeatletter
    
    \@addtoreset{equation}{section}
  \makeatother

\providecommand{\abs}[1]{\left\lvert#1\right\rvert}

\newcommand\cala{\mathcal{A}}
\newcommand\calb{\mathcal{B}}
\newcommand\calc{\mathcal{C}}
\newcommand\cald{\mathcal{D}}

\renewcommand\hom{\operatorname{Hom}}
\newcommand\id{\mathrm{id}}
\newcommand\one{\mathbf{1}}

\newcommand\z{\mathbb{Z}}
\DeclareMathOperator{\ad}{Ad}
\DeclareMathOperator{\bimod}{Bimod}

\DeclareMathOperator{\en}{End}
\DeclareMathOperator{\frob}{Frob}
\DeclareMathOperator{\homog}{Homog}
\DeclareMathOperator{\ind}{Ind}

\DeclareMathOperator{\obj}{Obj}
\DeclareMathOperator{\rep}{Rep}
\DeclareMathOperator{\aut}{Aut}
\DeclareMathOperator{\tr}{Tr}
\DeclareMathOperator{\trl}{Tr^L}

\begin{document}
\tikzset{>=stealth}
\tikzset{cross/.style={preaction={-,draw=white,line width=6pt}}}
\tikzset{block/.style={draw, fill=white, rectangle, minimum width=0.5}}
\maketitle

\emph{Dedicated to Yasuyuki Kawahigashi on the occasion of his sixtieth birthday}

\begin{abstract}
Let $G$ be a group. We give a categorical definition of the $G$-equivariant $\alpha$-induction associated with a given $G$-equivariant Frobenius algebra in a $G$-braided multitensor category, which generalizes the $\alpha$-induction for $G$-twisted representations of conformal nets. For a given $G$-equivariant Frobenius algebra in a spherical $G$-braided fusion category, we construct a $G$-equivariant Frobenius algebra, which we call a $G$-equivariant $\alpha$-induction Frobenius algebra, in a suitably defined category called neutral double. This construction generalizes Rehren's construction of $\alpha$-induction Q-systems. Finally, we define the notion of the $G$-equivariant full center of a $G$-equivariant Frobenius algebra in a spherical $G$-braided fusion category and show that it indeed coincides with the corresponding $G$-equivariant $\alpha$-induction Frobenius algebra, which generalizes a theorem of Bischoff, Kawahigashi and Longo.
\end{abstract}

\tableofcontents

\section{Introduction}

In rational conformal field theory, we can construct two-dimensional conformal field theories from a chiral conformal field theory. Mathematically, it corresponds to finding nonnegative integral combinations of some functions on the complex upper half plane which is invariant under the action of the modular group $SL(2,\mathbb{Z})$. Such combinations are called modular invariants. It is well-known that the modular invariants for the $SU(2)$ Wess--Zumino--Witten models admit an A-D-E classification i.e. they correspond to the simply-laced Dynkin diagrams \cite{cappelliade}\cite{MR906000}.

In the framework of algebraic quantum field theory, a chiral conformal field theory is modeled by a conformal net. When a conformal net $\cala$ satisfies the property called \emph{complete rationality} \cite[Definition 8]{MR1838752}, it produces a modular tensor category as its Doplicher--Haag--Roberts (DHR) category \cite{MR297259}\cite{MR334742}, which is denoted by $\rep \cala$ in this article, see \cite[Corollary 37]{MR1838752}. In this case, modular invariants can be formulated as matrices with nonnegative integer coefficients indexed by the simple objects of $\rep \mathcal{A}$ that commute with the $S$ and $T$ matrices of the modular tensor category.

It was shown by B\"{o}ckenhauer, Evans and Kawahigashi \cite{MR1652746}\cite{MR1671970}\cite{MR1706884}\cite{MR1729094}\cite{MR1777347} that \emph{$\alpha$-induction} \cite[Proposition 3.9]{MR1332979} is a powerful categorical tool for producing modular invariants. Namely, when a finite index standard extension $\mathcal{A} \subset \mathcal{B}$ of a completely rational conformal net $\mathcal{A}$ \cite[Section 3]{MR1332979} is given, the $\alpha$-induction for the extension assigns to each DHR endomorphism $\lambda$ of $\mathcal{A}$ some (not necessarily DHR) endomorphisms $\alpha^\pm (\lambda)$ of $\mathcal{B}$ and indeed $\langle \alpha^+(\lambda),\alpha^-(\mu) \rangle \coloneqq \dim \hom (\alpha^+(\lambda),\alpha^-(\mu))$ for simple objects $\lambda$ and $\mu$ of $\rep \mathcal{A}$ is a modular invariant matrix \cite[Theorem 5.7]{MR1729094}. Moreover, Rehren \cite{MR1754521} showed that the matrix is physical i.e. he constructed a commutative Q-system $\bigoplus_{\lambda, \mu} \langle \alpha^+(\lambda),\alpha^-(\mu) \rangle \lambda \boxtimes \overline{\mu}$ in $\rep \mathcal{A} \boxtimes (\rep \mathcal{A})^{\mathrm{rev}}$, and it was shown by Kong and Runkel \cite[Theorem 3.22]{MR2551797} and Bischoff, Kawahigashi and Longo \cite[Proposition 5.2]{MR3424476} that any \emph{physical} modular invariant matrix can be obtained in this fashion.

When a group $G$ acts on $\mathcal{A}$, we can consider the notion of a \emph{$G$-twisted representation} of $\mathcal{A}$, which generalizes the notion of a DHR endomorphism of $\mathcal{A}$ \cite{MR2183964}. Recently, the notion of the $\alpha$-induction for $G$-twisted representations, which we call \emph{$G$-equivariant $\alpha$-induction}, was introduced by Nojima \cite{MR4153896}. Then, it is natural to expect that $G$-equivariant $\alpha$-induction is also a powerful tool for producing modular invariants in some sense. In this article, we generalize some of the classical results introduced above, which gives some evidence for this expectation. Indeed, in some cases, our construction can produce modular invariants of the fixed point net $\cala^G$ (Theorem \ref{theorem_modular}).  

The main results of this article are twofold. First, we generalize Rehren's construction of $\alpha$-induction Q-systems \cite{MR1754521} to $G$-equivariant $\alpha$-induction (Theorem \ref{mainthm1}). Indeed, this theorem can be stated in a purely algebraic setting. Namely, we construct the \emph{$G$-equivariant $\alpha$-induction Frobenius algebra} for a neutral symmetric special simple $G$-equivariant Frobenius algebra in a split spherical $G$-braided fusion category (see Subsections \ref{subsection_gbraided_cat} and \ref{subsection_equiv_frob} for terminologies). For this, we have to introduce the notion of the \emph{neutral double} of a split semisimple $G$-braided multitensor category as an ambient category (Definition \ref{def_neutral_double}). Second, we generalize \cite[Proposition 4.18]{MR3424476}, which is an important ingredient in the characterization of physical modular invariants in terms of $\alpha$-induction \cite[Proposition 5.2]{MR3424476}. Namely, we introduce the notion of the $G$-equivariant full center of a neutral symmetric special $G$-equivariant Frobenius algebra in a split spherical $G$-braided fusion category with nonzero dimension (Definition \ref{def_fullcenter}) and show that it indeed coincides with the corresponding $G$-equivariant $\alpha$-induction Frobenius algebra when both of them are well-defined (Theorem \ref{mainthm2}).

This article is structured as follows. In Section \ref{section_preliminaries}, we review some preliminaries on bicategories and multitensor categories equipped with group actions. In particular, we recall the notions of a $G$-braided multitensor category and a $G$-equivariant Frobenius algebra for a group $G$. We also introduce an appropriate notion of equivalence between bicategories with group actions. In Section \ref{section_equiv_alpha}, we give the categorical definition of $G$-equivariant $\alpha$-induction. For this, we define an action of $G$ on the bicategory of $G$-equivariant Frobenius algebras. In Section \ref{section_equiv_frob}, we define the notion of the neutral double of a split semisimple $G$-braided multitensor category and construct the $G$-equivariant $\alpha$-induction Frobenius algebras. Finally, in Section \ref{section_fullcenter}, we define the notion of a $G$-equivariant full center and show the coincidence of $G$-equivariant $\alpha$-induction Frobenius algebras and $G$-equivariant full centers.

\section{Preliminaries}
\label{section_preliminaries}

\subsection{Group actions on bicategories}

In this subsection, we recall the notion of a group action on a monoidal category and that on a bicategory and give an appropriate notion of equivalence between bicategories with group actions, which turns out to be equivalent to the notion introduced in \cite{MR3936135}. Moreover, we prove a coherence theorem for $G$-actions on bicategories, see Theorem \ref{theorem_equivcoh}, cf. \cite[Theorem 3.1]{MR3936135}.

An \emph{action} of a group $G$ on a monoidal category $\calc$ is a monoidal functor from the monoidal category $\underline{G}$ of the elements of $G$ with only identity morphisms to the monoidal category of monoidal endofunctors on $\calc$. 

We also have to consider group actions on bicategories in this article (see Subsections \ref{subsection_action_bicat_frob} and \ref{subsection_mainthm2} below). For basic notions of bicategories, see e.g. \cite[Section 1.5]{MR2094071}. Here we fix our notation: the composition of 1-cells of $\calc$ is denoted by $\otimes_\calc$ or $\otimes$, for we regard bicategories as generalizations of monoidal categories, and often omitted. The associativity constraint is denoted by $a^\calc$. The unit is denoted by $\mathbf{1}^\calc$, and the left and right unit constraints are denoted respectively by $l^\calc$ and $r^\calc$. A pseudofunctor $F: \calc \to \cald$ consists of a map $F$ between 0-cells, functors $F_{A,B}: \hom_\calc(A,B) \to \hom_\cald(F(A),F(B))$ for $A,B \in \obj(\calc)$, natural isomorphisms $J^F_{A,B,C}: F_{B,C}(-) \otimes F_{A,B}(-) \cong F_{A,C}(- \otimes -)$ for $A,B,C \in \obj(\calc)$ and $\varphi^F_A:\one^\cald_{F(A)} \cong F_{A,A}(\one^\calc_A)$ for $A \in \obj(\calc)$. A pseudonatural transformation $\tau : F \to F'$ consists of 1-cells $\tau^0_C: F(C) \to F'(C)$ for $C \in \obj(\calc)$ and invertible 2-cells $\tau_\lambda : \tau^0_{C'} F(\lambda) \cong F'(\lambda) \tau^0_{C}$ for 1-cells $\lambda: C \to C'$. The vertical and horizontal compositions of pseudonatural transformations (see also \cite[Section 1.1]{MR3936135}) are denoted respectively by $\circ$ and $\ast$. 

By stating the coherence theorem for pseudofunctors \cite[Subsection 2.3.3]{MR3076451} in the form of ``all diagrams commute", we can see that the 2-cells obtained by vertically and horizontally composing components of $a^\calc$, $l^\calc$ and $r^\calc$ of bicategories $\calc$ and $J^F$'s and $\varphi^F$'s of pseudofunctors $F$ are indeed canonical, which allows us to suppress these 2-cells. In particular, we may suppress the constraints of the tricategory of bicategories. We can also introduce some graphical representations in the tricategory of bicategories, which are used only in this subsection. A pseudonatural transformation is represented as an arrow from top to bottom, and the vertical and horizontal compositions of pseudonatural transformations are represented respectively by the vertical and horizontal concatenation of arrows. Note that we do not have to be careful of relative vertical positions when taking horizontal concatenations since the comparison constraints between pseudonatural transformations (see \cite[Section 1.1]{MR3936135}) are canonical by coherence.

An \emph{action} of a group $G$ on a bicategory $\calc$ is a monoidal pseudofunctor (i.e. a trihomomorphism between monoidal bicategories, see \cite[Chapters 2 and 3]{MR1261589}) from $\uuline{G}$ to $\en(\calc)$, where $\uuline{G}$ denotes $G$ regarded as a monoidal bicategory and $\en(\calc)$ denotes the monoidal bicategory of pseudofunctors from $\calc$ to itself, see \cite[Sections 1.1 and 2]{MR3936135}. More explicitly, an action $\gamma$ of $G$ on $\calc$ consists of pseudofunctors (biequivalences) $\gamma(g)$ from $\calc$ to itself, pseudonatural equivalences $\chi^{\gamma}_{g,h}: \gamma(g) \gamma(h) \simeq \gamma(gh)$, a pseudonatural equivalence $\iota^\gamma: \id_\calc \simeq \gamma(e)$ and invertible modifications $\omega^{\gamma}_{g,h,k}: \chi^{\gamma}_{g h,f} \circ (\chi^{\gamma}_{g,h} \ast \id_{\gamma(f)}) \cong \chi^{\gamma}_{g,hf} \circ (\id_{\gamma(g)} \ast \chi^{\gamma}_{h,f})$, $\kappa^{\gamma}_{g}: \chi^{\gamma}_{e,g} \circ (\iota^\gamma \ast \id_{\gamma(g)}) \cong \id_{\gamma(g)}$ and $\zeta^{\gamma}_{g}: \chi^{\gamma}_{g,e} \circ (\id_{\gamma(g)} \ast \iota^\gamma) \cong \id_{\gamma(g)}$ for $g,h,k \in G$. Graphically, let a fork with two inputs $\gamma(g)$ and $\gamma(h)$ and one output $\gamma(g h)$ denote $\chi^\gamma_{g,h}$, and let a small circle with one output $\gamma(e)$ denote $\iota^\gamma$. Then, $\omega^\gamma$, $\kappa^\gamma$ and $\zeta^\gamma$ correspond respectively to associativity and left and right unit properties of an algebra in a weak sense, which satisfy the pentagon and triangle axioms \cite[Definition 3.1]{MR1261589}.

For a 1-cell $\lambda$ of $\calc$ and $g \in G$, the 1-cell $\gamma(g) (\lambda)$ is often denoted by ${}^g \lambda$. Moreover, ${}^g \lambda \mu$ denotes ${}^g (\lambda) \otimes \mu$ for $g \in G$ and 1-cells $\lambda$ and $\mu$, and ${}^g f \otimes f'$ denotes ${}^g (f) \otimes f'$ for $g \in G$ and 2-cells $f$ and $f'$ in this article.

An appropriate notion of morphisms between (strict) 2-categories with unital (see \cite[Definition 2.1]{MR3936135}) group actions is given in \cite[Definition 2.3]{MR3936135}. Here, we define the notion of an equivalence in a general setting.

First, recall that for a biequivalence $F: \calc \to \cald$, we can take a pseudofunctor $F^{-1}: \cald \to \calc$ with pseudonatural equivalences $\mathrm{ev}^F : F^{-1}F \simeq \id_\calc$ and $\mathrm{coev}^F: \id_\cald \simeq FF^{-1}$ by fixing data that consist of 0-cells $F^{-1}(D) \in \obj(\calc)$, equivalence 1-cells $\mathrm{coev}^{F,0}_D: D \to F(F^{-1}(D))$, left adjoint inverses $(\mathrm{coev}^{F,0}_D)^\vee$ of $\mathrm{coev}^{F,0}_D$ (i.e. left duals of $\mathrm{coev}^{F,0}_D$ in the bicategory $\cald$ with invertible evaluation and coevaluation maps) for $D \in \obj(\cald)$ and left adjoint functors $F_{C,C'}^\vee$ of $F_{C,C'}$ for $C,C' \in \obj(\calc)$, see e.g. \cite[Proposition 1.5.13]{MR2094071}. The set of these data is also denoted by $F^{-1}$ and referred to as an \emph{adjoint inverse} of $F$. Another choice of an adjoint inverse only yields a pseudonaturally equivalent pseudofunctor $F^{-1}$ by a standard duality argument. The pseudonatural equivalences $\mathrm{ev}^F$ and $\mathrm{coev}^F$ are graphically represented by arcs as in the case of duality in bicategories. Their adjoint inverses are represented by opposite arcs. 

We also recall that we have a natural isomorphism $J^{F^\vee}_{A, B, C}: F_{B, C}^\vee(-) \otimes F_{A, B}^\vee(-) \cong F_{A, C}^\vee(- \otimes -)$ for $A, B, C \in \obj(\calc)$ and an isomorphism $\varphi^{F^\vee}_A: \mathbf{1}^{\calc}_A \cong F^\vee_{A, A}(\mathbf{1}^\cald_{F(A)})$ for $A \in \obj(\calc)$ with the coherence conditions as that for pseudofunctors. Indeed, they are defined by putting $(J^{F^\vee}_{A, B, C})_{\lambda,\mu} \coloneqq F^\vee_{A, C}((\widetilde{\mathrm{coev}}^{F}_\lambda)^{-1} \otimes (\widetilde{\mathrm{coev}}^F_\mu)^{-1}) \circ (\widetilde{\mathrm{ev}}^F_{F^\vee_{B, C}(\lambda) F^\vee_{A, B}(\mu)})^{-1}$ for $\lambda: F(B) \to F(C)$ and $\mu: F(A) \to F(B)$, and $\varphi^{F^\vee}_A \coloneqq (\widetilde{\mathrm{ev}}^F_{\mathbf{1}^\calc_A})^{-1}$, where $\widetilde{\mathrm{ev}}^F_\nu$ and $\widetilde{\mathrm{coev}}^F_\nu$ respectively denote the components of the evaluation and coevaluation maps of $F_{A,B}$ for $\nu: F(A) \to F(B)$. Then, as in the case of pseudofunctors, we may suppress $J^{F^\vee}$'s and $\varphi^{F^\vee}$'s. By the definition of $J^{F^\vee}$'s and $\varphi^{F^\vee}$'s, the naturality of $\widetilde{\mathrm{ev}}^F$ and $\widetilde{\mathrm{coev}}^F$ and conjugate equations for $F$, we can see that $\widetilde{\mathrm{ev}}^F$ and $\widetilde{\mathrm{coev}}^F$ are monoidal i.e. $\widetilde{\mathrm{coev}}^F_{\lambda \mu} = \widetilde{\mathrm{coev}}^F_{\lambda} \otimes \widetilde{\mathrm{coev}}^F_{\mu}$ for $\lambda :F(B) \to F(C)$ and $\mu :F(A) \to F(B)$, and $\widetilde{\mathrm{coev}}^F_{\mathbf{1}_{F(A)}^\cald} = \id_{\mathbf{1}_{F(A)}^\cald}$ (and similar statements for $\widetilde{\mathrm{ev}}^F$), where $(J^{F^\vee}_{A, B, C})_{\lambda, \mu}$ and $\varphi^{F^\vee}_A$ are suppressed.

When a biequivalence $F:\calc \to \cald$ is given, we can transport $H \in \en(\calc)$ to $\en(\cald)$ by fixing an adjoint inverse $F^{-1}$ of $F$ and putting $\ad(F)(H) \coloneqq F H F^{-1} \in \en(\cald)$.

\begin{lem}
\label{adismonoidallem}
Let $F: \calc \to \cald$ be a biequivalence between bicategories. Then $\ad (F): \en(\calc) \to \en(\cald)$ can be regarded as a monoidal pseudofunctor. Another choice of an adjoint inverse of $F$ only yields a monoidally equivalent one (i.e. there exists a triequivalence with an identical 1-cell).

\begin{proof}
Let $H, K, L, H', K' \in \obj(\en(\calc))$ and let $\sigma \in \hom(H, H')$, $\tau \in \hom(K, K')$ and $\rho \in \hom(H', H'')$. Define functors $\ad(F)_{H, H'}$ to be
\begin{align*}
	\id_{F} \ast - \ast \id_{F^{-1}}: \hom(H,H') \to \hom(\ad(F)(H),\ad(F)(H')).
\end{align*}
Define $J^{\ad(F)}_{\rho,\sigma}: \ad(F)(\rho) \circ \ad(F)(\sigma) \cong \ad(F)(\rho \circ \sigma)$ and $\varphi^{\ad(F)}_H : \id_{\ad(F)(H)} \cong \ad(F)(\id_H)$ to be canonical isomorphisms, which makes $\ad(F)$ into a pseudofunctor by coherence. Define pseudonatural equivalences $\chi_{H,K}^{\ad(F),0} \coloneqq \id_{F H} \ast \mathrm{ev}^F \ast \id_{KF^{-1}}: \ad(F)(H) \ad(F)(K) \simeq \ad(F)(H K)$ and define invertible modifications $\chi^{\ad(F)}_{\sigma,\tau}: \chi_{H',K'}^{\ad(F),0} \circ (\ad(F)(\sigma) \ast \ad(F)(\tau)) \cong \ad(F)(\sigma \ast \tau) \circ \chi_{H,K}^{\ad(F),0}$ to be canonical isomorphisms, which give a pseudonatural equivalence $\chi^{\ad(F)} \coloneqq (\chi^{\ad(F),0},\chi^{\ad(F)}): \ad(F)(-) \circ \ad(F)(-) \simeq \ad(F)(- \circ -)$ by coherence. Put $\iota^{\ad(F)} \coloneqq \mathrm{coev}^F: \id_\cald \simeq \ad(F)(\id_\calc) = FF^{-1}$. Define modifications $\omega^{\ad(F)}_{H,K,L}: \chi_{H K,L}^{\ad(F),0} \circ (\chi_{H,K}^{\ad(F),0} \ast \id_{\ad(F)(L)}) \cong \chi_{H,K L}^{\ad(F),0} \circ (\id_{\ad(F)(H)} \ast \chi_{K,L}^{\ad(F),0})$ to be canonical isomorphisms, which give an invertible modification $\omega^{\ad(F)}: \chi^{\ad(F)} \circ (\chi^{\ad(F)} \ast \id_{\ad(F)}) \cong \chi^{\ad(F)} \circ (\id_{\ad(F)} \ast \chi^{\ad(F)})$ by coherence. Then, Axiom (HTA1) in \cite[Definition 3.1]{MR1261589} holds by coherence. 

Put $f_D \coloneqq \mathrm{coev}^{F,0}_D$ for $D \in \obj(\cald)$. Let $F^\vee(\lambda)$ denote $F_{C,C'}^\vee(\lambda)$ for $\lambda \in \hom_\cald(F(C),F(C'))$ and $C,C' \in \obj(\calc)$. Since $\mathrm{ev}^{F,0}_C = F^\vee (f_{F(C)}^\vee)$ for $C \in \obj(\calc)$ by construction, we can define an invertible 2-cell
\begin{align*}
	\xi_C^F: ((\id_F \ast \mathrm{ev}^{F}) \circ (\mathrm{coev}^{F} \ast \id_F))^0_C = F(\mathrm{ev}^{F,0}_C) \mathrm{coev}^{F,0}_{F(C)} = FF^\vee(f_{F(C)}^\vee) f_{F(C)} \to (\id_F)^0_C = \mathbf{1}^\cald_{F(C)}
\end{align*}
to be $\mathrm{ev}_{f_{F(C)}}((\widetilde{\mathrm{coev}}^F_{f^\vee_{F(C)}})^{-1} \otimes \id_{f_{F(C)}})$. Since $F^{-1}(\lambda) = F^\vee (f_{D'} \lambda f_D^\vee)$ for $\lambda \in \hom_\cald(D,D')$ by construction, we can also define an invertible 2-cell
\begin{align*}
	&\tilde{\xi}^F_D: ((\mathrm{ev}^{F} \ast \id_{F^{-1}}) \circ (\id_{F^{-1}} \ast \mathrm{coev}^{F}))_D^0 = \mathrm{ev}^{F,0}_{F^{-1}(D)} F^{-1}(\mathrm{coev}^{F,0}_D) = F^\vee(f^\vee_{FF^{-1}(D)} f_{FF^{-1}(D)} f_D f_D^\vee) \\ 
	&\quad \quad \to (\id_{F^{-1}})^0_D = \mathbf{1}^\calc_{F^{-1}(D)}
\end{align*}
to be $F^\vee(\mathrm{ev}_{f_{FF^{-1}(D)}} \otimes \mathrm{coev}_{f_D}^{-1})$ for $D \in \obj(\cald)$. We show that $\xi^F$ and $\tilde{\xi}^F$ are indeed modifications. For this, recall that $\mathrm{coev}^F_\nu = (\widetilde{\mathrm{coev}}^F_{\nu} \otimes \id_{f_D}) (\id_{f_{D'} \lambda} \otimes \mathrm{ev}^{-1}_{f_D})$ for $\nu \in \hom_\cald(D,D')$ and $\mathrm{ev}^F_\rho = (\widetilde{\mathrm{ev}}^F_\rho \otimes \id_{F^\vee(f^\vee_{F(C)})}) F^\vee(\mathrm{ev}_{f_{F(C')}} \otimes \id_{F(\rho) f^\vee_{F(C)}})$ for $\rho \in \hom_\calc(C,C')$ by construction. Then, we can see by direct computations that $\xi^F$ and $\tilde{\xi}^F$ are modifications.

Then, put $\kappa^{\ad(F)}_H \coloneqq \xi^F \ast \id_{HF^{-1}}$ and $\zeta^{\ad(F)}_H \coloneqq \id_{FH} \ast \tilde{\xi}^F$ for $H \in \obj(\en(\calc))$, which define modifications $\kappa^{\ad(F)}: \chi^{\ad(F)} \circ (\iota^{\ad(F)} \ast \id_{\ad(F)}) \cong \id_{\ad(F)}$ and $\zeta^{\ad(F)}: \chi^{\ad(F)} \circ (\id_{\ad(F)} \ast \iota^{\ad(F)}) \cong \id_{\ad(F)}$ by coherence. To check that these modifications satisfy Axiom (HTA2) in \cite[Definition 3.1]{MR1261589}, it is enough to show $\id_{\mathrm{ev}^F} \otimes (\id_{F^{-1}} \ast \xi^F) = \id_{\mathrm{ev}^F} \otimes (\tilde{\xi}^F \ast \id_F)$. For this, note that the component of the comparison constant $\mathrm{ev}^F \circ (\id_{F^{-1}F} \ast \mathrm{ev}^F) \cong \mathrm{ev}^F \circ (\mathrm{ev}^F \ast \id_{F^{-1}F})$ for $C \in \obj(\calc)$ is given by $F^\vee(\mathrm{ev}_{f_{F(C)}}\otimes (\widetilde{\mathrm{coev}}^F_{f^\vee_{F(C)}})^{-1} \otimes \id_{FF^{-1}F(C)})$. Then, we can find that it is enough to show $\mathrm{ev}_{f_{F(C)}} \otimes \id_{f_{F(C)}^\vee f_{F(C)}} = \id_{f_{F(C)}^\vee f_{F(C)}} \otimes \mathrm{ev}_{f_{F(C)}}$, which follows since $f_{F(C)}^\vee$ is an adjoint inverse. Thus, $\ad(F)$ is a monoidal pseudofunctor.

Finally, we show that another choice of an adjoint inverse $F_2^{-1}$ yields an equivalent monoidal pseudofunctor $\ad_2 (F)$. Let $\tilde{f}_D$ for $D \in \obj(\cald)$ denote $\mathrm{coev}^{F,0}_D$ for $F_2^{-1}$, and let $F^\vee_2$ denote $F^\vee$ for $F^{-1}_2$. Also, let $\mathrm{coev}_2^F$ and $\mathrm{ev}^F_2$ respectively denote $\mathrm{coev}^F$ and $\mathrm{ev}^F$ for $F_2^{-1}$. For $H \in \obj(\en(\calc))$, put $\eta_H^0 \coloneqq \id_{FH} \ast \tau: \ad(F)(H) \simeq \ad_2(F)(H)$, where $\tau: F^{-1} \simeq F_2^{-1}$ denotes the pseudonatural equivalence obtained from $\mathrm{coev}^F_2 (\mathrm{coev}^F)^\vee: FF^{-1} \simeq FF_2^{-1}$ and ${}^\vee$ denotes a left adjoint inverse in $\en(\cald)$, by duality. For $H, H' \in \obj(\en(\calc))$ and $\sigma \in \hom(H, H')$, define $\eta_\sigma: \eta^0_{H'} \circ \ad(F)(\sigma) \cong \ad_2(F)(\sigma) \circ \eta^0_H$ to be the canonical isomorphism, which gives a pseudonatural equivalence $\eta=(\eta^0, \eta): \ad(F) \simeq \ad_2(F)$ by coherence. Then, for $C \in \obj(\calc)$, define an invertible 2-cell 
\begin{align*}
	&\alpha_C: (\mathrm{ev}^F_2 \circ (\tau \ast \id_F))_C^0 = F^\vee_2(\tilde{f}^\vee_{F(C)}) F^\vee(f_{FF_2^{-1}F(C)}^\vee f_{FF_2^{-1}F(C)} \tilde{f}_{F(C)} f^\vee_{F(C)} f^\vee_{FF^{-1}F(C)} f_{FF^{-1}F(C)}) \\
	&\quad \quad \to \mathrm{ev}^{F,0}_C = F^\vee(f_{F(C)}^\vee)
\end{align*}
to be 
\begin{align*}
	F^\vee((\mathrm{ev}_{\tilde{f}_{F(C)}} \otimes \id_{f^\vee_{F(C)}})(\id_{\tilde{f}^\vee_{F(C)}} \otimes \mathrm{ev}_{f_{FF_2^{-1}F(C)}} \otimes \id_{\tilde{f}_{F(C)} f^\vee_{F(C)}} \otimes \mathrm{ev}_{f_{FF^{-1}F(C)}}))(\beta_{\tilde{f}^\vee_{F(C)}} \otimes \id_{\tau^0_{F(C)}}),
\end{align*}
where $\beta$ denotes the canonical isomorphism $F^\vee_2 \cong F^\vee$. It is routine to check that $\alpha$ is indeed a modification $\mathrm{ev}^F_2 \circ (\tau \ast \id_F) \cong \mathrm{ev}^F$. For $H,K \in \obj(\en(\calc))$, define $\Pi_{H,K}: \chi^{\ad_2(F),0}_{H,K} \circ (\eta_H^0 \ast \eta_K^0) \cong \eta_{HK}^0 \circ \chi^{\ad(F),0}_{H,K}$ to be $\id_{FH} \ast \alpha \ast \id_K \ast \id_{\tau}$, which gives a modification $\Pi: \chi^{\ad_2(F)} \circ (\eta \ast \eta) \cong \eta \circ \chi^{\ad(F)}$ by coherence. Then, $\eta$ and $\Pi$ satisfy the condition in \cite[pp. 21--22]{MR1261589} by coherence. Finally, define $M: \eta_{\id_\calc}^0 \circ \iota^{\ad(F)} \cong \iota^{\ad_2(F)}$ to be the composition of the modifications in Figure \ref{graphical_lemma_adismonoidal_proof}, where $\tilde{\xi}_2^F$ denotes $\tilde{\xi}^F$ for $F_2^{-1}$. Note that e.g. $\xi^F$ indeed denotes $(\xi^F \ast \id_{F_2^{-1}}) \otimes \id_{\mathrm{coev}^F_2}$ etc. Then, $\eta, \Pi$ and $M$ satisfy the conditions in \cite[pp. 23--24]{MR1261589} by $\id_{\mathrm{ev}^F} \otimes (\id_{F^{-1}} \ast \xi^F_2) = \id_{\mathrm{ev}^F} \otimes (\tilde{\xi}^F_2 \ast \id_F)$. Thus, $(\eta,\Pi,M)$ is a monoidal pseudonatural equivalence $\ad(F) \simeq \ad_2(F)$.
\begin{figure}[htb]
	\centering
	\begin{tikzpicture}
		\draw[<-] (0,0) -- (0,1) arc (180:0:0.25) -- (0.5,0);
		\draw[->] (1.5,0.5) -- (2.5,0.5);
		\node[block] at (0.5,0.5){$\tau$};
		\node at (2,1){$\tilde{\xi}^{F-1}_2$};
		\node at (0,-0.25){$F$};
		\node at (0.5,-0.25){$F_2^{-1}$};
		\node at (1,1){$F^{-1}$};
		\begin{scope}[shift={(3,0)}]
			\draw[<-] (0,0) -- (0,1) arc (180:0:0.25) -- (0.5,0.25) arc (180:360:0.25) arc (180:0:0.25) -- (1.5,0);
			\draw[->] (2,0.5) -- (3,0.5);
			\node[block] at (0.5,0.5){$\tau$};
			\node at (0,-0.25){$F$};
			\node at (1.5,-0.25){$F_2^{-1}$};
			\node at (1,1){$F^{-1}$};
			\node at (2.5,1){$\alpha$};
		\end{scope}
		\begin{scope}[shift={(6.5,0)}]
			\draw[<-] (0,0) -- (0,1) arc (180:0:0.25) -- (0.5,0.25) arc (180:360:0.25) arc (180:0:0.25) -- (1.5,0);
			\draw[->] (2,0.5) -- (3,0.5);
			\node at (0,-0.25){$F$};
			\node at (1.5,-0.25){$F_2^{-1}$};
			\node at (1,1){$F^{-1}$};
			\node at (2.5,1){$\xi^{F}$};
		\end{scope}
		\begin{scope}[shift={(10,0.5)}]
			\draw[<-] (0,0) arc (180:0:0.25);
			\node at (0,-0.25){$F$};
			\node at (0.5,-0.25){$F_2^{-1}$};
		\end{scope}
	\end{tikzpicture}
	\caption{A modification $M$}
	\label{graphical_lemma_adismonoidal_proof}
\end{figure}
\end{proof}
\end{lem}

\begin{rem}
	The proof above shows that a biequivalence is indeed a 2-dualizable pseudofunctor in the sense of the tricategory version of \cite[Definitions 5.1 and 6.1]{campbell2022riemann}.
\end{rem}

\begin{defi}
	\label{definition_g_biequivalence}
Let $(\calc,\gamma^\calc)$ and $(\cald,\gamma^\cald)$ be pairs of bicategories and actions of a group $G$. A \emph{$G$-biequivalence} between $(\calc,\gamma^\calc)$ and $(\cald,\gamma^\cald)$ is the pair $F = (F,\eta^F)$ of a biequivalence $F: \calc \to \cald$ and a monoidal equivalence $\eta^F: \ad(F) \circ \gamma^\calc \simeq \gamma^\cald$.
\end{defi}

The existence of a $G$-biequivalence between bicategories does not depend on a choice of $F^{-1}$ by Lemma \ref{adismonoidallem}.

\begin{lem}
	$G$-biequivalence is indeed an equivalence relation.
	
	\begin{proof}
		Reflexivity follows since $\ad(\id_\calc) \simeq \id_{\en(\calc)}$ as monoidal pseudofunctors by coherence. We show transitivity. Let $F: \calc \to \cald$ and $H: \cald \to \mathcal{E}$ be biequivalences. We fix $F^{-1}$ and $H^{-1}$, set $(HF)^{-1} (E) \coloneqq F^{-1} H^{-1}(E)$, $\mathrm{coev}^{HF,0}_E \coloneqq H(f_{H^{-1}(E)}) h_E$ and $(\mathrm{coev}^{HF,0}_E)^\vee \coloneqq h_E^\vee H(f_{H^{-1}(E)}^\vee)$ for $E \in \obj(\mathcal{E})$, where $f \coloneqq \mathrm{coev}^{F,0}$ and $h \coloneqq \mathrm{coev}^{H,0}$, and set $(HF)^\vee \coloneqq F^\vee H^\vee$. Then, for $\lambda \in \hom_{\mathcal{E}}(E,E')$, we have an invertible 2-cell
		\begin{align*}
			F^\vee(\widetilde{\mathrm{ev}}^H_{f_{H^{-1}(E)}} \otimes \id_{H^{-1}(\lambda)} \otimes \widetilde{\mathrm{ev}}^H_{f_{H^{-1}(E)}^\vee}): \ &(HF)^{-1}(\lambda) = F^\vee(H^\vee H(f_{H^{-1}(E')}) H^{-1}(\lambda) H^\vee H(f^\vee_{H^{-1}(E)})) \\
			&\to F^{-1} H^{-1} (\lambda) = F^\vee(f_{H^{-1}(E')} H^{-1}(\lambda) f^\vee_{H^{-1}(E)}),
		\end{align*}
		which defines a pseudonatural equivalence $\tau: (HF)^{-1} \simeq F^{-1} H^{-1}$. Moreover, we define an invertible 2-cell
		\begin{align*}
			&(\mathrm{ev}^H \circ \ad(H) (\mathrm{ev}^F) \circ (\tau \ast \id_{HF}))_C^0 = F^\vee (f^\vee_{F(C)} f_{F(C)} H^\vee(h^\vee_{HF(C)}) f^\vee_{H^{-1}HF(C)}) \\
			&\to \mathrm{ev}^{HF,0}_C = F^\vee (H^\vee(h^\vee_{HF(C)}) H^\vee H (f^\vee_{H^{-1}HF(C)}))
		\end{align*}
		for $C \in \obj(\calc)$ to be $F^\vee((\id_{H^\vee(h^\vee_{HF(C)})} \otimes (\widetilde{\mathrm{ev}}^H_{f^\vee_{H^{-1}HF(C)}})^{-1})(\mathrm{ev}_{f_{F(C)}} \otimes \id_{H^\vee(h^\vee_{HF(C)}) f^\vee_{H^{-1}HF(C)}}))$, which indeed defines a modification $\mathrm{ev}^H \circ \ad(H) (\mathrm{ev}^F) \circ (\tau \ast \id_{HF}) \cong \mathrm{ev}^{HF}$. Then, we obtain a monoidal equivalence $\ad(HF) \simeq \ad(H) \ad(F)$ as in the final part of the proof of Lemma \ref{adismonoidallem}, and therefore $G$-biequivalence is transitive. Finally, we show symmetry. Since $\ad(F^{-1} F) \simeq \ad(F^{-1}) \ad(F)$ as monoidal pseudofunctors for a biequivalence $F$ by the argument so far, it is enough to show that for biequivalences $F, \tilde{F}: \calc \to \cald$ with a pseudonatural equivalence $\tau: F \simeq \tilde{F}$, we have a monoidal equivalence $\ad(F) \simeq \ad(\tilde{F})$. Define $\tilde{\tau}: F^{-1} \simeq \tilde{F}^{-1}$ to be $(\mathrm{ev}^F \ast \id_{\tilde{F}^{-1}}) \circ (\id_{F^{-1}} \ast \tau^\vee \ast \id_{\tilde{F}^{-1}}) \circ (\id_{F^{-1}} \ast \mathrm{coev}^{\tilde{F}})$. Also, define an invertible modification $\alpha: \mathrm{ev}^{\tilde{F}} \circ (\tilde{\tau} \ast \tau) \cong \mathrm{ev}^F$ by Figure \ref{graphical_lemma_equivalence_proof_alpha}. Then, by putting $\eta_H \coloneqq \tau \ast \id_H \ast \tilde{\tau}$ and $\Pi_{H,K} \coloneqq \id_{FH} \ast \alpha \ast \id_{KF^{-1}}$ for $H, K \in \obj(\en(\calc))$, we obtain a pseudonatural equivalence $\eta$ and an invertible modification $\Pi$ with the condition in \cite[pp. 21--22]{MR1261589} as in the proof of Lemma \ref{adismonoidallem}. Define a modification $M$ by Figure \ref{graphical_lemma_equivalence_proof_m}. We can see that $(\eta, \Pi, M)$ is a monoidal equivalence $\ad (F) \simeq \ad(\tilde{F})$ by the duality of $\tau$. 
		\begin{figure}[htb]
			\begin{tabular}{cc}
				\begin{minipage}[b]{0.45\hsize}
					\centering
					\begin{tikzpicture}
						\draw[<-] (0,0.25) -- (0,-1) arc (180:360:0.25) -- (0.5,0) arc (180:0:0.25) -- (1,-1) arc (180:360:0.25) -- (1.5,0.25);
						\draw[->] (2,-0.5) -- (2.5,-0.5);
						\node[block] at (0.5,-0.5){$\tau^\vee$}; 
						\node[block] at (1.5,-0.5){$\tau$}; 
						\node at (0,0.5){$F^{-1}$};
						\node at (1.5,0.5){$F$};
						\node at (1.75,-1){$\tilde{F}$};
						\node at (2.25,0){$\xi^{\tilde{F}}$};
						\begin{scope}[shift={(3,0)}]
							\draw[<-] (0,0.25) -- (0,-1) arc (180:360:0.25) -- (0.5,0.25);
							\draw[->] (1,-0.5) -- (1.5,-0.5);
							\node[block] at (0.5,-0.75){$\tau^\vee$}; 
							\node[block] at (0.5,-0.25){$\tau$};
							\node at (0,0.5){$F^{-1}$};
							\node at (0.5,0.5){$F$};
							\node at (1.25,0){$\mathrm{ev}_{\tau}$};
						\end{scope}
						\begin{scope}[shift={(5,-0.75)}]
							\draw[<-] (0,0.25) arc (180:360:0.25);
							\node at (0,0.5){$F^{-1}$};
							\node at (0.5,0.5){$F$};
						\end{scope}
					\end{tikzpicture}
					\caption{A modification $\alpha$}
					\label{graphical_lemma_equivalence_proof_alpha}
				\end{minipage}
				\begin{minipage}[b]{0.45\hsize}
					\centering
					\begin{tikzpicture}
						\draw[<-] (0,-0.25) -- (0,1) arc (180:0:0.25) -- (0.5,0) arc (180:360:0.25) -- (1,1) arc (180:0:0.25) -- (1.5,-0.25);
						\draw[->] (2,0.5) -- (2.5,0.5);
						\node[block] at (1,0.5){$\tau^\vee$}; 
						\node[block] at (0,0.5){$\tau$}; 
						\node at (1.5,-0.5){$\tilde{F}^{-1}$};
						\node at (0,-0.5){$\tilde{F}$};
						\node at (-0.25,1){$F$};
						\node at (2.25,1){$\xi^{F}$};
						\begin{scope}[shift={(3,0)}]
							\draw[<-] (0,-0.25) -- (0,1) arc (180:0:0.25) -- (0.5,-0.25);
							\draw[->] (1,0.5) -- (1.5,0.5);
							\node[block] at (0,0.75){$\tau^\vee$}; 
							\node[block] at (0,0.25){$\tau$}; 
							\node at (0.5,-0.5){$\tilde{F}^{-1}$};
							\node at (0,-0.5){$\tilde{F}$};
							\node at (1.25,1){$\mathrm{coev}^{-1}_\tau$};
						\end{scope}
						\begin{scope}[shift={(5,0.5)}]
							\draw[<-] (0.5,0) arc (0:180:0.25);
							\node at (0.5,-0.25){$\tilde{F}^{-1}$};
							\node at (0,-0.25){$\tilde{F}$};
						\end{scope}
					\end{tikzpicture}
					\caption{A modification $M$}
					\label{graphical_lemma_equivalence_proof_m}	
				\end{minipage}
			\end{tabular}
		\end{figure}
	\end{proof}
\end{lem}

By definition \cite[Section 3.3]{MR1261589}, for a $G$-biequivalence $F : (\calc,\gamma^\calc) \to (\cald,\gamma^\cald)$, the monoidal equivalence $\eta^F$ consists of pseudonatural equivalences $\eta^F_g : \ad(F) (\gamma^\calc(g)) \simeq \gamma^\cald(g)$, invertible modifications $\Pi^{F}_{g,h}: \chi^{\gamma^\cald}_{g,h} \circ (\eta_g^F \ast \eta_h^F) \cong \eta^F_{gh} \circ \ad(F)(\chi^{\gamma^\calc}_{g,h}) \circ \chi^{\ad(F),0}_{\gamma^\calc(g),\gamma^\cald(h)}$ for $g,h \in G$ and an invertible modification $M^F: \eta^F_e \circ \ad(F)(\iota^{\gamma^\calc}) \circ \iota^{\ad(F)} \cong \iota^{\gamma^\cald}$.  Let $\eta^F_g$ be graphically represented by a fork with three inputs $F$, $\gamma^\calc(g)$ and $F^{-1}$ and one output $\gamma^\cald(g)$. 

Next, we compare our definition with that in \cite[Definition 2.3]{MR3936135}. Suppose we are given a $G$-biequivalence $F:(\calc, \gamma^\calc) \to (\cald, \gamma^\cald)$. Then, we obtain $\tilde{\eta}^F_g : F \gamma^\calc(g) \simeq \gamma^\cald(g) F$ for $g \in G$ by duality: namely, we put $\tilde{\eta}^F_g \coloneqq (\eta^F_g \ast \id_F) \circ (\id_{F \gamma^\calc(g)} \ast (\mathrm{ev}^{F})^{\vee})$. Moreover, we can define invertible modifications $\tilde{\Pi}^F_{g,h}: (\chi^{\gamma^\cald}_{g,h} \ast \id_F) \circ (\id_{\gamma^\cald(g)} \ast \tilde{\eta}_h^F) \circ (\tilde{\eta}^F_g \ast \id_{\gamma^\calc(h)}) \cong \tilde{\eta}^F_{gh} \circ (\id_F \ast \chi^{\gamma^\calc}_{g,h})$ for $g,h \in G$ and $\tilde{M}^F :\tilde{\eta}^F_e \circ (\id_F \ast \iota^{\gamma^\calc}) \cong \iota^{\gamma^\cald} \ast \id_F$ by Figures \ref{graphical_pi_to_tildepi} and \ref{graphical_m_to_tildem}. We can check that they satisfy the conditions in Figures \ref{graphical_coherence_omega}, \ref{graphical_coherence_kappa} and \ref{graphical_coherence_zeta}, where a crossing from $F \gamma^\calc(g)$ to $\gamma^\cald(g) F$ denotes $\tilde{\eta}^F_g$ for $g \in G$, with some standard computations, but note that we use $\id_{\mathrm{ev}^F} \otimes (\id_{F^{-1}} \ast \xi^F) = \id_{\mathrm{ev}^F} \otimes (\tilde{\xi}^F \ast \id_F)$ in the check for Figure \ref{graphical_coherence_zeta}. Conversely, when a biequivalence $F$ and a triple $(\tilde{\eta}^F,\tilde{\Pi}^F,\tilde{M}^F)$ with the conditions in Figures \ref{graphical_coherence_omega}, \ref{graphical_coherence_kappa} and \ref{graphical_coherence_zeta} are given, we can construct a triple $(\eta^F, \Pi^F, M^F)$ by putting $\eta^F_g \coloneqq (\id_{\gamma^\cald(g)} \ast (\mathrm{coev}^{F})^{\vee}) \circ (\tilde{\eta}^F_g \ast \id_{F^{-1}})$ for $g \in G$ and defining $\Pi_{g,h}^F$ for $g, h \in G$ and $M^F$ by Figures \ref{graphical_tildepi_to_pi} and \ref{graphical_tildem_to_m}. We can see that $(F, \eta^F, \Pi^F, M^F)$ is a $G$-biequivalence, using the following lemma to check the condition in \cite[p. 24]{MR1261589}. 
\begin{figure}[htb]
	\centering
	\begin{tikzpicture}
		\draw[->] (0,0) -- (0,-0.25) arc (180:360:0.25) arc (180:0:0.25) arc (180:360:0.25) arc (180:0:0.25) -- (2,-1.5);
		\draw (0.25,0) -- (0.25,-0.75) arc (180:360:0.5) -- (1.25,0);
		\draw[->] (0.75,-1.25) -- (0.75,-1.5);
		\draw[->] (2.5,-0.75) -- (3.5,-0.75);
		\node at (0.125,0.25){$F\gamma^\calc(g)$};
		\node at (1.25,0.25){$\gamma^\calc(h)$};
		\node at (0.75,-1.75){$\gamma^\cald(gh)$};
		\node at (2,-1.75){$F$};
		\node at (3,-0.25){$\Pi^F_{g,h}$};
		\begin{scope}[shift={(4,0)}]
			\draw[->] (0,0) -- (0,-0.5) arc (180:360:0.75)-- (1.5,-0.25) arc (180:0:0.25) -- (2,-1.5);
			\draw (0.25,0) -- (0.25,-0.5) arc (180:360:0.5) -- (1.25,0);
			\draw[->] (0.75,-1) -- (0.75,-1.5);
			\draw[->] (2.5,-0.75) -- (3.5,-0.75);
			\draw (0.75,-0.25) circle [radius=0.25];
			\node at (0.125,0.25){$F\gamma^\calc(g)$};
			\node at (1.25,0.25){$\gamma^\calc(h)$};
			\node at (0.75,-1.75){$\gamma^\cald(gh)$};
			\node at (2,-1.75){$F$};
			\node at (3,-0.25){$\mathrm{coev}^{-1}_{\mathrm{ev}^F}$};
		\end{scope}
		\begin{scope}[shift={(8,0)}]
			\draw[->] (0,0) -- (0,-0.5) arc (180:360:0.75)-- (1.5,-0.25) arc (180:0:0.25) -- (2,-1.5);
			\draw (0.25,0) -- (0.25,-0.5) arc (180:360:0.5) -- (1.25,0);
			\draw[->] (0.75,-1) -- (0.75,-1.5);
			\node at (0.125,0.25){$F\gamma^\calc(g)$};
			\node at (1.25,0.25){$\gamma^\calc(h)$};
			\node at (0.75,-1.75){$\gamma^\cald(gh)$};
			\node at (2,-1.75){$F$};	
		\end{scope}
	\end{tikzpicture}
	\caption{Construction of $\tilde{\Pi}^F$ from $\Pi^F$}
	\label{graphical_pi_to_tildepi}
\end{figure}
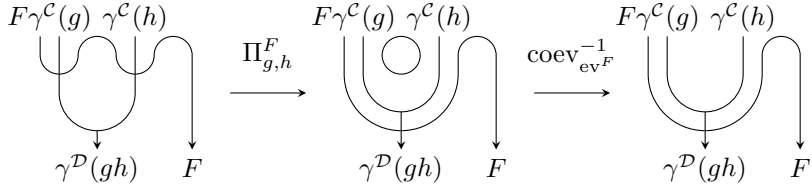
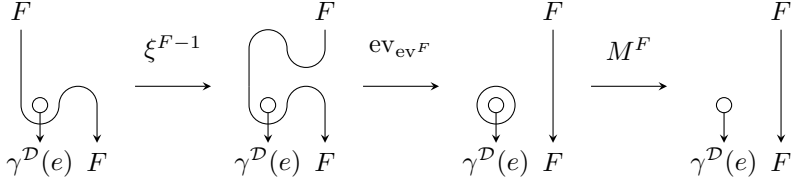
\begin{figure}[htb]
	\centering
	\begin{tikzpicture}
		\draw[->] (0,0.75) -- (0,-0.25) arc (180:360:0.25) arc (180:0:0.25) -- (1,-0.75);
		\draw[->] (0.25,-0.25) -- (0.25,-0.75);
		\draw[->] (1.5,0) -- (2.5,0);
		\draw[fill=white] (0.25,-0.25) circle [radius=0.1];
		\node at (0,1){$F$};
		\node at (0.25,-1){$\gamma^\cald(e)$};
		\node at (1,-1){$F$};
		\node at (2,0.5){$\xi^{F-1}$};
		\begin{scope}[shift={(3,0)}]
			\draw[->] (0,0) -- (0,-0.25) arc (180:360:0.25) arc (180:0:0.25) -- (1,-0.75);
			\draw (0,0) -- (0,0.5) arc (180:0:0.25) arc (180:360:0.25) -- (1,0.75);
			\draw[->] (0.25,-0.25) -- (0.25,-0.75);
			\draw[->] (1.5,0) -- (2.5,0);
			\draw[fill=white] (0.25,-0.25) circle [radius=0.1];
			\node at (1,1){$F$};
			\node at (0.25,-1){$\gamma^\cald(e)$};
			\node at (1,-1){$F$};
			\node at (2,0.5){$\mathrm{ev}_{\mathrm{ev}^F}$};
		\end{scope}
		\begin{scope}[shift={(6,0)}]
			\draw (0.25,-0.25) circle [radius=0.25];
			\draw[->] (1,0.75) -- (1,-0.75);
			\draw[->] (0.25,-0.25) -- (0.25,-0.75);
			\draw[->] (1.5,0) -- (2.5,0);
			\draw[fill=white] (0.25,-0.25) circle [radius=0.1];
			\node at (1,1){$F$};
			\node at (0.25,-1){$\gamma^\cald(e)$};
			\node at (1,-1){$F$};
			\node at (2,0.5){$M^F$};
		\end{scope}
		\begin{scope}[shift={(9,0)}]
			\draw[->] (1,0.75) -- (1,-0.75);
			\draw[->] (0.25,-0.25) -- (0.25,-0.75);
			\draw[fill=white] (0.25,-0.25) circle [radius=0.1];
			\node at (1,1){$F$};
			\node at (0.25,-1){$\gamma^\cald(e)$};
			\node at (1,-1){$F$};
		\end{scope}
	\end{tikzpicture}
	\caption{Construction of $\tilde{M}^F$ from $M^F$}
	\label{graphical_m_to_tildem}
\end{figure}

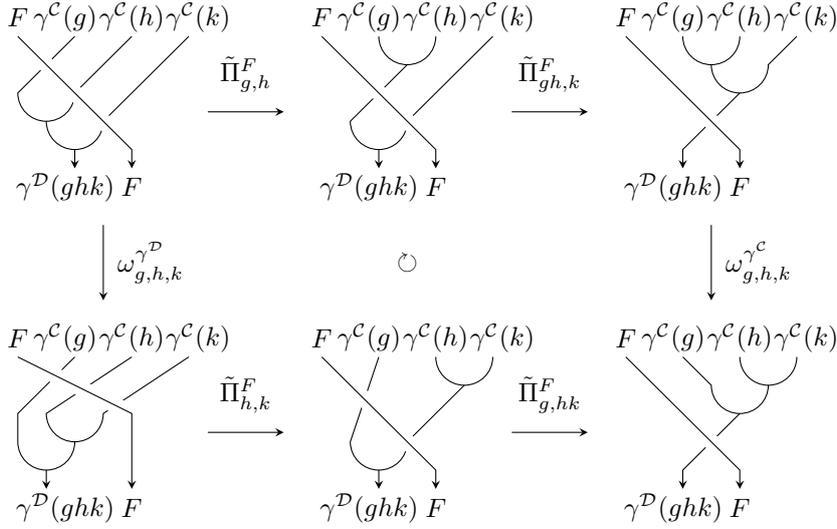
\begin{figure}[htb]
	\centering
	\begin{tikzpicture}
		\draw (0.75,0) -- (0,-0.75) arc (180:360:0.375) -- (1.5,0);
		\draw (0.375,-1.125) arc (180:360:0.375) -- (2.25,0);
		\draw[cross,->] (0,0) -- (1.5,-1.5) -- (1.5,-1.75);
		\draw[->] (0.75,-1.5) -- (0.75,-1.75);
		\draw[->] (2.5,-1) -- (3.5,-1);
		\draw[->] (1.125,-2.5) -- (1.125,-3.5);
		\node at (0,0.25){$F$};
		\node at (0.625,0.25){$\gamma^\calc(g)$};
		\node at (1.5,0.25){$\gamma^\calc(h)$};
		\node at (2.375,0.25){$\gamma^\calc(k)$};
		\node at (0.625,-2){$\gamma^\cald(g h k)$};
		\node at (1.5,-2){$F$};
		\node at (3,-0.5){$\tilde{\Pi}^F_{g, h}$};
		\node at (1.75,-3){$\omega^{\gamma^\cald}_{g, h, k}$};			
		\begin{scope}[shift={(4,0)}]
			\draw (0.75,0) arc (180:360:0.375);
			\draw (1.125,-0.375) -- (0.375,-1.125) arc (180:360:0.375) -- (2.25,0);
			\draw[cross,->] (0,0) -- (1.5,-1.5) -- (1.5,-1.75);
			\draw[->] (0.75,-1.5) -- (0.75,-1.75);
			\draw[->] (2.5,-1) -- (3.5,-1);
			\node at (0,0.25){$F$};
			\node at (0.625,0.25){$\gamma^\calc(g)$};
			\node at (1.5,0.25){$\gamma^\calc(h)$};
			\node at (2.375,0.25){$\gamma^\calc(k)$};
			\node at (0.625,-2){$\gamma^\cald(g h k)$};
			\node at (1.5,-2){$F$};
			\node at (3,-0.5){$\tilde{\Pi}^F_{g h, k}$};
			\node at (1.125,-3){$\circlearrowright$};
		\end{scope}
		\begin{scope}[shift={(8,0)}]
			\draw (0.75,0) arc (180:360:0.375);
			\draw (1.125,-0.375) arc (180:360:0.375) -- (2.25,0);
			\draw[->] (1.5,-0.75) -- (0.75,-1.5) -- (0.75,-1.75);
			\draw[cross,->] (0,0) -- (1.5,-1.5) -- (1.5,-1.75);
			\draw[->] (1.125,-2.5) -- (1.125,-3.5);
			\node at (0,0.25){$F$};
			\node at (0.625,0.25){$\gamma^\calc(g)$};
			\node at (1.5,0.25){$\gamma^\calc(h)$};
			\node at (2.375,0.25){$\gamma^\calc(k)$};
			\node at (0.625,-2){$\gamma^\cald(g h k)$};
			\node at (1.5,-2){$F$};
			\node at (1.75,-3){$\omega^{\gamma^\calc}_{g, h, k}$};			
		\end{scope}
		\begin{scope}[shift={(0,-4.25)}]
			\draw (0.75,0) -- (0,-0.75) -- (0,-1.125) arc (180:360:0.375);
			\draw (1.5,0) -- (0.375,-0.75) arc (180:360:0.375) -- (2.25,0);
			\draw[cross,->] (0,0) -- (1.5,-0.75) -- (1.5,-1.75);
			\draw[->] (0.375,-1.5) -- (0.375,-1.75);
			\draw[->] (2.5,-1) -- (3.5,-1);
			\node at (0,0.25){$F$};
			\node at (0.625,0.25){$\gamma^\calc(g)$};
			\node at (1.5,0.25){$\gamma^\calc(h)$};
			\node at (2.375,0.25){$\gamma^\calc(k)$};
			\node at (0.625,-2){$\gamma^\cald(g h k)$};
			\node at (1.5,-2){$F$};
			\node at (3,-0.5){$\tilde{\Pi}^F_{h, k}$};
		\end{scope}
		\begin{scope}[shift={(4,-4.25)}]
			\draw (1.5,0) arc (180:360:0.375);
			\draw (0.75,0) -- (0.375,-1.125) arc (180:360:0.375) -- (1.875,-0.375);
			\draw[cross,->] (0,0) -- (1.5,-1.5) -- (1.5,-1.75);
			\draw[->] (0.75,-1.5) -- (0.75,-1.75);
			\draw[->] (2.5,-1) -- (3.5,-1);
			\node at (0,0.25){$F$};
			\node at (0.625,0.25){$\gamma^\calc(g)$};
			\node at (1.5,0.25){$\gamma^\calc(h)$};
			\node at (2.375,0.25){$\gamma^\calc(k)$};
			\node at (0.625,-2){$\gamma^\cald(g h k)$};
			\node at (1.5,-2){$F$};
			\node at (3,-0.5){$\tilde{\Pi}^F_{g, h k}$};			
		\end{scope}
		\begin{scope}[shift={(8,-4.25)}]
			\draw (1.5,0) arc (180:360:0.375);
			\draw (0.75,0) -- (1.125,-0.375) arc (180:360:0.375);
			\draw[->] (1.5,-0.75) -- (0.75,-1.5) -- (0.75,-1.75);
			\draw[cross,->] (0,0) -- (1.5,-1.5) -- (1.5,-1.75);
			\node at (0,0.25){$F$};
			\node at (0.625,0.25){$\gamma^\calc(g)$};
			\node at (1.5,0.25){$\gamma^\calc(h)$};
			\node at (2.375,0.25){$\gamma^\calc(k)$};
			\node at (0.625,-2){$\gamma^\cald(g h k)$};
			\node at (1.5,-2){$F$};
		\end{scope}
	\end{tikzpicture}
	\caption{The coherence for $\omega$'s}
	\label{graphical_coherence_omega}
\end{figure}

\begin{figure}[htb]
	\begin{tabular}{cc}
		\begin{minipage}[b]{0.45\hsize}
			\centering
			\begin{tikzpicture}
				\draw (0.5,0) -- (0,-0.5) -- (0,-1) arc (180:360:0.25) -- (1,-0.5) -- (1,0);
				\draw[->] (0.25,-1.25) -- (0.25,-1.5);
				\draw[cross,->] (0,0) -- (1,-1) -- (1,-1.5); 
				\draw[fill=white] (0.5,0) circle [radius=0.1];
				\draw[->] (1.5,-0.75) -- (2.5,-0.75);
				\draw[->] (0.5,-2.25) -- (0.5,-3.25);
				\node at (0,0.25){$F$};
				\node at (1,0.25){$\gamma^\calc(g)$};
				\node at (0.25,-1.75){$\gamma^\cald(g)$};
				\node at (1,-1.75){$F$};
				\node at (2,-0.25){$\tilde{\Pi}^F_{e, g}$};
				\node at (1,-2.75){$\tilde{M}^F$};
				\node at (2,-2.75){$\circlearrowright$};	
				\begin{scope}[shift={(3,0)}]
					\draw (0.5,0) arc (180:360:0.25);
					\draw[->] (0.75,-0.25) -- (0.25,-0.75) -- (0.25,-1.5);
					\draw[cross,->] (0,0) -- (1,-1) -- (1,-1.5); 
					\draw[fill=white] (0.5,0) circle [radius=0.1];
					\draw[->] (0.5,-2.25) -- (0.5,-3.25);
					\node at (0,0.25){$F$};
					\node at (1,0.25){$\gamma^\calc(g)$};
					\node at (0.25,-1.75){$\gamma^\cald(g)$};
					\node at (1,-1.75){$F$};
					\node at (1,-2.75){$\kappa^{\gamma^\calc}_g$};	
				\end{scope}
				\begin{scope}[shift={(0,-4)}]
					\draw (0,-1) arc (180:360:0.25) -- (1,-0.5) -- (1,0);
					\draw[->] (0.25,-1.25) -- (0.25,-1.5);
					\draw[cross,->] (0,0) -- (1,-1) -- (1,-1.5); 
					\draw[fill=white] (0,-1) circle [radius=0.1];
					\draw[->] (1.5,-0.75) -- (2.5,-0.75);
					\node at (0,0.25){$F$};
					\node at (1,0.25){$\gamma^\calc(g)$};
					\node at (0.25,-1.75){$\gamma^\cald(g)$};
					\node at (1,-1.75){$F$};
					\node at (2,-0.25){$\kappa^{\gamma^\cald}_{g}$};
				\end{scope}
				\begin{scope}[shift={(3,-4)}]
					\draw[->] (1,0) -- (0,-1);
					\draw[cross,->] (0,0) -- (1,-1);
					\node at (0,0.25){$F$};
					\node at (1,0.25){$\gamma^\calc(g)$};
					\node at (0,-1.25){$\gamma^\cald(g)$};
					\node at (1,-1.25){$F$};
				\end{scope}
			\end{tikzpicture}
			\caption{The coherence for $\kappa$'s}
			\label{graphical_coherence_kappa}
		\end{minipage}
		\begin{minipage}[b]{0.45\hsize}
			\centering
			\begin{tikzpicture}
				\draw (0.5,0) -- (0,-0.5) -- (0,-1) arc (180:360:0.25) -- (1,-0.5) -- (1,0);
				\draw[->] (0.25,-1.25) -- (0.25,-1.5);
				\draw[cross,->] (0,0) -- (1,-1) -- (1,-1.5); 
				\draw[fill=white] (1,0) circle [radius=0.1];
				\draw[->] (1.5,-0.75) -- (2.5,-0.75);
				\draw[->] (0.5,-2.25) -- (0.5,-3.25);
				\node at (-0.125,0.25){$F$};
				\node at (0.5,0.25){$\gamma^\calc(g)$};
				\node at (0.25,-1.75){$\gamma^\cald(g)$};
				\node at (1,-1.75){$F$};
				\node at (2,-0.25){$\tilde{\Pi}^F_{g, e}$};
				\node at (1,-2.75){$\tilde{M}^F$};
				\node at (2,-2.75){$\circlearrowright$};	
				\begin{scope}[shift={(3,0)}]
					\draw (0.5,0) arc (180:360:0.25);
					\draw[->] (0.75,-0.25) -- (0.25,-0.75) -- (0.25,-1.5);
					\draw[cross,->] (0,0) -- (1,-1) -- (1,-1.5); 
					\draw[fill=white] (1,0) circle [radius=0.1];
					\draw[->] (0.5,-2.25) -- (0.5,-3.25);
					\node at (-0.125,0.25){$F$};
					\node at (0.5,0.25){$\gamma^\calc(g)$};
					\node at (0.25,-1.75){$\gamma^\cald(g)$};
					\node at (1,-1.75){$F$};
					\node at (1,-2.75){$\zeta^{\gamma^\calc}_g$};	
				\end{scope}
				\begin{scope}[shift={(0,-4)}]
					\draw (0.5,0) -- (0,-0.5) -- (0,-1) arc (180:360:0.25);
					\draw[->] (0.25,-1.25) -- (0.25,-1.5);
					\draw[cross,->] (0,0) -- (1,-1) -- (1,-1.5); 
					\draw[fill=white] (0.5,-1) circle [radius=0.1];
					\draw[->] (1.5,-0.75) -- (2.5,-0.75);
					\node at (-0.125,0.25){$F$};
					\node at (0.5,0.25){$\gamma^\calc(g)$};
					\node at (0.25,-1.75){$\gamma^\cald(g)$};
					\node at (1,-1.75){$F$};
					\node at (2,-0.25){$\zeta^{\gamma^\cald}_g$};
				\end{scope}
				\begin{scope}[shift={(3,-4)}]
					\draw[->] (1,0) -- (0,-1);
					\draw[cross,->] (0,0) -- (1,-1);
					\node at (0,0.25){$F$};
					\node at (1,0.25){$\gamma^\calc(g)$};
					\node at (0,-1.25){$\gamma^\cald(g)$};
					\node at (1,-1.25){$F$};
				\end{scope}
			\end{tikzpicture}
			\caption{The coherence for $\zeta$'s}
			\label{graphical_coherence_zeta}
		\end{minipage}
	\end{tabular}
\end{figure}

\begin{figure}[htb]
	\centering
	\begin{tikzpicture}
		\draw (0.5,0) -- (0,-0.5) -- (0.5,-1) arc (180:360:0.25) -- (2,0);
		\draw[cross,->] (0,0) -- (0.5,-0.5) arc (180:360:0.25) -- (1,0);
		\draw[cross,->] (1.5,0) -- (2,-0.5) arc (180:360:0.25) -- (2.5,0);
		\draw[->] (0.75,-1.25) -- (0.75,-1.5);
		\draw[->] (2.75,-0.75) -- (3.25,-0.75);
		\node at (0.5,0.25){$F \gamma^\calc(g) F^{-1}$};
		\node at (2.25,0.25){$F \gamma^\calc(h) F^{-1}$};
		\node at (0.75,-1.75){$\gamma^\cald(gh)$};
		\node at (3,-1.25){$\xi^{F-1}$};
		\begin{scope}[shift={(3.5,0)}]
			\draw (0.5,1) -- (0,0.5) -- (0,-0.5) -- (0.5,-1) arc (180:360:0.25) -- (2,0) -- (2,1);
			\draw[cross,->] (0,1) -- (0.5,0.5) arc (180:360:0.25) -- (1,1);
			\draw[cross,->] (0.5,-0.25) -- (1.5,-1.25) arc (180:360:0.25) -- (2.5,-0.75) -- (2.5,1);
			\draw[->] (0.75,-1.25) -- (0.75,-1.5);
			\draw[->] (2.75,-0.75) -- (3.75,-0.75);
			\draw (0.5,-0.25) arc (180:0:0.25) arc (180:360:0.25) -- (1.5,1);
			\node at (0.5,1.25){$F \gamma^\calc(g) F^{-1}$};
			\node at (2.25,1.25){$F \gamma^\calc(h) F^{-1}$};
			\node at (0.75,-1.75){$\gamma^\cald(gh)$};
			\node at (3.25,-1.25){$\mathrm{coev}^{-1}_{\mathrm{coev}^F}$};
		\end{scope}
		\begin{scope}[shift={(7.5,0)}]
			\draw (0.5,0) -- (0,-0.5) -- (0.5,-1) arc (180:360:0.25) -- (2,0);
			\draw[cross,->] (0,0) -- (1,-1) arc (180:360:0.25) -- (2.5,0);
			\draw[->] (1.5,0) arc (360:180:0.25);
			\draw[->] (0.75,-1.25) -- (0.75,-1.5);
			\draw[->] (2.5,-0.75) -- (3.5,-0.75);
			\node at (0.5,0.25){$F \gamma^\calc(g) F^{-1}$};
			\node at (2.25,0.25){$F \gamma^\calc(h) F^{-1}$};
			\node at (0.75,-1.75){$\gamma^\cald(gh)$};
			\node at (3,-1.25){$\tilde{\Pi}^F_{g, h}$};
		\end{scope}
		\begin{scope}[shift={(11.25,0)}]
			\draw (0.5,0) -- (0,-0.5) -- (0.5,-1) arc (180:360:0.25) -- (2,0);
			\draw[cross,->] (0,0) -- (1,-1) arc (180:360:0.25) -- (2.5,0);
			\draw[->] (1.5,0) arc (360:180:0.25);
			\draw[->] (0.75,-1.25) -- (0.75,-1.5);
			\node at (0.5,0.25){$F \gamma^\calc(g) F^{-1}$};
			\node at (2.25,0.25){$F \gamma^\calc(h) F^{-1}$};
			\node at (0.75,-1.75){$\gamma^\cald(gh)$};
		\end{scope}
	\end{tikzpicture}
	\caption{Construction of $\Pi^F$ from $\tilde{\Pi}^F$}
	\label{graphical_tildepi_to_pi}
\end{figure}
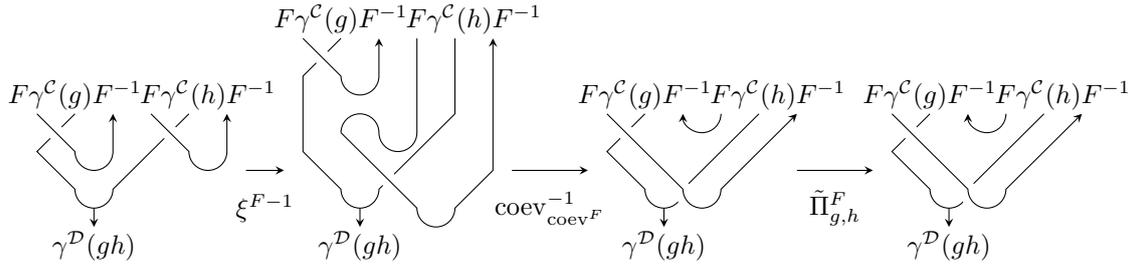

\begin{figure}[htb]
	\centering
	\begin{tikzpicture}
		\draw[->] (0.5,0) -- (0,-0.5) -- (0,-0.75);
		\draw[fill=white] (0.5,0) circle [radius=0.1];
		\draw[cross] (0,0) -- (0.5,-0.5) arc (180:360:0.25) -- (1,0) arc (0:180:0.5);
		\draw[->] (1.5,-0.25) -- (2.5,-0.25);
		\node at (0,-1){$\gamma^\cald(e)$};
		\node at (-0.25,0){$F$};
		\node at (2,0.25){$\tilde{M}^F$};
		\begin{scope}[shift={(3,0)}]
			\draw[->] (0,-0.25) -- (0,-0.75);
			\draw[fill=white] (0,-0.25) circle [radius=0.1];
			\draw[cross] (0.5,0) -- (0.5,-0.5) arc (180:360:0.25) -- (1,0) arc (0:180:0.25);
			\draw[->] (1.5,-0.25) -- (2.5,-0.25);
			\node at (0,-1){$\gamma^\cald(e)$};
			\node at (0.25,0){$F$};
			\node at (2,0.25){$\mathrm{ev}_{\mathrm{coev}^F}$};
		\end{scope}
		\begin{scope}[shift={(6,0)}]
			\draw[->] (0,-0.25) -- (0,-0.75);
			\draw[fill=white] (0,-0.25) circle [radius=0.1];
			\node at (0,-1){$\gamma^\cald(e)$};
		\end{scope}
	\end{tikzpicture}
	\caption{Construction of $M^F$ from $\tilde{M}^F$}
	\label{graphical_tildem_to_m}
\end{figure}
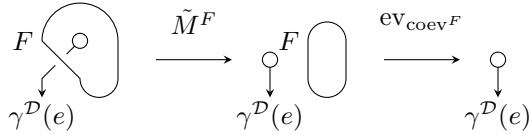

\begin{lem}
	For a biequivalence $F$, we have $(\xi^F \ast \id_{F^{-1}}) \otimes \id_{\mathrm{coev}^F} = (\id_F \ast \tilde{\xi}^F) \otimes \id_{\mathrm{coev}^F}$. 
	
	\begin{proof}
		We already proved $\id_{\mathrm{ev}^F} \otimes (\id_{F^{-1}} \ast \xi^F) = \id_{\mathrm{ev}^F} \otimes (\tilde{\xi}^F \ast \id_F)$ in the proof of Lemma \ref{adismonoidallem}. Then, from $(\id_{\mathrm{ev}^F} \otimes (\id_{F^{-1}} \ast \xi^F))(\id_{\mathrm{ev}^F \ast \mathrm{ev}^F} \otimes (\tilde{\xi}^F \ast \id_{\mathrm{coev}^F} \ast \id_F)) = (\id_{\mathrm{ev}^F} \otimes (\tilde{\xi}^F \ast \id_F)) (\id_{\mathrm{ev}^F \ast \mathrm{ev}^F} \otimes (\id_{F^{-1}} \ast \id_{\mathrm{coev}^F} \xi^F))$ we obtain $\id_{\mathrm{ev}^F \ast \mathrm{ev}^F} \otimes (\tilde{\xi}^F \ast \id_{\mathrm{coev}^F} \ast \id_F) = \id_{\mathrm{ev}^F \ast \mathrm{ev}^F} \otimes (\id_{F^{-1}} \ast \id_{\mathrm{coev}^F} \xi^F)$ and therefore $\id_{\mathrm{ev}^F \ast \mathrm{ev}^F} \otimes (\id_{F^{-1}} \ast \xi^F \ast \id_{F^{-1}F}) \otimes (\id_{F^-1} \ast \id_{\mathrm{coev}^F} \ast \id_{F}) = \id_{\mathrm{ev}^F \ast \mathrm{ev}^F} \otimes (\id_{F^{-1}F} \ast \tilde{\xi}^F \ast \id_F) \otimes (\id_{F^-1} \ast \id_{\mathrm{coev}^F} \ast \id_{F})$. Thus, we obtain the statement since we may compose $\xi^{F-1} \ast \tilde{\xi}^{F-1}$ with the modifications in the statement.
	\end{proof}
\end{lem}

Thus, we can equivalently define a $G$-equivalence to be a tuple $(F, \tilde{\eta}^F, \tilde{\Pi}^F, \tilde{M}^F)$ with the conditions in Figures \ref{graphical_coherence_omega}, \ref{graphical_coherence_kappa} and \ref{graphical_coherence_zeta}, which reduces to \cite[Definition 2.3]{MR3936135} when the action is unital.

Then, we state the coherence theorem for group actions on bicategories \cite[Theorem 3.1]{MR3936135} in the form of ``all diagrams commute'' for our convenience.

\begin{thm}
	\label{theorem_equivcoh}
	Let $\calc$ be a bicategory with an action $\gamma^\calc$ of a group $G$. Define a set $W = \bigsqcup_{n \ge 0} W_n$ of words recursively by the following rules: $\mathbf{1}_C^\calc, (\chi^{\gamma^\calc}_{g,h})^0_C, (\iota^{\gamma^\calc})_C^0 \in W_0$ for any 0-cell $C$ and $g,h \in G$, $- \in W_1$, $\otimes \in W_{2}$, $\gamma^\calc(g) \in W_1$ for $g \in G$ and $w'((w_i)_i) \in W_{\sum_i n_i}$ for $w' \in W_n$ and a family $(w_i)_{i=1}^n$ with $w_i \in W_{n_i}$. Let the functor $\calc^n \to \calc$, where $\calc^0$ denotes the category with only one object and its identity morphism, corresponding naturally to a word $w \in W_n$ be denoted again by $w$. Define a set $I$ of morphisms recursively by the following rules: the components of $a^\calc,l^\calc,r^\calc, J^{\gamma^\calc(g)}$ and $\varphi^{\gamma^\calc(g)}$ are in $I$, $(\chi^{\gamma^\calc}_{g,h})_\lambda$ and $(\iota^{\gamma^\calc})_\lambda$ for any 1-cell $\lambda$ and $g,h \in G$ are in $I$, $\omega^{\gamma^\calc}_{g,h,k}, \kappa^{\gamma^\calc}_g, \zeta^{\gamma^\calc}_g \in I$ for $g,h,k \in G$, $f^{-1} \in I$ for $f \in I$, $f \circ f' \in I$ for $f,f' \in I$ if it is well-defined and $w((f_i)_i) \in I$ for $w \in W_n$ and $(f_i)_{i=1}^n \subset I$. Then, for any $w,w' \in W_n$ and 1-cells $(\lambda_i)_{i=1}^n$, a 2-cell $w((\lambda_i)_i) \to w'((\lambda_i)_i)$ in $I$ is unique if it exists.

	\begin{proof}
		By the proof of \cite[Theorem 3.1]{MR3936135}, there exists a $G$-biequivalence $F$ from $\calc$ to a 2-category $\cald$ with a strict action $\gamma^\cald$ (see \cite[Definition 2.2]{MR3936135}) such that every equivalence 1-cell is indeed an isomorphism. We regard $F$ as data $(F, \tilde{\eta}^F, \tilde{\Pi}^F, \tilde{M}^F)$ as above. We define 2-cells $C^v_{(\lambda_i)_i}$ for $v \in W_n$ and 1-cells $(\lambda_i)_{i=1}^n$ recursively by the following rules: $C^{\mathbf{1}_C^\calc} \coloneqq (\varphi^F_C)^{-1}: F(\mathbf{1}^\calc_C) \cong \mathbf{1}_{F(C)}^\cald$ for $C \in \obj(\calc)$, $C^{\otimes}_{\lambda,\mu} \coloneqq (J^F_{\lambda,\mu})^{-1} : F(\lambda \mu) \cong F(\lambda) F(\mu)$ for 1-cells $\lambda$ and $\mu$ of $\calc$, $C^{\gamma^\calc(g)}_\lambda \coloneqq \id \otimes (\tilde{\eta}^F_g)_\lambda : F\gamma^\calc(g)(\lambda) \cong ((\tilde{\eta}^F_g)_{C'}^0)^{-1} \gamma^\cald(g) F (\lambda) (\tilde{\eta}^F_g)_{C}^0$ for $\lambda \in \hom_\calc(C, C')$ and $g \in G$, $C^{(\chi_{g,h}^{\gamma^\calc})_C^0} \coloneqq \id \otimes (\tilde{\Pi}_{g,h}^F)_C^{-1}: F((\chi_{g,h}^{\gamma^\calc})_C^0) \cong ((\tilde{\eta}^{F}_{gh})^0_C)^{-1} \gamma^\cald(g)((\tilde{\eta}^F_h)^0_C) (\tilde{\eta}^F_g)^0_{\gamma^\calc(h)(C)}$ for $C \in \obj(\calc)$ and $g,h \in G$, $C^{(\iota^{\gamma^\calc})^0_C} \coloneqq \id \otimes \tilde{M}^F_C : F((\iota^{\gamma^\calc})_C^0) \cong ((\tilde{\eta}_e^F)_C^0)^{-1}$ for $C \in \obj(\calc)$ and $C^{v((w_i)_i)}_{(\lambda_i)_i} \coloneqq v((C^{w_i}_{\lambda_i})_i) \circ C^v_{(w_i(\lambda_i))_i}$ for $v \in W_n$, $(w_i)_{i=1}^n \subset W$ and 1-cells $(\lambda_i)_{i=1}^n$.

		We show $C^w_{(\lambda_i)_i} = C^{w'}_{(\lambda_i)_i} \circ F(f)$ for any 2-cell $f:w((\lambda_i)_i) \to w'((\lambda_i)_i)$, which uniquely determines $f$. It is enough to prove this when $f$ is a generator of $I$. The statement for $a^\calc,l^\calc,r^\calc,J^{\gamma^\calc(g)}$ and $\varphi^{\gamma^\calc(g)}$ follows since $F$ is a pseudofunctor and $\tilde{\eta}^F_g$ is a pseudonatural transformation for any $g \in G$, which is standard. The statement for $(\chi^{\gamma^\calc}_{g,h})_{\lambda}$ and $(\iota^{\gamma^\calc})_\lambda$ follows since $\tilde{\Pi}^F_{g,h}$ and $\tilde{M}^F$ are modifications. The statement for $\omega^{\gamma^\calc}_{g,h,k}, \kappa^{\gamma^\calc}_g$ and $\zeta^{\gamma^\calc}_g$ follows from Figures \ref{graphical_coherence_omega}, \ref{graphical_coherence_kappa} and \ref{graphical_coherence_zeta}.
	\end{proof}
\end{thm}

Thanks to this theorem, we may hereafter suppress the 2-cells in $I$ in this article.

Finally, we return to the case of monoidal categories. 

\begin{defi}
	Let $(\calc, \gamma^\calc)$ and $(\cald, \gamma^\cald)$ be monoidal categories with actions of a group $G$. A \emph{monoidal $G$-equivalence} is a $G$-biequivalence with monoidal natural isomorphisms $\eta^F_g : \ad(F) \circ \gamma^\calc(g) \cong \gamma^\cald (g)$ (see \cite[Definition 2.4.8]{egno}) and identical $\Pi^F$ and $M^F$.
\end{defi}

As we have already seen, we can give an equivalent definition of a monoidal $G$-equivalence with $\tilde{\eta}^F$, where $\tilde{\Pi}^F$ and $\tilde{M}^F$ are identical. Namely, we can define a $G$-biequivalence to be the pair $(F, \tilde{\eta}^F)$ with $\tilde{\eta}^F_{gh} = (\id \ast \tilde{\eta}^F_h) \circ (\tilde{\eta}^F_g \ast \id)$ for $g,h \in G$ and $\tilde{\eta}^F_e = \id$, where we suppressed some natural isomorphisms by Theorem \ref{theorem_equivcoh}. Since these equations do not include an adjoint inverse $F^{-1}$, we can give the following definition, which recovers that in \cite[Section 3.1]{MR3671186} for unital actions. 

\begin{comment}
\begin{figure}[htb]
	\begin{tabular}{cc}
		\begin{minipage}[c]{0.45\hsize}
			\centering
			\begin{tikzcd}
				F \gamma^\calc(g) \gamma^\calc(h) \dar{\id \ast \chi^{\gamma^\calc}_{g,h}} \rar{\tilde{\eta}^F_g \ast \id} \drar[phantom]{\circlearrowright}
					& \gamma^\cald (g) F \gamma^\calc (h) \dar{\id \ast \tilde{\eta}^F_h} \\
				F \gamma^\calc(gh) \drar{\tilde{\eta}^F_{gh}}
					&\gamma^\cald(g) \gamma^\cald(h) F \dar{\chi^{\gamma^\cald}_{g,h} \ast \id} \\
					&\gamma^\cald(g h) F
			\end{tikzcd}
			\caption{The coherence for $\chi$'s}
			\label{diagram_cohcond_chi}
		\end{minipage}
		\begin{minipage}[c]{0.45\hsize}
			\centering
			\begin{tikzcd}
				F \drar[swap]{\iota^{\gamma^\cald} \ast \id} \rar{\id \ast \iota^{\gamma^\calc}} \drar[phantom,bend left=20]{\circlearrowright}
					& F \gamma^\calc(e) \dar{\tilde{\eta}^F_e} \\
				& \gamma^\cald(e) F
			\end{tikzcd}
			\caption{The coherence for $\iota$'s}
			\label{diagram_cohcond_iota}
		\end{minipage}
	\end{tabular}
\end{figure}
\end{comment}

\begin{defi}
	\label{definition_monoidal_gfunctor}
	Let $(\calc, \gamma^\calc)$ and $(\cald, \gamma^\cald)$ be monoidal categories with actions of a group $G$. A \emph{monoidal $G$-functor} is the pair $(F, \tilde{\eta}^F)$ of a monoidal functor $F$ and a family $\tilde{\eta}^F$ of monoidal isomorphisms $\tilde{\eta}_g^F: F \gamma^\calc(g) \cong \gamma^\cald(g) F$ for $g \in G$ with $\tilde{\eta}^F_{gh} = (\id \ast \tilde{\eta}^F_h) \circ (\tilde{\eta}^F_g \ast \id)$ for $g,h \in G$ and $\tilde{\eta}^F_e = \id$.
\end{defi}

By definition, a monoidal $G$-functor is a monoidal $G$-equivalence if and only if it is an equivalence as a functor. One can also define the notion of a \emph{$G$-pseudofunctor}, recovering \cite[Definition 2.3]{MR3936135}, with Figures \ref{graphical_coherence_omega}, \ref{graphical_coherence_kappa} and \ref{graphical_coherence_zeta} by replacing $F$ by a general pseudofunctor, which is not necessary for this article.

\subsection{Equivariantly braided tensor categories}
\label{subsection_gbraided_cat}

In this subsection, we recall the notions of a $G$-crossed multitensor category and a $G$-braided multitensor category.

Our terminologies on tensor category theory follow those in \cite{egno}. In particular, a \emph{multitensor category} over a field $k$ is a locally finite $k$-linear abelian rigid monoidal category with a bilinear monoidal product as in \cite[Definition 4.1.1]{egno}. A \emph{tensor functor} is an exact faithful $k$-linear monoidal functor as in \cite[Definition 4.2.5]{egno}.

\begin{rem}
Although tensor categories that arise in algebraic quantum field theory do not have zero objects as remarked in \cite[Remark 2.11]{MR2183964}, because we are interested only in their rigid subcategories, which are semisimple, we regard such a subcategory as an abelian category by adding a zero object. Note that a linear functor on such a category is automatically exact.
\end{rem}

When $\calc$ is moreover a multitensor category, an \emph{action} of $G$ is a monoidal functor from $\underline{G}$ to the monoidal category of the tensor autoequivalences on $\calc$. When $\calc$ is moreover pivotal (see e.g. \cite[Definition 4.7.8]{egno}), a \emph{pivotal action} of $G$ is a monoidal functor from $\underline{G}$ to the monoidal category of pivotal tensor autoequivalences on $\calc$. Recall that a monoidal functor $F:\calc \to \cald$ is \emph{pivotal} when $\delta^\cald F = F\delta^\calc$, where $\delta^\calc$ denotes the pivotal structure of $\calc$.

First, we recall the notion of a \emph{grading} on a multitensor category. Let $\{ \calc_i \}_{i \in I}$ be a family of additive categories. Its \emph{direct sum} $\calc = \bigoplus_{i \in I} \calc_i$ is an additive category with a family of additive functors $\{ I_i : \calc_i \to \calc \}_{i \in I}$ which has the universal property that for any additive category $\cald$ and a family of additive functors $\{ F_i: \calc_i \to \cald \}_{i \in I}$ there exists an additive functor $F: \calc \to \cald$ with a family of natural isomorphisms $\sigma^F = \{ \sigma^{F,i}: F I_i \cong F_i \}_{i \in I}$ such that for any other such pair $(\tilde{F},\tilde{\sigma}^F)$ and a family of natural transformations $\{ \tau^i : F_i \to \tilde{F}_i \}_{i \in I}$ we have a unique natural transformation $\tau : F \to \tilde{F}$ with $\tau^i \circ \sigma^{F,i} = \sigma^{\tilde{F},i} \circ \tau I_i$ for any $i \in I$. Such a category $\calc$ always exists: indeed, we can explicitly give $\calc$ as the category of families of objects whose all but finite components are zero, see e.g. \cite[Section 1.3]{egno}. The category $\calc$ is unique up to a unique equivalence by universality. We can indeed take $\sigma^{F, i}$'s to be identities, and in this case, we refer to the natural transformation $\tau$ obtained by universality as the \emph{extension} of the original family of natural transformations $\{ \tau^i \}_{i \in I}$.

We can equivalently define the direct sum $\calc = \bigoplus_{i \in I} \calc_i$ to be an additive category which has $\calc_i$'s as its subcategories such that every object is a direct sum of objects of $\calc_i$'s and $\hom_\calc(\lambda, \mu) = \{ 0 \}$ for $\lambda \in \obj(\calc_i)$ and $\mu \in \obj(\calc_j)$ with $i \neq j$. Indeed, the explicit construction above satisfies this property, and conversely by decomposing the objects of $\calc$ we obtain an additive equivalence from $\calc$ to $\bigoplus_{i \in I} \calc_i$.

The objects in the set $\homog(\calc) \coloneqq \bigcup_{i \in I} \obj(\calc_i)$ are called \emph{homogeneous}. We write $\partial_\calc \lambda = i$ or simply $\partial \lambda = i$ if $\lambda \in \obj(\calc_i)$. Since $\obj(\calc_i) \cap \obj(\calc_j)$ consists of the zero objects of $\calc$ for $i \neq j$, $\partial$ is single-valued on nonzero objects. 

By definition, we can always decompose an object into a direct sum of homogeneous objects. The following lemma shows that such a decomposition, which we call a \emph{homogeneous decomposition}, is essentially unique.

\begin{lem}
\label{homogdecompunique}
Let $\lambda$ be an object of an $I$-graded multitensor category and let $\lambda = \bigoplus_{i \in F} \lambda_i = \bigoplus_{j \in F'} \lambda_j'$ be decompositions in two ways with finite subsets $F$ and $F'$ of $I$. Then, for every $i \in F \cap F'$ we have $\lambda_i \cong \lambda_{i}'$ as subobjects (i.e. there exists an isomorphism $f \in \hom (\lambda_i, \lambda_{i}')$ such that $\iota_i = \iota'_i \circ f$, where $\iota_i$ and $\iota'_i$ denote the morphisms that embeds respectively $\lambda_i$ and $\lambda_i'$ into $\lambda$) by a unique subobject isomorphism. Moreover, $\lambda_i \cong 0$ for every $i \in F \setminus F'$ and $\lambda_j' \cong 0$ for every $j \in F' \setminus F$.

\begin{proof}
Let $p_i$ and $p_i'$ denote the projections of the decompositions $\lambda = \bigoplus_i \lambda_i$ and $\lambda = \bigoplus_j \lambda_j'$ respectively, and let $\iota_i$ and $\iota'_j$ denote the corresponding inclusions. Since $p_i \circ \iota_j' = 0$ and $p_j' \circ \iota_i =0$ for $i \in F$ and $j \in F' \setminus \{i\}$, we have $(p_i \circ \iota_i') \circ (p_i' \circ \iota_i) =\sum_{j \in F'} p_i \circ \iota_j' \circ p_j' \circ \iota_i = \id_{\lambda_i}$ and similarly $ (p_i' \circ \iota_i) \circ (p_i \circ \iota_i') = \id_{\lambda_i'}$ for every $i \in F \cap F'$. Moreover, $\iota_i' \circ p_i' \circ \iota_i = \sum_{j \in I} \iota_j' \circ p_j' \circ \iota_i = \iota_i$, which shows that $p_i' \circ \iota_i \in \hom(\lambda_i, \lambda_i')$ is a subobject isomorphism for every $i \in F \cap F'$. Uniqueness follows from $p_i' \circ \iota_i' = \id_{\lambda_i'}$. When $i \in F \setminus F'$, we have $p_i = \sum_{j \in F'} p_i \circ \iota'_j \circ p'_j = 0$. Then $\id_{\lambda_i} = p_i \circ \iota_i = 0$ and therefore $\lambda_i \cong 0$. A similar argument gives $\lambda_j' \cong 0$ for $j \in F' \setminus F$.
\end{proof}
\end{lem}

The following notions, which were introduced in \cite[Section 2]{turaev2000homotopy} and \cite[Definition 2.1]{kirillov2004g}, are what we are interested in.

\begin{defi}
	\label{definition_crossed}
Let $G$ be a group. A \emph{$G$-grading} on a multitensor category $\calc$ is its decomposition into a direct sum $\calc = \bigoplus_{g \in G} \calc_g$ of a family $\{ \calc_g \}_{g \in G}$ of abelian full subcategories of $\calc$ with $\obj(\calc_g) \otimes \obj(\calc_h) \subset \obj(\calc_{gh})$ for any $g,h \in G$. A multitensor category with a $G$-grading is called a \emph{$G$-graded multitensor category}. A \emph{$G$-graded tensor functor} between $G$-graded multitensor categories is a tensor functor $F: \calc \to \cald$ with $F(\obj(\calc_g)) \subset \obj(\cald_g)$ for every $g \in G$.

A \emph{$G$-crossed multitensor category} is the pair $\calc = (\calc, \gamma^\calc)$ of a pivotal $G$-graded multitensor category $\calc$ and a pivotal action $\gamma^\calc$ of $G$ on $\calc$ such that $\gamma^\calc(g)(\obj(\calc_h)) \subset \obj(\calc_{ghg^{-1}})$ for $g,h \in G$. A \emph{$G$-crossed functor} between $G$-crossed multitensor categories is a tensor $G$-functor (i.e. a monoidal $G$-functor which is a tensor functor as a monoidal functor, see Definition \ref{definition_monoidal_gfunctor}) which is $G$-graded and pivotal. A $G$-crossed functor is a \emph{$G$-crossed equivalence} if it is an equivalence as a functor.

A \emph{$G$-braided multitensor category} is the pair $\calc = (\calc, b^\calc)$ of a $G$-crossed multitensor category $\calc$ and a family $b^\calc$ of isomorphisms $b^\calc_{\lambda, \mu}: \lambda \mu \cong {}^g \mu \lambda $ for $\lambda \in \obj(\calc_g), \mu \in \obj(\calc)$ and $g \in G$ which is natural in $\lambda$ and $\mu$ such that ${}^h b^\calc_{\lambda, \mu} = b_{{}^h \lambda, {}^h \mu}^\calc$, $b^{\calc}_{\lambda_1 \lambda_2, \mu} = (b^\calc_{\lambda_1, {}^{g_2} \mu} \otimes \id_{\lambda_2}) \circ (\id_{\lambda_1} \otimes b^\calc_{\lambda_2,\mu})$ and $b^\calc_{\lambda, \mu_1 \mu_2} = (\id_{{}^g \mu_1} \otimes b^\calc_{\lambda,\mu_2}) \circ (b^\calc_{\lambda, \mu_1} \otimes \id_{\mu_2})$ for $\lambda \in \obj(\calc_g)$, $\lambda_1 \in \obj(\calc_{g_1})$, $\lambda_2 \in \obj(\calc_{g_2})$, $\mu,\mu_1,\mu_2 \in \obj(\calc)$ and $g,h,g_1,g_2 \in G$, where we suppressed some isomorphisms by Theorem \ref{theorem_equivcoh}. We call $b^\calc$ the \emph{$G$-braiding} of $\calc$. A \emph{$G$-braided functor} between $G$-braided multitensor categories $\calc$ and $\cald$ is a $G$-crossed functor $F: \calc \to \cald$ such that $F(b^\calc_{\lambda, \mu}) = ((\tilde{\eta}^F_g)_\mu \otimes \id_{F(\lambda)}) \circ b^\cald_{F(\lambda),F(\mu)}$ for $\lambda \in \obj(\calc_g), \mu \in \obj(\calc)$ and $g \in G$. A $G$-braided functor is a \emph{$G$-braided equivalence} if it is a $G$-crossed equivalence.

For a $G$-braided multitensor category $\calc$, we can define a family of isomorphisms $b^{\calc-}$ by putting $b_{\mu,\lambda}^{\calc-} \coloneqq b^{\calc-1}_{\lambda, \mu}: {}^g \mu \lambda \cong \lambda \mu$ for $\lambda \in \obj(\calc_g), \mu \in \obj(\calc)$ and $g \in G$. We call $b^{\calc-}$ the \emph{reverse} of $b^\calc$.
\end{defi}

Note that we automatically have $\mathbf{1}_\calc \in \obj(\calc_e)$ for a $G$-graded multitensor category $\calc$. We also have $\lambda^\vee, {}^\vee \lambda \in \obj(\calc_{g^{-1}})$ for $\lambda \in \obj(\calc_g)$. The objects in $\calc_e$ are said to be \emph{neutral}.

For a $G$-crossed multitensor category $\calc$, $\lambda \in \homog(\calc) \setminus \{ 0 \}$ and $\mu \in \obj(\calc)$, we let ${}^\lambda \mu$ denote ${}^{\partial \lambda} \mu$. Although $\partial \lambda$ is ill-defined when $\lambda$ is zero, we often write ${}^\lambda \mu$ even when $\lambda$ is possibly zero because we mostly use this notation to consider morphisms that become zero when $\lambda$ is zero. Then, we have ${}^{\lambda_1 \lambda_2}\mu = {}^{\lambda_1}({}^{\lambda_2}\mu)$ for $\lambda_1, \lambda_2 \in \homog(\calc)$ and $\mu \in \obj(\calc)$ since $\partial (\lambda_1 \lambda_2) = \partial \lambda_1 \partial \lambda_2$. We let ${}^{\overline{\lambda}} \mu$ denote ${}^{(\partial \lambda)^{-1}} \mu$ for $\lambda \in \homog(\calc)$ and $\mu \in \obj(\calc)$. These notations also apply to morphisms.

A multitensor category that arises in algebraic quantum field theory is indeed a \emph{${}^\ast$-multitensor category} i.e. a multitensor category over $\mathbb{C}$ with a contravariant strict monoidal antilinear involution endofunctor ${}^\ast$ that is the identity on objects and has positivity, see \cite[Section 2.4]{MR1966524}. Note that we assume rigidity in our definition. For an action $\gamma$ of a group $G$ on a ${}^\ast$-multitensor category, we always assume that $\gamma(g)$ commutes with the ${}^\ast$-involution for any $g \in G$. For a $G$-braided ${}^\ast$-multitensor category, we always assume that the components of the $G$-braiding are unitary.

Here, we give our motivating example of a $G$-braided (${}^\ast$-)tensor category. Let $\cala$ be an irreducible M\"{o}bius covariant net (on $S^1$) with Haag duality on $\mathbb{R}$, see e.g. \cite[Section 2]{MR2123931} and \cite[Section 2.2]{MR1652746} for definitions. It is well-known \cite[Sections II and IV]{MR297259} \cite{MR1016869} \cite{MR1104414} (see also \cite[Section 2]{MR1652746}) that the category $\rep \cala$ of finite index Doplicher--Haag--Roberts (DHR) endomorphisms of $\cala$ is a braided ${}^\ast$-tensor category.

\begin{defi}
	Let $\cala$ be a M\"{o}bius covariant net on $S^1$. The \emph{automorphism group} of $\cala$ is defined to be
	\begin{align*}
		\aut(\cala) \coloneqq \{ g \in \mathcal{U}(H_\cala) \mid g \cala(I) g^{-1} = \cala(I) \ \text{for} \ \forall I \in \mathcal{I} \ \text{and} \ g \Omega_\cala = \Omega_\cala \},
	\end{align*}
	where $\mathcal{U}$ denotes the group of unitary operators, $H_\cala$ denotes the Hilbert space of $\cala$, $\mathcal{I}$ denotes the set of intervals in $S^1$,  and $\Omega_\cala$ denotes the vacuum of $\cala$. We do not consider topologies. A group homomorphism $G \to \aut(\cala)$ is called an \emph{action} of $G$ on $\cala$. An action is \emph{faithful} if it is injective. 
\end{defi}

\begin{rem}
	For any $I \in \mathcal{I}$, the adjoint action $\ad: \aut(\cala) \to \aut(\cala(I))$ is injective by Reeh--Schlieder property.
\end{rem}

The following notions were defined in \cite[Definition 2.8]{MR2183964}.

\begin{defi}
	\label{eg_twisted_rep}
	Let $\cala$ be a M\"{o}bius covariant net on $S^1$ with an action $\beta$ of a group $G$. For $g \in G$, an endomorphism $\lambda$ of $\cala_\infty \coloneqq \bigcup_{I \in \mathcal{I}_{\mathbb{R}}} \cala(I)$, where $\mathcal{I}_{\mathbb{R}}$ denotes the set of bounded intervals in $\mathbb{R}$, is \emph{$g$-localized} in $I \in \mathcal{I}_{\mathbb{R}}$ if $\lambda|_{\cala(I'_{\mathrm{L}})} = \id_{\cala(I'_{\mathrm{L}})}$ and $\lambda|_{\cala(I'_{\mathrm{R}})} = \ad \beta(g)|_{\cala(I'_{\mathrm{R}})}$, where $I'_{\mathrm{L}}$ (resp. $I'_{\mathrm{R}}$) denotes the left (resp. right) connected component of the open complement $I'$ of $I$. A $g$-localized endomorphism $\lambda$ is a \emph{$g$-twisted DHR endomorphism} of $\cala$ if for any $\tilde{I} \in \mathcal{I}_{\mathbb{R}}$, there exists a unitary $u \in \cala_\infty$ such that $\ad u \circ \lambda$ is localized in $\tilde{I}$. The ${}^\ast$-category of rigid $g$-twisted DHR endomorphisms is denoted by $g\text{-}\rep \cala$, and an object of the ${}^\ast$-tensor category $G\text{-}\rep \cala \coloneqq \bigoplus_{g \in G} g\text{-}\rep \cala$ is called a \emph{$G$-twisted DHR endomorphism} of $\cala$.
\end{defi}

\begin{rem}
	\label{remark_cross_not_faithful}
	In \cite[Definition 2.8]{MR2183964}, the category $G\text{-}\rep \cala$ is defined to be the category generated by $g\text{-}\rep \cala$'s in $\en \cala_\infty$, and the faithfulness of an action $\beta$ is assumed to assure that $g\text{-}\rep \cala$'s are mutually disjoint. Here, we have defined $G\text{-}\rep \cala$ as a direct sum category from the beginning, and therefore the faithfulness of an action $\beta$ is not needed.
\end{rem}

An action $\beta$ of $G$ on $\cala$ induces the adjoint action on the ${}^\ast$-tensor category $G\text{-}\rep \cala$: ${}^g \lambda \coloneqq \ad \beta(g) \circ \lambda \circ \ad \beta(g^{-1})$ for $\lambda \in G\text{-}\rep \cala$ and $g \in G$, which makes $G\text{-}\rep \cala$ into a $G$-crossed ${}^\ast$-tensor category. If moreover $\cala$ satisfies Haag duality on $\mathbb{R}$, the category $G\text{-}\rep \cala$ turns into a $G$-braided ${}^\ast$-tensor category by \cite[Theorem 2.21]{MR2183964}, see also \cite[Section 3]{MR4153896}. The argument recovers the classical DHR theory when $G$ is trivial.

We go back to general theory. We collect here some graphical calculi for $G$-braided multitensor categories. Our graphical notation follows that in \cite{MR1729094} i.e. a morphism in a monoidal category is represented as an arrow from top to bottom. Now, let $\calc$ be a $G$-braided multitensor category. We do not draw the isomorphisms in Theorem \ref{theorem_equivcoh}. For $\lambda \in \homog(\calc)$ and $\mu \in \obj(\calc)$, the component $b^\calc_{\lambda, \mu}$ of the $G$-braiding is represented by the crossing in Figure \ref{graphicalbraiding}. The component $b^{\calc-}_{\mu, \lambda}$ of the reverse is represented by the reverse crossing in Figure \ref{graphicalreversebraiding}.

\begin{figure}[htb]
\begin{minipage}{0.45\linewidth}
\centering
\begin{tikzpicture}
\draw[->] (0,0) -- (-1,-1);
\draw[->, cross] (-1,0) -- (0,-1);
\node at (-1, 0.25){$\lambda$};
\node at (0, 0.25){$\mu$};
\node at (-1, -1.25){${}^\lambda \mu$};
\node at (0, -1.25){$\lambda$};
\end{tikzpicture}
\caption{The $G$-braiding $b^\calc_{\lambda, \mu}$}
\label{graphicalbraiding}
\end{minipage}
\begin{minipage}{0.45\linewidth}
\centering
\begin{tikzpicture}
\draw[->] (-1,0) -- (0,-1);
\draw[->,cross] (0,0) -- (-1,-1);
\node at (-1, -1.25){$\lambda$};
\node at (0, -1.25){$\mu$};
\node at (-1, 0.25){${}^\lambda \mu$};
\node at (0, 0.25){$\lambda$};
\end{tikzpicture}
\caption{The reverse $b^{\calc-}_{\mu, \lambda}$}
\label{graphicalreversebraiding}
\end{minipage}
\end{figure}

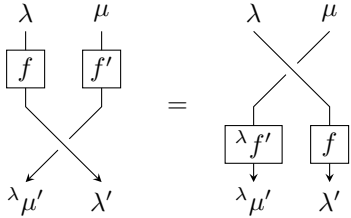
\begin{figure}[htb]
\centering
\begin{tikzpicture}
\draw[->] (0,1) -- (0,0) -- (-1,-1);
\draw[->, cross] (-1,1) -- (-1,0) -- (0,-1);
\draw[fill=white] (-0.25, 0.75) rectangle (0.25, 0.25) node[midway]{$f'$};
\draw[fill=white] (-1.25, 0.75) rectangle (-0.75, 0.25) node[midway]{$f$};
\node at (-1, 1.25){$\lambda$};
\node at (0, 1.25){$\mu$};
\node at (-1, -1.25){${}^{\lambda} \mu'$};
\node at (0, -1.25){$\lambda'$};
\node at (1, 0){$=$};
\begin{scope}[shift={(3,0)}]
\draw[->] (0,1) -- (-1,0) -- (-1,-1);
\draw[->, cross] (-1,1) -- (0,0) -- (0,-1);
\draw[fill=white] (-0.25, -0.25) rectangle (0.25, -0.75) node[midway]{$f$};
\draw[fill=white] (-1.375, -0.25) rectangle (-0.625, -0.75) node[midway]{${}^\lambda f'$};
\node at (-1, 1.25){$\lambda$};
\node at (0, 1.25){$\mu$};
\node at (-1, -1.25){${}^\lambda \mu'$};
\node at (0, -1.25){$\lambda'$};
\end{scope}
\end{tikzpicture}
\caption{The naturality of a $G$-braiding}
\label{bfe}
\end{figure}

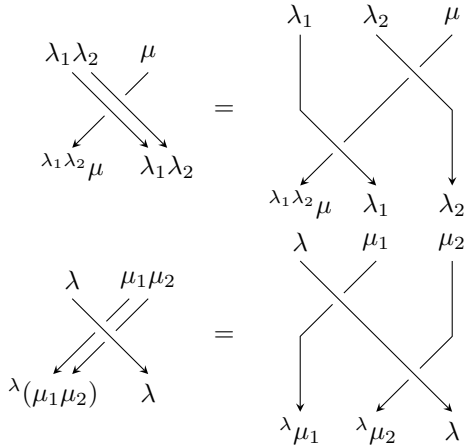
\begin{figure}[htb]
\centering
\begin{tikzpicture}
\draw[->] (0,0) -- (-1,-1);
\draw[->, cross] (-1,0) -- (0,-1);
\draw[->, cross] (-0.75,0) -- (0.25,-1);
\node at (-1, 0.25){$\lambda_1 \lambda_2$};
\node at (0, 0.25){$\mu$};
\node at (-1, -1.25){${}^{\lambda_1 \lambda_2} \mu$};
\node at (0.25, -1.25){$\lambda_1 \lambda_2$};
\node at (1,-0.5){$=$};
\begin{scope}[shift={(4,0.5)}]
\draw[->] (0,0) -- (-2,-2);
\draw[->, cross] (-1,0) -- (0,-1) -- (0,-2);
\draw[->, cross] (-2,0) -- (-2,-1) -- (-1,-2);
\node at (-2, 0.25){$\lambda_1$};
\node at (-1, 0.25){$\lambda_2$};
\node at (0, 0.25){$\mu$};
\node at (-2, -2.25){${}^{\lambda_1 \lambda_2} \mu$};
\node at (0, -2.25){$\lambda_2$};
\node at (-1, -2.25){$\lambda_1$};
\end{scope}
\begin{scope}[shift={(0,-3)}]
\draw[->] (0,0) -- (-1,-1);
\draw[->] (-0.25,0) -- (-1.25,-1);
\draw[->, cross] (-1,0) -- (0,-1);
\node at (-1, 0.25){$\lambda$};
\node at (0, 0.25){$\mu_1 \mu_2$};
\node at (-1.25, -1.25){${}^{\lambda} (\mu_1 \mu_2)$};
\node at (0, -1.25){$\lambda$};
\node at (1,-0.5){$=$};
\begin{scope}[shift={(4,0.5)}]
\draw[->] (0,0) -- (0,-1) -- (-1,-2);
\draw[->] (-1,0) -- (-2,-1) -- (-2,-2);
\draw[->, cross] (-2,0) -- (0,-2);
\node at (-2, 0.25){$\lambda$};
\node at (-1, 0.25){$\mu_1$};
\node at (0, 0.25){$\mu_2$};
\node at (-2, -2.25){${}^{\lambda} \mu_1$};
\node at (0, -2.25){$\lambda$};
\node at (-1, -2.25){${}^\lambda \mu_2$};
\end{scope}
\end{scope}
\end{tikzpicture}
\caption{Axioms for the $G$-braiding $b^\calc$}
\label{graphicalbraidingeq}
\end{figure}

The naturality of $b^\calc$ is represented in Figure \ref{bfe}, where we may assume $\partial \lambda = \partial \lambda'$ since otherwise the morphisms in the equation are zero. We also have a similar representation for $b^{\calc-}$. The first axiom of a $G$-braiding says that ${}^h b^\calc_{\lambda, \mu}$ and $b^\calc_{{}^h \lambda, {}^h \mu}$ are represented by the same diagram. The second and third axioms for $b^\calc$ are represented in Figure \ref{graphicalbraidingeq}.

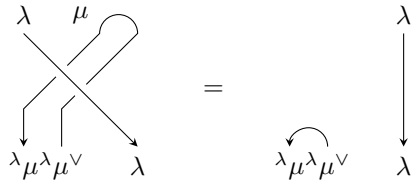
\begin{figure}[H]
\centering
\begin{tikzpicture}
\draw[<-] (0, -2) -- (0,-1.5) -- (1, -0.5) arc (180:0:0.25) -- (0.5, -1.5) -- (0.5,-2);
\draw[->, cross] (0, -0.5) -- (1.5,-2);
\node at (0, -0.25){$\lambda$};
\node at (0, -2.25){${}^\lambda \mu$};
\node at (0.5, -2.25){${}^\lambda \mu^\vee$};
\node at (0.75, -0.25){$\mu$};
\node at (1.5, -2.25){$\lambda$};
\node at (2.5,-1.25){$=$};
\begin{scope}[shift={(3.5,0)}]
\draw[<-] (0, -2) arc (180:0:0.25);
\draw[->, cross] (1.5, -0.5) -- (1.5,-2);
\node at (1.5, -0.25){$\lambda$};
\node at (0, -2.25){${}^\lambda \mu$};
\node at (0.5, -2.25){${}^\lambda \mu^\vee$};
\node at (1.5, -2.25){$\lambda$};
\end{scope}
\end{tikzpicture}
\caption{Moving an arc}
\label{graphicalmovearc}
\end{figure}

Since $\mathrm{coev}_{{}^g \mu}$ and ${}^g \mathrm{coev}_{\mu}$ (resp. $\mathrm{ev}_{{}^g \mu}$ and ${}^g \mathrm{ev}_{\mu}$, and similar for right duals) are represented by the same diagram for $\mu \in \obj(\calc)$ and $g \in G$ (see e.g. \cite[Exercise 2.10.6]{egno}), for example, we can perform the calculation in Figure \ref{graphicalmovearc} by Figures \ref{bfe} and \ref{graphicalbraidingeq}.

The following lemma is the $G$-equivariant version of \cite[Fig. 27]{MR1729094}.

\begin{lem}
	\label{lem_crossing_rotation}
Let $\calc$ be a $G$-braided multitensor category. Let $\lambda \in \homog(\calc)$ and let $\mu \in \obj(\calc)$. Then, the equations in Figure \ref{graphicalrotation} hold.

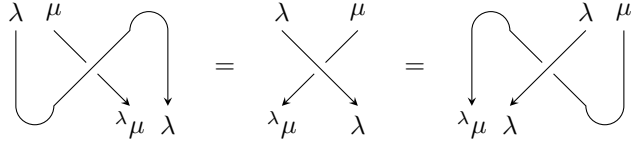
\begin{figure}[H]
\centering
\begin{tikzpicture}
\draw[->] (0,0) -- (-1,-1);
\draw[->, cross] (-1,0) -- (0,-1);
\node at (-1, 0.25){$\lambda$};
\node at (0, 0.25){$\mu$};
\node at (-1, -1.25){${}^\lambda \mu$};
\node at (0, -1.25){$\lambda$};
\node at (-1.75,-0.5) {$=$};
\node at (0.75,-0.5) {$=$};
\begin{scope}[shift={(-3,0)}]
\draw[->] (-1,0) -- (0,-1);
\draw[->, cross] (-1.5,0) -- (-1.5,-1) arc (180:360:0.25) -- (0,0) arc (180:0:0.25) -- (0.5, -1);
\node at (0.5, -1.25){$\lambda$};
\node at (-1, 0.25){$\mu$};
\node at (0, -1.25){${}^\lambda \mu$};
\node at (-1.5, 0.25){$\lambda$};
\end{scope}
\begin{scope}[shift={(3,0)}]
\draw[->] (0.5,0) -- (0.5,-1) arc (360:180:0.25) -- (-1,0) arc (0:180:0.25) -- (-1.5,-1);
\draw[->, cross] (0,0) -- (-1,-1);
\node at (0, 0.25){$\lambda$};
\node at (0.5, 0.25){$\mu$};
\node at (-1.5, -1.25){${}^\lambda \mu$};
\node at (-1, -1.25){$\lambda$};
\end{scope}
\end{tikzpicture}
\caption{Rotations of a crossing}
\label{graphicalrotation}
\end{figure}

\begin{proof}
By the graphical calculation in Figure \ref{graphical_rotation_proof}, the leftmost diagram in Figure \ref{graphicalrotation} is the right inverse of $b^{\calc-}_{\mu, \lambda}$. One can similarly show that it is also the left inverse and therefore the left equality in Figure \ref{graphicalrotation}. The proof of the right equality is similar.
\begin{figure}[H]
\centering
\begin{tikzpicture}
\draw[->] (-1,0) -- (0.5,-1.5);
\draw[->, cross] (-1.5,0) -- (-1.5,-1) arc (180:360:0.25) -- (0,0) arc (180:0:0.25) -- (0.5, -1) -- (0, -1.5);
\node at (0, -1.75){$\lambda$};
\node at (-1, 0.25){$\mu$};
\node at (0.5, -1.75){$\mu$};
\node at (-1.5, 0.25){$\lambda$};
\node at (1.5, -0.75){$=$};
\begin{scope}[shift={(4,0)}]
\draw[->] (0,0) -- (0,-1.5);
\draw[->, cross] (-1.5,0) -- (-1.5,-1) arc (180:360:0.25) arc (180:0:0.25) -- (-0.5, -1.5);
\node at (-0.5, -1.75){$\lambda$};
\node at (0, 0.25){$\mu$};
\node at (0, -1.75){$\mu$};
\node at (-1.5, 0.25){$\lambda$};
\node at (1, -0.75){$=$};
\begin{scope}[shift={(2.5,0)}]
\draw[->] (0,0) -- (0,-1.5);
\draw[->, cross] (-0.5,0) -- (-0.5, -1.5);
\node at (-0.5, -1.75){$\lambda$};
\node at (0, 0.25){$\mu$};
\node at (0, -1.75){$\mu$};
\node at (-0.5, 0.25){$\lambda$};
\end{scope}
\end{scope}
\end{tikzpicture}
\caption{A part of the proof of rotations}
\label{graphical_rotation_proof}
\end{figure}
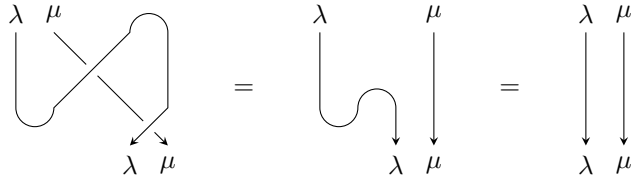
\end{proof}
\end{lem}

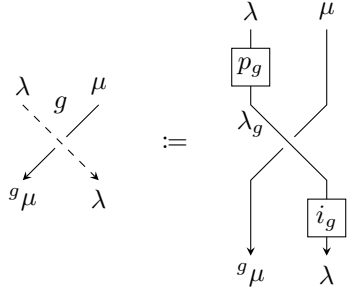
\begin{figure}[H]
	\centering
	\begin{tikzpicture}
	\draw[->] (1,0) -- (0,-1);
	\draw[->,dashed,cross] (0,0) -- (1,-1);
	\node at (0,0.25) {$\lambda$};
	\node at (1,0.25) {$\mu$};
	\node at (0,-1.25) {${}^g \mu$};
	\node at (1,-1.25) {$\lambda$};
	\node at (0.5,0) {$g$};
	\node at (2,-0.5) {$\coloneqq$};
	\begin{scope}[shift={(3,0)}]
	\draw[->] (1,1) -- (1,0) -- (0,-1) -- (0,-2);
	\draw[->,cross] (0,1) -- (0,0) -- (1,-1) -- (1,-2);
	\node at (0,1.25) {$\lambda$};
	\node at (1,1.25) {$\mu$};
	\node at (0,-2.25) {${}^g \mu$};
	\node at (1,-2.25) {$\lambda$};
	\node at (0,-0.25) {$\lambda_g$};
	\draw[fill=white] (-0.25, 0.75) rectangle (0.25,0.25) node[midway]{$p_g$};
	\draw[fill=white] (0.75, -1.25) rectangle (1.25,-1.75) node[midway]{$i_g$};
	\end{scope}
	\end{tikzpicture}
	\caption{A dashed crossing}
	\label{graphicaldashedcrossing}
	\end{figure}

Let $\calc$ be a $G$-braided multitensor category. For any $\lambda, \mu \in \obj(\calc)$ and $g \in G$, $\lambda$ not necessarily homogeneous, let a dashed crossing labeled by $g$ from $\lambda \mu$ to ${}^g \mu \lambda$ denote the morphism in Figure \ref{graphicaldashedcrossing}, where $i_g$ and $p_g$ respectively denote the inclusion and projection of $\lambda_g$ in a homogeneous decomposition $\lambda = \bigoplus_g \lambda_g$. Then, this morphism does not depend on the choice of $i_g$ and $p_g$ by Lemma \ref{homogdecompunique} and by the naturality of $b^\calc$. Note also that a dashed crossing is natural. 
	
Finally, we consider the notion of the $G$-equivariant version of a ribbon structure. Let $\calc$ be a $G$-braided multitensor category. For $\lambda \in \homog(\calc)$, define $\theta_\lambda^\calc : \lambda \cong {}^\lambda \lambda$ by putting $\theta_\lambda^\calc \coloneqq (\id_{{}^\lambda \lambda} \otimes \mathrm{ev}_{\lambda^\vee}) b^\calc_{\lambda^{\vee \vee}, \lambda} (\delta^\calc_\lambda \otimes \mathrm{coev}_\lambda)$. If $\calc$ is a split spherical fusion category, then $\theta^\calc$ satisfies $\theta_{\mathbf{1}^\calc}^\calc = \id_{\mathbf{1}^\calc}$, $\theta_{\lambda \mu} = (\theta_{{}^{\lambda \mu \overline{\lambda} \lambda}} \otimes \theta_{{}^{\lambda} \mu})$, ${}^g \theta^\calc_\lambda = \theta^\calc_{{}^g \lambda}$ and $(\theta_\lambda^\calc)^\vee = {}^\lambda \theta_{\lambda^\vee}^\calc$ for $\lambda, \mu \in \homog(\calc)$ and $g \in G$ by \cite[Lemma 2.3]{kirillov2004g}. 

\begin{defi}
	A $G$-braided multitensor category is a \emph{$G$-ribbon} multitensor category if $\theta^\calc$ satisfies $\theta_{\mathbf{1}^\calc}^\calc = \id_{\mathbf{1}^\calc}$, $\theta^\calc_{\lambda \mu} = (\theta^\calc_{{}^{\lambda \mu \overline{\lambda} \lambda}} \otimes \theta^\calc_{{}^{\lambda} \mu})$, ${}^g \theta^\calc_\lambda = \theta^\calc_{{}^g \lambda}$ and $(\theta_\lambda^\calc)^\vee = {}^\lambda \theta_{\lambda^\vee}^\calc$ for $\lambda, \mu \in \homog(\calc)$ and $g \in G$. 
\end{defi}

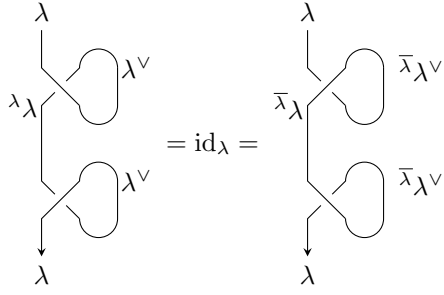
\begin{figure}[htb]
		\centering
		\begin{tikzpicture}
			\draw (1,-0.5) arc (0:180:0.25) -- (0,-1) -- (0,-2) -- (0.5,-2.5);
			\draw[cross] (0,0) -- (0,-0.5) -- (0.5,-1);
			\draw (0.5,-1) arc (180:360:0.25) -- (1,-0.5);
			\draw[cross,->] (1,-2) arc (0:180:0.25) -- (0,-2.5) -- (0,-3);
			\draw (0.5,-2.5) arc (180:360:0.25) -- (1,-2);
			\node at (0,0.25){$\lambda$};
			\node at (1.25,-0.5){$\lambda^\vee$};
			\node at (-0.25,-1){${}^\lambda \lambda$};
			\node at (0,-3.25){$\lambda$};
			\node at (1.25,-2){$\lambda^\vee$};
			\node at (2.25,-1.5){$= \id_\lambda =$};
			\begin{scope}[shift={(3.5,0)}]
				\draw (0,0) -- (0,-0.5) -- (0.5,-1);
				\draw[->] (1,-2) arc (0:180:0.25) -- (0,-2.5) -- (0,-3);
				\draw[cross] (1,-0.5) arc (0:180:0.25) -- (0,-1) -- (0,-2) -- (0.5,-2.5);
				\draw (0.5,-2.5) arc (180:360:0.25) -- (1,-2);
				\draw (0.5,-1) arc (180:360:0.25) -- (1,-0.5);
				\node at (0,0.25){$\lambda$};
				\node at (1.5,-0.5){${}^{\overline{\lambda}} \lambda^\vee$};
				\node at (-0.25,-1){${}^{\overline{\lambda}} \lambda$};
				\node at (0,-3.25){$\lambda$};
				\node at (1.5,-2){${}^{\overline{\lambda}} \lambda^\vee$};
			\end{scope}
		\end{tikzpicture}
		\caption{A Reidemeister move for a framed tangle}
		\label{graphical_reidemeister}
	\end{figure}

For a $G$-ribbon multitensor category $\calc$ and $\lambda \in \homog(\calc)$, we can do the graphical calculations in Figure \ref{graphical_reidemeister}, where we identify the left dual $\lambda^\vee$ and right dual ${}^\vee \lambda$ of $\lambda$ using the pivotal structure $\delta^\calc$. Indeed, the equalities follow from $(\theta_\lambda^\calc)^\vee = {}^\lambda \theta_{\lambda^\vee}^\calc$.

\subsection{Equivariant Frobenius algebras}
\label{subsection_equiv_frob}

In this subsection, we recall the notion of a $G$-equivariant Frobenius algebra.

A \emph{Frobenius algebra} in a monoidal category is an algebra object with a compatible coalgebra structure, see e.g. \cite[Section 2.3]{MR2187404} for the precise definition and more terminologies. The product, unit, coproduct, and counit of a Frobenius algebra $A$ are denoted by $m_A,\eta_A,\Delta_A$, and $\varepsilon_A$ respectively. 

For a pivotal multitensor category $\calc$, we say a Frobenius algebra $A$ in $\calc$ is \emph{symmetric} if $\delta^\calc_A = \id_A$ when $A^{\vee}$ is taken to be $A$, see \cite[Definition 2.22]{MR2187404}. When $A$ is moreover \emph{special} (i.e. $m_A \circ \Delta_A$ and $\varepsilon_A \circ \eta_A$ are invertible scalars, see \cite[Definition 2.22(iii)]{MR2187404}), the product of the invertible scalars $m_A \circ \Delta_A$ and $\varepsilon_A \circ \eta_A$ is equal to $\dim A$. Since any symmetric special Frobenius algebra $A$ has a unique corresponding coalgebra structure $(\tilde{\Delta}_A,\tilde{\varepsilon}_A)$ with $m_A \tilde{\Delta}_A = \id_A$ and $\tilde{\varepsilon}_A \eta_A = \dim A$ (see e.g. \cite[Remark 3.14]{MR1966524}), we assume from the beginning that $m_A \Delta_A = \id_A$ and $\varepsilon_A \eta_A = \dim A$.

A \emph{Q-system} in a ${}^\ast$-multitensor category is a triple $(\theta, w, x)$ such that $(\theta, x^\ast, w,x, (\dim \theta) w^\ast)$ is a special Frobenius algebra, see e.g. \cite[Definition 3.8]{bklr}. A Q-system is automatically symmetric by definition. 

In this article, we have to consider \emph{equivariant} Frobenius algebras. To define this notion, note that if $A$ is an algebra (resp. Frobenius algebra) in a monoidal category with an action of a group $G$, then ${}^g A = ({}^g A, {}^g m_A, {}^g \eta_A)$ (resp. ${}^g A = ({}^g A, {}^g m_A, {}^g \eta_A, {}^g \Delta_A, {}^g \varepsilon_A)$) is again an algebra (resp. Frobenius algebra).

\begin{defi}
Let $\calc$ be a monoidal category with an action of a group $G$. A \emph{$G$-equivariant algebra} in $\calc$ is a pair $(A,z^A)$ of an algebra $A$ in $\calc$ and a $G$-equivariant structure $z^A$ on $A$ (i.e. $z^A_{gh} = z_g^A {}^g z_h^A$ for $g, h \in G$, where we suppress the canonical isomorphism $(\chi_{g,h})_A: {}^{g}({}^h A) \cong {}^{g h} A$ in Theorem \ref{theorem_equivcoh}, see e.g. \cite[Definition 2.7.2]{egno}) such that $z^A_g:{}^g A \cong A$ is an algebra isomorphism for every $g \in G$. Similarly, a \emph{$G$-equivariant Frobenius algebra} in $\calc$ is a pair $A = (A,z^A)$ of a Frobenius algebra $A$ in $\calc$ and a $G$-equivariant structure $z^A$ on $A$ such that $z^A_g$ is a Frobenius algebra isomorphism for every $g \in G$. A \emph{$G$-equivariant Q-system} in a ${}^\ast$-multitensor category is a $G$-equivariant Frobenius algebra $A$ such that $A$ is a Q-system and $z_g^A$ is unitary for every $g \in G$. A \emph{homomorphism} between $G$-equivariant (Frobenius) algebras is a (Frobenius) algebra homomorphism that is $G$-equivariant. By an \emph{isomorphism} between $G$-equivariant Q-systems, we always mean a unitary isomorphism.
\end{defi}

\begin{eg}
It is well-known \cite{MR1257245} (see also \cite[Theorem 3.11]{bklr}) that a finite index extension of a type III factor $N$ gives a Q-system in the ${}^\ast$-tensor category $\en_0(N)$ of finite index unital ${}^\ast$-endomorphisms of $N$ (see e.g. \cite[Chapter 2]{bklr}) up to isomorphisms and vice versa. Indeed, when a group $G$ acts on $N$, we have the adjoint action of $G$ on $\en(N)$, and a finite index extension of $N$ with an extension of the action of $G$ corresponds to a $G$-equivariant Q-system in $\en_0(N)$, see \cite[Section 2.4]{MR4153896}.
\end{eg}

Here, we give examples coming from algebraic quantum field theory.

\begin{defi}
	Let $\mathcal{J}$ be a directed set and let $\cala = \{ \cala(I) \}_{I \in \mathcal{J}}$ be a family of type III von Neumann algebras on a Hilbert space $H_\cala$ with a common cyclic separating vector $\Omega_\cala \in H_\cala$ such that $\cala(I) \subset \cala(J)$ if $I \le J$. A \emph{standard extension} $\cala \subset \calb$ of $\cala$ is a pair $(\calb, E^\calb)$ of a family $\calb = \{ \calb(I) \}_{I \in \mathcal{J}}$ of von Neumann algebras on a Hilbert space $H_\calb$ that extends $H_\cala$ and a family of faithful normal conditional expectations $E^\calb = \{ E^\calb_I : \calb(I) \to \cala(I) \}_{I \in \mathcal{J}}$ such that $\cala(I) \subset \calb(I)$, $\Omega_\cala \in H_\cala \subset H_\calb$ is cyclic and separating for $\calb(I)$ and $\varphi_{\Omega_\cala} \circ E_I^\calb = \varphi_{\Omega_\cala}$ for any $I \in \mathcal{J}$, where $\varphi_{\Omega_\cala}$ denotes the vector state of $\Omega_\cala$, and $\calb(I) \subset \calb(J)$ and $E^\calb_J|_{\calb(I)} = E^\calb_I$ if $I \le J$.
	
	A standard extension $\cala \subset \calb$ is of \emph{finite index} if $E^\calb_I$ is of finite index for any $I \in \mathcal{J}$.

	An \emph{isomorphism} between standard extensions $\cala \subset \calb_1$ and $\cala \subset \calb_2$ of $\cala$ is a unitary $u: H_{\calb_1} \to H_{\calb_2}$ such that $u$ commutes with $\cala_\infty \coloneqq \bigcup_{I \in \mathcal{J}} \cala(I)$, $u \calb_1(I) u^\ast = \calb_2(I)$ for any $I \in \mathcal{J}$ and $u \Omega_{\cala} = \Omega_{\cala}$.
\end{defi}

Let $\cala$ be an irreducible local M\"{o}bius covariant net. Then, the restriction of $\cala$ on $\mathbb{R}$, which is again denoted by $\cala$, is a family of type III factors indexed by the directed set of intervals in $\mathbb{R}$ for which the vacuum $\Omega_\cala$ is a common cyclic separating vector. Moreover, for an inclusion $\cala \subset \calb$ of irreducible local M\"{o}bius covariant nets, there exists a unique family of faithful normal conditional expectations that makes $\cala \subset \calb$ into a standard extension of $\cala$ by Bisognano--Wichmann property and Takesaki's theorem \cite{MR0303307}.

By \cite[Theorem 4.9]{MR1332979}, a finite index standard extension of $\cala$ bijectively corresponds to a Q-system in $\rep \cala$ up to isomorphisms. If moreover $\cala$ satisfies Haag duality on $\mathbb{R}$, then a finite index local standard extension of $\cala$ (i.e. a finite index standard extension $\calb$ of $\cala$ such that $\calb(I)$ and $\calb(J)$ commute if $I \cap J = \emptyset$) bijectively corresponds to a commutative Q-system in $\rep \cala$ up to isomorphisms.

\begin{rem}
	We consider not only local extensions but also \emph{nonlocal} extensions i.e. $\calb$ is not necessarily a local M\"{o}bius covariant net in the example above, which is crucial for the theory of $\alpha$-induction, see \cite[Section 5]{MR1777347}. Note also that we do not assume even the relative locality of $\calb$ with respect to $\cala$ (i.e. $\cala(I)$ and $\calb(J)$ commute if $I \cap J = \emptyset$), which is indeed automatic by \cite[Theorem 4.9]{MR1332979}.
\end{rem}

\begin{defi}
Let $\cala$ be an irreducible local M\"{o}bius covariant net on $S^1$. Suppose we have an action $\beta_\cala: G \to \aut(\cala)$ of a group $G$ on $\cala$. For a standard extension $\cala \subset \calb$, let $\aut(\calb)$ denote the group which consists of the unitaries preserving the local algebras by adjoint action and the vacuum as that for the M\"{o}bius covariant nets. We say an action $\beta_\calb: G \to \aut (\calb)$ of $G$ on $\calb$ \emph{extends} the action $\beta_\cala$ on $\cala$ if $\ad \beta_\calb(g)(a) = \ad \beta_\cala(g)(a)$ for any interval $I$ in $\mathbb{R}$, $a \in \cala(I)$ and $g \in G$, see \cite[Definition 6.2]{MR3933137}. We refer to the pair $\cala \subset \calb = (\cala \subset \calb, \beta_\calb)$ of a standard extension and an extension of $\beta_\cala$ as a \emph{$G$-equivariant standard extension} of $(\cala, \beta_\cala)$. An \emph{isomorphism} between two $G$-equivariant extensions $\cala \subset \calb_1$ and $\cala \subset \calb_2$ is a unitary $u$ between the Hilbert spaces $H_{\calb_1}$ and $H_{\calb_2}$ on which $\calb_1$ and $\calb_2$ act respectively such that $u$ is an isomorphism of standard extensions and $\ad u \circ \beta_{\calb_1} = \beta_{\calb_2}$.
\end{defi}

\begin{rem}
	\begin{enumerate}
		\item When $\calb$ is the restriction of a local M\"{o}bius covariant net with strong additivity, the local algebras on $S^1$ are generated by the local algebras on $\mathbb{R}$ and therefore $\aut(\calb)$ defined above coincides with that of a local M\"{o}bius covariant net.
		\item If $\ad \beta_\calb(g)(a) = \ad \beta_\cala(g)(a)$ for any interval $I$ and $a \in \cala(I)$, then $\ad \beta_\calb(g) (\cala(I)) \subset \cala(I)$ and therefore $\beta_\calb(g)$ commutes with the Jones projection $e_\cala$ of $\cala \subset \calb$. Hence indeed $\ad \beta_\calb(g) \circ E_I = E_I \circ \ad \beta_\calb(g)$ for any interval $I$ and $g \in G$ since $\beta_\calb(g) E_I(b) \beta_\calb(g)^{-1} \Omega_\cala = \beta_\calb(g) e_\cala b e_\cala \beta_\calb(g)^{-1} \Omega_\cala = e_\cala \beta_\calb(g) b \beta_\calb(g)^{-1} e_\cala \Omega_\cala = E_I(\beta_\calb(g) b \beta_\calb(g)^{-1}) \Omega_\cala$ for $b \in \calb(I)$. Similarly, for an isomorphism $u$ between $G$-equivariant extensions, we have $\ad u \circ E_I^{\calb_1} = E_I^{\calb_2} \circ \ad u$ automatically for any interval $I$.
	\end{enumerate}
\end{rem}

\begin{prop}
	Let $\cala$ be an irreducible local M\"{o}bius covariant net on $S^1$ with Haag duality on $\mathbb{R}$. Then, a finite index $G$-equivariant standard extension $\cala \subset \calb$ bijectively corresponds to a $G$-equivariant Q-system in $\rep \cala$ up to isomorphisms.

	\begin{proof}
		Suppose a finite index $G$-equivariant standard extension $\cala \subset \calb$ is given. Then, a corresponding Q-system $\theta$ is localized in some interval $I$ and can be restricted to a dual canonical endomorphism of $\cala(J) \subset \calb(J)$ for any interval $J$ including $I$.  By applying an argument in \cite[Lemma 2.8]{MR4153896} for the inclusion $\cala(I) \subset \calb(I)$, we can give $z_g \in \hom_{\en(\cala(I))}({}^g \theta|_{\cala(I)}, \theta|_{\cala(I)})$ for every $g \in G$, which is what we want since the argument for $\cala(J) \subset \calb(J)$ gives the same morphism $z_g$ and therefore $z_g \in \hom_{\en \cala_\infty}({}^g \theta, \theta)$. If we change $\theta$ for an isomorphic one, then we only get isomorphic $G$-equivariant Q-system by \cite[Lemma 2.8]{MR4153896}. If we change $\calb$ for an isomorphic $\tilde{\calb}$ with a unitary $u$, then $U \coloneqq \ad u \circ - \circ \ad u^\ast$ is a strict 2-functor between the 2-categories of morphisms of $\cala(I) \subset \calb(I)$ and $\cala(I) \subset \tilde{\calb}(I)$, see e.g. \cite[Section 1.3]{MR1966524} for the 2-category of morphisms. Since $u$ and therefore $U$ intertwine group actions, $U(\theta) = \theta$, which follows since $u$ commutes with $\cala_\infty$, is a $G$-equivariant Q-system corresponding to $\tilde{\calb}$. Thus, changing $\calb$ does not affect the resulting $G$-equivariant Q-systems.

Conversely, suppose we are given a $G$-equivariant Q-system $(\theta,w,x,z)$ in $\rep \cala$ that is localized in $I$. Then, $w,x, z_g \in \cala(I)$ by the Haag duality assumption on $\mathbb{R}$ of $\cala$, and we can construct a finite index standard extension $\iota: \cala \subset \calb$ on the GNS Hilbert space $L^2(\calb(I))$ of $\calb(I)$ associated with the vacuum state, see the proof of \cite[Theorem 4.9]{MR1332979}. By \cite[Lemma 2.10]{MR4153896}, we have an extension $\beta^\calb_I: G \to \aut(\calb(I))$ of an action $\ad \beta_\cala|_{\cala(I)}: G \to \aut(\cala(I))$. Then, for every $g \in G$, define a linear operator $\beta_\calb(g)$ on $L^2(\calb(I))$ by putting $\beta_\calb(g)[b] \coloneqq [\beta^\calb_I(g)(b)]$ for $[b] \in \calb(I) \subset L^2(\calb(I))$, which satisfies ${}^g b \coloneqq \beta^\calb_I(g)(b) = \beta_\calb(g) b \beta_\calb(g)^{-1}$ for $b \in \calb(I)$, in particular ${}^g a = \beta_\calb(g) a \beta_\calb(g)^{-1} = \ad \beta_\cala(g)(a)$ for $a \in \cala(I)$, by definition. We show that $\beta_\calb(g) \in \aut(\calb)$ and therefore $\beta_\calb$ is an extension of $\beta$. First, $\beta_\calb(g)$ is unitary since it is an automorphism that preserves the vacuum. Next, we check that $\ad \beta_\calb(g)(\calb(J)) \subset \calb(J)$ for an arbitrary interval $J$. For this, we show $\ad \beta_\calb(g)(\tilde{a}) = \ad \beta_\cala(g)(\tilde{a})$ for $\tilde{a} \in \cala_\infty$. Let $K$ be an interval containing $\tilde{a}$. By construction, we have an element $v \in \calb(I)$ with $L^2(\calb(I)) = v^\ast L^2(\cala(I))$ as a representation of $\cala(K)$, see the proof of \cite[Theorem 4.9]{MR1332979}. We also have ${}^g v = z_g^\ast v$ by \cite[Lemma 2.8]{MR4153896}. Since $\beta_\calb(g)|_{L^2(\cala(I))} = \beta_\cala(g)$ under the identification $L^2(\cala(I)) \cong H_\cala$ as cyclic representations of $\cala(I)$, we have
\begin{align*}
\beta_\calb(g) \tilde{a} \beta_\calb(g)^{-1} v^\ast [a] &= \beta_\calb(g) \tilde{a} {}^{g^{-1}} v^\ast [{}^{g^{-1}}a] = \beta_\calb(g) v^\ast \theta (\tilde{a}) z_{g^{-1}} [{}^{g^{-1}}a] = {}^g v^\ast \beta_\cala(g) \theta(\tilde{a})z_{g^{-1}}[{}^{g^{-1}}a] \\
&= v^\ast z_g ({}^g \theta)(\ad \beta_\cala(g)(\tilde{a})) \beta_\cala(g) z_{g^{-1}}[{}^{g^{-1}}a] \\
&= v^\ast \theta(\ad \beta_\cala(g)(\tilde{a})) z_g \ad \beta_\cala(g) (z_{g^{-1}})[a] = v^\ast \theta(\ad \beta_\cala(g)(\tilde{a}))[a] \\ 
&= \ad \beta_\cala(g)(\tilde{a})v^\ast [a]
\end{align*}
for every $a \in \cala(I)$ and therefore ${}^g \tilde{a} \coloneqq \ad \beta_\calb(g) (\tilde{a}) = \ad \beta_\cala(g) (\tilde{a})$. Note that $\beta_\calb(g)|_{L^2(\cala(I))} = \beta_\cala(g)$ is used in the third equality and $z_e = {}^e z_e = \id_\theta$ in the penultimate equality. Hence in particular ${}^g \cala(J) \subset \cala(J)$. Then, recall that we have a Q-system isomorphism $u:\theta \cong \tilde{\theta}$, where $\tilde{\theta}$ is localized in $J$, with $\calb(J) = \cala(J)uv$ by construction. We can make $u$ into a $G$-equivariant Q-system isomorphism by putting $z_g^{\tilde{\theta}} \coloneqq u z_g {}^g u^\ast$. Then, $z_g^{\tilde{\theta}} \in \cala(J)$ by Haag duality and therefore
\begin{align*}
{}^g \calb(J) = {}^g(\cala(J)) {}^gu {}^g v \subset \cala(J) {}^g u z_g^\ast v = \cala(J) {}^g u z_g^\ast v = \cala(J) (z_g^{\tilde{\theta}})^\ast u v = \cala(J)uv = \calb(J).
\end{align*}
Thus $\beta_\calb$ is an extension of $\beta_\cala$. If we change $\theta$ for an isomorphic $\tilde{\theta}$ localized in an interval $J$ with a unitary $u$ then we construct $\cala \subset \tilde{\calb}$ on $L^2(\tilde{\calb}(J))$. We have the counterpart $\tilde{v} \in \tilde{\calb}(J)$ of $v$ and $U: L^2(\calb(J)) \to L^2(\calb(I)); [a \tilde{v}] \mapsto a u[v]$ ($a \in \cala(J)$) defines an isomorphism of standard extensions, see the proof of \cite[Theorem 4.9]{MR1332979}. Moreover, $U$ intertwines group actions and therefore is an isomorphism of $G$-equivariant standard extensions since 
\begin{align*}
U^\ast \beta_{\calb}(g) U[a \tilde{v}] &= U^\ast \beta_{\calb}(g) a u [v] = {}^g (au) U^\ast [{}^g v] = {}^g(au) z_g^\ast U^\ast [v] = {}^g a (z^{\tilde{\theta}}_g)^\ast U^\ast u[v] = [{}^g a (z^{\tilde{\theta}}_g)^\ast \tilde{v}] \\
&= \beta_{\tilde{\calb}}(g)[a\tilde{v}]
\end{align*}
for any $g \in G$ and $a \in \cala(J)$. Thus, changing $\theta$ only yields an isomorphic $G$-equivariant extension. Note that in particular, we can harmlessly replace an interval in which a given $G$-equivariant Q-system is localized by another (say larger) one. It is easy to see that our constructions are mutually inverse.
	\end{proof}
\end{prop}

\begin{rem}
	\label{nonfaithrem}
In the proposition above, the action $\beta_\cala$ does not need to be faithful. That is, $\beta_\cala$ can have a nontrivial kernel. We encounter such a situation e.g. when $\cala = \calb^G$ for a M\"{o}bius covariant net $\calb$.
\end{rem}

See also a result \cite[Proposition 6.3]{MR3933137} for finite group actions on completely rational conformal nets. Note that our proof does not need either the finiteness of a group or the locality of a net. 

We go back to general theory. We collect here some graphical calculi for equivariant algebras. We follow the graphical notations for algebras and Frobenius algebras in \cite[Equation 2.22]{MR2187404}. Namely, the product and unit of an algebra are represented respectively by a fork and a small circle. 

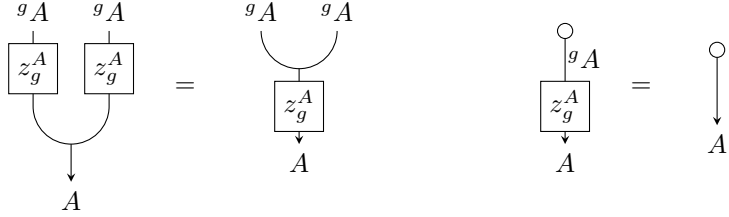
\begin{figure}[htb]
	\centering
	\begin{tikzpicture}
		\draw (0,1) arc (180:360:0.5);
			\draw[->] (0.5,0.5) -- (0.5,-0.5);
			\node[block] at (0.5,0){$z_g^A$};
			\node at (0,1.25){${}^g A$};
			\node at (1,1.25){${}^g A$};
			\node at (0.5,-0.75){$A$};
		\node at (-1,0.25){$=$};
		\begin{scope}[shift={(-3,0)}]
		\draw (0,1) -- (0,0) arc (180:360:0.5) -- (1,1);
		\draw[->] (0.5,-0.5) -- (0.5,-1);
		\node[block] at (0,0.5){$z_g^A$};
		\node[block] at (1,0.5){$z_g^A$};
		\node at (0.5,-1.25){$A$};
		\node at (0,1.25){${}^g A$};
		\node at (1,1.25){${}^g A$};
		\end{scope}
		\begin{scope}[shift={(4,-0.5)}]
			\draw[<-] (0,0) -- (0,1.5);
			\draw[fill=white] (0,1.5) circle [radius=0.1];
			\node[block] at (0,0.5){$z_g^A$};
			\node at (0,-0.25){$A$};
			\node at (0.25,1.125){${}^g A$};
			\node at (1,0.75){$=$};
			\begin{scope}[shift={(2,0.25)}]
				\draw[<-] (0,0) -- (0,1);
				\draw[fill=white] (0,1) circle [radius=0.1];
				\node at (0,-0.25){$A$};
			\end{scope}
		\end{scope}
	\end{tikzpicture}
	\caption{$z_g^A$ is an algebra homomorphism}
	\label{graphicalequivstrhom}
\end{figure}

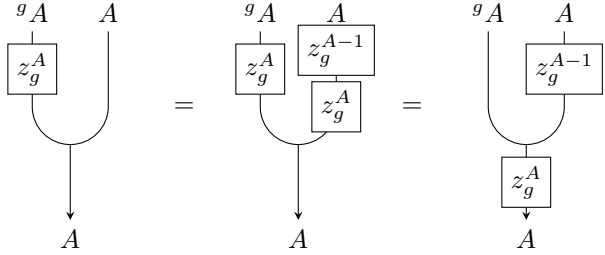
\begin{figure}[htb]
	\centering
	\begin{tikzpicture}
		\draw (0,1) -- (0,0) arc (180:360:0.5) -- (1,1);
		\draw[->] (0.5,-0.5) -- (0.5,-1.5);
		\node[block] at (0,0.5){$z_g^A$};
		\node at (0.5,-1.75){$A$};
		\node at (0,1.25){${}^g A$};
		\node at (1,1.25){$A$};
		\node at (2,0){$=$};
		\begin{scope}[shift={(3,0)}]
			\draw (0,1) -- (0,0) arc (180:360:0.5) -- (1,1);
		\draw[->] (0.5,-0.5) -- (0.5,-1.5);
		\node[block] at (0,0.5){$z_g^A$};
		\node[block] at (1,0){$z_g^A$};
		\node[block] at (1,0.75){$z_g^{A-1}$};
		\node at (0.5,-1.75){$A$};
		\node at (0,1.25){${}^g A$};
		\node at (1,1.25){$A$};
		\node at (2,0){$=$};		
		\end{scope}
		\begin{scope}[shift={(6,0)}]
			\draw (0,1) -- (0,0) arc (180:360:0.5) -- (1,1);
		\draw[->] (0.5,-0.5) -- (0.5,-1.5);
		\node[block] at (0.5,-1){$z_g^A$};
		\node[block] at (1,0.5){$z_g^{A-1}$};
		\node at (0.5,-1.75){$A$};
		\node at (0,1.25){${}^g A$};
		\node at (1,1.25){$A$};
		\end{scope}
	\end{tikzpicture}
	\caption{A useful calculus for equivariant algebras}
	\label{graphical_equiv_str_hom_calculi}
\end{figure}

Let $A$ be a $G$-equivariant algebra. Then, it is graphically represented in Figure \ref{graphicalequivstrhom} that $z^A_g$ is an algebra homomorphism for $g \in G$. The first equality in Figure \ref{graphicalequivstrhom} is often used in the form of Figure \ref{graphical_equiv_str_hom_calculi}. When $A$ is moreover a $G$-equivariant Frobenius algebra, we have similar representations of the coassociativity and counit property of $A$. It is often used in the form of Figure \ref{graphical_equiv_str} that $z^A$ is a $G$-equivariant structure. In particular, by putting $h = g^{-1}$, we obtain Figure \ref{graphical_equiv_str_inverse} since $z^A_g {}^g z^A_{g^{-1}} = z_e^A = \id_A$.

\begin{figure}[H]
	\begin{tabular}{cc}
		\begin{minipage}[b]{0.45\hsize}
			\centering
			\begin{tikzpicture}
				\draw[->] (-0.5,-0.5) -- (-0.5,-1.5) -- (-1,-2) -- (-1,-3);
				\draw[dashed,cross,->] (-1,-0.5) -- (-1,-1.5) -- (-0.5,-2) -- (-0.5,-3);
				\node[block] at (-0.5,-1){$z_h^A$};
				\node[block] at (-1,-2.5){$z_g^A$};
				\node at (-1,-0.25){$\lambda$};
				\node at (-0.5,-3.25){$\lambda$};
				\node at (-0.5,-0.25){${}^h A$};
				\node at (-1,-3.25){$A$};
				\node at (-1.25,-1.75){$g$};
				\node at (0,-1.75){$=$};
				\begin{scope}[shift={(2,0)}]
					\draw[->] (-0.5,-0.5) -- (-0.5,-1.5) -- (-1,-2) -- (-1,-3);
					\draw[dashed,cross,->] (-1,-0.5) -- (-1,-1.5) -- (-0.5,-2) -- (-0.5,-3);
					\node[block] at (-1,-2.5){$z_{gh}^A$};
					\node at (-1,-0.25){$\lambda$};
					\node at (-0.5,-3.25){$\lambda$};
					\node at (-0.5,-0.25){${}^h A$};
					\node at (-1,-3.25){$A$};
					\node at (-1.25,-1.75){$g$};
				\end{scope}
			\end{tikzpicture}
			\caption{$z^A$ is a $G$-equivariant structure}
			\label{graphical_equiv_str}
		\end{minipage}
		\begin{minipage}[b]{0.45\hsize}
			\centering
			\begin{tikzpicture}
				\draw[->] (-0.5,-0.5) -- (-0.5,-1.5) -- (-1,-2) -- (-1,-3);
				\draw[dashed,cross,->] (-1,-0.5) -- (-1,-1.5) -- (-0.5,-2) -- (-0.5,-3);
				\node[block] at (-0.5,-1){$z_{g^{-1}}^A$};
				\node[block] at (-1,-2.5){$z_g^A$};
				\node at (-1,-0.25){$\lambda$};
				\node at (-0.5,-3.25){$\lambda$};
				\node at (-0.5,-0.25){${}^{g^{-1}} A$};
				\node at (-1,-3.25){$A$};
				\node at (-1.25,-1.75){$g$};
				\node at (0,-1.75){$=$};
				\begin{scope}[shift={(2,0)}]
					\draw[->] (-0.5,-0.5) -- (-0.5,-1.5) -- (-1,-2) -- (-1,-3);
					\draw[dashed,cross,->] (-1,-0.5) -- (-1,-1.5) -- (-0.5,-2) -- (-0.5,-3);
					\node at (-1,-0.25){$\lambda$};
					\node at (-0.5,-3.25){$\lambda$};
					\node at (-0.5,-0.25){${}^{g^{-1}} A$};
					\node at (-1,-3.25){$A$};
					\node at (-1.25,-1.75){$g$};
				\end{scope}
			\end{tikzpicture}
			\caption{$z^A_g {}^g z^A_{g^{-1}} = \id_A$}
			\label{graphical_equiv_str_inverse}
		\end{minipage}
	\end{tabular}
\end{figure}

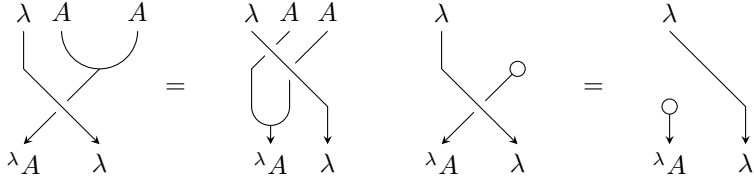
\begin{figure}[H]
	\centering
	\begin{tikzpicture}
		\draw (0,0) arc (180:360:0.5);
		\draw[->] (0.5,-0.5) -- (-0.5,-1.5);
		\draw[cross,->] (-0.5,0) -- (-0.5,-0.5) -- (0.5,-1.5);
		\node at (0,0.25){$A$};
		\node at (1,0.25){$A$};
		\node at (-0.5,0.25){$\lambda$};
		\node at (-0.5,-1.75){${}^\lambda A$};
		\node at (0.5,-1.75){$\lambda$};
		\node at (1.5,-0.75){$=$};
		\begin{scope}[shift={(3,0)}]
			\draw (0,0) -- (-0.5,-0.5) -- (-0.5,-1) arc (180:360:0.25) -- (0,-0.5) -- (0.5,0);
			\draw[->] (-0.25,-1.25) -- (-0.25,-1.5);
			\draw[cross,->] (-0.5,0) -- (0.5,-1) -- (0.5,-1.5);
			\node at (0,0.25){$A$};
			\node at (0.5,0.25){$A$};
			\node at (-0.5,0.25){$\lambda$};
			\node at (-0.25,-1.75){${}^\lambda A$};
			\node at (0.5,-1.75){$\lambda$};
		\end{scope}
		\begin{scope}[shift={(5.5,0)}]
			\draw[->] (0.5,-0.5) -- (-0.5,-1.5);
			\draw[fill=white] (0.5,-0.5) circle [radius=0.1];
			\draw[cross,->] (-0.5,0) -- (-0.5,-0.5) -- (0.5,-1.5);
			\node at (-0.5,0.25){$\lambda$};
			\node at (-0.5,-1.75){${}^\lambda A$};
			\node at (0.5,-1.75){$\lambda$};
			\node at (1.5,-0.75){$=$};
			\begin{scope}[shift={(3,0)}]
				\draw[->] (-0.5,-1) -- (-0.5,-1.5);
				\draw[fill=white] (-0.5,-1) circle [radius=0.1];
				\draw[cross,->] (-0.5,0) -- (0.5,-1) -- (0.5,-1.5);
				\node at (-0.5,0.25){$\lambda$};
				\node at (-0.5,-1.75){${}^\lambda A$};
				\node at (0.5,-1.75){$\lambda$};
			\end{scope}
		\end{scope}
	\end{tikzpicture}
	\caption{Algebra structures and crossings}
	\label{graphical_product_bfe}
\end{figure}

Finally, note that when $A$ is an algebra in a $G$-braided multitensor category, we can move the algebra structures through crossings as in Figure \ref{graphical_product_bfe} by the definition of ${}^\lambda A$.

\section{Equivariant $\alpha$-induction}
\label{section_equiv_alpha}

The $\alpha$-induction for twisted representations of M\"{o}bius covariant nets (see Example \ref{eg_twisted_rep}) was introduced by Nojima in \cite{MR4153896}. In this section, we reformulate this notion in terms of tensor categories (Definition \ref{def_alpha_ind} and Remark \ref{rem_nojima}), which is an equivariant generalization of Ostrik's work \cite[Section 5.1]{MR1976459}, see also \cite[Section 4.6]{bklr}. For this, we begin with crossed structures on bimodule categories.

\subsection{Group actions on the bicategories of equivariant Frobenius algebras}
\label{subsection_action_bicat_frob}

In this subsection, we induce a group action on the bicategory of equivariant Frobenius algebras from that on the ambient multitensor category (Proposition \ref{bimodgroupactiondef}). We also see that this indeed makes bimodule categories $G$-crossed (Proposition \ref{prop_bimod_crossed}).

It is known \cite[Sections 4 and 5]{MR2075605} that the special Frobenius algebras in a multitensor category $\calc$ form a rigid bicategory, whose 1-cells are bimodules as in ordinary ring theory. Note that this bicategory is defined only up to equivalences since we have to fix relative tensor products of bimodules to obtain a composition of 1-cells. In particular, for a special Frobenius algebra $A$ in $\calc$, the category $\bimod_\calc(A)$ of $A$-bimodules in $\calc$ is a multitensor category. A special Frobenius algebra $A$ is called \emph{simple} if $\bimod_\calc(A)$ is moreover a tensor category i.e. $\en_{\bimod_\calc(A)}(A)$ is one-dimensional, see \cite[Definition 2.26]{MR2187404}.

\begin{figure}[htb]
	\centering
	\begin{tikzpicture}
		\draw[rounded corners] (0,0) -- (0,-0.5) -- (0.5,-0.5);
		\draw[->] (0.5,0) -- (0.5,-1) -- (-0.5,-2);
		\draw[cross,->] (-0.5,0) -- (-0.5,-1) -- (0.5,-2);
		\draw (0.5,-0.4) arc (90:270:0.1);
		\node at (-0.5,0.25){$\lambda$};
		\node at (0.5,0.25){$\mu$};
		\node at (0,0.25){$A$};
		\node at (-0.5,-2.25){${}^\lambda \mu$};
		\node at (0.5,-2.25){$\lambda$};
		\node at (1.25,-1){$=$};
		\begin{scope}[shift={(2.5,0)}]
			\draw[rounded corners] (0,0) -- (-0.5,-0.5) -- (-0.5,-1.5) -- (0,-1.5);
			\draw[->] (0.5,0) -- (0.5,-0.5) -- (0,-1) -- (0,-2);
			\draw[cross,->] (-0.5,0) -- (0.5,-1) -- (0.5,-2);
			\draw (0,-1.4) arc (90:270:0.1);
			\node at (-0.5,0.25){$\lambda$};
			\node at (0.5,0.25){$\mu$};
			\node at (0,0.25){$A$};
			\node at (0,-2.25){${}^\lambda \mu$};
			\node at (0.5,-2.25){$\lambda$};
			\node at (-0.5,-1.75){${}^\lambda A$};
		\end{scope}
		\begin{scope}[shift={(5.5,0)}]
			\draw[rounded corners] (0.5,0) -- (0.5,-0.5) -- (0,-0.5);
			\draw[->] (0,0) -- (0,-1) -- (-0.5,-1.5) -- (-0.5,-2);
			\draw[cross,->] (-0.5,0) -- (-0.5,-1) -- (0.5,-2);
			\draw (0,-0.4) arc (90:-90:0.1);
			\node at (-0.5,0.25){$\lambda$};
			\node at (0.5,0.25){$A$};
			\node at (0,0.25){$\mu$};
			\node at (-0.5,-2.25){${}^\lambda \mu$};
			\node at (0.5,-2.25){$\lambda$};
			\node at (1.25,-1){$=$};
			\begin{scope}[shift={(2.5,0)}]
				\draw (0.5,0) -- (0.5,-0.5) -- (0,-1);
				\draw[rounded corners] (0,-1) -- (0,-1.5) -- (-0.5,-1.5);
				\draw[->] (0,0) -- (-0.5,-0.5) -- (-0.5,-2);
				\draw[cross,->] (-0.5,0) -- (0.5,-1) -- (0.5,-2);
				\draw (-0.5,-1.4) arc (90:-90:0.1);
				\node at (-0.5,0.25){$\lambda$};
				\node at (0.5,0.25){$A$};
				\node at (0,0.25){$\mu$};
				\node at (-0.5,-2.25){${}^\lambda \mu$};
				\node at (0.5,-2.25){$\lambda$};
				\node at (0,-1.75){${}^\lambda A$};
			\end{scope}	
		\end{scope}
	\end{tikzpicture}
	\caption{Module products and crossings}
	\label{graphical_move_module}
\end{figure}
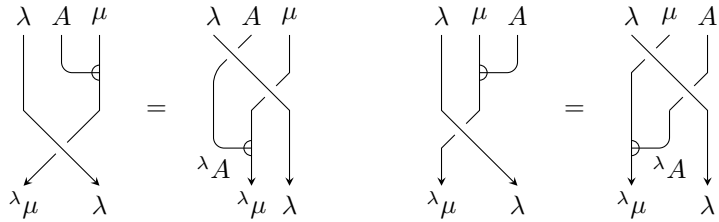

Let $\calc$ be a multitensor category with an action of $G$ and let $A$ and $B$ be algebras in $\calc$. For an $A$-$B$-bimodule $\lambda = (\lambda, m^{\mathrm{L}}_\lambda, m^{\mathrm{R}}_\lambda)$ in $\calc$, ${}^g \lambda = ({}^g \lambda, {}^g m^{\mathrm{L}}_\lambda, {}^g m^{\mathrm{R}}_\lambda)$ is an ${}^g A$-${}^g B$-bimodule. Note that we can move module products through crossings as in Figure \ref{graphical_move_module} for $\lambda \in \homog(\calc)$ and an $A$-$B$-bimodule $\mu$ by this definition. We modify this procedure to obtain an action that restricts to the multitensor categories of bimodules.

\begin{prop}
	\label{bimodgroupactiondef}
	Let $\calc$ be a multitensor category equipped with an action $\gamma^\calc$ of a group $G$. Let $\frob^G(\calc)$ denote the bicategory of special $G$-equivariant Frobenius algebras in $\calc$.  Then, for a 1-cell $\lambda = (\lambda,m^{\mathrm{L}}_\lambda,m^{\mathrm{R}}_\lambda) : A \to B$ in $\frob^G(\calc)$ and $g \in G$, we can give a 1-cell ${}^g \lambda \coloneqq ({}^g \lambda, {}^g m^{\mathrm{L}}_\lambda \circ (z^{A-1}_g \otimes \id_\lambda), {}^g m^{\mathrm{R}}_{ \lambda}(\id_{{}^g \lambda} \otimes z^{B-1}_g)): A \to B$. This assignment naturally defines an action of $G$ on $\frob^G(\calc)$. Another choice of relative tensor products only yields a $G$-biequivalent action.

	\begin{proof}
		The assignment $\gamma(g): \lambda \mapsto {}^g \lambda$ defined above is well-defined since the graphical calculations in Figures \ref{graphicalbimodgroupactionassoc} and \ref{graphicalbimodgroupactionunit}, where $z \coloneqq z^B$ and small half circles denote module products as in \cite[Section 3.6]{bklr}, shows that ${}^g \lambda$ is indeed a left $B$-module, and the proof for right $A$-modularity is similar. We regard $\gamma(g)$ as a family of functors by putting $\gamma(g)(f) \coloneqq \gamma^\calc(g)(f)$ for a 2-cell $f$. 
		
		\begin{figure}[htb]
			\centering
			\begin{tikzpicture}
					\draw (0,0)--(0,2);
					\draw[rounded corners] (0,0.5) -- (-1.5,0.5) -- (-1.5,2);
					\draw[rounded corners] (0,1) -- (-0.5,1) -- (-0.5,2);
					\draw (0,0.6) arc (90:270:0.1);
					\draw (0,1.1) arc (90:270:0.1);
					\node at (0,-0.25){${}^g \lambda$};
					\node at (0,2.25){${}^g \lambda$};
					\node at (-0.5,2.25){$B$};
					\node at (-1.5,2.25){$B$};
					\node[block] at (-0.5,1.5){$z_g^{-1}$};
					\node[block] at (-1.5,1.5){$z_g^{-1}$};
					\node at (0.5,1){$=$};
					\node at (-1.5,0.25){${}^g B$};
					\node at (-0.75,0.8){${}^g B$};
				\begin{scope}[shift={(3,0)}]
					\draw (0,0)--(0,2);
					\draw[rounded corners] (-1.5,1.5) -- (-1.5,2);
					\draw[rounded corners] (-0.5,1.5) -- (-0.5,2);
					\draw (-1.5,1.5) arc (180:360:0.5);
					\draw[rounded corners] (-1,1) -- (-1,0.5)--(0,0.5);
					\draw (0,0.6) arc (90:270:0.1);
					\node at (0,-0.25){${}^g \lambda$};
					\node at (0,2.25){${}^g \lambda$};
					\node at (-0.5,2.25){$B$};
					\node at (-1.5,2.25){$B$};
					\node[block] at (-0.5,1.5){$z_g^{-1}$};
					\node[block] at (-1.5,1.5){$z_g^{-1}$};
					\node at (0.5,1){$=$};
					\node at (-1,0.25){${}^g B$};
				\end{scope}
				\begin{scope}[shift={(6,0)}]
					\draw (0,0)--(0,2);
					\draw (-1.5,2) arc (180:360:0.5);
					\draw[rounded corners] (-1,1.5) -- (-1,0.5)--(0,0.5);
					\draw (0,0.6) arc (90:270:0.1);
					\node at (0,-0.25){${}^g \lambda$};
					\node at (0,2.25){${}^g \lambda$};
					\node at (-0.5,2.25){$B$};
					\node at (-1.5,2.25){$B$};
					\node[block] at (-1,1){$z_g^{-1}$};
					\node at (-1,0.25){${}^g B$};
				\end{scope}
			\end{tikzpicture}
			\caption{The left $B$-modularity of ${}^g \lambda$}
			\label{graphicalbimodgroupactionassoc}
		\end{figure}
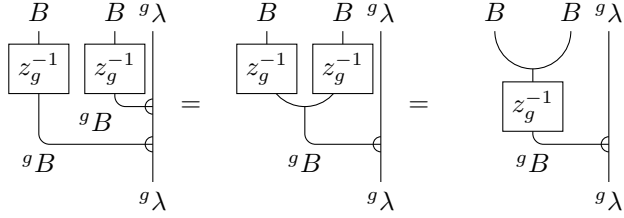

		\begin{figure}[htb]
			\centering
			\begin{tikzpicture}
				\draw (0,0.5)--(0,2.5);
				\draw[rounded corners] (0,1) -- (-0.5,1) -- (-0.5,2.5);
				\draw[fill=white] (-0.5,2.5) circle [radius=0.1];
				\draw (0,1.1) arc (90:270:0.1);
				\node at (0,0.25){${}^g \lambda$};
				\node at (0,2.75){${}^g \lambda$};
				\node[block] at (-0.5,1.5){$z_g^{-1}$};
				\node at (0.5,1.5){$=$};
				\node at (-0.75,0.8){${}^g B$};
				\node at (-0.75,2.25){$B$};
				\begin{scope}[shift={(2,0)}]
					\draw (0,0.5)--(0,2.5);
					\draw[rounded corners] (0,1) -- (-0.5,1) -- (-0.5,2.5);
					\draw[fill=white] (-0.5,2.5) circle [radius=0.1];
					\draw (0,1.1) arc (90:270:0.1);
					\node at (0,0.25){${}^g \lambda$};
					\node at (0,2.75){${}^g \lambda$};
					\node at (0.5,1.5){$=$};
					\node at (-0.75,0.8){${}^g B$};
				\end{scope}
				\begin{scope}[shift={(3.25,0)}]
					\draw (0,0.5)--(0,2.5);
					\node at (0,0.25){${}^g \lambda$};
					\node at (0,2.75){${}^g \lambda$};
				\end{scope}
			\end{tikzpicture}
			\caption{The left unit property of ${}^g \lambda$}
			\label{graphicalbimodgroupactionunit}
		\end{figure}
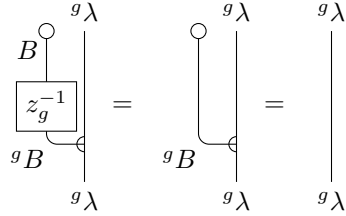
		
		We make $\gamma(g)$ into a pseudofunctor. For 1-cells $\lambda: A \to B$ and $\mu: B \to C$, let $(\mu \otimes_B \lambda, s_{\mu,\lambda}, r_{\mu,\lambda})$ denote a retract of $\mu \otimes_\calc \lambda$ with the idempotent $e_{\mu,\lambda} \coloneqq s_{\mu,\lambda} \circ r_{\mu,\lambda}$. Note that $s_{\mu,\lambda}$ and $r_{\mu,\lambda}$ are natural in $\mu$ and $\lambda$ by the definition of the bifunctor $\otimes_B$. By the graphical calculation in Figure \ref{bimodgroupactionmonoidalproof}, we have $e_{{}^g \mu, {}^g \lambda} = {}^g e_{\mu,\lambda}$ and therefore a unique subobject isomorphism $J^{\gamma(g)}_{\mu, \lambda} \coloneqq {}^g (r_{\mu, \lambda}) \circ s_{{}^g \mu, {}^g \lambda} : {}^g \mu \otimes_B {}^g \lambda \cong {}^g(\mu \otimes_B \lambda)$. We can show that the bicategory version of \cite[Diagram 2.23]{egno} commutes by $e_{{}^g \mu, {}^g \lambda} = {}^g e_{\mu,\lambda}$ and the naturality of $s$ and $r$. We also put $\varphi^{\gamma(g)}_A \coloneqq (z_g^A)^{-1}$ for $A \in \obj(\frob^G(\calc))$, which is an $A$-bimodule morphism since $z_g^A$ is an algebra homomorphism. Then, by the definition of the left and right unit isomorphisms of $\frob^G(\calc)$, the commutativity of the bicategory versions of \cite[Diagrams 2.25 and 2.26]{egno} follows from the definition of $\gamma$ and the naturality of $s$. Thus, $(\gamma(g),J^{\gamma(g)},\varphi^{\gamma(g)})$ is a pseudofunctor.

		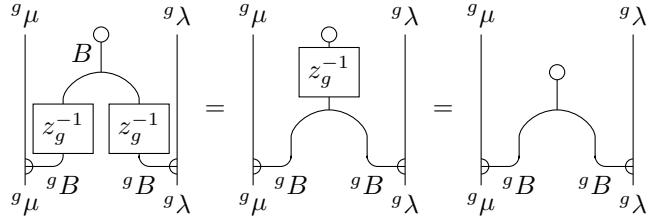
\begin{figure}[htb]
			\centering
			\begin{tikzpicture}
				\draw (0,1.25) -- (0,-0.75);
				\draw[rounded corners] (0,-0.5) -- (0.5,-0.5) -- (0.5,0.25) arc (180:0:0.5) -- (1.5,-0.25) -- (1.5,-0.5) -- (2,-0.5);
				\draw (2,1.25) -- (2,-0.75);
				\draw (1,0.75) -- (1,1.25);
				\draw[fill=white] (1,1.25) circle [radius=0.1];
				\node[block] at (0.5,0) {$z_g^{-1}$};
				\node[block] at (1.5,0) {$z_g^{-1}$};
				\draw (0,-0.4) arc (90:-90:0.1);
				\draw (2,-0.4) arc (90:270:0.1);
				\node at (2.5,0.25){$=$};
				\node at (0,1.5){${}^g \mu$};
				\node at (2,1.5){${}^g \lambda$};
				\node at (0.75,1){$B$};
				\node at (0,-1){${}^g \mu$};
				\node at (2,-1){${}^g \lambda$};
				\node at (0.5,-0.75){${}^g B$};
				\node at (1.5,-0.75){${}^g B$};
				\begin{scope}[shift={(3,0)}]
					\draw (0,1.25) -- (0,-0.75);
					\draw[rounded corners] (0,-0.5) -- (0.5,-0.5) -- (0.5,-0.25) arc (180:0:0.5) -- (1.5,-0.5) -- (2,-0.5);
					\draw (2,1.25) -- (2,-0.75);
					\draw (1,0.25) -- (1,1.25);
					\draw[fill=white] (1,1.25) circle [radius=0.1];
					\node[block] at (1,0.75) {$z_g^{-1}$};
					\draw (0,-0.4) arc (90:-90:0.1);
					\draw (2,-0.4) arc (90:270:0.1);
					\node at (2.5,0.25){$=$};
					\node at (0,1.5){${}^g \mu$};
					\node at (2,1.5){${}^g \lambda$};
					\node at (0,-1){${}^g \mu$};
					\node at (2,-1){${}^g \lambda$};
					\node at (0.5,-0.75){${}^g B$};
					\node at (1.5,-0.75){${}^g B$};
				\end{scope}
				\begin{scope}[shift={(6,0)}]
					\draw (0,1.25) -- (0,-0.75);
					\draw[rounded corners] (0,-0.5) -- (0.5,-0.5) -- (0.5,-0.25) arc (180:0:0.5) -- (1.5,-0.5) -- (2,-0.5);
					\draw (2,1.25) -- (2,-0.75);
					\draw (1,0.25) -- (1,0.75);
					\draw[fill=white] (1,0.75) circle [radius=0.1];
					\draw (0,-0.4) arc (90:-90:0.1);
					\draw (2,-0.4) arc (90:270:0.1);
					\node at (0,1.5){${}^g \mu$};
					\node at (2,1.5){${}^g \lambda$};
					\node at (0,-1){${}^g \mu$};
					\node at (2,-1){${}^g \lambda$};
					\node at (0.5,-0.75){${}^g B$};
					\node at (1.5,-0.75){${}^g B$};
				\end{scope}
			\end{tikzpicture}
			\caption{$e_{{}^g \mu, {}^g \lambda} = {}^g e_{\mu,\lambda}$}
			\label{bimodgroupactionmonoidalproof}
		\end{figure}

		Since $z_{gh} = z_g {}^g z_h$, we have a canonical invertible 2-cell $(\chi^{\gamma}_{g, h})_\lambda: {}^g ({}^h \lambda) \cong {}^{gh} \lambda$ by the coherence in $\calc$, which gives a pseudonatural equivalence $\chi^{\gamma}_{g, h}: \gamma(g) \gamma(h) \simeq \gamma(gh)$ with $\chi^{\gamma,0}_{g, h} \coloneqq \mathbf{1}$. Similarly, $z_e = \id_B$ gives a pseudonatural equivalence $\iota^{\gamma}: \id \simeq \gamma(e)$ with $\iota^{\gamma,0} = \mathbf{1}$. We can define all the remaining modifications to be canonical invertible 2-cells, and then $\gamma$ is an action of $G$ by coherence. Another choice of relative tensor products only yields a $G$-biequivalent action since we get the same $\chi^\gamma$ and $\iota^\gamma$.
	\end{proof}
\end{prop}

	It is known \cite[Proposition 3.33]{bklr} that for a type III factor $N$, there is a biequivalence between the bicategory of Q-systems in $\en_0(N)$ and the 2-category of finite index extensions of $N$ whose 1-cells from $M_1$ to $M_2$ are the subobjects of $\iota_2 \lambda \overline{\iota_1}$ for $\lambda \in \en_0(N)$, where $\iota_i: N \subset M_i$ for $i=1,2$. When a group $G$ acts on $N$, in the same way, we have a biequivalence between the bicategory $\operatorname{Q}^G(\en_0(N))$ of the $G$-equivariant Q-systems in $\en_0(N)$ and the 2-category $\operatorname{Ext}^G(N)$ of finite index extensions of $N$ with $G$-actions whose 1-cells are as above. 
	
\begin{prop}
	This biequivalence can be made into a $G$-biequivalence, where $G$ acts on $\operatorname{Q}^G(\en_0(N))$ by Proposition \ref{bimodgroupactiondef} and on $\operatorname{Ext}^G(N)$ by the adjoint action.
	
	\begin{proof}
		First, note that if a 1-cell $\lambda$ of $\operatorname{Q}^G(\en_0(N))$ is standard (see \cite[Section 3.6]{bklr}), then ${}^g \lambda$ is again standard for $g \in G$. Therefore, the $G$-action on $\operatorname{Q}^G(\en_0(N))$ restricts to the subbicategory $\operatorname{Q}^G_{\mathrm{st}}(\en_0(N))$ whose 1-cells are standard bimodules, and we have a $G$-biequivalence $\operatorname{Q}^G_{\mathrm{st}}(\en_0(N)) \simeq \operatorname{Q}^G(\en_0(N))$ by \cite[Proposition 3.32(i)]{bklr}.

		Then, let $F: \operatorname{Ext}^G(N) \to \operatorname{Q}^G_{\mathrm{st}}(\en_0(N))$ denote the biequivalence defined by putting $F(\varphi) \coloneqq \overline{\iota_1} \varphi \iota_2$ for $\varphi \in \hom_{\operatorname{Ext}^G(N)}(M_2,M_1)$, see \cite[Propositions 3.33 and 3.38]{bklr}. It is enough to make $F$ into a $G$-biequivalence. For this, note that for $M \in \obj(\operatorname{Ext}^G(N))$, we have $z^{F(M)}_g \in \hom({}^g \overline{\iota}, \overline{\iota})$ for $g \in G$, where $\iota: N \subset M$, since
		\begin{align*}
		z_g ({}^g \overline{\iota}) (nv) &= z_g ({}^g \theta)(n) {}^g(\theta(z_{g^{-1}}^{ \ast}) x) = z_g ({}^g \theta)(n)({}^g \theta)(z_g) {}^g x = \theta(n)z_g ({}^g \theta)(z_g){}^g x = \theta(n) x z_g \\
		&= \overline{\iota}(nv) z_g
		\end{align*}
		for $n \in N$ and $v \coloneqq \mathrm{coev}_{\overline{\iota}}$, where $\theta \coloneqq F(M) = \overline{\iota} \iota$, $x \coloneqq \overline{\iota}(v)$ and $z \coloneqq z^{\theta}$. Then, put $(\tilde{\eta}^F_g)_{\varphi} \coloneqq z_g^{F(M_1)\ast} \otimes \id_{\varphi \iota_2} : \overline{\iota_1} {}^g \varphi \iota_2 \cong {}^g(\overline{\iota_1} \varphi \iota_2)$ for $\varphi \in \hom_{\operatorname{Ext}^G(N)}(M_2,M_1)$ and $g \in G$, which is a bimodule isomorphism since ${}^g v^\ast = v^\ast z_g$. Since $s_{F(\varphi),F(\psi)} = \id_{\overline{\iota_1} \varphi} \otimes v_2 \otimes \id_{\psi \iota_3}$ and $r_{{}^g F(\varphi), {}^g F(\psi)} = \id_{{}^g\overline{\iota_1} {}^g \varphi} \otimes {}^g v_2^\ast \otimes \id_{{}^g \psi \iota_3}$ for $\psi \in \hom_{\operatorname{Ext}^G(N)}(M_3,M_2)$ and $g \in G$, where $v_2 \coloneqq \mathrm{coev}_{\overline{\iota_2}}$, when we take $F(\varphi) \otimes_{F(M_2)} F(\psi) = F(\varphi \psi)$ by Figure \ref{bimodgroupactionmonoidalproof} and the proof of \cite[Proposition 3.38]{bklr}, we obtain $(\tilde{\eta}^F_g)_{\varphi \psi} = (\tilde{\eta}^F_g)_{\varphi} \otimes_{F(M_2)} (\tilde{\eta}^F_g)_{\psi}$ from ${}^g v^\ast = v^\ast z_g$. Therefore, by putting $\tilde{\eta}^{F,0}_g \coloneqq \mathbf{1}$, we get a pseudonatural equivalence $\tilde{\eta}_g^F: F \gamma^{\operatorname{Ext}^G(N)}(g) \simeq \gamma^{\operatorname{Q}^G_{\mathrm{st}}(\en_0(N))}(g) F$. Since $(\tilde{\eta}_{g h}^F)_\varphi = {}^g (\tilde{\eta}^F_h)_{\varphi} (\tilde{\eta}^F_g)_{{}^h \varphi}$ by $z_{g h} = z_g {}^g z_h$, we can make $(F, \tilde{\eta}^F)$ into a $G$-biequivalence by coherence.
	\end{proof}
\end{prop}

By definition, the restriction of the action defined in Proposition \ref{bimodgroupactiondef} to $\calc$ coincides with $\gamma^\calc$, and that to $\bimod_\calc(A)$ for $A \in \obj(\frob^G(\calc))$ is an action on a multitensor category. We show that $\bimod_\calc(A)$ is indeed $G$-crossed (when $A$ is symmetric and neutral, see Propositions \ref{bimodpivotalprop} and \ref{prop_bimod_crossed}).

\begin{prop}
	\label{bimodpivotalprop}
	Let $A$ be a symmetric special Frobenius algebra in a pivotal multitensor category $\calc$. Then, $\bimod_\calc(A)$ is again pivotal. Another choice of relative tensor products only yields an isomorphic pivotal structure. If $\calc$ is moreover a tensor category and $A$ is simple, then $\tr^{\mathrm{L}}(\delta^{\bimod_\calc(A)}_\lambda f) = \tr^{\mathrm{L}} (\delta^\calc_\lambda f)/\dim_{\delta^\calc} (A)$ and $\tr^{\mathrm{R}}(f \delta^{\bimod_\calc(A)-1}_\lambda) = \tr^{\mathrm{R}} (f \delta^{\calc-1}_{\lambda})/\dim_{\delta^\calc} (A)$ for $f \in \en_{\bimod_\calc(A)}(\lambda)$. In particular, $\dim_{\delta^{\bimod_\calc(A)}}(\lambda) = \dim_{\delta^{\calc}}(\lambda)/\dim_{\delta^\calc}(A)$ for $\lambda \in \obj(\bimod_\calc(A))$, and if $\calc$ is spherical, then $\bimod_\calc(A)$ is again spherical. If a group $G$ acts pivotally on $\calc$ (not necessarily tensor) and $A \in \obj(\frob^G(\calc))$, then the induced action of $G$ on $\bimod_\calc(A)$ is again pivotal.

	\begin{proof}
		Let $\lambda \in \obj(\bimod_\calc(A))$ and let $m_\lambda$ denote either the left or right module product of $\lambda$. By the proof of the rigidity of $\bimod_\calc(A)$ in \cite[Section 5]{MR2075605}, we have $m_{\lambda^{\vee \vee}} = (m_\lambda)^{\vee \vee}$ in $\calc$ when $A^{\vee}$ is taken to be $A$. Therefore, $\delta^\calc_\lambda$ is an $A$-bimodule morphism by the naturality and monoidality of $\delta^\calc$ if $\delta_A^\calc = \id_A$. Let us put $\delta^{\bimod_\calc(A)}_\lambda \coloneqq \delta^\calc_\lambda$, which is natural in $\lambda$ by definition. The monoidality of $\delta^{\bimod_\calc(A)}$ follows from the naturality and monoidality of $\delta^\calc$. Thus, $\delta^{\bimod_\calc(A)}$ is a pivotal structure on $\bimod_\calc(A)$. Another choice of relative tensor products only yields an isomorphic pivotal structure since an equivalence from the new bimodule category to the original one is defined to be an identity as a functor.
		
		Now, suppose $\calc$ is a tensor category and $A$ is simple. By the proof of rigidity \cite[Section 5]{MR2075605}, evaluation and coevaluation maps of $\lambda \in \obj(\bimod_\calc(A))$ are respectively given in Figures \ref{graphicalbimodev} and \ref{graphicalbimodcoev} as morphisms in $\calc$. Then, we find that $e_{\lambda, \lambda^\vee}$, the idempotent for the subobject $\lambda \otimes_A \lambda^\vee$ of $\lambda \otimes_\calc \lambda^\vee$ (see the proof of Proposition \ref{bimodgroupactiondef}), in $\tr^{\mathrm{L}}(\delta^{\bimod_\calc(A)}_\lambda f)$ cancels for $f \in \en_{\bimod_\calc(A)}(\lambda)$ by the definition of the left module product of $\lambda^\vee$, see \cite[Section 5]{MR2075605}. Therefore, we have 
		\begin{align*}
			\tr^{\mathrm{L}}(\delta^{\bimod_\calc(A)}_\lambda f) &= (\varepsilon_A \eta_A)^{-1} \varepsilon_A \tr^{\mathrm{L}}(\delta^{\bimod_\calc(A)}_\lambda f) \eta_A \\
			&= \dim_{\delta^\calc}(A)^{-1} \varepsilon_A \mathrm{ev}^{\bimod_\calc(A)}_{\lambda^\vee} \delta^{\bimod_\calc(A)}_\lambda f \mathrm{coev}^{\bimod_\calc(A)}_{\lambda} \eta_A \\
			&= \dim_{\delta^\calc}(A)^{-1} \mathrm{ev}^{\calc}_{\lambda^\vee} \delta^{\calc}_\lambda f \mathrm{coev}^{\calc}_{\lambda} = \dim_{\delta^\calc}(A)^{-1} \tr^{\mathrm{L}} (\delta^\calc f).
		\end{align*}
		The proof for $\tr^{\mathrm{R}}$ is similar. Finally, the last statement follows since the group action on the morphisms of $\bimod_\calc(A)$ coincides with that of $\calc$.
		\begin{figure}[htb]
			\begin{tabular}{cc}
				\begin{minipage}[b]{0.45\hsize}
					\centering
					\begin{tikzpicture}
						\draw[rounded corners,<-] (1,0.25) -- (1,1.25) arc (0:180:0.25) -- (0.5,1) -- (0.25,1);
						\draw[rounded corners] (-0.25,1.75) -- (-0.25,0.5) arc (180:360:0.25) -- (0.25,1.75);
						\draw (0.25,1.1) arc (90:-90:0.1);
						\draw (0.75,1.5) -- (0.75,1.75);
						\draw[fill=white] (0.75,1.75) circle [radius=0.1];
						\draw (0,1.75) -- (0,2.25);
						\node at (1,0){$A$};
						\node at (-0.5,0.5){$\lambda^\vee$};
						\node at (0.5,0.5){$\lambda$};
						\node[block] at (0,1.75){$s_{\lambda^\vee,\lambda}$};
						\node at (0,2.5){$\lambda^\vee \otimes_A \lambda$};
					\end{tikzpicture}
					\caption{Evaluation in $\bimod_\calc(A)$}
					\label{graphicalbimodev}
				\end{minipage}
				\begin{minipage}[b]{0.45\hsize}
					\centering
					\begin{tikzpicture}
						\draw[rounded corners] (1,-0.25) -- (1,0.75) arc (180:0:0.25) -- (1.5,-0.25);
						\draw[rounded corners] (0.5,1) -- (0.5,0) -- (1,0);
						\draw (1,0.1) arc (90:270:0.1);
						\draw[->] (1.25,-0.5) -- (1.25,-1);
						\node at (0.5,1.25){$A$};
						\node at (0.75,0.5){$\lambda$};
						\node at (1.75,0.5){$\lambda^\vee$};
						\node[block] at (1.25,-0.5){$r_{\lambda, \lambda^\vee}$};
						\node at (1.25,-1.25){$\lambda \otimes_A \lambda^\vee$};
					\end{tikzpicture}
					\caption{Coevaluation in $\bimod_\calc(A)$}
					\label{graphicalbimodcoev}
				\end{minipage}
			\end{tabular}
		\end{figure}
	\end{proof}
\end{prop}

\begin{rem}
	When $A$ is moreover swap-commutative in the sense of \cite[Definition 5.5]{MR2029790}, each left $A$-module can be regarded as an $A$-bimodule, and the left $A$-modules in $\calc$ form a multitensor subcategory of $\bimod_\calc(A)$. Thus, the multitensor category of left $A$-modules in $\calc$ inherits the pivotal structure of $\bimod_\calc(A)$, and Proposition \ref{bimodpivotalprop} reproduces \cite[Proposition 5.18 and Lemma 5.19]{MR2029790} as its corollaries.
\end{rem}

\begin{prop}
	\label{prop_bimod_crossed}
	Let $\calc$ be a $G$-crossed multitensor category and let $A \in \obj(\frob^G(\calc))$ be symmetric neutral (recall that neutral means $A \in \obj(\calc_e)$, see the sentence right after Definition \ref{definition_crossed}). Then, we can define a $G$-grading on $\bimod_\calc(A)$ by putting $\obj(\bimod_\calc(A)_g) \coloneqq \{ \mu \in \obj(\bimod_\calc(A)) \mid F(\mu) \in \obj(\calc_g) \}$ for $g \in G$, where $F$ denotes the forgetful functor $\bimod_\calc(A) \to \calc$. Combined with the action defined in Proposition \ref{bimodgroupactiondef}, $\bimod_\calc(A)$ becomes a $G$-crossed multitensor category. Another choice of relative tensor products only yields an isomorphic $G$-crossed structure.

	\begin{proof}
		By Propositions \ref{bimodgroupactiondef} and \ref{bimodpivotalprop}, the category $\bimod_\calc(A)$ is equipped with a pivotal $G$-action. By the definition of the $G$-grading, the homogeneous decomposition in $\calc$ gives that in $\bimod_\calc(A)$. For $\lambda \in \obj(\bimod_\calc(A)_g)$ and $\mu \in \obj(\bimod_\calc(A)_h)$, since $A$ is neutral, $\lambda \otimes_A \mu$ is the cokernel of a morphism $\lambda \otimes_\calc A \otimes_\calc \mu \to \lambda \otimes_\calc \mu$ in $\calc_{gh}$, and therefore $\lambda \otimes_A \mu \in \obj(\calc_{gh})$. Finally, ${}^k (\obj(\bimod_\calc(A)_g)) \subset \obj(\bimod_\calc(A)_{kgk^{-1}})$ for $k \in G$ since $\calc$ is $G$-crossed. Thus, $\bimod_\calc(A)$ is $G$-crossed. Another choice of relative tensor products only yields an isomorphic $G$-crossed structure as in the proof of Proposition \ref{bimodpivotalprop}.
	\end{proof}
\end{prop}

\subsection{Categorical definition of equivariant $\alpha$-induction}
\label{subsection_def_alphaind}

In this subsection, we define $G$-equivariant $\alpha$-induction (Definition \ref{def_alpha_ind}).

\begin{prop}
Let $A$ be a neutral $G$-equivariant algebra in a $G$-braided multitensor category $\calc$. Put $m_{A\lambda}^{\mathrm{L}} \coloneqq m_A \otimes \id_\lambda \in \hom(AA\lambda, A\lambda)$ and define $m_{A\lambda}^{\mathrm{R}\pm} \in \hom(A \lambda A,A\lambda)$ to be the morphisms given in Figure \ref{graphicalalphainduction} for every $\lambda \in \obj(\calc)$. Then $\alpha^{G \pm}_A(\lambda) \coloneqq (A\lambda,m^{\mathrm{L}}_{A\lambda},m^{\mathrm{R}\pm}_{A\lambda})$ are $A$-bimodules.

\begin{figure}[htb]
\centering
\begin{tikzpicture} 
\draw[<-] (0,0) -- (0,0.5);
\draw (0.5,1) arc (360:180:0.5);
\draw (0.5,1) -- (0.5,2) -- (1.5,3);
\node[block] at (0.5,1.5){$z_g^A$};
\draw[dashed,cross] (0.5,3)--(1.5,2);
\draw[->] (1.5,2)--(1.5,0);
\draw (-0.5,3) -- (-0.5,1);
\node at (0,-0.25){$A$};
\node at (1.5,-0.25){$\lambda$};
\node at (0.5,3.25){$\lambda$};
\node at (1.5,3.25){$A$};
\node at (-0.5,3.25){$A$};
\node at (-1.75,1.5){$m_{A\lambda}^{\mathrm{R}+} \coloneqq \displaystyle{\sum_{g \in G}}$};
\node at (1,3){$g$};
\begin{scope}[shift={(5,0)}]
\draw[<-] (0,0) -- (0,0.5);
\draw (0.5,1) arc (360:180:0.5);
\draw (0.5,3)--(1.5,2);
\draw[cross] (0.5,1) -- (0.5,2) -- (1.5,3);
\draw[->] (1.5,2)--(1.5,0);
\draw (-0.5,3) -- (-0.5,1);
\node at (0,-0.25){$A$};
\node at (1.5,-0.25){$\lambda$};
\node at (0.5,3.25){$\lambda$};
\node at (1.5,3.25){$A$};
\node at (-0.5,3.25){$A$};
\node at (-1.5,1.5){$m_{A\lambda}^{\mathrm{R}-} \coloneqq$};
\end{scope}
\end{tikzpicture}
\caption{The right module structures of $\alpha_A^{G\pm}(\lambda)$}
\label{graphicalalphainduction}
\end{figure}
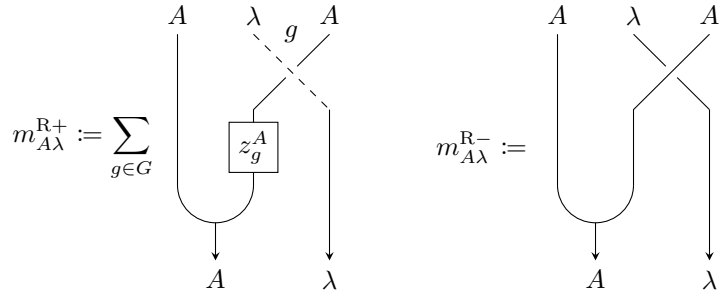

\begin{proof}
We only show the statement for $\alpha_A^{G+}(\lambda)$ because the proof for $\alpha_A^{G-}$ is similar. $m_{A\lambda}^\mathrm{L}$ and $m_{A\lambda}^{\mathrm{R}+}$ are denoted respectively by $m^\mathrm{L}$ and $m^\mathrm{R}$ in this proof. First, $(A\lambda, m^\mathrm{L})$ are left $A$-modules since $A$ is an algebra. Then, the right $A$-modularity of $(A\lambda, m^\mathrm{R})$ follows from the graphical calculation in Figure \ref{graphicalalpharightassoc}. Note that we have by definition two dashed crossings labeled by, say, $g,h \in G$, but only the components labeled by $g = h$ survive. Note also that we used Figure \ref{graphical_product_bfe} at the second equation. The right unit property can also be checked easily by a graphical calculation, and therefore $(A\lambda,m^\mathrm{R})$ is a right $A$-module. Finally, the bimodularity of $\alpha_A^{G+}(\lambda)$ follows from the associativity of $A$.
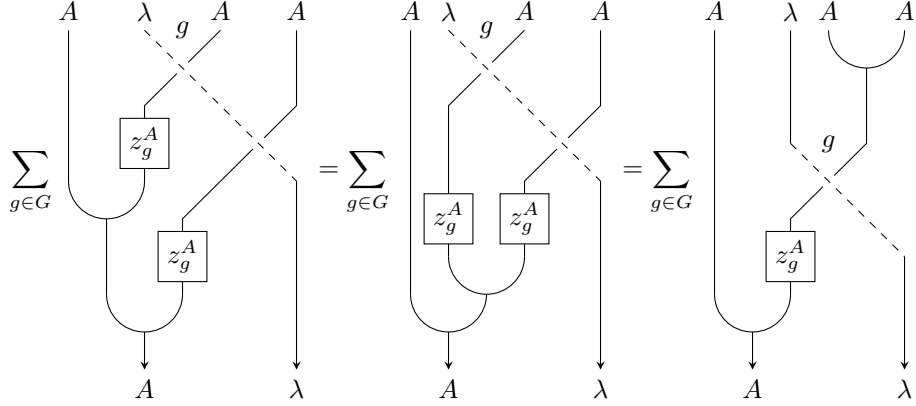
\begin{figure}[htb]
\centering
\begin{tikzpicture}
\draw (0,-0.5) -- (0,0.5);
\draw (0.5,1) arc (360:180:0.5);
\draw (0.5,1) -- (0.5,2) -- (1.5,3);
\node[block] at (0.5,1.5){$z_g^A$};
\draw (1,-0.5) arc (360:180:0.5);
\draw (1,-0.5) -- (1,0.5) -- (2.5,2) -- (2.5,3);
\node[block] at (1,0){$z_g^A$};
\draw[dashed,cross] (0.5,3)--(2.5,1);
\draw[->] (2.5,1)--(2.5,-1.5);
\draw (-0.5,3) -- (-0.5,1);
\draw[->] (0.5,-1) -- (0.5,-1.5);
\node at (0.5,-1.75){$A$};
\node at (2.5,-1.75){$\lambda$};
\node at (0.5,3.25){$\lambda$};
\node at (1.5,3.25){$A$};
\node at (2.5,3.25){$A$};
\node at (-0.5,3.25){$A$};
\node at (-1,1){$\displaystyle{\sum_{g \in G}}$};
\node at (1,3){$g$};
\begin{scope}[shift={(4,0)}]
\draw (0,3) -- (0,-0.5);
\draw (1.5,0) arc (360:180:0.5);
\draw (0.5,0) -- (0.5,2) -- (1.5,3);
\node[block] at (0.5,0.5){$z_g^A$};
\draw (1,-0.5) arc (360:180:0.5);
\draw (1.5,0) -- (1.5,1) -- (2.5,2) -- (2.5,3);
\node[block] at (1.5,0.5){$z_g^A$};
\draw[dashed,cross] (0.5,3)--(2.5,1);
\draw[->] (2.5,1)--(2.5,-1.5);
\draw[->] (0.5,-1) -- (0.5,-1.5);
\node at (0.5,-1.75){$A$};
\node at (2.5,-1.75){$\lambda$};
\node at (0.5,3.25){$\lambda$};
\node at (1.5,3.25){$A$};
\node at (2.5,3.25){$A$};
\node at (0,3.25){$A$};
\node at (-0.75,1){$=\displaystyle{\sum_{g \in G}}$};
\node at (1,3){$g$};
\end{scope}
\begin{scope}[shift={(8,0)}]
\draw (0,3) -- (0,-0.5);
\draw (1,-0.5) arc (360:180:0.5);
\draw (2.5,3) arc (360:180:0.5);
\draw (1,-0.5) -- (1,0.5) -- (2,1.5) -- (2,2.5);
\node[block] at (1,0){$z_g^A$};
\draw[dashed,cross] (1,1.5)--(2.5,0);
\draw (1,1.5) -- (1,3);
\draw[->] (2.5,0)--(2.5,-1.5);
\draw[->] (0.5,-1) -- (0.5,-1.5);
\node at (0.5,-1.75){$A$};
\node at (2.5,-1.75){$\lambda$};
\node at (1,3.25){$\lambda$};
\node at (1.5,3.25){$A$};
\node at (2.5,3.25){$A$};
\node at (0,3.25){$A$};
\node at (-0.75,1){$=\displaystyle{\sum_{g \in G}}$};
\node at (1.5,1.5){$g$};
\end{scope}
\end{tikzpicture}
\caption{The right $A$-modularity of $\alpha^{G+}_A(\lambda)$}
\label{graphicalalpharightassoc}
\end{figure}
\end{proof}
\end{prop}

\begin{rem}
	\label{rem_nojima}
When $A$ is a neutral $G$-equivariant Q-system in a $G$-braided ${}^\ast$-multitensor category, the bimodules $\alpha^{G \pm}_A(\lambda)$ defined above are standard. In particular, when $A \in \rep \cala$ and is localized in an interval $I$ for a local M\"{o}bius covariant net $\cala$ with Haag duality on $\mathbb{R}$ and a group $G$ acting on $\cala$, for every $\lambda \in \obj(G\text{-}\rep \cala)$ we can regard $\alpha^{G \pm}_A(\lambda)$ as standard $A$-bimodules in $\en \cala(I)$. Therefore, they correspond to objects in $\en \calb(I)$, where $\calb$ is the extension of $\cala$ corresponding to $A$, by \cite[Proposition 3.32(ii)]{bklr} via the formula (3.6.3) in \cite[Section 3.6]{bklr}, which gives the equations (4.7) in \cite[Section 4.3]{MR4153896}. Thus in this case our definition coincides with what is considered in \cite{MR4153896}. 

Note that here we do not assume that the action of $G$ on $\cala$ is faithful, see Remarks \ref{remark_cross_not_faithful} and \ref{nonfaithrem}. Such a situation was already considered in \cite{MR4153896}. Namely, a group $G$ acts on $\calb$, and the restriction of this action on $\cala$ can have a nontrivial kernel, by which the quotient of $G$ is denoted by $G'$. For $g_1, g_2 \in G$ with $p(g_1) = p(g_2) = g'$, where $p: G \to G'$ is the quotient map, we have two induced homomorphisms $\alpha^{g_1;+}(\lambda)$ and $\alpha^{g_2;+}(\lambda)$ for $\lambda \in g'\text{-}\rep \cala$. In our framework, they are just the induced homomorphisms of $\lambda \in g_1\text{-}\rep \cala$ and $\lambda \in g_2\text{-}\rep \cala$. Thus, the framework in \cite{MR4153896} is included in ours.
\end{rem}

We regard $\alpha_A^{G \pm}$ as functors by putting $\alpha_A^{G\pm} \coloneqq \id_A \otimes -$ on morphisms.

\begin{prop}
	\label{prop_alpha_is_crossed}
Let $A$ be a neutral special symmetric $G$-equivariant Frobenius algebra in a $G$-braided multitensor category $\calc$. Note that in this setting $\bimod_\calc(A)$ is a $G$-crossed multitensor category by Proposition \ref{prop_bimod_crossed}. Then, $\alpha_A^{G \pm}: \calc \to \bimod_\calc(A)$ can be regarded as $G$-crossed tensor functors.

\begin{proof}
We only show the statement for $\alpha \coloneqq \alpha_A^{G+}$. For $\lambda, \mu \in \obj(\calc)$, the idempotent $e_{\alpha(\lambda),\alpha(\mu)}$ is given in Figure \ref{graphicalalphaindidempotent} by definition. We can split it with the morphisms in Figure \ref{graphicalalphaindsplit}. Indeed, $s_{\alpha(\lambda),\alpha(\mu)} \circ r_{\alpha(\lambda),\alpha(\mu)} = e_{\alpha(\lambda),\alpha(\mu)}$ follows from the graphical calculation in Figure \ref{graphicalalphaindidempotent}, and $r_{\alpha(\lambda),\alpha(\mu)} \circ s_{\alpha(\lambda),\alpha(\mu)} = \id_{A \lambda \mu}$ follows from the graphical calculation in Figure \ref{graphicalalphaindsplitproof}. Note that we used Figures \ref{graphical_equiv_str_hom_calculi} and \ref{graphical_equiv_str_inverse}. Therefore, there exists a canonical isomorphism $\alpha(\lambda) \otimes_A \alpha(\mu) \cong A \lambda \mu$ in $\calc$. Since indeed $s_{\alpha(\lambda),\alpha(\mu)}$ is $A$-bimodular by the graphical calculation in Figure \ref{graphicalalphaindsplitproofs}, where its left $A$-modularity follows from the Frobenius property of $A$, we have a canonical isomorphism $J^\alpha_{\lambda,\mu}: \alpha(\lambda) \otimes_A \alpha(\mu) \cong \alpha(\lambda \mu)$, which makes $\alpha$ into a tensor functor as in the proof of Proposition \ref{bimodgroupactiondef}.

\begin{figure}[H]
	\centering
	\begin{tikzpicture}
		\draw[<-] (0.25,0) -- (0.25,0.25);
		\draw[<-] (1.75,0) -- (1.75,0.25);
		\draw (1.25,2.25) -- (1.25,2.5);
		\draw[fill=white] (1.25,2.5) circle [radius=0.1];
		\draw (0,2.5) -- (0,0.5) arc (180:360:0.25) -- (0.5,1.5) -- (1,2) arc (180:0:0.25) -- (1.5,0.5) arc (180:360:0.25) -- (2,2.5);
		\node[block] at (0.5,1){$z_g$};
		\draw[<-] (2.5,0) -- (2.5,2.5);
		\draw[dashed,cross] (0.5,2.5) -- (0.5,2) --(1,1.5);
		\draw[->] (1,1.5) -- (1,0);
		\node at (0,2.75){$A$};
		\node at (2,2.75){$A$};
		\node at (0.5,2.75){$\lambda$};
		\node at (2.5,2.75){$\mu$};
		\node at (0.25,-0.25){$A$};
		\node at (1,-0.25){$\lambda$};
		\node at (1.75,-0.25){$A$};
		\node at (2.5,-0.25){$\mu$};
		\node at (0.75,2.25){$g$};
		\node at (-0.5,1.25){$\displaystyle{\sum_g}$};
		\node at (3,1.25){$=$};
		\begin{scope}[shift={(4.25,0)}]
			\draw[<-] (0.25,0) -- (0.25,1.25);
			\draw[<-] (1.75,0) -- (1.75,0.25);
			\draw (1.25,2.25) -- (1.25,2.5);
			\draw[fill=white] (1.25,2.5) circle [radius=0.1];
			\draw (0,2.5) -- (0,1.5) arc (180:360:0.25) -- (1,2) arc (180:0:0.25) -- (1.5,0.5) arc (180:360:0.25) -- (2,2.5);
			\node[block] at (0.25,0.5){$z_g$};
			\node[block] at (-0.1,2.1){$z_g^{-1}$};
			\draw[<-] (2.5,0) -- (2.5,2.5);
			\draw[dashed,cross] (0.5,2.5) -- (0.5,2) --(1,1.5);
			\draw[->] (1,1.5) -- (1,0);
			\node at (0,2.75){$A$};
			\node at (2,2.75){$A$};
			\node at (0.5,2.75){$\lambda$};
			\node at (2.5,2.75){$\mu$};
			\node at (0.25,-0.25){$A$};
			\node at (1,-0.25){$\lambda$};
			\node at (1.75,-0.25){$A$};
			\node at (2.5,-0.25){$\mu$};
			\node at (0.75,2.25){$g$};
			\node at (-0.5,1.25){$\displaystyle{\sum_g}$};
			\node at (3,1.25){$=$};
		\end{scope}
		\begin{scope}[shift={(8.5,0)}]
			\draw (0.25,1) -- (0.25,1.25);
			\draw (0,2.5) -- (0,1.5) arc (180:360:0.25) -- (1,2) -- (1,2.5);
			\draw[<->] (0,0) -- (0,0.75) arc (180:0:0.25) -- (1,0.25) -- (1,0);
			\node[block] at (0,0.5){$z_g$};
			\node[block] at (-0.1,2.1){$z_g^{-1}$};
			\draw[<-] (1.5,0) -- (1.5,2.5);
			\draw[dashed,cross] (0.5,2.5) -- (0.5,2) --(1,1.5) -- (1,0.75) -- (0.5,0.25);
			\draw[->] (0.5,0.25) -- (0.5,0);
			\node at (0,2.75){$A$};
			\node at (1,2.75){$A$};
			\node at (0.5,2.75){$\lambda$};
			\node at (1.5,2.75){$\mu$};
			\node at (0,-0.25){$A$};
			\node at (0.5,-0.25){$\lambda$};
			\node at (1,-0.25){$A$};
			\node at (1.5,-0.25){$\mu$};
			\node at (0.75,2.25){$g$};
			\node at (-0.5,1.25){$\displaystyle{\sum_g}$};
			\node at (2,1.25){$=$};
		\end{scope}
		\begin{scope}[shift={(11.75,0)}]
			\draw (0.25,1) -- (0.25,1.25);
			\draw (0,2.5) -- (0,1.5) arc (180:360:0.25) -- (0.5,2) -- (1,2.5);
			\draw[<->] (0,0) -- (0,0.75) arc (180:0:0.25) -- (1,0.25) -- (1,0);
			\node[block] at (0.5,1.75){$z_g$};
			\node[block] at (1,0.5){$z_{g^{-1}}$};
			\draw[<-] (1.5,0) -- (1.5,2.5);
			\draw[dashed,cross] (0.5,2.5) -- (1,2) -- (1,1.5) -- (0.25,0.75);
			\draw[->] (0.25,0.75) -- (0.25,0);
			\node at (0,2.75){$A$};
			\node at (1,2.75){$A$};
			\node at (0.5,2.75){$\lambda$};
			\node at (1.5,2.75){$\mu$};
			\node at (0,-0.25){$A$};
			\node at (0.25,-0.25){$\lambda$};
			\node at (1,-0.25){$A$};
			\node at (1.5,-0.25){$\mu$};
			\node at (1,2.25){$g$};
			\node at (-0.5,1.25){$\displaystyle{\sum_g}$};
		\end{scope}
	\end{tikzpicture}
	\caption{The idempotent $e_{\alpha(\lambda),\alpha(\mu)}$}
	\label{graphicalalphaindidempotent}
\end{figure}
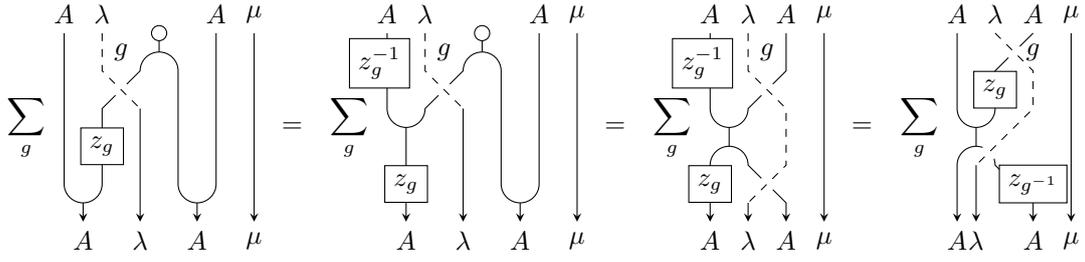

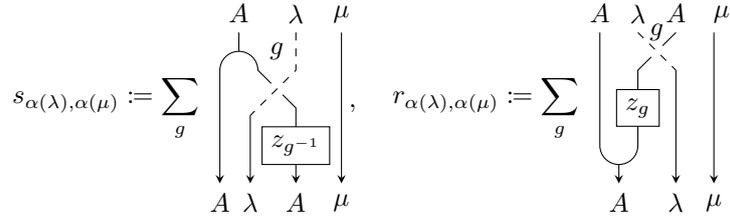
\begin{figure}[H]
	\centering
	\begin{tikzpicture}
		\draw (0,1) -- (0,0.75);
		\draw[<->] (-0.25,-1)-- (-0.25,0.5) arc (180:0:0.25) -- (0.75,0) -- (0.75,-1);
		\draw[dashed,cross] (0.75,1) -- (0.75,0.5) -- (0.15,-0.1);
		\draw[->] (0.15,-0.1) -- (0.15,-1);
		\draw[->] (1.35,1) -- (1.35,-1);
		\node[block] at (0.75,-0.5){$z_{g^{-1}}$};
		\node at (0.5,0.75){$g$};
		\node at (0,1.25){$A$};
		\node at (0.75,1.25){$\lambda$};
		\node at (1.35,1.25){$\mu$};
		\node at (-0.25,-1.25){$A$};
		\node at (0.15,-1.25){$\lambda$};
		\node at (0.75,-1.25){$A$};
		\node at (1.35,-1.25){$\mu$};
		\node at (-1.75,0){$s_{\alpha(\lambda),\alpha(\mu)} \coloneqq \displaystyle{\sum_g}$};
		\node at (1.5,0){,};
		\begin{scope}[shift={(5,0)}]
			\draw[->] (0,-0.75) -- (0,-1);
			\draw (0.75,1) -- (0.25,0.5) -- (0.25,-0.5) arc (360:180:0.25) -- (-0.25,1);
			\draw[dashed,cross] (0.25,1) -- (0.75,0.5);
			\draw[->] (0.75,0.5) -- (0.75,-1);
			\draw[->] (1.25,1) -- (1.25,-1);
			\node[block] at (0.25,0){$z_{g}$};
			\node at (0.5,1){$g$};
			\node at (-0.25,1.25){$A$};
			\node at (0.25,1.25){$\lambda$};
			\node at (1.35,1.25){$\mu$};
			\node at (0.75,1.25){$A$};
			\node at (0.75,-1.25){$\lambda$};
			\node at (0,-1.25){$A$};
			\node at (1.25,-1.25){$\mu$};
			\node at (-1.75,0){$r_{\alpha(\lambda),\alpha(\mu)} \coloneqq \displaystyle{\sum_g}$};	
		\end{scope}
	\end{tikzpicture}
	\caption{A splitting of $e_{\alpha(\lambda),\alpha(\mu)}$}
	\label{graphicalalphaindsplit}
\end{figure}

\begin{figure}[H]
	\centering
	\begin{tikzpicture}
		\draw (0,1) -- (0,0.75);
		\draw (-0.25,-1)-- (-0.25,0.5) arc (180:0:0.25) -- (0.75,0) -- (0.75,-1);
		\draw (1.35,1) -- (1.35,-1);
		\node[block] at (0.75,-0.5){$z_{g^{-1}}^A$};
		\node at (0.5,0.75){$g$};
		\node at (0,1.25){$A$};
		\node at (0.75,1.25){$\lambda$};
		\node at (1.35,1.25){$\mu$};
		\node at (-0.75,-1){$\displaystyle{\sum_g}$};
		\begin{scope}[shift={(0,-2)}]
			\draw[->] (0,-0.75) -- (0,-1);
			\draw (0.75,1) -- (0.25,0.5) -- (0.25,-0.5) arc (360:180:0.25) -- (-0.25,1);
			\draw[->] (0.75,0.5) -- (0.75,-1);
			\draw[->] (1.35,1) -- (1.35,-1);
			\node[block] at (0.25,0){$z_{g}^A$};
			\node at (0.75,-1.25){$\lambda$};
			\node at (0,-1.25){$A$};
			\node at (1.35,-1.25){$\mu$};	
		\end{scope}
		\draw[dashed,cross] (0.75,1) -- (0.75,0.5) -- (0.15,-0.1) -- (0.15,-0.9) -- (0.75,-1.5);
		\node at (1.75,-1){$=$};
		\begin{scope}[shift={(2.5,0)}]
			\draw (0,1) -- (0,0.75);
			\draw (-0.25,-1)-- (-0.25,0.5) arc (180:0:0.25) -- (0.25,0) -- (0.25,-1);
			\draw (1.35,1) -- (1.35,-1);
			\node at (0,1.25){$A$};
			\node at (0.75,1.25){$\lambda$};
			\node at (1.35,1.25){$\mu$};
			\begin{scope}[shift={(0,-2)}]
				\draw[->] (0,-0.75) -- (0,-1);
				\draw (0.25,1) -- (0.25,0.5) -- (0.25,-0.5) arc (360:180:0.25) -- (-0.25,1);
				\draw[->] (0.75,0.5) -- (0.75,-1);
				\draw[->] (1.35,1) -- (1.35,-1);
				\node at (0.75,-1.25){$\lambda$};
				\node at (0,-1.25){$A$};
				\node at (1.35,-1.25){$\mu$};	
			\end{scope}
			\draw (0.75,1) -- (0.75,-1.5);
			\node at (1.75,-1){$=$};
			\node at (2.25,-1){$\id$};
		\end{scope}
	\end{tikzpicture}
	\caption{$r_{\alpha(\lambda),\alpha(\mu)} \circ s_{\alpha(\lambda),\alpha(\mu)} = \id_{A \lambda \mu}$}
	\label{graphicalalphaindsplitproof}
\end{figure}
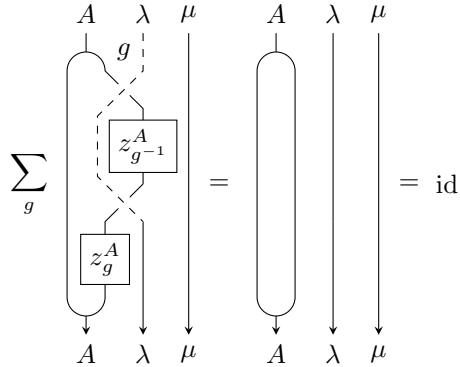

\begin{figure}[H]
	\centering
	\begin{tikzpicture}
		\draw (0,2.5) -- (0,0.75);
		\draw[->] (0.5,0.25) -- (0.75,0) -- (0.75,-1);
		\draw[<-] (-0.25,-1)-- (-0.25,0.5) arc (180:0:0.25) arc (180:360:0.25) -- (0.75,1.5) -- (1.75,2.5);
		\draw[dashed,cross] (0.75,2.5)--(0.75,2) -- (1.25,1.5) -- (1.25,0.5) -- (0.15,-0.1);
		\draw[dashed,cross] (1.25,2.5)--(1.75,2);
		\draw[->] (0.15,-0.1) -- (0.15,-1);
		\draw[->] (1.75,2) -- (1.75,-1);
		\node[block] at (0.75,-0.5){$z_{g^{-1}}$};
		\node[block] at (0.75,1){$z_{gh}$};
		\node at (0.5,1.75){$g$};
		\node at (2,2.25){$h$};
		\node at (0,2.75){$A$};
		\node at (1.75,2.75){$A$};
		\node at (0.75,2.75){$\lambda$};
		\node at (1.25,2.75){$\mu$};
		\node at (-0.25,-1.25){$A$};
		\node at (0.15,-1.25){$\lambda$};
		\node at (0.75,-1.25){$A$};
		\node at (1.75,-1.25){$\mu$};
		\node at (-0.75,0.75){$\displaystyle{\sum_{g,h}}$};
		\node at (2.25,0.75){$=$};
		\begin{scope}[shift={(3.75,0)}]
			\draw (0,2.5) -- (0,0.75);
			\draw[->] (0.5,0.25) -- (0.5,0) -- (1.25,-1);
			\draw[<-] (-0.25,-1)-- (-0.25,0.5) arc (180:0:0.25) arc (180:360:0.25) -- (0.75,1.5) -- (1.75,2.5);
			\draw[dashed,cross,->] (0.75,2.5)--(0.75,1.5) -- (1.25,1.25) -- (1.25,-0.5) -- (0.75,-1);
			\draw[dashed,cross] (1.25,2.5)--(1.75,2);
			\draw[->] (1.75,2) -- (1.75,-1);
			\node[block] at (0.5,-0.25){$z_{g}^{-1}$};
			\node[block] at (0.75,1){$z_{g}$};
			\node[block] at (1.25,1.75){$z_{h}$};
			\node at (0.5,1.75){$g$};
			\node at (2,2.25){$h$};
			\node at (0,2.75){$A$};
			\node at (1.75,2.75){$A$};
			\node at (0.75,2.75){$\lambda$};
			\node at (1.25,2.75){$\mu$};
			\node at (-0.25,-1.25){$A$};
			\node at (0.75,-1.25){$\lambda$};
			\node at (1.25,-1.25){$A$};
			\node at (1.75,-1.25){$\mu$};
			\node at (-0.75,0.75){$\displaystyle{\sum_{g,h}}$};
			\node at (2.25,0.75){$=$};
		\end{scope}
		\begin{scope}[shift={(7.5,0)}]
			\draw (0,2.5) -- (0,0.75);
			\draw[->] (0.5,-0.25) -- (1.25,-1);
			\draw[<-] (-0.25,-1)-- (-0.25,0.5) arc (180:0:0.25) -- (0.25,0) arc (180:360:0.25) -- (0.75,1.5) -- (1.75,2.5);
			\draw[dashed,cross,->] (0.75,2.5)--(0.75,1.5) -- (1.25,1.25) -- (1.25,-0.5) -- (0.75,-1);
			\draw[dashed,cross] (1.25,2.5)--(1.75,2);
			\draw[->] (1.75,2) -- (1.75,-1);
			\node[block] at (0.25,0.25){$z_{g}^{-1}$};
			\node[block] at (1.25,1.75){$z_{h}$};
			\node at (0.5,1.75){$g$};
			\node at (2,2.25){$h$};
			\node at (0,2.75){$A$};
			\node at (1.75,2.75){$A$};
			\node at (0.75,2.75){$\lambda$};
			\node at (1.25,2.75){$\mu$};
			\node at (-0.25,-1.25){$A$};
			\node at (0.75,-1.25){$\lambda$};
			\node at (1.25,-1.25){$A$};
			\node at (1.75,-1.25){$\mu$};
			\node at (-0.75,0.75){$\displaystyle{\sum_{g,h}}$};
			\node at (2.25,0.75){$=$};
		\end{scope}
		\begin{scope}[shift={(11.25,0)}]
			\draw (0,2.5) -- (0,2.25);
			\draw[->] (1,-0.75) -- (1,-1);
			\draw[<-] (-0.25,-1)-- (-0.25,2) arc (180:0:0.25) -- (0.75,1.5) -- (0.75,-0.5) arc (180:360:0.25) -- (2,0.5) -- (2,2.5);
			\draw[dashed,cross,->] (0.75,2.5) -- (0.75,2) --(0.25,1.5) -- (0.25,-1);
			\draw[dashed,cross] (1.5,2.5)-- (1.5,0.5) --(2,0);
			\draw[->] (2,0) -- (2,-1);
			\node[block] at (0.75,1){$z_{g^{-1}}$};
			\node[block] at (1.25,-0.25){$z_{h}$};
			\node at (0.5,2.25){$g$};
			\node at (1.75,0.75){$h$};
			\node at (0,2.75){$A$};
			\node at (2,2.75){$A$};
			\node at (0.75,2.75){$\lambda$};
			\node at (1.5,2.75){$\mu$};
			\node at (-0.25,-1.25){$A$};
			\node at (0.25,-1.25){$\lambda$};
			\node at (1,-1.25){$A$};
			\node at (2,-1.25){$\mu$};
			\node at (-0.75,0.75){$\displaystyle{\sum_{g,h}}$};
		\end{scope}
	\end{tikzpicture}
	\caption{The right $A$-modularity of $s_{\alpha(\lambda),\alpha(\mu)}$}
	\label{graphicalalphaindsplitproofs}
\end{figure}
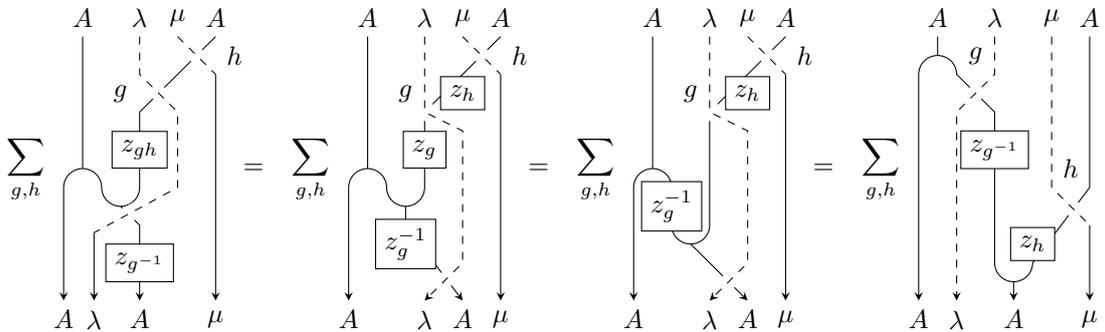

Put $(\tilde{\eta}^{\alpha}_g)_\lambda \coloneqq z_g \otimes \id_\lambda$ for every $g \in G$, which is a left $A$-module isomorphism ${}^g \alpha(\lambda) \cong \alpha({}^g \lambda)$ since $z_g$ is an algebra homomorphism. It is indeed also right $A$-modular by the graphical calculation in Figure \ref{graphicalalphaindiscrossed}, where we used $z_h {}^{hg^{-1}} z_g^{-1} = z_h {}^h z_{g^{-1}} = z_{hg^{-1}}$ for $g,h \in G$. The pair $(\alpha, \tilde{\eta}^\alpha)$ is a tensor $G$-functor by $z_{gh} = z_g {}^g z_h$. It preserves the gradings since $A$ is neutral.
\begin{figure}[htb]
	\centering
	\begin{tikzpicture}
		\draw[<-] (0,0) -- (0,0.5);
		\draw (0.5,1) arc (360:180:0.5);
		\draw (0.5,1) -- (0.5,2) -- (1,2.5) -- (1,3.5);
		\node[block] at (0.5,1.5){${}^g z_h$};
		\draw[dashed,cross] (0.5,2.5)--(1,2);
		\draw[->] (1,2)--(1,0);
		\draw (-0.5,3.5) -- (-0.5,1);
		\draw (0.5,3.5) -- (0.5,2.5);
		\node at (0,-0.25){$A$};
		\node at (1,-0.25){$\lambda$};
		\node at (0.5,3.75){$\lambda$};
		\node at (1,3.75){$A$};
		\node at (-0.5,3.75){$A$};
		\node at (-1,1.75){$\displaystyle{\sum_{g}}$};
		\node[block] at (-0.5,3){$z_g^{-1}$};
		\node[block] at (1,3){$z_g^{-1}$};
		\node at (1.5,2.25){$ghg^{-1}$};
		\begin{scope}[shift={(4,0)}]
			\draw[<-] (0,0) -- (0,1.5);
			\draw (0.5,2) arc (360:180:0.5);
			\draw (0.5,2) -- (0.5,3) -- (1,3.5);
			\node[block] at (0.25,2.5){$z_{g h g^{-1}}$};
			\draw[dashed,cross] (0.5,3.5)--(1,3);
			\draw[->] (1,3)--(1,0);
			\draw (-0.5,3.5) -- (-0.5,2);
			\node at (0,-0.25){$A$};
			\node at (1,-0.25){$\lambda$};
			\node at (0.5,3.75){$\lambda$};
			\node at (1,3.75){$A$};
			\node at (-0.5,3.75){$A$};
			\node at (-1,1.75){$\displaystyle{\sum_{g}}$};
			\node[block] at (0,0.75){$z_g^{-1}$};
			\node at (1.5,3.25){$ghg^{-1}$};
			\node at (-1.75,1.75){$=$};
		\end{scope}
	\end{tikzpicture}
	\caption{The right $A$-modularity of $(\tilde{\eta}^{\alpha}_g)_\lambda$}
	\label{graphicalalphaindiscrossed}
\end{figure}
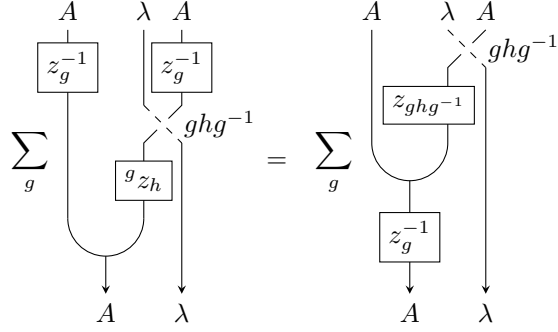
\end{proof}
\end{prop}

\begin{defi}
	\label{def_alpha_ind}
Let $A$ be a neutral special symmetric $G$-equivariant Frobenius algebra in a $G$-braided multitensor category $\calc$. We refer to the $G$-crossed tensor functors $\alpha^{G\pm}_A: \calc \to \bimod_\calc(A)$ as the \emph{$G$-equivariant $\alpha$-induction functors} for $A$. 
\end{defi}

\section{Equivariant $\alpha$-induction Frobenius algebras}
\label{section_equiv_frob}

In this section, we construct equivariant $\alpha$-induction Frobenius algebras (Theorem \ref{mainthm1}), which generalizes Rehren's work \cite{MR1754521} (see also \cite[Definition 4.17]{MR3424476}) and is one of our main theorems. To give their ambient categories, we introduce the \emph{neutral double} construction (Definition \ref{def_neutral_double}) for $G$-braided multitensor categories.

\subsection{The neutral double of a $G$-braided tensor category}

In this subsection, we introduce the neutral double construction (Definition \ref{def_neutral_double}). For this, we begin with crossed products.

\begin{prop}
	\label{prop_semidir_def}
Let $G$ be a group. Let $\calc$ be a $G$-graded multitensor category and let $\cald$ be a multitensor category with an action $\gamma^\cald$ of $G$. Then, we can define a bilinear bifunctor $\otimes: (\calc \boxtimes \cald) \times (\calc \boxtimes \cald) \to \calc \boxtimes \cald$ by putting $(c_1 \boxtimes d_1) \otimes (c_2 \boxtimes d_2) \coloneqq c_1 c_2 \boxtimes {}^{\overline{c_2}} d_1 d_2$ for $c_1 \in \obj(\calc)$, $c_2 \in \homog(\calc)$ and $d_1,d_2 \in \obj(\cald)$ and make $\calc \boxtimes \cald$ into a multitensor category, which is denoted by $\calc \ltimes_{\gamma^\cald} \cald$ or simply $\calc \ltimes \cald$. Another choice of the Deligne tensor product $\calc \boxtimes \cald$ only yields a strictly isomorphic multitensor category. If $\calc$ and $\cald$ are tensor (resp. multifusion, resp. fusion) categories, then so is $\calc \ltimes \cald$.

\begin{proof}
Since $\calc \times \cald \times \calc \times \cald$ has a $G$-grading $\bigoplus_{g \in G} \calc \times \cald \times \calc_g \times \cald$ and $\otimes$ is linear exact in each variable (see \cite[Proposition 4.2.1]{egno}), $\otimes$ can be extended to a linear exact functor on $\calc \times \cald \times \calc \times \cald$ and therefore on $(\calc \boxtimes \cald) \times (\calc \boxtimes \cald)$. Since
\begin{align*}
(( c_1 \boxtimes d_1) \otimes ( c_2 \boxtimes d_2 )) \otimes ( c_3 \boxtimes d_3 ) &= (c_1 c_2) c_3 \boxtimes {}^{\overline{c_3}}({}^{\overline{c_2}} d_1 d_2) d_3 \\
( c_1 \boxtimes d_1) \otimes (( c_2 \boxtimes d_2 ) \otimes ( c_3 \boxtimes d_3 )) &= c_1 (c_2 c_3) \boxtimes {}^{\overline{c_2c_3}} d_1 ({}^{\overline{c_3}}d_2 d_3)
\end{align*}
for $c_1 \in \obj(\calc)$, $c_2,c_3 \in \homog(\calc)$ and $d_1,d_2,d_3 \in \obj(\cald)$, we can define a natural isomorphism $a: (- \otimes -) \otimes - \cong - \otimes (- \otimes -)$ to be the extension of the canonical isomorphism given by coherence (Theorem \ref{theorem_equivcoh}), which satisfies the pentagon axiom by coherence. We put $\mathbf{1} \coloneqq \mathbf{1}_\calc \boxtimes \mathbf{1}_\cald$ and similarly define left and right unit isomorphisms by coherence. Thus, $\calc \ltimes \cald$ turns into a monoidal category. Another choice of $\calc \boxtimes \cald$ only yields a strictly isomorphic multitensor category by universality since we define the tensor product by extension.  

Next, we show that $\calc \ltimes \cald$ is rigid and therefore a multitensor category. As in the proof of \cite[Proposition 5.17]{MR1106898}, it is enough to check the rigidity of the objects of $\calc_g \times \cald$ for $g \in G$. We show that a left dual of $c \boxtimes d$ for $c \in \homog(\calc)$ and $d \in \obj(\cald)$ is given by $c^\vee \boxtimes {}^{c} d^\vee$. Indeed, since $(c^\vee \boxtimes {}^{c} d^\vee) (c \boxtimes d) = c^\vee c \boxtimes d^\vee d$ and $(c \boxtimes d) (c^\vee \boxtimes {}^{c} d^\vee) = c c^\vee \boxtimes {}^c (d d^\vee)$, where we suppressed some isomorphisms by Theorem \ref{theorem_equivcoh}, we can define evaluation and coevaluation maps by putting $\mathrm{ev}_{c \boxtimes d} \coloneqq \mathrm{ev}_c \boxtimes \mathrm{ev}_d$ and $\mathrm{coev} \coloneqq \mathrm{coev}_c \boxtimes {}^c \mathrm{coev}_d$. The conjugate equations follow from those for $c$ and $d$ and therefore $(c \boxtimes d)^\vee = c^\vee \boxtimes {}^{c} d^\vee$. We can also show that ${}^\vee (c \boxtimes d) = {}^\vee c \boxtimes {}^{c \vee} d$ and therefore $\calc \ltimes \cald$ is rigid.
\end{proof}
\end{prop}

We call $\calc \ltimes \cald$ the \emph{crossed product} of $\calc$ and $\cald$. Let us see that $\calc \ltimes \cald$ inherits structures on $\calc$ and $\cald$.

\begin{prop}
\label{semidirpivotalprop}
Let $G$ be a group, let $\calc$ be a $G$-graded pivotal multitensor category, and let $\cald$ be a pivotal multitensor category with a pivotal action of $G$. Then, the pivotal structures on $\calc$ and $\cald$ induce a pivotal structure $\delta^{\calc \ltimes \cald}$ on $\calc \ltimes \cald$. Another choice of $\calc \boxtimes \cald$ only yields an isomorphic pivotal structure. 

\begin{proof}
Since $(c \boxtimes d)^{\vee \vee} = c^{\vee \vee} \boxtimes d^{\vee \vee}$ for $c \in \obj(\calc)$ and $d \in \obj(\cald)$ by the proof of Proposition \ref{prop_semidir_def}, we can define a natural isomorphism $\delta^{\calc \ltimes \cald}: \id \cong (-)^{\vee \vee}$ by putting $\delta_{c \boxtimes d}^{\calc \ltimes \cald} \coloneqq \delta^\calc_c \boxtimes \delta^\cald_d$. Since the $G$-action on $\cald$ is pivotal, we have
\begin{align*}
	\delta_{(c_1 \boxtimes d_1) \otimes (c_2 \boxtimes d_2)}^{\calc \ltimes \cald} &= \delta_{c_1 c_2 \boxtimes{}^{\overline{c_2}}d_1d_2}^{\calc \ltimes \cald} = \delta^\calc_{c_1 c_2} \boxtimes \delta^\cald_{{}^{\overline{c_2}}d_1d_2} = (\delta^\calc_{c_1} \otimes \delta^\calc_{c_2}) \boxtimes ({}^{\overline{c_2}} \delta^\cald_{d_1} \otimes \delta^\cald_{d_2}) \\
	&= (\delta^\calc_{c_1} \boxtimes \delta^\cald_{d_1}) \otimes  (\delta^\calc_{c_2} \boxtimes \delta^\cald_{d_2}) = \delta^{\calc \ltimes \cald}_{c_1 \boxtimes d_1} \otimes \delta^{\calc \ltimes \cald}_{c_2 \boxtimes d_2}
\end{align*}
for $c_1 \in \obj(\calc)$, $c_2 \in \homog(\calc)$ and $d_1,d_2 \in \obj(\cald)$ by the monoidality of $\delta^\calc$ and $\delta^\cald$. Thus, $\delta^{\calc \ltimes \cald}$ is monoidal and therefore a pivotal structure on $\calc \ltimes \cald$. Another choice of $\calc \boxtimes \cald$ only yields an isomorphic pivotal structure since we define the pivotal structure by extension.
\end{proof}
\end{prop}

\begin{rem}
	When $\calc$ and $\cald$ in the proposition above are moreover tensor categories, we have $\dim_{\delta^{\calc \ltimes \cald}}(c \boxtimes d) = \dim_{\delta^\calc}(c) \dim_{\delta^\cald}(d)$ for $c \in \obj(\calc)$ and $d \in \obj(\cald)$ by definition. In particular, when $\calc$ and $\cald$ are split (i.e. $\en(\lambda) \cong k$ for any simple object $\lambda$, see \cite[Section 4.16]{egno}) semisimple spherical tensor categories, since every simple object of $\calc \ltimes \cald$ is of the form $c \boxtimes d$ for some simple objects $c$ and $d$ respectively of $\calc$ and $\cald$ by \cite[Theorem 27]{MR3108080}, we have $\dim(\lambda) = \dim(\lambda^\vee)$ for every simple object $\lambda = c \boxtimes d$ of $\calc \ltimes \cald$ by $(c \boxtimes d)^\vee = c^\vee \boxtimes {}^c d^\vee$, see the proof of Proposition \ref{prop_semidir_def}. Therefore, in this case, $\calc \ltimes \cald$ is again a split semisimple spherical tensor category.
\end{rem}

\begin{prop}
\label{prop_semidir_action}
Let $\calc$ be a $G$-crossed multitensor category and let $\cald$ be a multitensor category with an action $\gamma^\cald$ of $G$. Then, an action of $G$ is naturally induced on $\calc \ltimes \cald$. Another choice of $\calc \boxtimes \cald$ only yields a $G$-tensor isomorphic one. When $\cald$ is moreover pivotal and $\gamma^\calc$ and $\gamma^\cald$ are pivotal actions, the induced action on $\calc \ltimes \cald$ is again pivotal.

\begin{proof}
Define a linear exact functor $\gamma(g)$ for every $g \in G$ by putting $\gamma(g) (c \boxtimes d) \coloneqq {}^g c \boxtimes {}^g d$ for $c \in \obj(\calc)$ and $d \in \obj(\cald)$. Since
\begin{align*}
\gamma (g) (c_1 \boxtimes d_1) \otimes \gamma(g) (c_2 \boxtimes d_2) &= {}^g c_1 {}^g c_2 \boxtimes {}^{g \overline{c_2} g^{-1}} ({}^g d_1) {}^g d_2 \\
\gamma (g) ((c_1 \boxtimes d_1) \otimes (c_2 \boxtimes d_2)) &= {}^g (c_1 c_2) \boxtimes {}^{g} ({}^{\overline{c_2}} d_1 d_2)
\end{align*}
for $c_1 \in \obj(\calc)$, $c_2 \in \homog(\calc)$ and $d_1,d_2 \in \obj(\cald)$ by $\partial {}^g c_2 = g \partial c_2 g^{-1}$, we can define a canonical natural isomorphism $J^{\gamma(g)}$ by coherence. We similarly have canonical natural isomorphisms $\chi_{g,h}^\gamma$ for $g,h \in G$ and $\iota^\gamma$ by coherence. Thus, $\gamma$ is an action of $G$. Another choice of $\calc \boxtimes \cald$ only yields a $G$-tensor isomorphic action since we define the action by extension. The last statement follows from the definitions of $\gamma$ and the pivotal structure $\delta^{\calc \ltimes \cald}$, see the proof of Proposition \ref{semidirpivotalprop}.
\end{proof}
\end{prop}

\begin{prop}
	\label{prop_semidir_crossed}
	Let $\calc$ and $\cald$ be split semisimple $G$-crossed multitensor categories. Then, a $G$-crossed structure is naturally induced on $\calc \ltimes \cald$. Another choice of $\calc \boxtimes \cald$ only yields a $G$-crossed isomorphic one.

	\begin{proof}
		We already know that $\calc \ltimes \cald$ is a pivotal category with a pivotal action of $G$ by Proposition \ref{prop_semidir_action}.  Since $\calc \boxtimes \cald = \bigoplus_{g} \bigoplus_h \calc_h \otimes \cald_{h^{-1} g}$ follows from \cite[Theorem 27]{MR3108080} by the assumption of split semisimplicity, we put $(\calc \ltimes \cald)_g \coloneqq \bigoplus_{h} \calc_h \boxtimes \cald_{h^{-1} g}$ for $g \in G$. Then, we have $(\calc \ltimes \cald)_g \otimes (\calc \ltimes \cald)_h \subset (\calc \ltimes \cald)_{g h}$ for $g,h \in G$ since $\partial ((c_1 \boxtimes d_1) (c_2 \boxtimes d_2)) = \partial (c_1 c_2 \boxtimes {}^{\overline{c_2}} d_1 d_2) = \partial c_1 \partial d_1 \partial c_2 \partial d_2$ for $c_1,c_2,d_1,d_2 \in \homog (\calc)$ by definition. Thus, $\calc \ltimes \cald$ is a $G$-graded multitensor category. Moreover, $\calc \ltimes \cald$ is $G$-crossed since $\partial^{\calc \ltimes \cald}({}^g c \boxtimes {}^g d) = g \partial c g^{-1} g \partial d g^{-1} = g \partial c \partial d g^{-1}$ for $c \in \homog(\calc), d \in \homog(\cald)$ and $g \in G$. Another choice of $\calc \boxtimes \cald$ only yields an isomorphic $G$-crossed structure since changing $\calc \boxtimes \cald$ preserves $\calc_h \boxtimes \cald_{h^{-1}g}$ and therefore yields an isomorphic $G$-grading.
	\end{proof}
\end{prop}

Next, we show that if $\calc$ and $\cald$ are moreover $G$-braided, then we can obtain a $G$-braiding on a subcategory of $\calc \ltimes \cald$.

\begin{thm}
	For a $G$-braided multitensor category $\calc$, define $\calc^{\mathrm{rev}}$ to be $\calc$ as an abelian category with an action of $G$. Put $\calc^{\mathrm{rev}}_g \coloneqq \calc_{g^{-1}}$ for $g \in G$. Then, by putting $\lambda \otimes_{\calc^{\mathrm{rev}}} \mu \coloneqq {}^{\partial_\calc \mu} \lambda \otimes_\calc \mu$ for $\lambda \in \obj(\calc^{\mathrm{rev}})$ and $\mu \in \homog(\calc^{\mathrm{rev}}) = \homog(\calc)$, and putting $\mathbf{1}^{\calc^{\mathrm{rev}}} \coloneqq \mathbf{1}^\calc$, we obtain a multitensor structure on $\calc^{\mathrm{rev}}$. Moreover, by putting $b^{\calc^{\mathrm{rev}}}_{\lambda, \mu} \coloneqq b^{\calc-}_{\lambda, \mu}$ (see Definition \ref{definition_crossed}), we obtain a $G$-braiding on $\calc^{\mathrm{rev}}$. We call $\calc^{\mathrm{rev}}$ the \emph{reverse} of $\calc$.

	\begin{proof}
		For $\lambda \in \obj(\calc)$ and $\mu, \nu \in \homog(\calc)$, we have $(\lambda \otimes_{\calc^{\mathrm{rev}}} \mu ) \otimes_{\calc^{\mathrm{rev}}} \nu = {}^\nu({}^\mu \lambda \mu) \nu$ and $\lambda \otimes_{\calc^{\mathrm{rev}}} (\mu \otimes_{\calc^{\mathrm{rev}}} \nu) = {}^{\nu \mu} \lambda ({}^\nu \mu \nu)$, where the objects on the right-hand sides are regarded as objects of $\calc$, since $\partial_\calc ({}^\nu \mu \nu) = \partial_\calc \nu \partial_\calc \mu$. Thus, we obtain an associative bilinear monoidal product by coherence and universality. We can similarly see that $\mathbf{1}^{\calc^{\mathrm{rev}}}$ is indeed a unit. 

		Next, we show that $\calc^{\mathrm{rev}}$ is rigid and therefore a multitensor category. For $\lambda \in \homog(\calc^{\mathrm{rev}})$, put $\lambda^\vee \coloneqq {}^{\overline{\lambda}} \lambda^\vee_\calc$, where $\lambda^\vee_\calc$ denotes the left dual of $\lambda$ in $\calc$. Then, $\mathrm{ev}_\lambda \coloneqq \mathrm{ev}^\calc_\lambda : \lambda^\vee \otimes_{\calc^{\mathrm{rev}}} \lambda =  \lambda_\calc^\vee \lambda \to \mathbf{1}$ and $\mathrm{coev}_\lambda \coloneqq {}^{\overline{\lambda}} \mathrm{coev}^\calc_\lambda : \mathbf{1} \to \lambda \otimes_{\calc^{\mathrm{rev}}} \lambda^\vee = {}^{\overline{\lambda}} (\lambda \lambda_\calc^\vee)$ give desired duality, where $\mathrm{ev}^\calc_\lambda$ and $\mathrm{coev}^\calc_\lambda$ denote evaluation and coevaluation maps of $\lambda$ in $\calc$. The proof of right duality is similar. Moreover, since $\lambda^{\vee \vee} = (\lambda_\calc^\vee)_\calc^\vee$, we can put $\delta^{\calc^{\mathrm{rev}}} \coloneqq \delta^\calc$, which indeed defines a pivotal structure on $\calc^{\mathrm{rev}}$ since the action of $G$ on $\calc$ is pivotal.
		
		The relation $\partial_{\calc^{\mathrm{rev}}}(\mu \otimes_{\calc^{\mathrm{rev}}} \nu)= \partial_\calc ({}^\nu \mu \nu)^{-1} = (\partial_\calc \nu \partial_\calc \mu)^{-1} = \partial_{\calc^{\mathrm{rev}}} \mu \partial_{\calc^{\mathrm{rev}}} \nu$ shows that $\calc^{\mathrm{rev}}$ is a $G$-graded multitensor category and therefore a $G$-crossed category since $\partial_{\calc^{\mathrm{rev}}} {}^g \mu = (g \partial_\calc \mu g^{-1})^{-1} = g \partial_{\calc^{\mathrm{rev}}} \mu g^{-1}$ for $g \in G$. Finally, it follows from the axioms for $b^{\calc}$ that $b^{\calc^{\mathrm{rev}}}$ is a $G$-braiding.
	\end{proof}
\end{thm}

Let $\calc_1$ (resp. $\calc_2$) be a split semisimple $G_1$-braided (resp. $G_2$-braided) multitensor category. Then, $\calc_1 \boxtimes \calc_2$ naturally turns into a $G_1 \times G_2$-braided multitensor category.

\begin{defi}
	\label{def_neutral_double}
	Let $\calc$ and $\cald$ be split semisimple $G$-crossed multitensor categories. Let $D(\calc, \cald)$ denote the full multitensor subcategory of $\calc \boxtimes \cald^{\mathrm{rev}}$ with $G$-grading restricted through the diagonal embedding $G \subset G \times G$. We regard $D(\calc,\cald)$ as a $G$-braided multitensor category by restricting the $G \times G$-braiding of $\calc \boxtimes \cald^{\mathrm{rev}}$. We refer to the $G$-braided multitensor category $D(\calc) \coloneqq D(\calc, \calc)$ as the \emph{neutral double} of $\calc$. 
\end{defi}

$D$ is for Double. Note that $D(\calc, \cald) = (\calc \ltimes \cald)_e$ as a pivotal multitensor category with an action of $G$. Note also that when $G$ is trivial, $D(\calc)$ is equal to $\calc \boxtimes \calc^{\mathrm{rev}}$ as an (ordinary) braided multitensor category, where $\calc^{\mathrm{rev}}$ is the reverse of the (ordinary) braided multitensor category $\calc$.

Finally, we give an application of these constructions to algebraic quantum field theory. Let $\cala_1$ and $\cala_2$ be local M\"{o}bius covariant nets. Then, we can obtain a local M\"{o}bius covariant net $\cala_1 \otimes \cala_2$ on the (1+1)-dimensional Minkowski space $\mathbb{R}^{1,1}$ in the sense of \cite[Section 2]{MR2029950} by putting $\cala_1 \otimes \cala_2 (I \times J) \coloneqq \cala_1(I) \otimes \cala_2(J)$ for an interval $I$ in the positive half line in $\mathbb{R}$ and an interval $J$ in the negative half line in $\mathbb{R}$. Note that $\cala_1 \otimes \cala_2$ does not satisfy Haag duality on $\mathbb{R}^{1,1}$ in the sense of \cite[Subsection 2.1]{MR2029950}.

For a local M\"{o}bius covariant net $\cala$ on $\mathbb{R}^{1,1}$, we can define its automorphism group $\aut \cala$ and the notion of a group action on $\cala$ as in the case of local M\"{o}bius covariant nets on $S^1$. 

\begin{defi}
	Let $\cala$ be an irreducible local M\"{o}bius covariant net on $\mathbb{R}^{1,1}$ with an action $\beta$ of a group $G$. For $g \in G$, an endomorphism $\lambda$ of $\cala_\infty \coloneqq \bigcup_{O \in \mathcal{DC}} \cala(O)$, where $\mathcal{DC}$ denotes the set of double cones in $\mathbb{R}^{1,1}$, is \emph{$g$-localized} in $O \in \mathcal{DC}$ if $\lambda|_{\cala(O'_{\mathrm{L}})} = \id_{\cala(O'_{\mathrm{L}})}$ and $\lambda|_{\cala(O'_{\mathrm{R}})} = \ad \beta(g)|_{\cala(O'_{\mathrm{R}})}$, where $O'_{\mathrm{L}}$ (resp. $O'_{\mathrm{R}}$) denotes the left (resp. right) connected component of the causal complement $O'$ of $O$. A $g$-localized endomorphism $\lambda$ is a \emph{$g$-twisted DHR endomorphism} of $\cala$ if for any $\tilde{O} \in \mathcal{DC}$, there exists a unitary $u \in \cala_\infty$ such that $\ad u \circ \lambda$ is localized in $\tilde{O}$. The ${}^\ast$-category of rigid $g$-twisted DHR endomorphisms is denoted by $g\text{-}\rep \cala$, and an object of the ${}^\ast$-tensor category $G\text{-}\rep \cala \coloneqq \bigoplus_{g \in G} g\text{-}\rep \cala$ is called a \emph{$G$-twisted DHR endomorphism} of $\cala$.
\end{defi}

The following proposition gives a physical interpretation of $\calc \boxtimes \cald^{\mathrm{rev}}$ and $D(\calc, \cald)$.

\begin{prop}
	Let $\cala_1$ (resp. $\cala_2$) be a completely rational irreducible local M\"{o}bius covariant net on $S^1$ with an action of a group $G_1$ (resp. $G_2$). Then, $(G_1 \times G_2)\text{-}\rep (\cala_1 \otimes \cala_2) \simeq (G_1\text{-}\rep \cala_1) \boxtimes (G_2\text{-}\rep \cala_2)^{\mathrm{rev}}$ as $G_1 \times G_2$-crossed ${}^\ast$-tensor categories, where $G_1 \times G_2$ acts on $\cala_1 \otimes \cala_2$ by the tensor product representation. In particular, when $G_1 = G_2 = G$, we have an equivalence $G\text{-}\rep (\cala_1 \otimes \cala_2) \simeq D(G\text{-}\rep \cala_1, G\text{-}\rep \cala_2)$, where $G$ acts on $\cala_1 \otimes \cala_2$ as the diagonal subgroup of $G \times G$.

	\begin{proof}
		For $\lambda \boxtimes \mu \in \obj((g\text{-}\rep \cala_1) \boxtimes (G_2\text{-}\rep \cala_2)^{\mathrm{rev}}_{h}) = \obj((g\text{-}\rep \cala_1) \boxtimes (h^{-1}\text{-}\rep \cala_2))$, we can define $\lambda \tilde{\boxtimes} \mu \in \obj((g, h)\text{-}\rep (\cala_1 \otimes \cala_2))$ by putting $\lambda \tilde{\boxtimes} \mu(a \otimes b) \coloneqq \lambda(a) \otimes {}^{h} (\mu(b))$ for $a \in \cala_1(I)$ and $b \in \cala_2(J)$ and $I \in \mathcal{I}_+, J \in \mathcal{I}_-$, where $\mathcal{I}_+$ (resp. $\mathcal{I}_-$) denotes the set of intervals in the positive (resp. negative) half line. This assignment defines a fully faithful ${}^\ast$-functor $F: D(G\text{-}\rep \cala_1, G\text{-}\rep \cala_2) \to G\text{-}\rep (\cala_1 \otimes \cala_2)$ by universality. This functor $F$ is a strict tensor functor since
		\begin{align*}
			(\lambda_1 \tilde{\boxtimes} \mu_1) \circ (\lambda_2 \tilde{\boxtimes} \mu_2) (a \otimes b) &= \lambda_1 \lambda_2(a) \otimes {}^{h_1} (\mu_1 ({}^{h_2}(\mu_2(b)))) = \lambda_1 \lambda_2(a) \otimes {}^{h_1 h_2} (({}^{h_2^{-1}} \mu_1) \mu_2(b)) \\
			&= ((\lambda_1 \boxtimes \mu_1) \otimes (\lambda_2 \boxtimes \mu_2)) (a \otimes b)
		\end{align*}
		for $\lambda_i \boxtimes \mu_i \in \obj((g_i\text{-}\rep \cala_1) \boxtimes (h_i^{-1}\text{-}\rep \cala_2))$ ($i=1,2$). Moreover, $F$ is a strict $G$-equivariant functor and therefore a strict $G$-crossed functor since
		\begin{align*}
			 ({}^{g'} \lambda \tilde{\boxtimes} {}^{h'} \mu)(a \otimes b) &= {}^{g'} \lambda ({}^{g'^{-1}}a) \otimes {}^{h' h h'^{-1}} ({}^{h'} (\mu({}^{h'^{-1}} b))) = {}^{g'} \lambda ({}^{g'^{-1}}a) \otimes {}^{h' h} (\mu({}^{h'^{-1}} b)) \\
			 &= ({}^{(g',h')}(\lambda \tilde{\boxtimes} \mu))(a \otimes b)
		\end{align*} 
		for $g',h' \in G$. 

		Finally, we show that any irreducible object $\rho \in \obj((G_1 \times G_2)\text{-}\rep (\cala_1 \otimes \cala_2))$ there exists $\lambda \boxtimes \mu \in \obj((G_1\text{-}\rep \cala_1) \boxtimes (G_2\text{-}\rep \cala_2)^{\mathrm{rev}})$ such that $\rho \cong \lambda \tilde{\boxtimes} \mu$ and therefore $F$ is essentially surjective. By the complete rationality assumption of $\cala_1$ and $\cala_2$, their vacuum Hilbert spaces $H_{\cala_1}$ and $H_{\cala_2}$ are separable by \cite[Proposition 15]{MR1838752} and Bisognano--Wichmann property. By the proof of \cite[Lemma 27]{MR1838752}, the factoriality of $\rho|_{\cala_1 \otimes \mathbb{C}}$ and $\rho|_{\mathbb{C} \otimes \cala_2}$ follows from the simplicity of $\rho$. Suppose $\rho$ is $(g,h)$-localized in $I \times J \in \mathcal{DC}$. Then, by Reeh--Schlieder property, $\rho|_{\cala_1(I') \otimes \mathbb{C}}$ is a faithful normal representation. Since $\cala_1(I') = \cala_1(I)'$ by Haag duality on $\mathbb{R}$, which is included in the complete rationality assumption, $\cala_1(I')$ is a type III factor and therefore $\rho|_{\cala_1(I') \otimes \mathbb{C}}$ is unitarily equivalent to a DHR endomorphism of $\cala_1$ by \cite[Corollary V.3.2]{MR1873025}, which is of type I by \cite[Corollary 14]{MR1838752}. Then, by the proof of \cite[Lemma 27]{MR1838752}, there are faithful normal representations $\pi_1, \pi_2$ of $(\cala_1)_\infty$ and $(\cala_2)_\infty$ on $H_{\cala_1}$ and $H_{\cala_2}$ respectively such that $\rho$ is unitarily equivalent to $\pi_1 \otimes \pi_2$. By replacing $\pi_1$ and $\pi_2$ with endomorphisms, $\rho$ is unitarily equivalent to $\lambda \tilde{\boxtimes} \mu$ for some $\lambda \otimes \mu \in \obj((g\text{-}\rep \cala_1) \boxtimes (h^{-1}\text{-}\rep \cala_2))$. This unitary is indeed in $\cala_1(I) \otimes \cala_2(J)$ by Haag duality on $\mathbb{R}$.
 	\end{proof}
\end{prop}

Thanks to this proposition, we can regard $(G_1 \times G_2)\text{-}\rep (\cala_1 \times \cala_2)$ as a $G_1 \times G_2$-braided ${}^\ast$-tensor category. Note that we cannot apply \cite[Proposition 2.17]{MR2183964} to this category due to the lack of Haag duality on $\mathbb{R}^{1,1}$.

\subsection{Equivariant $\alpha$-induction Frobenius algebras}

In this subsection, we construct equivariant $\alpha$-induction Frobenius algebras (Theorem \ref{mainthm1}) and show that it is symmetric and has the equivariant version of commutativity.  

First, we introduce the notion of the \emph{conjugate} of a morphism of a $G$-braided multitensor category.

\begin{defi}
Let $\calc$ be a $G$-braided multitensor category. Let $\{ \lambda_i \}_{i=1}^n, \{ \mu_j \}_{j=1}^m \subset \homog(\calc)$ and let $f \in \hom(\lambda_1 \cdots \lambda_n,\mu_1 \cdots \mu_m)$. Then, the \emph{conjugate} $\overline{f}$ of $f$ is defined to be the morphism in Figure \ref{morphismconjdef}. 
\begin{figure}[htb]
\centering
\begin{tikzpicture}
\draw[-] (0,0.5) -- (1.25,1.75) arc (180:0:0.125) -- (2.75,0.5) -- (2.75, 0) -- (1,-1.75);
\draw[-, cross] (0.25,0.5) -- (1.5,1.75) arc (180:0:0.125) --(3, 0.5) -- (3, -1.75);
\node at (1,0.75){$\cdots$};
\draw[-, cross] (1.25, 0.5) -- (2.5,1.75) arc (180:0:0.125) -- (4,0.5) -- (4,-1.75);
\node at (2,0.75){$\cdots$};
\draw[-,cross] (2.25, 0.5) -- (3.5,1.75) arc (180:0:0.125) -- (5,0.5) -- (5,-1.75);
\draw[->] (0,0) -- (-1.25,-1.25) arc (360:180:0.125) -- (-2.75,0) -- (-2.75, 2.25);
\draw[->, cross] (0.25,0) -- (-1,-1.25) arc (360:180:0.125) -- (-2.5,0) -- (-2.5,0.5) -- (-0.75, 2.25);
\draw[->, cross] (1.25,0) -- (0,-1.25) arc (360:180:0.125) -- (-1.5,0) -- (-1.5,0.5) -- (0.25, 2.25);
\draw[->,cross] (2.25,0) -- (1,-1.25) arc (360:180:0.125) -- (-0.5,0) -- (-0.5, 0.5) -- (1.25, 2.25);
\draw[fill=white] (-0.25,0.5) rectangle (2.5, 0) node[midway]{$f$};
\node at (0.5,-0.25){$\cdots$};
\node at (1.5,-0.25){$\cdots$};
\node at (0.75,-2){${}^{\overline{\lambda_2 \cdots \lambda_n}} \lambda_1^\vee$};
\node at (2.375,-2){${}^{\overline{\lambda_3 \cdots \lambda_n}} \lambda_2^\vee$};
\node at (4,-2){${}^{\overline{\lambda_{i+1} \cdots \lambda_n}} \lambda_i^\vee$};
\node at (5,-2){$\lambda_n^\vee$};
\node at (-3,2.5) {${}^{\overline{\mu_2 \cdots \mu_m}} \mu_1^\vee$};
\node at (-1.5,2.5){${}^{\overline{\mu_3 \cdots \mu_m}}\mu_2^\vee$};
\node at (0.25,2.5){${}^{\overline{\mu_{j+1} \cdots \mu_m}} \mu_j^\vee$};
\node at (1.5,2.5){$\mu_m^\vee$};
\node at (-1,0.25){$\cdots$};
\node at (-2,0.25){$\cdots$};
\node at (3.5,0.25){$\cdots$};
\node at (4.5,0.25){$\cdots$};
\end{tikzpicture}
\caption{The conjugate $\overline{f}$ of $f \in \hom(\lambda_1 \cdots \lambda_n,\mu_1 \cdots \mu_m)$}
\label{morphismconjdef}
\end{figure}
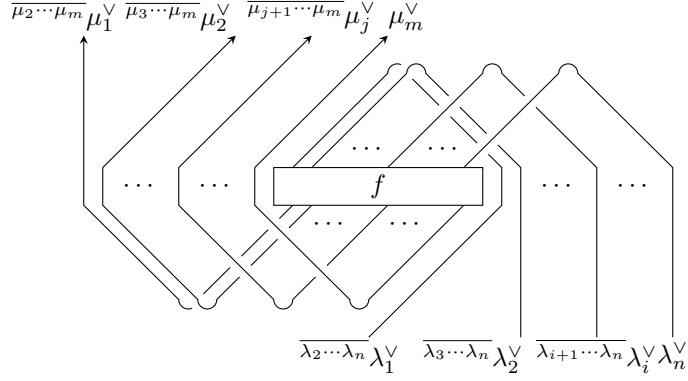
\end{defi}

\begin{rem}
	When $\calc$ is an (ordinary) braided ${}^\ast$-multitensor category, $\overline{f}$ for $f \in \hom(\lambda, \mu_1 \mu_2)$ in the sense of \cite[Section 4.1]{MR3424476} is $\overline{f^\ast}$ in our sense.
\end{rem}

\begin{lem}
\label{conjugatecomposition}
Let $\calc$ be a $G$-braided multitensor category. Let $\{ \lambda_i \}_{i=1}^n, \{ \mu_j \}_{j=1}^m, \{ \nu_k \}_{k=1}^l \subset \homog(\calc)$, $f \in \hom(\lambda_1 \cdots \lambda_n,\mu_1 \cdots \mu_m)$ and $f' \in \hom(\mu_1 \cdots \mu_m, \nu_1 \cdots \nu_l)$. Then, $\overline{f' \circ f} = \overline{f} \circ \overline{f'}$.

\begin{proof}
It follows from the graphical calculation in Figure \ref{graphical_conjugation_contravariant}.
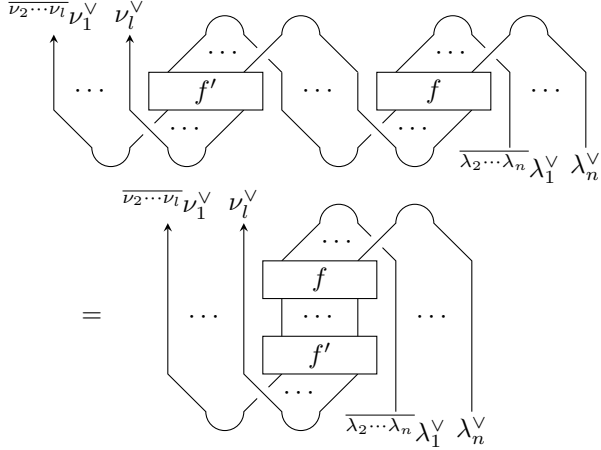
\begin{figure}[htb]
\centering
\begin{tikzpicture}
\draw[-] (0,0) -- (0.5,0.5) arc (180:0:0.25) -- (1.5,0);
\draw[-, cross] (1,0) -- (1.5,0.5) arc (180:0:0.25) -- (2.5,0);
\draw[->] (0,-0.5) -- (-0.5,-1) arc (360:180:0.25) -- (-1.5,-0.5) -- (-1.5, 0.5);
\draw[->, cross] (1,-0.5) -- (0.5,-1) arc (360:180:0.25) -- (-0.5,-0.5) -- (-0.5, 0.5);
\draw[fill=white] (-0.25, 0) rectangle (1.25, -0.5) node[midway]{$f'$};
\node at (-1.5, 0.75) {${}^{\overline{\nu_2 \cdots \nu_l}} \nu_1^\vee$};
\node at (-0.5, 0.75) {$\nu_l^\vee$};
\node at (0.25, -0.75) {$\cdots$};
\node at (0.75, 0.25) {$\cdots$};
\node at (-1, -0.25) {$\cdots$};
\begin{scope}[shift={(3,0)}]
\draw[-] (0,0) -- (0.5,0.5) arc (180:0:0.25) -- (1.5,0) -- (1.5, -1);
\draw[-, cross] (1,0) -- (1.5,0.5) arc (180:0:0.25) -- (2.5,0) -- (2.5, -1);
\draw[-] (0,-0.5) -- (-0.5,-1) arc (360:180:0.25) -- (-1.5,-0.5) -- (-1.5, 0);
\draw[-, cross] (1,-0.5) -- (0.5,-1) arc (360:180:0.25) -- (-0.5,-0.5) -- (-0.5,0);
\draw[fill=white] (-0.25, 0) rectangle (1.25, -0.5) node[midway]{$f$};
\node at (1.5, -1.25) {${}^{\overline{\lambda_2 \cdots \lambda_n}} \lambda_1^\vee$};
\node at (2.5, -1.25) {$\lambda_n^\vee$};
\node at (0.25, -0.75) {$\cdots$};
\node at (0.75, 0.25) {$\cdots$};
\node at (-1, -0.25) {$\cdots$};
\node at (2, -0.25) {$\cdots$};
\end{scope}
\begin{scope}[shift={(1.5,-3.5)}]
\draw[->] (0,-0.5) -- (-0.5,-1) arc (360:180:0.25) -- (-1.5,-0.5) -- (-1.5, 1.5);
\draw[->, cross] (1,-0.5) -- (0.5,-1) arc (360:180:0.25) -- (-0.5,-0.5) -- (-0.5, 1.5);
\draw[fill=white] (-0.25, 0) rectangle (1.25, -0.5) node[midway]{$f'$};
\node at (-1.5, 1.75) {${}^{\overline{\nu_2 \cdots \nu_l}} \nu_1^\vee$};
\node at (-0.5, 1.75) {$\nu_l^\vee$};
\node at (0.25, -0.75) {$\cdots$};
\node at (2, 0.25) {$\cdots$};
\begin{scope}[shift={(0,1)}]
\draw[-] (0,0) -- (0.5,0.5) arc (180:0:0.25) -- (1.5,0) -- (1.5,-2);
\draw[-, cross] (1,0) -- (1.5,0.5) arc (180:0:0.25) -- (2.5,0) -- (2.5,-2);
\draw[-] (0,-0.5) -- (0,-1);
\draw[-] (1,-0.5) -- (1,-1);
\draw[fill=white] (-0.25, 0) rectangle (1.25, -0.5) node[midway]{$f$};
\node at (1.5, -2.25) {${}^{\overline{\lambda_2 \cdots \lambda_n}} \lambda_1^\vee$};
\node at (2.5, -2.25) {$\lambda_n^\vee$};
\node at (0.5, -0.75) {$\cdots$};
\node at (0.75, 0.25) {$\cdots$};
\node at (-1, -0.75) {$\cdots$};
\node at (-2.5, -0.75) {$=$};
\end{scope}
\end{scope}
\end{tikzpicture}
\caption{Conjugation is contravariant}
\label{graphical_conjugation_contravariant}
\end{figure}
\end{proof}
\end{lem}

\begin{lem}
\label{conjugatemonoidalproduct}
Let $\calc$ be a $G$-braided multitensor category. Let $\{ \lambda_i \}_{i=1}^n, \{ \mu_j \}_{j=1}^m, \{ \nu_k \}_{k=1}^l, \{\rho_p\}_{p=1}^q$ be families in $\homog(\calc)$ and let $f \in \hom(\lambda_1 \cdots \lambda_n,\mu_1 \cdots \mu_m)$ and $f' \in \hom(\nu_1 \cdots \nu_l, \rho_1 \cdots \rho_q)$. Then, $\overline{f \otimes f'} = {}^{\overline{\nu_1 \cdots \nu_l}} (\overline{f}) \otimes \overline{f'}$.

\begin{proof}
We may assume $\partial (\nu_1 \cdots \nu_l) = \partial (\rho_1 \cdots \rho_q)$ since otherwise $f'=0$. Then, the statement follows from the graphical calculation in Figure \ref{graphical_conjugation_tensor}, where we used the assumption $\partial (\nu_1 \cdots \nu_l) = \partial (\rho_1 \cdots \rho_q)$ at the final equality. 
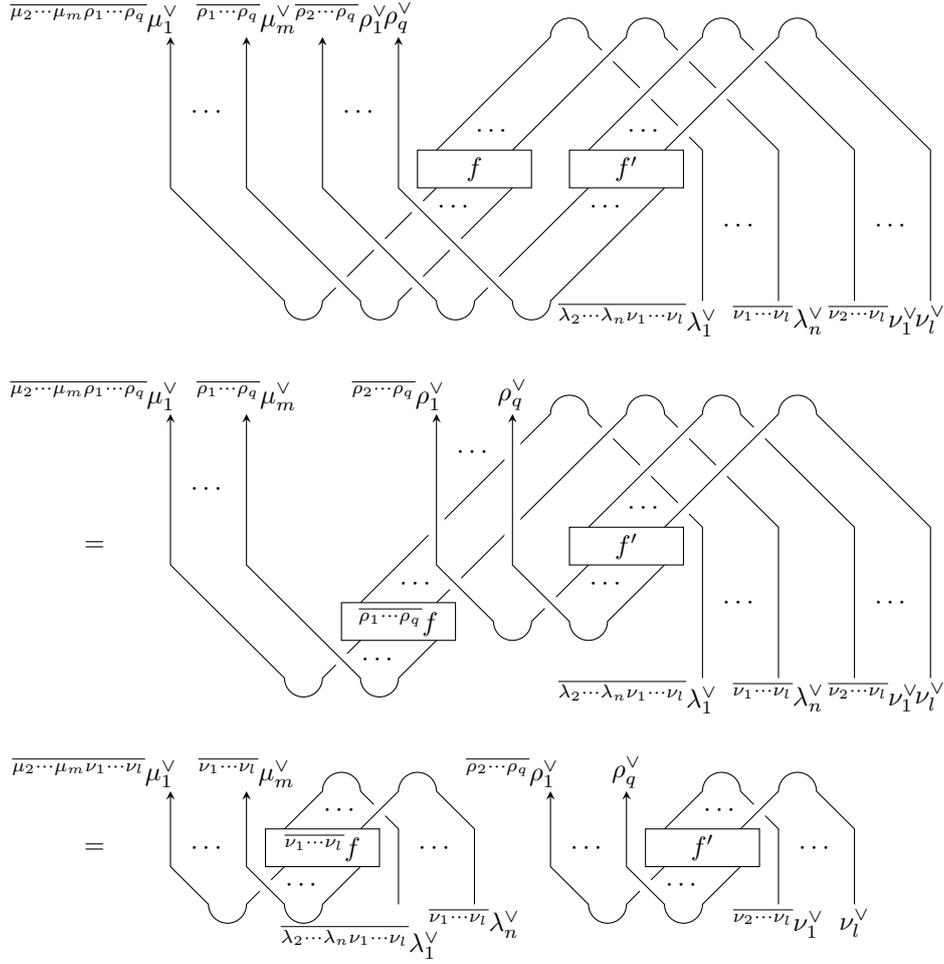
\begin{figure}[htb]
\centering
\begin{tikzpicture}
\foreach \x [evaluate=\x as \y using \x+1.5, evaluate=\x as \z using \x+3.5] in {0,1,2,3}
{
	\draw[-, cross] (\x,0) -- (\y,1.5) arc (180:0:0.25) -- (\z,0) -- (\z, -2);
}%the for loop describing upper curves
\foreach \x [evaluate=\x as \y using \x-1.5, evaluate=\x as \z using \x-3.5] in {0,1,2,3}
{
	\draw[->, cross] (\x,-0.5) -- (\y,-2) arc (360:180:0.25) -- (\z,-0.5) -- (\z, 1.5);
}%the for loop describing lower curves
\draw[fill=white] (-0.25, 0) rectangle (1.25, -0.5) node[midway]{$f$};
\draw[fill=white] (1.75, 0) rectangle (3.25, -0.5) node[midway]{$f'$};
\node at (-3,0.5){$\cdots$};
\node at (-1,0.5){$\cdots$};
\node at (6,-1){$\cdots$};
\node at (4,-1){$\cdots$};
\node at (0.25,-0.75){$\cdots$};
\node at (0.75,0.25){$\cdots$};
\node at (2.25,-0.75){$\cdots$};
\node at (2.75,0.25){$\cdots$};
\node at (-4.5, 1.75){${}^{\overline{\mu_2 \cdots \mu_m \rho_1 \cdots \rho_q}} \mu_1^\vee$};
\node at (-2.5, 1.75){${}^{\overline{\rho_1 \cdots \rho_q}} \mu_m^\vee$};
\node at (-1.25, 1.75){${}^{\overline{\rho_2 \cdots \rho_q}} \rho_1^\vee$};
\node at (-0.5, 1.75){$\rho_q^\vee$};
\node at (2.65, -2.25){${}^{\overline{\lambda_2 \cdots \lambda_n \nu_1 \cdots \nu_l}} \lambda_1^\vee$};
\node at (4.5, -2.25){${}^{\overline{\nu_1 \cdots \nu_l}} \lambda_n^\vee$};
\node at (5.75, -2.25){${}^{\overline{\nu_2 \cdots \nu_l}} \nu_1^\vee$};
\node at (6.5, -2.25){$\nu_l^\vee$};
\begin{scope}[shift={(0,-5)}]
\foreach \x [evaluate=\x as \y using \x+1.5, evaluate=\x as \z using \x+3.5, evaluate=\w using \x-1] in {0,1}
{
	\draw[-, cross] (\w,-1) -- (\y,1.5) arc (180:0:0.25) -- (\z,0) -- (\z, -2);
}%the for loop describing upper curves
\foreach \x [evaluate=\x as \y using \x+1.5, evaluate=\x as \z using \x+3.5] in {2,3}
{
	\draw[-, cross] (\x,0) -- (\y,1.5) arc (180:0:0.25) -- (\z,0) -- (\z, -2);
}%the for loop describing upper curves
\foreach \x [evaluate=\x as \y using \x-1.5, evaluate=\x as \z using \x-3.5, evaluate=\w using \x-1] in {0,1}
{
	\draw[->, cross] (\w,-1.5) -- (\y,-2) arc (360:180:0.25) -- (\z,-0.5) -- (\z, 1.5);
}%the for loop describing lower curves
\foreach \x [evaluate=\x as \y using \x-0.75, evaluate=\x as \z using \x-2] in {2,3}
{
	\draw[->, cross] (\x,-0.5) -- (\y,-1.25) arc (360:180:0.25) -- (\z,-0.5) -- (\z, 1.5);
}%the for loop describing lower curves
\draw[fill=white] (-1.25, -1) rectangle (0.25, -1.5) node[midway]{${}^{\overline{\rho_1 \cdots \rho_q}}f$};
\draw[fill=white] (1.75, 0) rectangle (3.25, -0.5) node[midway]{$f'$};
\node at (-3,0.5){$\cdots$};
\node at (0.5,1){$\cdots$};
\node at (6,-1){$\cdots$};
\node at (4,-1){$\cdots$};
\node at (-0.75,-1.75){$\cdots$};
\node at (-0.25,-0.75){$\cdots$};
\node at (2.25,-0.75){$\cdots$};
\node at (2.75,0.25){$\cdots$};
\node at (-4.5, 1.75){${}^{\overline{\mu_2 \cdots \mu_m \rho_1 \cdots \rho_q}} \mu_1^\vee$};
\node at (-2.5, 1.75){${}^{\overline{\rho_1 \cdots \rho_q}} \mu_m^\vee$};
\node at (-0.5, 1.75){${}^{\overline{\rho_2 \cdots \rho_q}} \rho_1^\vee$};
\node at (1, 1.75){$\rho_q^\vee$};
\node at (2.65, -2.25){${}^{\overline{\lambda_2 \cdots \lambda_n \nu_1 \cdots \nu_l}} \lambda_1^\vee$};
\node at (4.5, -2.25){${}^{\overline{\nu_1 \cdots \nu_l}} \lambda_n^\vee$};
\node at (5.75, -2.25){${}^{\overline{\nu_2 \cdots \nu_l}} \nu_1^\vee$};
\node at (6.5, -2.25){$\nu_l^\vee$};
\node at (-4.5,-0.25){$=$};
\begin{scope}[shift={(-2,-4)}]
\draw[-] (0,0) -- (0.5,0.5) arc (180:0:0.25) -- (1.5,0) -- (1.5, -1);
\draw[-, cross] (1,0) -- (1.5,0.5) arc (180:0:0.25) -- (2.5,0) -- (2.5, -1);
\draw[->] (0,-0.5) -- (-0.5,-1) arc (360:180:0.25) -- (-1.5,-0.5) -- (-1.5, 0.5);
\draw[->, cross] (1,-0.5) -- (0.5,-1) arc (360:180:0.25) -- (-0.5,-0.5) -- (-0.5, 0.5);
\draw[fill=white] (-0.25, 0) rectangle (1.25, -0.5) node[midway]{${}^{\overline{\nu_1 \cdots \nu_l}}f$};
\node at (-2.5, 0.75) {${}^{\overline{\mu_2 \cdots \mu_m \nu_1 \cdots \nu_l}} \mu_1^\vee$};
\node at (-0.5, 0.75) {${}^{\overline{\nu_1 \cdots \nu_l}} \mu_m^\vee$};
\node at (0.25, -0.75) {$\cdots$};
\node at (0.75, 0.25) {$\cdots$};
\node at (-1, -0.25) {$\cdots$};
\node at (2, -0.25) {$\cdots$};
\node at (1, -1.5){${}^{\overline{\lambda_2 \cdots \lambda_n \nu_1 \cdots \nu_l}} \lambda_1^\vee$};
\node at (2.5, -1.25){${}^{\overline{\nu_1 \cdots \nu_l}} \lambda_n^\vee$};
\node at (-2.5, -0.25) {$=$};
\end{scope}
\begin{scope}[shift={(3,-4)}]
\draw[-] (0,0) -- (0.5,0.5) arc (180:0:0.25) -- (1.5,0) -- (1.5, -1);
\draw[-, cross] (1,0) -- (1.5,0.5) arc (180:0:0.25) -- (2.5,0) -- (2.5, -1);
\draw[->] (0,-0.5) -- (-0.5,-1) arc (360:180:0.25) -- (-1.5,-0.5) -- (-1.5, 0.5);
\draw[->, cross] (1,-0.5) -- (0.5,-1) arc (360:180:0.25) -- (-0.5,-0.5) -- (-0.5, 0.5);
\draw[fill=white] (-0.25, 0) rectangle (1.25, -0.5) node[midway]{$f'$};
\node at (-2, 0.75) {${}^{\overline{\rho_2 \cdots \rho_q}} \rho_1^\vee$};
\node at (-0.5, 0.75) {$\rho_q^\vee$};
\node at (0.25, -0.75) {$\cdots$};
\node at (0.75, 0.25) {$\cdots$};
\node at (-1, -0.25) {$\cdots$};
\node at (2, -0.25) {$\cdots$};
\node at (1.5, -1.25){${}^{\overline{\nu_2 \cdots \nu_l}} \nu_1^\vee$};
\node at (2.5, -1.25){$\nu_l^\vee$};
\end{scope}
\end{scope}
\end{tikzpicture}
\caption{Conjugation and the monoidal product}
\label{graphical_conjugation_tensor}
\end{figure}
\end{proof}
\end{lem}

\begin{lem}
\label{conjugatestar}
Let $\calc$ be a $G$-braided ${}^\ast$-multitensor category. Let $\{ \lambda_i \}_{i=1}^n, \{ \mu_j \}_{j=1}^m \subset \homog(\calc)$ and let $f \in \hom(\lambda_1 \cdots \lambda_n,\mu_1 \cdots \mu_m)$. Then $\overline{f^\ast} = (\overline{f})^\ast$ if we define $\overline{f}$ using standard solutions of conjugate equations (for the definition of a standard solution in a ${}^\ast$-multitensor category, see \cite[Definition 8.29]{MR3994584}).

\begin{proof}
Since standard left and right duals of a morphism in a ${}^\ast$-multitensor category coincide by \cite[Proposition 8.33]{MR3994584}, the equation in Figure \ref{graphical_conjugation_star} holds by Figure \ref{graphicalrotation}. Therefore, the statement follows by taking ${}^\ast$-involution.
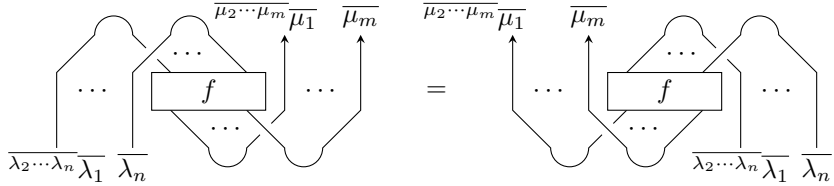
\begin{figure}[H]
\centering
\begin{tikzpicture}
\begin{scope}[xscale=1, yscale=-1]
\draw[->] (0,0) -- (0.5,0.5) arc (180:0:0.25) -- (1.5,0) -- (1.5, -1);
\draw[->, cross] (1,0) -- (1.5,0.5) arc (180:0:0.25) -- (2.5,0) -- (2.5, -1);
\draw[-] (0,-0.5) -- (-0.5,-1) arc (360:180:0.25) -- (-1.5,-0.5) -- (-1.5, 0.5);
\draw[-, cross] (1,-0.5) -- (0.5,-1) arc (360:180:0.25) -- (-0.5,-0.5) -- (-0.5, 0.5);
\node at (1.25, -1.25){${}^{\overline{\mu_2 \cdots \mu_m}} \overline{\mu_1}$};
\node at (2.5, -1.25){$\overline{\mu_m}$};
\node at (-1.5, 0.75) {${}^{\overline{\lambda_2 \cdots \lambda_n}} \overline{\lambda_1}$};
\node at (-0.5, 0.75) {$\overline{\lambda_n}$};
\node at (0.25, -0.75) {$\cdots$};
\node at (0.75, 0.25) {$\cdots$};
\node at (-1, -0.25) {$\cdots$};
\node at (2, -0.25) {$\cdots$};
\end{scope}
\draw[fill=white] (-0.25, 0.5) rectangle (1.25, 0) node[midway]{$f$};
\node at (3.5,0.25) {$=$};
\begin{scope}[shift={(6,0.5)}]
\draw[-] (0,0) -- (0.5,0.5) arc (180:0:0.25) -- (1.5,0) -- (1.5, -1);
\draw[-, cross] (1,0) -- (1.5,0.5) arc (180:0:0.25) -- (2.5,0) -- (2.5, -1);
\draw[->] (0,-0.5) -- (-0.5,-1) arc (360:180:0.25) -- (-1.5,-0.5) -- (-1.5, 0.5);
\draw[->, cross] (1,-0.5) -- (0.5,-1) arc (360:180:0.25) -- (-0.5,-0.5) -- (-0.5, 0.5);
\draw[fill=white] (-0.25, 0) rectangle (1.25, -0.5) node[midway]{$f$};
\node at (-2, 0.75) {${}^{\overline{\mu_2 \cdots \mu_m}} \overline{\mu_1}$};
\node at (-0.5, 0.75) {$\overline{\mu_m}$};
\node at (0.25, -0.75) {$\cdots$};
\node at (0.75, 0.25) {$\cdots$};
\node at (-1, -0.25) {$\cdots$};
\node at (2, -0.25) {$\cdots$};
\node at (1.5, -1.25){${}^{\overline{\lambda_2 \cdots \lambda_n}} \overline{\lambda_1}$};
\node at (2.5, -1.25){$\overline{\lambda_n}$};
\end{scope}
\end{tikzpicture}
\caption{Conjugation and the ${}^\ast$-structure}
\label{graphical_conjugation_star}
\end{figure}

\end{proof}
\end{lem}

Next, let $\calc$ be a semisimple pivotal tensor category and let $\tr$ denote $\trl(\delta^\calc -)$. For $\lambda, \mu \in \obj(\calc)$ with $d_{\lambda} \coloneqq \dim \lambda \neq 0$, the bilinear form given by putting $(f',f) \coloneqq d_\lambda^{-1} \tr(f' \circ f)$ for $f \in \hom(\lambda, \mu)$ and $f' \in \hom(\mu,\lambda)$ is nondegenerate by the proof of \cite[Lemma II.4.2.3]{MR3617439}. Therefore, $\hom(\mu, \lambda)$ can be identified with the dual space $\hom(\lambda, \mu)^\vee$ of $\hom(\lambda, \mu)$, and for a basis $\{ e_i \}_{i=1}^{\langle \lambda, \mu \rangle}$ of $\hom(\lambda,\mu)$, where $\langle \lambda, \mu \rangle \coloneqq \dim \hom (\lambda, \mu)$, we can take the corresponding dual basis $\{\tilde{e}_i \}_{i=1}^{\langle \lambda, \mu \rangle}$ of $\hom(\mu, \lambda)$. Note that if $\calc$ is split and $\lambda$ is simple, then $\tilde{e}_i e_j = (\tilde{e}_i,e_j) \id_\lambda = \delta_{ij} \id_\lambda$.

In particular, when $\calc$ is a ${}^\ast$-tensor category, we can make $\hom(\lambda, \mu)$ into a Hilbert space by putting $(f_1|f_2) \coloneqq d_\lambda^{-1} \tr(f_1^\ast f_2)$ for $f_1, f_2 \in \hom(\lambda, \mu)$. By the Gram--Schmidt orthonormalization, we can take an orthonormal basis $\{ e_i \}_i$ of $\hom(\lambda, \mu)$. Note that if $\lambda$ is simple, then $f_1^\ast f_2 = (f_1|f_2) \id_\lambda$ and therefore $\{ e_i \}_i$ is a family of orthogonal isometries with the dual basis $\{ e_i^\ast \}_i$.

Now, we can prove one of our main theorems. Let $A$ be a neutral symmetric special simple $G$-equivariant Frobenius algebra in a split $G$-braided fusion category $\calc$. Fix a complete system $I$ of representatives of the isomorphism classes of the simple objects of $\calc$ and define $\theta \in \obj(D(\calc))$ by putting 
\begin{align*}
\theta &\coloneqq \bigoplus_{\substack{\lambda_1, \lambda_2 \in I \\ \partial \lambda_1 = \partial \lambda_2}} \langle \alpha_{A}^{G+}(\lambda_1), \alpha_{A}^{G-}(\lambda_2) \rangle \lambda_1 \boxtimes \lambda_2^\vee.
\end{align*}

Let $\eta \in \hom(\mathbf{1}_{D(\calc)},\theta)$ denote the inclusion and let $\varepsilon \in \hom(\theta, \mathbf{1}_{D(\calc)})$ denote the projection multiplied by $\dim \theta$. 

For any $\lambda_1, \lambda_2 \in I$ with $\partial \lambda_1 = \partial \lambda_2$, fix a basis $\{ \phi^L_l \}_{l}$ of $\hom(\alpha_A^{G+}(\lambda_1), \alpha_A^{G-}(\lambda_2))$, where $L$ is a multi-index $L=(\lambda_1,\lambda_2,l)$ which consists of $\lambda_1,\lambda_2 \in I$ with $\partial \lambda_1 = \partial \lambda_2$ and an integer $1 \le l \le \langle \alpha_{A}^{G+}(\lambda_1), \alpha_{A}^{G-}(\lambda_2) \rangle$. For any $\lambda, \mu, \nu \in I$, fix a basis $\{e^{\nu, \lambda \mu}_i \}_i$ of $\hom(\nu, \lambda \mu)$. Define $\tilde{m} \in \hom (\theta^2, \theta)$ and $\Delta \in \hom (\theta, \theta^2)$ by putting
\begin{align*}
\tilde{m} &\coloneqq \bigoplus_{L,M,N} \sum_{i,j} \frac{1}{d_{\nu_1}} \tr(\tilde{\phi}^N_n \alpha_A^{G-}(\tilde{e}_{j}^{\nu_2, \lambda_2 \mu_2}) (\phi_l^{L} \otimes_A \phi_m^{M}) \alpha_A^{G+}(e_{i}^{\nu_1, \lambda_1 \mu_1})) (\tilde{e}_{i}^{\nu_1, \lambda_1 \mu_1} \boxtimes \overline{e_{j}^{\nu_2, \lambda_2 \mu_2}}) \\
\Delta &\coloneqq \bigoplus_{L,M,N} \sum_{i,j} \frac{d_{\lambda_2}d_{\mu_2}}{d_{\nu_2} d_{\nu_1}} \tr(\alpha_A^{G+}(\tilde{e}_{i}^{\nu_1, \lambda_1 \mu_1}) (\tilde{\phi}_l^{L} \otimes_A \tilde{\phi}_m^{M}) \alpha_A^{G-}(e_{j}^{\nu_2, \lambda_2 \mu_2}) \phi^N_n) (e_{i}^{\nu_1, \lambda_1 \mu_1} \boxtimes \overline{\tilde{e}_{j}^{\nu_2, \lambda_2 \mu_2}}),
\end{align*}
where $M=(\mu_1,\mu_2,m)$ and $N=(\nu_1, \nu_2, n)$ are multi-indices like $L$. 

For any $g \in G$ and $\lambda \in I$, fix an isomorphism $u_g^\lambda : {}^g \lambda \cong \lambda (g)$, where $\lambda(g)$ denotes the object in $I$ which is isomorphic to ${}^g \lambda$. For every $g \in G$, define $z_g \in \hom({}^g \theta, \theta)$ by putting
\begin{align*}
	z_g &\coloneqq \bigoplus_{L,l'} \frac{1}{d_{\lambda_1}} \tr(\alpha^{G+}_A(\tilde{u}_g^{\lambda_1}) \tilde{\phi}_{l'}^{L(g)} \alpha^{G-}_A(u_g^{\lambda_2}) \tilde{\eta}^{g-1}_{\lambda_2} {}^g \phi_l^L \tilde{\eta}^g_{\lambda_1}) (u_g^{\lambda_1} \boxtimes \overline{\tilde{u}_g^{\lambda_2}}),
\end{align*}
where $L(g)$ is the multi-index $(\lambda_1(g), \lambda_2(g),m)$ and $\tilde{\eta}^g_\lambda \coloneqq z_g^{A-1} \otimes_\calc \id_{{}^g \lambda}$ for $\lambda \in I$.

When $A$ is a Q-system in a ${}^\ast$-fusion category $\calc$, let $w \in \hom(\mathbf{1}_{D(\calc)}, \theta)$ denote the isometric inclusion. Take $\{ \phi^L_l \}_l$'s and $\{ e^{\nu, \lambda \mu}_i \}_i$'s to be orthonormal bases and define $x \in \hom(\theta, \theta^2)$ to be
\begin{align*}
	\bigoplus_{L,M,N} \sum_{i,j} \sqrt{\frac{d_{\lambda_2}d_{\mu_2}}{d_{\theta} d_{\nu_2}}} \frac{1}{d_{\nu_1}} \tr(\alpha^{G+}_A(e_{i}^{\nu_1, \lambda_1 \mu_1})^\ast (\phi_l^{L} \otimes_A \phi_m^{M})^\ast \alpha^{G-}_A(e_{j}^{\nu_2, \lambda_2 \mu_2})\phi^N_n)(e_{i}^{\nu_1, \lambda_1 \mu_1} \boxtimes \overline{e_{j}^{\nu_2, \lambda_2 \mu_2 \ast}}).
\end{align*}

Take $u_g^\lambda$'s to be unitaries and define $z \in \hom({}^g \theta, \theta)$ by putting
\begin{align*}
	z_g &\coloneqq \bigoplus_{L,l'} \frac{1}{d_{\lambda_1}} \tr(\alpha^{G+}_A(u_g^{\lambda_1})^\ast \phi_{l'}^{L(g) \ast} \alpha^{G-}_A(u_g^{\lambda_2}) \tilde{\eta}^{g \ast}_{\lambda_2} {}^g \phi_l^L \tilde{\eta}^g_{\lambda_1}) (u_g^{\lambda_1} \boxtimes \overline{u_g^{\lambda_2 \ast}}).
\end{align*}

\begin{thm}
	\label{mainthm1}
Let $A$ be a neutral symmetric special simple $G$-equivariant Frobenius algebra in a split $G$-braided fusion category $\calc$. Then, the tuple $\Theta^G_{\alpha} (A) \coloneqq (\theta, \tilde{m}, \eta, \Delta, \varepsilon, z)$ is a $G$-equivariant Frobenius algebra in $D(\calc)$. 

If $\calc$ is moreover spherical, then $\tilde{m} \circ \Delta = \dim \Theta_\alpha^G(A)$. Another choice of $\{ e^{\nu, \lambda \mu}_i\}_i$'s and $u_g^\lambda$'s yields the same $G$-equivariant Frobenius algebra. Another choice of $I$, a direct sum and $\lambda_2^\vee$'s in the definition of $\theta$ and $\{ \phi_l^L \}_l$'s yields an isomorphic $G$-equivariant Frobenius algebra. We call $\Theta^G_\alpha(A)$ the \emph{$G$-equivariant $\alpha$-induction Frobenius algebra} associated with $A$. 

When $A$ is a Q-system in a ${}^\ast$-fusion category $\calc$, $\Theta_\alpha^G(A) \coloneqq (\theta, w, x, z)$ is a $G$-equivariant Q-system. Another choice of $\{ e^{\nu, \lambda \mu}_i\}_i$'s and $u_g^\lambda$'s yields the same $G$-equivariant Q-system. Another choice of $I$, a direct sum and $\lambda_2^\vee$'s in the definition of $\theta$ and $\{ \phi_l^L \}_l$'s yields an isomorphic $G$-equivariant Q-system. We call $\Theta^G_\alpha(A)$ the \emph{$G$-equivariant $\alpha$-induction Q-system} associated with $A$.

\begin{proof}
First, we note that $\tilde{m}, \Delta$ and $z$ are well-defined. Indeed,
\begin{align*}
\hom(\theta, \theta^2) &= \bigoplus_{L,M,N} \hom(\nu_1 \boxtimes \nu_2^\vee, (\lambda_1 \boxtimes \lambda_2^\vee) \otimes (\mu_1 \boxtimes \mu_2^\vee)) \\
&= \bigoplus_{L,M,N} \hom(\nu_1 \boxtimes \nu_2^\vee, \lambda_1 \mu_1 \boxtimes {}^{\overline{\mu_2}} \lambda_2^\vee \mu_2^\vee)
\end{align*}
by $\partial \mu_1 = \partial \mu_2$ and $e_{i}^{\nu_1, \lambda_1 \mu_1} \boxtimes \overline{\tilde{e}_{j}^{\nu_2, \lambda_2 \mu_2}} \in \hom(\nu_1 \boxtimes \nu_2^\vee, \lambda_1 \mu_1 \boxtimes {}^{\overline{\mu_2}} \lambda_2^\vee \mu_2^\vee)$ by definition. The coefficient of $\Delta$ is also well-defined. Indeed, we can take the dual bases $\{ \tilde{e}^{\nu, \lambda \mu}_i \}_i$'s since $d_\nu \neq 0$ by \cite[Proposition 4.8.4]{egno}. The category $\bimod_\calc(A)$ is a tensor category by assumption, pivotal by Proposition \ref{bimodpivotalprop} and semisimple by \cite[Proposition 6.1]{MR2075605}. Since $\dim \alpha^{G \pm}_A (\lambda) = d_{\lambda} \neq 0$ for $\lambda \in I$ by Proposition \ref{semidirpivotalprop}, we can take the dual bases $\{ \tilde{\phi}^L_l \}_l$'s. Thus, $\Delta$ is well-defined. Similarly, $m$ and $z$ are well-defined.

Next, we note that $\Theta_\alpha^G(A)$ does not depend on a choice of $\{ e^{\nu, \lambda \mu}_i\}_i$'s since the effect of changing $\{ e^{\nu, \lambda \mu}_i\}_i$'s in the definition of $\tilde{m}$ and $\Delta$ is canceled by the corresponding change of $\{ \tilde{e}^{\nu, \lambda \mu}_i\}_i$'s. For the same reason, it does not depend on a choice of $u_g^\lambda$'s. Note also that another choice of $\{ \phi_l^L\}_l$'s only yields another direct sum in the definition of $\theta$. Once we have shown that $\Theta_\alpha^G(A)$ is a $G$-equivariant Frobenius algebra, we can see that another direct sum yields a canonically isomorphic one since we define structures to be direct sums. Also, another choice of left duals yields a canonically isomorphic one by the definition of conjugation. Another choice of $I$ yields an isomorphic one by fixing isomorphisms between simple objects.

Then, we show that $\Theta_\alpha^G(A)$ is a Frobenius algebra. Because the proof is similar to that of \cite[Theorem 1.4]{MR1754521}, we only show coassociativity here. We have 
\begin{align*}
&(\Delta \otimes \id_\theta) \Delta \\ 
&= \bigoplus_{L,M,M',N'} \sum_{N,i,j,i',j'} \frac{d_{\lambda_2} d_{\mu_2} d_{\mu_2'}}{d_{\nu_2'}d_{\nu_1} d_{\nu_1'}} \tr(\alpha_A^{G+}(\tilde{e}_{i}^{\nu_1, \lambda_1 \mu_1}) (\tilde{\phi}_l^{L} \otimes_A \tilde{\phi}_m^{M}) \alpha_A^{G-}(e_{j}^{\nu_2, \lambda_2 \mu_2}) \phi^N_n) \\
&\hspace{3.5cm} \times \tr(\alpha_A^{G+}(\tilde{e}_{i'}^{\nu'_1, \nu_1 \mu'_1}) (\tilde{\phi}_n^{N} \otimes_A \tilde{\phi}_{m'}^{M'}) \alpha_A^{G-}(e_{j'}^{\nu'_2, \nu_2 \mu'_2}) \phi^{N'}_{n'}) \\ 
&\hspace{3.5cm} \times (e_{i}^{\nu_1, \lambda_1 \mu_1} \otimes \id_{\mu_1'}) e_{i'}^{\nu_1', \nu_1 \mu'_1} \boxtimes ({}^{\overline{\mu_1'}} (\overline{\tilde{e}_{j}^{\nu_2, \lambda_2 \mu_2}}) \otimes \id_{\mu_2'}) \overline{\tilde{e}_{j'}^{\nu_2', \nu_2 \mu'_2}} \\
&= \bigoplus_{L,M,M',N'} \sum_{N,i,j,i',j'} \frac{d_{\lambda_2} d_{\mu_2} d_{\mu_2'}}{d_{\nu_2'}d_{\nu_1} d_{\nu_1'}} \tr(\alpha_A^{G+}(\tilde{e}_{i}^{\nu_1, \lambda_1 \mu_1}) (\tilde{\phi}_l^{L} \otimes_A \tilde{\phi}_m^{M}) \alpha_A^{G-}(e_{j}^{\nu_2, \lambda_2 \mu_2}) \phi^N_n) \\
&\hspace{3.5cm} \times \tr(\alpha_A^{G+}(\tilde{e}_{i'}^{\nu'_1, \nu_1 \mu'_1}) (\tilde{\phi}_n^{N} \otimes_A \tilde{\phi}_{m'}^{M'}) \alpha_A^{G-}(e_{j'}^{\nu'_2, \nu_2 \mu'_2}) \phi^{N'}_{n'}) \\ 
&\hspace{3.5cm} \times (e_{i}^{\nu_1, \lambda_1 \mu_1} \otimes \id_{\mu_1'}) e_{i'}^{\nu_1', \nu_1 \mu'_1} \boxtimes \overline{\tilde{e}_{j'}^{\nu_2', \nu_2 \mu'_2}(\tilde{e}_{j}^{\nu_2, \lambda_2 \mu_2} \otimes \id_{\mu_2'})}
\end{align*}
by $\partial \mu'_1 = \partial \mu'_2$ and Lemmata \ref{conjugatecomposition} and \ref{conjugatemonoidalproduct}. By the cyclicity of the trace, we have
\begin{align*}
	&\frac{1}{d_{\nu_1}} \tr(\alpha_A^{G+}(\tilde{e}_{i}^{\nu_1, \lambda_1 \mu_1}) (\tilde{\phi}_l^{L} \otimes_A \tilde{\phi}_m^{M}) \alpha_A^{G-}(e_{j}^{\nu_2, \lambda_2 \mu_2}) \phi^N_n) \\
	&= \frac{d_{\nu_2}}{d_{\nu_1}} \frac{1}{d_{\nu_2}} \tr(\phi^N_n \alpha_A^{G+}(\tilde{e}_{i}^{\nu_1, \lambda_1 \mu_1}) (\tilde{\phi}_l^{L} \otimes_A \tilde{\phi}_m^{M}) \alpha_A^{G-}(e_{j}^{\nu_2, \lambda_2 \mu_2})) \\
	&= \left( \frac{d_{\nu_2}}{d_{\nu_1}} \phi_n^N, \alpha_A^{G+}(\tilde{e}_{i}^{\nu_1, \lambda_1 \mu_1}) (\tilde{\phi}_l^{L} \otimes_A \tilde{\phi}_m^{M}) \alpha_A^{G-}(e_{j}^{\nu_2, \lambda_2 \mu_2}) \right).
\end{align*}
Note that $\{ d_{\nu_2}/d_{\nu_1} \phi_n^N \}_n$ is the dual basis of $\{ \tilde{\phi}_n^N \}_n$ by the cyclicity of the trace since $\dim \alpha_A^{G+}(\nu_1) = d_{\nu_1}$ and $\dim \alpha_A^{G-}(\nu_2) = d_{\nu_2}$ by Proposition \ref{bimodpivotalprop}. Then, by the Fourier expansion in the basis $\{ \tilde{\phi}_n^N \}$ summing over the index $n$ in $N$, we obtain
\begin{align*}
	&(\Delta \otimes \id_\theta) \Delta \\
	&= \bigoplus_{L,M,M',N'} \sum_{\nu_1,\nu_2,i,j,i',j'} \frac{d_{\lambda_2} d_{\mu_2} d_{\mu_2'}}{d_{\nu_2'}d_{\nu_1'}} \\
	&\quad \times \tr(\alpha_A^{G+}(\tilde{e}_{i'}^{\nu'_1, \nu_1 \mu'_1} (\tilde{e}_i^{\nu_1, \lambda_1 \mu_1} \otimes \id_{\mu_1'})) (\tilde{\phi}_l^{L} \otimes_A \tilde{\phi}_m^M \otimes_A \tilde{\phi}_{m'}^{M'}) \alpha_A^{G-}((e^{\nu_2,\lambda_2\mu_2} \otimes \id_{\mu_2'})e_{j'}^{\nu'_2, \nu_2 \mu'_2}) \phi^{N'}_{n'}) \\
	&\quad \times (e_{i}^{\nu_1, \lambda_1 \mu_1} \otimes \id_{\mu_1'}) e_{i'}^{\nu_1', \nu_1 \mu'_1} \boxtimes \overline{\tilde{e}_{j'}^{\nu_2', \nu_2 \mu'_2}(\tilde{e}_{j}^{\nu_2, \lambda_2 \mu_2} \otimes \id_{\mu_2'})}.
\end{align*}
By a similar calculation, we obtain
\begin{align*}
	&(\id_\theta \otimes \Delta) \Delta \\
	&= \bigoplus_{L,M,M',N'} \sum_{\nu_1,\nu_2,i,j,i',j'} \frac{d_{\lambda_2} d_{\mu_2} d_{\mu_2'}}{d_{\nu_2'}d_{\nu_1'}} \\
	&\quad \times \tr(\alpha_A^{G+}(\tilde{e}_{i'}^{\nu'_1, \lambda_1 \nu_1} (\id_{\lambda_1} \otimes \tilde{e}_i^{\nu_1, \mu_1 \mu'_1})) (\tilde{\phi}_l^{L} \otimes_A \tilde{\phi}_m^M \otimes_A \tilde{\phi}_{m'}^{M'}) \alpha_A^{G-}((\id_{\lambda_2} \otimes e^{\nu_2,\mu_2\mu'_2})e_{j'}^{\nu'_2, \lambda_2 \nu_2}) \phi^{N'}_{n'}) \\
	&\quad \times (\id_{\lambda_1} \otimes e_i^{\nu_1, \mu_1 \mu'_1}) e_{i'}^{\nu'_1, \lambda_1 \nu_1} \boxtimes \overline{\tilde{e}_{j'}^{\nu'_2, \lambda_2 \nu_2}(\id_{\lambda_2} \otimes \tilde{e}^{\nu_2,\mu_2\mu'_2})}.
\end{align*}
Since
\begin{align*}
	\{ (e_{i}^{\nu_k, \lambda_k \mu_k} \otimes \id_{\mu_k'}) e_{i'}^{\nu_k', \nu_k \mu'_k} \}_{\nu_k,i,i'} \quad \text{and} \quad \{ (\id_{\lambda_k} \otimes e_i^{\nu_k, \mu_k \mu'_k}) e_{i'}^{\nu'_k, \lambda_k \nu_k} \}_{\nu_k, i,i'}
\end{align*}
are both bases of $\hom(\nu'_k, \lambda_k \mu_k \mu'_k)$ and their dual bases are respectively given by
\begin{align*}
	\{ \tilde{e}_{i'}^{\nu_k', \nu_k \mu'_k} (\tilde{e}_{i}^{\nu_k, \lambda_k \mu_k} \otimes \id_{\mu_k'}) \}_{\nu_k,i,i'} \quad \text{and} \quad \{ \tilde{e}_{i'}^{\nu'_k, \lambda_k \nu_k}(\id_{\lambda_k} \otimes \tilde{e}_i^{\nu_k, \mu_k \mu'_k}) \}_{\nu_k, i,i'}
\end{align*}
for $k =1,2$, we can interchange these bases in the representations of $(\Delta \otimes \id_\theta) \Delta$ and $(\id_\theta \otimes \Delta) \Delta$ and therefore obtain $(\Delta \otimes \id_\theta) \Delta = (\id_\theta \otimes \Delta) \Delta$. We can also show $\tilde{m} \circ \Delta = d_\theta$ when $\calc$ is spherical as in the proof of \cite[Theorem 1.4]{MR1754521} by Proposition \ref{bimodpivotalprop}.

Finally, we show that $\Theta_\alpha^G(A)$ is a $G$-equivariant Frobenius algebra, which is essentially new. First, we show that $z_g$ is invertible for every $g \in G$. Indeed, we show that the inverse is given by
\begin{align*}
	z_g^{-1} = \bigoplus_{L,l'} \frac{1}{d_{\lambda_1}} \tr(\tilde{\eta}^{g-1}_{\lambda_1} {}^g \tilde{\phi}_l^L \tilde{\eta}^g_{\lambda_2} \alpha^{G-}_A(\tilde{u}_g^{\lambda_2}) \phi_{l'}^{L(g)} \alpha^{G+}_A(u_g^{\lambda_1})) (\tilde{u}_g^{\lambda_1} \boxtimes \overline{u_g^{\lambda_2}}) 
\end{align*}
By Lemma \ref{conjugatecomposition}, we have
\begin{align*}
z_g^{-1} z_g &= \bigoplus_L \sum_{l'} \frac{1}{d_{\lambda_1}^2} \tr(\tilde{\eta}^{g-1}_{\lambda_1} {}^g \tilde{\phi}_l^L \tilde{\eta}^g_{\lambda_2} \alpha^{G-}_A(\tilde{u}_g^{\lambda_2}) \phi_{l'}^{L(g)} \alpha^{G+}_A(u_g^{\lambda_1})) \\
&\hspace{1.5cm} \times \tr(\alpha^{G+}_A(\tilde{u}_g^{\lambda_1}) \tilde{\phi}_{l'}^{L(g)} \alpha^{G-}_A(u_g^{\lambda_2}) \tilde{\eta}^{g-1}_{\lambda_2} {}^g \phi_l^L \tilde{\eta}^g_{\lambda_2}).
\end{align*}
By the cyclicity of the trace, we have
\begin{align*}
&\frac{1}{d_{\lambda_1}} \tr(\alpha^{G+}_A(\tilde{u}_g^{\lambda_1}) \tilde{\phi}_{l'}^{L(g)} \alpha^{G-}_A(u_g^{\lambda_2}) \tilde{\eta}^{g-1}_{\lambda_2} {}^g \phi_l^L \tilde{\eta}^g_{\lambda_2}) = (\tilde{\phi}_{l'}^{L(g)}, \alpha^{G-}_A(u_g^{\lambda_2}) \tilde{\eta}^{g-1}_{\lambda_2} {}^g \phi_l^L \tilde{\eta}^g_{\lambda_2} \alpha^{G+}_A(\tilde{u}_g^{\lambda_1})),
\end{align*}
so that we obtain
\begin{align*}
z_g^{-1} z_g &= \bigoplus_L \frac{1}{d_{\lambda_1}} \tr({}^g(\tilde{\phi}_l^L \phi_l^L)) = \bigoplus_L (\tilde{\phi}_l^L, \phi_l^L) = \id_{{}^g \theta}
\end{align*}
by the Fourier expansion in the basis $\{ \phi_{l'}^{L(g)} \}_{l'}$ since the group action on $\bimod_\calc(A)$ is pivotal by Proposition \ref{bimodpivotalprop} and therefore preserves traces. We can similarly check $z_g z_g^{-1} = \id_{\theta}$. Thus, $z_g$ is invertible.

Next, we show that $(\theta, z)$ is a $G$-equivariant object. For any $g, h \in G$, we have
\begin{align*}
z_g {}^g z_h &= \bigoplus_{L,l''} \sum_{l'} \frac{1}{d_{\lambda_1}^2} \tr(\alpha^{G+}_A(\tilde{u}_g^{\lambda_1(h)}) \tilde{\phi}_{l''}^{L(gh)} \alpha^{G-}_A(u_g^{\lambda_2(h)})\tilde{\eta}^{g-1}_{\lambda_2(h)} {}^g \phi_{l'}^{L(h)} \tilde{\eta}^g_{\lambda_1(h)}) \\
&\hspace{1.5cm} \times \tr(\alpha^{G+}_A(\tilde{u}_h^{\lambda_1}) \tilde{\phi}_{l'}^{L(h)} \alpha^{G-}_A(u_h^{\lambda_2}) \tilde{\eta}^{h-1}_{\lambda_2} {}^h \phi_{l}^L \tilde{\eta}^h_{\lambda_1}) (u_g^{\lambda_1(h)} {}^g u_h^{\lambda_1} \boxtimes \overline{{}^g \tilde{u}_h^{\lambda_2} \tilde{u}_g^{\lambda_2(h)}})
\end{align*}
by Lemma \ref{conjugatecomposition}. Since $\tilde{\eta}^g_{{}^{h} \lambda_1} \alpha^{G+}_A({}^g \tilde{u}_h^{\lambda_1}) = {}^g \alpha^{G+}_A(\tilde{u}_h^{\lambda_1}) \tilde{\eta}^g_{\lambda_1(h)}$ and $\tilde{\eta}^g_{\lambda_2(h)} {}^g \alpha^{G-}_A(\tilde{u}_h^{\lambda_2}) = \alpha_A^{G-}({}^g \tilde{u}_h^{\lambda_2}) \tilde{\eta}^g_{{}^{h} \lambda_2}$ (see the proof of Proposition \ref{prop_alpha_is_crossed}), we have
\begin{align*}
	z_g {}^g z_h &= \bigoplus_{L,l''} \frac{1}{d_{\lambda_1}} \tr(\alpha_A^{G+}({}^g \tilde{u}_h^{\lambda_1} \tilde{u}_g^{\lambda_1(h)}) \tilde{\phi}_{l''}^{L(gh)} \alpha^{G-}_A(u_g^{\lambda_2(h)}{}^g u_h^{\lambda_2})\tilde{\eta}^{gh-1}_{\lambda_2} {}^{gh} \phi_{l}^L \tilde{\eta}^{gh}_{\lambda_1}) \\
	&\hspace{1.5cm} \times (u_g^{\lambda_1(h)} {}^g u_h^{\lambda_1} \boxtimes \overline{{}^g \tilde{u}_h^{\lambda_2} \tilde{u}_g^{\lambda_2(h)}})
\end{align*}
by the cyclicity of the trace and the Fourier expansion in the basis $\{ \phi_{l'}^{L(h)} \}_{l'}$. The right hand side is equal to $z_{gh}$ since $u_g^{\lambda_k(h)}{}^g u_h^{\lambda_k}$ is an isomorphism ${}^{gh} \lambda_k \cong \lambda_k(gh)$ for $k=1,2$ and $z$ does not depend on a choice of isomorphisms. Therefore, $(\theta, z)$ is a $G$-equivariant object.

Finally, we show that $z_g$ is a Frobenius algebra homomorphism for every $g \in G$. Since we may take a basis of $\hom(\alpha^{G+}_A(\mathbf{1}_\calc), \alpha^{G-}_A(\mathbf{1}_\calc))$ to be the identity and take $u_g^{\mathbf{1}_\calc}$'s to be canonical isomorphisms, we have $z_g {}^g \eta = \eta$ and $\varepsilon z_g = {}^g \varepsilon$ for every $g \in G$. We show $\Delta z_g = (z_g \otimes z_g) {}^g \Delta$. We have
\begin{align*}
\Delta z_g &= \bigoplus_{L,M,N} \sum_{n',i,j} \frac{d_{\lambda_2} d_{\mu_2}}{d_{\nu_2} d_{\nu_1^2}} \tr(\alpha^{G+}_A(\tilde{e}_i^{\nu_1(g), \lambda_1 \mu_1})(\tilde{\phi}_l^{L} \otimes_A \tilde{\phi}_m^{M}) \alpha^{G-}_A(e_j^{\nu_2(g), \lambda_2 \mu_2})\phi^{N(g)}_{n'})\\
&\hspace{2cm} \times \tr(\alpha^{G+}_A(\tilde{u}_g^{\nu_1}) \tilde{\phi}_{n'}^{N(g)} \alpha^{G-}_A(u_g^{\nu_2}) \tilde{\eta}^{g-1}_{\nu_2} {}^g \phi^N_{n} \tilde{\eta}^g_{\nu_1})(e_i^{\nu_1(g), \lambda_1 \mu_1} u_g^{\nu_1} \boxtimes \overline{\tilde{u}_g^{\nu_2} \tilde{e}_j^{\nu_2(g), \lambda_2 \mu_2}}) \\
&= \bigoplus_{L,M,N} \sum_{i,j} \frac{d_{\lambda_2} d_{\mu_2}}{d_{\nu_2} d_{\nu_1}} \tr(\alpha^{G+}_A(\tilde{u}_g^{\nu_1} \tilde{e}_i^{\nu_1(g), \lambda_1 \mu_1})(\tilde{\phi}_l^{L} \otimes_A \tilde{\phi}_m^{M}) \alpha^{G-}_A(e_j^{\nu_2(g), \lambda_2 \mu_2}u_g^{\nu_2}) \tilde{\eta}^{g-1}_{\nu_2} {}^g \phi^N_{n} \tilde{\eta}^g_{\nu_1})\\
&\hspace{2cm} \times (e_i^{\nu_1(g), \lambda_1 \mu_1} u_g^{\nu_1} \boxtimes \overline{\tilde{u}_g^{\nu_2} \tilde{e}_j^{\nu_2(g), \lambda_2 \mu_2}})
\end{align*}
by Lemma \ref{conjugatecomposition}, the cyclicity of the trace and the Fourier expansion in the basis $\{ \phi_{n'}^{N(g)} \}_{n'}$. On the other hand, we have 
\begin{align*}
&(z_g \otimes z_g){}^g \Delta \\
&= \bigoplus_{L,M,N} \sum_{l',m',i,j} \frac{d_{\lambda_2} d_{\mu_2}}{d_{\nu_2}d_{\lambda_1} d_{\mu_1} d_{\nu_1}} \\
&\quad \times \tr(\alpha^{G+}_A(\tilde{u}_g^{\lambda_1(g^{-1})}) \tilde{\phi}_{l}^{L} \alpha^{G-}_A(u_g^{\lambda_2(g^{-1})})\tilde{\eta}^{g-1}_{\lambda_2(g^{-1})} {}^g \phi_{l'}^{L(g^{-1})} \tilde{\eta}^g_{\lambda_1(g^{-1})}) \\
&\quad \times \tr(\alpha^{G+}_A(\tilde{u}_g^{\mu_1(g^{-1})}) \tilde{\phi}_{m}^{M} \alpha^{G-}_A(u_g^{\mu_2(g^{-1})}) \tilde{\eta}^{g-1}_{\mu_2(g^{-1})} {}^g \phi_{m'}^{M(g^{-1})} \tilde{\eta}^g_{\mu_1(g^{-1})}) \\
&\quad \times \tr(\alpha^{G+}_A(\tilde{e}_j^{\nu_1, \lambda_1(g^{-1}) \mu_1(g^{-1})}) (\tilde{\phi}_{l'}^{L(g^{-1})} \otimes_A \tilde{\phi}_{m'}^{M(g^{-1})}) \alpha^{G-}_A(e_j^{\nu_2, \lambda_2(g^{-1}) \mu_2(g^{-1})}) \phi^N_n) \\
&\quad \times (u_g^{\lambda_1(g^{-1})} \otimes u_g^{\mu_1(g^{-1})}) {}^g e_i^{\nu_1, \lambda_1(g^{-1}) \mu_1(g^{-1})} \boxtimes \overline{{}^g \tilde{e}_j^{\nu_2, \lambda_2(g^{-1}) \mu_2(g^{-1})} (\tilde{u}_g^{\lambda_2(g^{-1})} \otimes \tilde{u}_g^{\mu_2(g^{-1})})}
\end{align*}
by $\partial \mu_1 = \partial \mu_2$ and Lemmata \ref{conjugatecomposition} and \ref{conjugatemonoidalproduct}. Since $\tilde{\eta}^g$ is monoidal by the proof of Proposition \ref{prop_alpha_is_crossed} and
\begin{align*}
	{}^g \alpha^{G+}_A(\tilde{e}_j^{\nu_1, \lambda_1(g^{-1}) \mu_1(g^{-1})}) \tilde{\eta}^g_{\lambda_1(g^{-1}) \mu_1(g^{-1})} &= \tilde{\eta}^g_{\nu_1} \alpha^{G+}_A({}^g \tilde{e}_j^{\nu_1, \lambda_1(g^{-1}) \mu_1(g^{-1})}) \\
	\tilde{\eta}^{g-1}_{\lambda_2(g^{-1}) \mu_2(g^{-1})} {}^g \alpha^{G-}_A(e_j^{\nu_2, \lambda_2(g^{-1}) \mu_2(g^{-1})}) &= \alpha^{G-}_A({}^g e_j^{\nu_2, \lambda_2(g^{-1}) \mu_2(g^{-1})}) \tilde{\eta}^{g-1}_{\nu_2}
\end{align*}
by definition, we obtain
\begin{align*}
	&(z_g \otimes z_g){}^g \Delta \\
	&= \bigoplus_{L,M,N} \sum_{i,j} \frac{d_{\lambda_2} d_{\mu_2}}{d_{\nu_2} d_{\nu_1}} \\
	&\hspace{2cm} \times \tr(\alpha^{G+}_A({}^g \tilde{e}_j^{\nu_1, \lambda_1(g^{-1}) \mu_1(g^{-1})}(\tilde{u}_g^{\lambda_1(g^{-1})} \otimes \tilde{u}_g^{\mu_1(g^{-1})}))(\tilde{\phi}_{l}^{L} \otimes_A \tilde{\phi}_{m}^{M}) \\
	&\hspace{3cm} \alpha^{G-}_A((u_g^{\lambda_2(g^{-1})} \otimes u_g^{\mu_2(g^{-1})}) {}^g e_j^{\nu_2, \lambda_2(g^{-1}) \mu_2(g^{-1})}) \tilde{\eta}^{g-1}_{\nu_2} {}^g \phi^N_n \tilde{\eta}^g_{\nu_1}) \\
	&\hspace{2cm} \times (u_g^{\lambda_1(g^{-1})} \otimes u_g^{\mu_1(g^{-1})}) {}^g e_i^{\nu_1, \lambda_1(g^{-1}) \mu_1(g^{-1})} \boxtimes \overline{{}^g \tilde{e}_j^{\nu_2, \lambda_2(g^{-1}) \mu_2(g^{-1})} (\tilde{u}_g^{\lambda_2(g^{-1})} \otimes \tilde{u}_g^{\mu_2(g^{-1})})}
\end{align*}
by the cyclicity of the trace and the Fourier expansion in these bases $\{ \phi^{L(g^{-1})}_{l'} \}_{l'}$ and $\{ \phi^{M(g^{-1})}_{m'} \}_{m'}$. Then, comparing the bases
\begin{align*}
	\{ e_i^{\nu_k(g), \lambda_k \mu_k} u_g^{\nu_k} \}_i \quad \text{and} \quad \{(u_g^{\lambda_k(g^{-1})} \otimes u_g^{\mu_k(g^{-1})}) {}^g e_i^{\nu_k, \lambda_k(g^{-1}) \mu_k(g^{-1})} \}_i
\end{align*}
of $\hom({}^g \nu_k, \lambda_k \mu_k)$ for $k=1,2$, we obtain the conclusion. We can similarly show $z_g m = m (z_g \otimes z_g)$. Thus, $\Theta^G_\alpha(A)$ is a $G$-equivariant Frobenius algebra. We can similarly prove the statements for Q-systems just by taking bases to be orthonormal and replacing dual bases by ${}^\ast$-conjugated bases.
\end{proof}
\end{thm}

Next, we show that $\Theta_\alpha^G(A)$ satisfies the $G$-equivariant version of commutativity. 

\begin{defi}
	Let $\calc$ be a $G$-braided multitensor category $\calc$. Then, a $G$-equivariant algebra $A$ in $\calc$ is \emph{$G$-commutative} if the equation in Figure \ref{graphical_gcomm} holds. Similarly, when $A$ is a $G$-equivariant Frobenius algebra, it is \emph{$G$-cocommutative} if the equation in Figure \ref{graphicalgcocommutativity} holds.
	
	\begin{figure}[htb]
		\begin{tabular}{cc}
			\begin{minipage}[b]{0.45\hsize}
				\centering
				\begin{tikzpicture}
					\draw (0.5,3) -- (-0.5,2) -- (-0.5,1) arc (180:360:0.5);
					\draw[dashed,cross] (0.5,2) -- (-0.5,3);
					\draw (0.5,2)--(0.5,1);
					\draw[<-] (0,0) -- (0,0.5);
					\node[block] at (-0.5,1.5){$z_g^A$};
					\node at (0,-0.25){$A$};
					\node at (-0.5,3.25){$A$};
					\node at (0.5,3.25){$A$};
					\node at (0,3){$g$};
					\node at (-1.5,1.5){$\displaystyle \sum_{g \in G}$};
					\node at (-0.75,2.25){${}^g A$};
					\node at (1,1.5){$=$};
					\begin{scope}[shift={(2,1)}]
						\draw[<-] (0,0) -- (0,0.5);
						\draw (-0.5,1) arc (180:360:0.5);
						\node at (0,-0.25){$A$};
						\node at (0.5,1.25){$A$};
						\node at (-0.5,1.25){$A$};
					\end{scope}
				\end{tikzpicture}
				\caption{$G$-commutativity}
				\label{graphical_gcomm}
			\end{minipage}
			\begin{minipage}[b]{0.45\hsize}
				\centering
				\begin{tikzpicture}
				\draw[->] (1,0) -- (0,-1) -- (0,-2);
				\draw[-,dashed,cross] (0,0) -- (1,-1);
				\draw[-] (1,0) arc (0:180:0.5);
				\draw[-] (0.5,1) -- (0.5,0.5);
				\draw[->] (1,-1) -- (1,-2);
				\node[block] at (0,-1.5){$z^A_g$};
				\node at (0.5,1.25) {$A$};
				\node at (-0.25,-0.75) {${}^g A$};
				\node at (0.5,0) {$g$};
				\node at (0,-2.25) {$A$};
				\node at (1,-2.25) {$A$};
				\node at (-1,-0.5){$\displaystyle \sum_{g \in G}$};
				\node at (1.5,-0.5) {$=$};
				\begin{scope}[shift={(2,-1)}]
				\draw[<->] (1,0) arc (0:180:0.5);
				\draw[-] (0.5,1) -- (0.5,0.5);
				\node at (0.5,1.25) {$A$};
				\node at (1,-0.25) {$A$};
				\node at (0,-0.25) {$A$};
				\end{scope}
				\end{tikzpicture}
				\caption{$G$-cocommutativity}
				\label{graphicalgcocommutativity}
			\end{minipage}
		\end{tabular}
	\end{figure}
\end{defi}	

\begin{lem}
	\label{lem_alpha_braiding}
	Let $A$ be a neutral special $G$-equivariant Frobenius algebra in a $G$-braided multitensor category $\calc$. Let $\lambda, \lambda' \in \obj(\calc_g)$ and $\mu,\mu' \in \obj(\calc)$. Then, for $f \in \hom(\alpha^{G+}_A(\lambda),\alpha^{G-}_A(\lambda'))$ and $f' \in \hom(\alpha^{G+}_A(\mu), \alpha^{G-}_A(\mu'))$, we have
	\begin{align*}
		(\id_A \otimes_\calc b^\calc_{\lambda',\mu'}) (f \otimes_A f') &= (\tilde{\eta}_{\mu'}^{g-1} {}^g f' \tilde{\eta}_\mu^g \otimes_A f) (\id_A \otimes_\calc b^\calc_{\lambda,\mu}) \\
		(f \otimes_A f') (\id_A \otimes_\calc b^{\calc-}_{\mu,\lambda}) &= (\id_A \otimes_\calc b^{\calc-}_{\mu',\lambda'}) (\tilde{\eta}^{g-1}_{\mu'} {}^g f' \tilde{\eta}^g_\mu \otimes_A f),
	\end{align*}
	where $\tilde{\eta}^g_\lambda \coloneqq z_g^{A-1} \otimes_\calc \id_{{}^g \lambda}$.

	\begin{proof}
		The statement for $b^{\calc-}$ follows from the graphical calculations in Figure \ref{graphical_alpha_braiding}, where $z \coloneqq z^A$. We used the right $A$-modularity of $f$ and the left $A$-modularity of $\tilde{\eta}_{\mu'}^{g-1} {}^g f' \tilde{\eta}_\mu^g$ respectively at the second and third equalities in the upper equation. We used the left $A$-modularity of $f$ and the right $A$-modularity of $\tilde{\eta}_{\mu'}^{g-1} {}^g f' \tilde{\eta}_\mu^g$ respectively at the second and third equalities in the lower equation. The proof for $b^\calc$ is similar.
	\end{proof}

	\begin{figure}[htb]
		\centering
		\begin{tikzpicture}
		\draw (0,1) -- (0,0.75);
		\draw (-0.25,-1.5)-- (-0.25,0.5) arc (180:0:0.25) -- (0.75,0) -- (0.75,-1.5);
		\draw (0.75,1) -- (1.35,0.4) -- (1.35,-1.5);
		\draw[cross] (1.25,1) -- (0.15,-0.1);
		\draw (0.15,-0.1) -- (0.15,-1.5);
		\node[block] at (0.75,-0.5){$z_{g^{-1}}$};
		\node at (0,1.25){$A$};
		\node at (0.75,1.25){${}^g \mu$};
		\node at (1.25,1.25){$\lambda$};
			\draw[->] (0.25,-1.5) -- (0.25,-2) -- (0.75,-2.5) -- (0.75,-3);
			\draw[cross] (0.75,-1.5) -- (0.75,-2) -- (0.25,-2.5) arc (360:180:0.25) -- (-0.25,-1.5);
			\draw[->] (0,-2.75) -- (0,-3);
			\draw[->] (1.25,-1.5) -- (1.25,-3);
			\node at (1.25,-3.25){$\mu'$};
			\node at (0,-3.25){$A$};
			\node at (0.75,-3.25){$\lambda'$};
		\draw[fill=white] (-0.5,-1.25) rectangle (0.35,-1.75) node[midway]{$f$};
		\draw[fill=white] (0.5,-1.25) rectangle (1.5,-1.75) node[midway]{$f'$};
		\node at (2,-1){$=$};
		\begin{scope}[shift={(3,0)}]
			\draw (0,1) -- (0,0.75);
			\draw (-0.25,-1.5)-- (-0.25,0.5) arc (180:0:0.25) -- (0.25,-0.25);
			\draw (0.75,1) -- (0.75,-0.5);
			\draw (0.25,-0.75) -- (0.75,-1.25);
			\draw[->] (0.75,-0.75) -- (1.25,-1.25) --(1.25,-1.5) -- (1.25,-3);
			\draw[cross] (1.25,1) -- (1.25,-0.75) -- (0.25,-1.25);
			\node[block] at (0.25,0.25){$z_{g}^{-1}$};
			\node at (0,1.25){$A$};
			\node at (0.75,1.25){${}^g \mu$};
			\node at (1.25,1.25){$\lambda$};
				\draw[->] (0.25,-1.5) -- (0.25,-2) -- (0.75,-2.5) -- (0.75,-3);
				\draw[cross] (0.75,-1.5) -- (0.75,-2) -- (0.25,-2.5) arc (360:180:0.25) -- (-0.25,-1.5);
				\draw (0.75,-1.25) -- (0.75,-1.5);
				\draw[->] (0,-2.75) -- (0,-3);
				\node at (1.25,-3.25){$\mu'$};
				\node at (0,-3.25){$A$};
				\node at (0.75,-3.25){$\lambda'$};
			\draw[fill=white] (-0.5,-1.25) rectangle (0.35,-1.75) node[midway]{$f$};
			\draw[fill=white] (0,-0.25) rectangle (1,-0.75) node[midway]{${}^g f'$};
			\node at (2,-1){$=$};
		\end{scope}
		\begin{scope}[shift={(6,0)}]
			\draw (0,1) -- (0,0.75);
			\draw (-0.25,-1.75)-- (-0.25,0.5) arc (180:0:0.25) -- (0.25,-0.25) -- (0.25,-1.75) arc (360:180:0.25);
			\draw (0.75,1) -- (0.75,-0.5);
			\draw[->] (0.75,-0.75) -- (1.25,-1.25) --(1.25,-1.5) -- (1.25,-3);
			\draw[cross,->] (1.25,1) -- (1.25,-0.75) -- (0.75,-1.25) -- (0.75,-3);
			\node[block] at (0.25,0.25){$z_{g}^{-1}$};
			\node[block] at (0.25,-1.5){$z_g$};
			\node at (0,1.25){$A$};
			\node at (0.75,1.25){${}^g \mu$};
			\node at (1.25,1.25){$\lambda$};
				\draw (0.75,-1.25) -- (0.75,-1.5);
				\draw[->] (0,-2) -- (0,-3);
				\node at (1.25,-3.25){$\mu'$};
				\node at (0,-3.25){$A$};
				\node at (0.75,-3.25){$\lambda'$};
			\draw[fill=white] (-0.25,-2.25) rectangle (1,-2.75) node[midway]{$f$};
			\draw[fill=white] (0,-0.25) rectangle (1,-0.75) node[midway]{${}^g f'$};
			\node at (2,-1){$=$};
		\end{scope}
		\begin{scope}[shift={(9,0)}]
			\draw[->] (0,1) -- (0,-3);
			\draw (0.75,1) -- (0.75,-0.5);
			\draw[->] (0.75,-0.75) -- (1.25,-1.25) --(1.25,-1.5) -- (1.25,-3);
			\draw[cross,->] (1.25,1) -- (1.25,-0.75) -- (0.75,-1.25) -- (0.75,-3);
			\node[block] at (0,0.5){$z_{g}^{-1}$};
			\node at (0,1.25){$A$};
			\node at (0.75,1.25){${}^g \mu$};
			\node at (1.25,1.25){$\lambda$};
				\draw (0.75,-1.25) -- (0.75,-1.5);
				\node at (1.25,-3.25){$\mu'$};
				\node at (0,-3.25){$A$};
				\node at (0.75,-3.25){$\lambda'$};
			\draw[fill=white] (-0.25,-2.25) rectangle (1,-2.75) node[midway]{$f$};
			\draw[fill=white] (-0.25,-0.25) rectangle (1,-0.75) node[midway]{${}^g f'$};
			\node[block] at (0,-1.5){$z_g$};
		\end{scope}
		\begin{scope}[shift={(0,-5)}]
			\draw (0,1) -- (0,0.75);
			\draw (-0.25,-1.5)-- (-0.25,0.5) arc (180:0:0.25) -- (0.75,0) -- (0.75,-1.5);
			\draw (1.35,1) -- (1.35,-1.5);
			\draw[dashed,cross] (0.75,1) -- (0.75,0.5) -- (0.15,-0.1);
			\draw (0.15,-0.1) -- (0.15,-1.5);
			\node[block] at (0.75,-0.5){$z_{h^{-1}}$};
			\node at (0,1.25){$A$};
			\node at (0.5,0.75){$h$};
			\node at (0.75,1.25){${}^g \mu$};
			\node at (1.35,1.25){$\lambda$};
				\draw[->] (0.25,-1.5) -- (0.25,-2) -- (1.25,-3);
				\draw[cross] (0.75,-1.5) -- (0.75,-2) -- (0.25,-2.5) arc (360:180:0.25) -- (-0.25,-1.5);
				\draw[->] (0,-2.75) -- (0,-3);
				\draw[->,cross] (1.25,-1.5) -- (1.25,-2.5) -- (0.75,-3);
				\node at (1.25,-3.25){$\mu'$};
				\node at (0,-3.25){$A$};
				\node at (0.75,-3.25){$\lambda'$};
				\draw[fill=white] (-0.5,-1.25) rectangle (0.35,-1.75) node[midway]{${}^g f'$};
				\draw[fill=white] (0.5,-1.25) rectangle (1.5,-1.75) node[midway]{$f$};
				\node at (2,-1){$=$};
				\node[block] at (-0.35,0){$z_g^{-1}$};
				\node[block] at (-0.25,-2.25){$z_g$};
				\node at (-1,-1){$\displaystyle \sum_h$};
		\end{scope}
		\begin{scope}[shift={(3.75,-5)}]
			\draw (0,1) -- (0,0.75);
			\draw (-0.25,-1.5)-- (-0.25,0.5) arc (180:0:0.25) -- (0.75,0) -- (0.75,-0.5);
			\draw[dashed,cross] (0.75,1) -- (0.75,0.5) -- (0,-0.25);
			\node at (0,1.25){$A$};
			\node at (0.5,0.75){$h$};
			\node at (0.75,1.25){${}^g \mu$};
			\node at (1.35,1.25){$\lambda$};
				\draw[->] (0,-0.75) -- (1.25,-2) -- (1.25,-3);
				\draw[cross] (0.75,-0.5) -- (0.75,-1) -- (0.25,-1.5) arc (360:180:0.25) -- (-0.25,-1.5);
				\draw[->] (0,-1.75) -- (0,-3);
				\draw[->,cross] (1.35,1) -- (1.35,-1.15) -- (0.75,-1.75) -- (0.75,-3);
				\node at (1.25,-3.25){$\mu'$};
				\node at (0,-3.25){$A$};
				\node at (0.75,-3.25){$\lambda'$};
				\draw[fill=white] (-0.5,-0.25) rectangle (0.25,-0.75) node[midway]{${}^g f'$};
				\draw[fill=white] (-0.25,-2) rectangle (1,-2.5) node[midway]{$f$};
				\node at (2,-1){$=$};
				\node[block] at (0.75,-0.5){$z_{h^{-1}}$};
				\node[block] at (-0.35,0.25){$z_g^{-1}$};
				\node[block] at (-0.25,-1.25){$z_g$};
				\node at (-1,-1){$\displaystyle \sum_h$};
		\end{scope}
		\begin{scope}[shift={(7.5,-6.5)}]
			\draw (0,2.5) -- (0,2.25) arc (90:180:0.25) -- (-0.25,0.5) arc (180:360:0.25) -- (0.25,1) -- (0.75,1.25) -- (0.75,1.75) -- (0.25,2) arc (0:90:0.25);
			\draw[cross,dashed] (0.75,2.5) -- (0.75,2) -- (0,1.75) -- (0,1.25) -- (0.75,1);
			\draw[->] (0,0.25) -- (0,-3);
			\draw[->] (0.75,1) -- (0.75,-1.5) -- (1.25,-2) -- (1.25,-3);
			\draw[cross,->] (1.25,2.5) -- (1.25,-1.5) -- (0.75,-2) -- (0.75,-3);
			\node[block] at (0,-0.25){$z_{g}^{-1}$};
			\node at (0,2.75){$A$};
			\node at (0.75,2.75){${}^g \mu$};
			\node at (1.25,2.75){$\lambda$};
				\draw (0.75,-1.25) -- (0.75,-1.5);
				\node at (1.25,-3.25){$\mu'$};
				\node at (0,-3.25){$A$};
				\node at (0.75,-3.25){$\lambda'$};
			\draw[fill=white] (-0.25,-2.25) rectangle (1,-2.75) node[midway]{$f$};
			\draw[fill=white] (-0.25,-0.75) rectangle (1,-1.25) node[midway]{${}^g f'$};
			\node[block] at (0,-1.75){$z_g$};
			\node[block] at (0.25,0.75){$z_h$};
			\node[block] at (0.75,1.5){$z_{h^{-1}}$};
			\node at (0.5,2.25){$h$};
			\node at (-1,0.5){$\displaystyle \sum_h$};
			\node at (2,0.5){$=$};
		\end{scope}
		\begin{scope}[shift={(10.5,-5)}]
			\draw[->] (0,1) -- (0,-3);
			\draw (0.75,1) -- (0.75,-0.5);
			\draw[->] (0.75,-0.75) -- (1.25,-1.25) --(1.25,-1.5) -- (1.25,-3);
			\draw[cross,->] (1.25,1) -- (1.25,-0.75) -- (0.75,-1.25) -- (0.75,-3);
			\node[block] at (0,0.5){$z_{g}^{-1}$};
			\node at (0,1.25){$A$};
			\node at (0.75,1.25){${}^g \mu$};
			\node at (1.25,1.25){$\lambda$};
				\draw (0.75,-1.25) -- (0.75,-1.5);
				\node at (1.25,-3.25){$\mu'$};
				\node at (0,-3.25){$A$};
				\node at (0.75,-3.25){$\lambda'$};
			\draw[fill=white] (-0.25,-2.25) rectangle (1,-2.75) node[midway]{$f$};
			\draw[fill=white] (-0.25,-0.25) rectangle (1,-0.75) node[midway]{${}^g f'$};
			\node[block] at (0,-1.5){$z_g$};
		\end{scope}
		\end{tikzpicture}
		\caption{A crossing and equivariant $\alpha$-induction}
		\label{graphical_alpha_braiding}
	\end{figure}
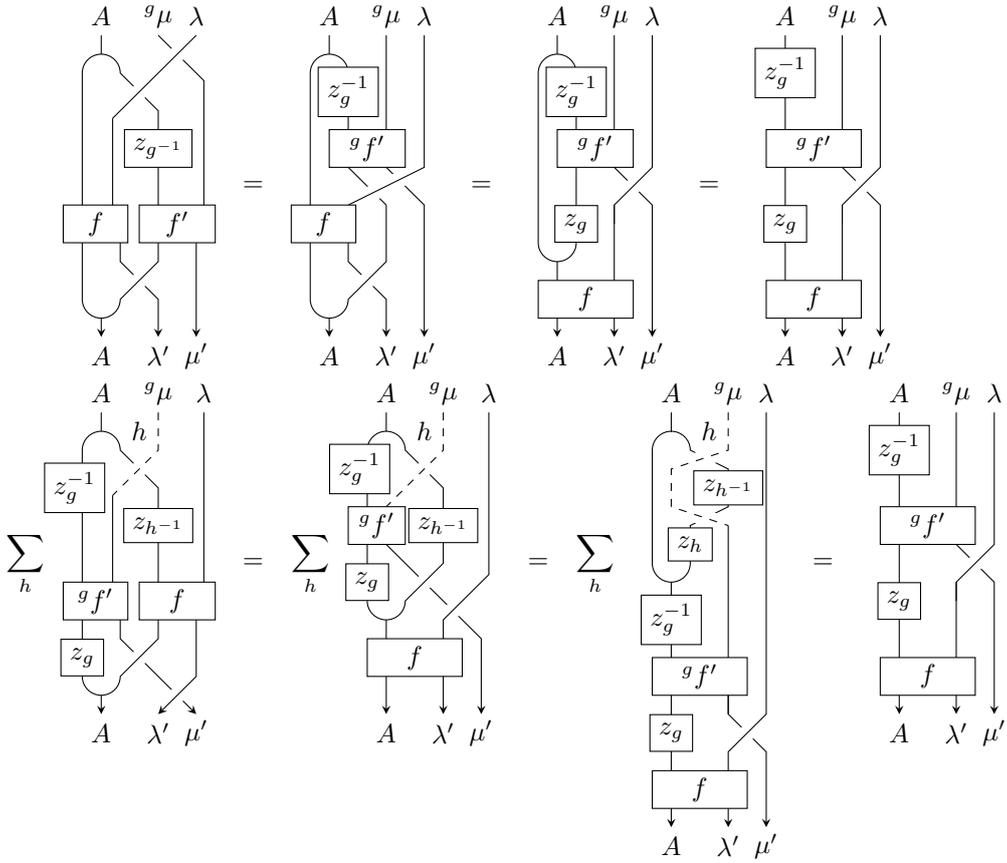

\end{lem}

\begin{lem}
\label{conjugatebraiding}
Let $\calc$ be a $G$-braided multitensor category. Then, $\overline{b^{\calc-}_{\mu, \lambda}} = b^{\calc-}_{{}^{\lambda^\vee, \mu^\vee}}$ for any $\lambda, \mu \in \homog(\calc)$.

\begin{proof}
The statement follows from the graphical calculation in Figure \ref{graphical_conj_revbr}.
\begin{figure}[htb]
\centering
\begin{tikzpicture}
\draw[-] (0,-1) -- (-1,0) -- (-0.5,0.5) arc (180:0:0.25) -- (0.5,0) -- (0.5,-1.5);
\draw[-,cross] (1,-1.5) -- (1,0.5) arc (0:180:0.25) -- (-1.5,-1.5);
\draw[->] (-1.5,-1.5) arc (360:180:0.25) -- (-2,0.5);
\draw[->,cross] (-0.5,-1.5) arc (360:180:0.25) --(-1.5,-1) -- (-1.5,0.5);
\draw[-] (-0.5,-1.5) -- (0,-1);
\node at (-2,0.75){${}^{\overline{\mu}}\lambda^\vee$};
\node at (-1.5,0.75){$\mu^\vee$};
\node at (0.5,-1.75){$\mu^\vee$};
\node at (1,-1.75){$\lambda^\vee$};
\node at (2,-0.5){$=$};
\begin{scope}[shift={(4,0)}]
\draw[<-] (-1,0) -- (0,-1);
\draw[<-,cross] (0,0) -- (-1,-1);
\node at (-1,-1.25){$\mu^\vee$};
\node at (0,-1.25){$\lambda^\vee$};
\node at (-1,0.25){${}^{\overline{\mu}}\lambda^\vee$};
\node at (0,0.25){$\mu^\vee$};
\end{scope}
\end{tikzpicture}
\caption{The conjugation of a reverse crossing}
\label{graphical_conj_revbr}
\end{figure}
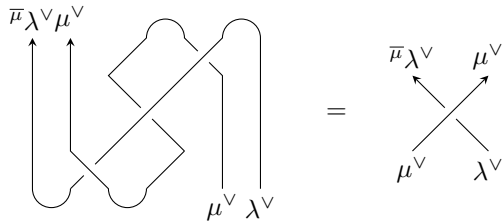
\end{proof}
\end{lem}

\begin{prop}
\label{prop_alphafrob_gcomm}
Let $A$ be a neutral symmetric special simple $G$-equivariant Frobenius algebra in a split $G$-braided fusion category $\calc$. Then, the $G$-equivariant $\alpha$-induction Frobenius algebra $\Theta_\alpha^G(A)$ associated with $A$ is $G$-commutative and $G$-cocommutative.

\begin{proof}
Let us follow the notation in Theorem \ref{mainthm1}. We only show cocommutativity because the proof of commutativity is similar. First, note that $\theta = \bigoplus_{g \in G} \bigoplus_{L, \partial \lambda_1 = g} \lambda_1 \boxtimes \lambda_2^\vee$ is a homogeneous decomposition. Therefore, the left-hand side of Figure \ref{graphicalgcocommutativity} for $\theta$ is given by
\begin{align*}
&\bigoplus_{L,M,N} \sum_{l',i,j} \frac{d_{\lambda_2} d_{\mu_2}}{d_{\nu_2} d_{\lambda_1} d_{\nu_1}} \tr(\alpha^{G+}_A(\tilde{u}_{\partial \mu_1}^{\lambda_1(\partial \overline{\mu_1})}) \tilde{\phi}_l^{L} \alpha^{G-}_A(u_{\partial \mu_1}^{\lambda_2(\partial \overline{\mu_1})})\tilde{\eta}^{\partial \mu_1-1}_{\lambda_2(\partial \overline{\mu_1})} {}^{\mu_1} \phi_{l'}^{L(\partial \overline{\mu_1})} \tilde{\eta}^{\partial \mu_1}_{\lambda_1(\partial \overline{\mu_1})})\\
&\times \tr(\alpha^{G+}_A(\tilde{e}_i^{\nu_1, \mu_1 (\lambda_1(\partial \overline{\mu_1}))}) (\tilde{\phi}_m^{M} \otimes_A \tilde{\phi}_{l'}^{L(\partial \overline{\mu_1})}) \alpha^{G-}_A(e_j^{\nu_2, \mu_2 (\lambda_2(\partial \overline{\mu_1}))}) \phi_n^N) \\
&\times ((u_{\partial \mu_1}^{\lambda_1(\partial \overline{\mu_1})} \otimes \id_{\mu_1}) b^\calc_{\mu_1, \lambda_1(\partial \overline{\mu_1})} e_i^{\nu_1, \mu_1 (\lambda_1(\partial \overline{\mu_1}))} \boxtimes \overline{\tilde{e}_j^{\nu_2, \mu_2 (\lambda_2(\partial \overline{\mu_1}))} b^{\calc-}_{\lambda_2(\partial \overline{\mu_1}),\mu_2} (\tilde{u}_{\partial \mu_1}^{\lambda_2(\partial \overline{\mu_1})} \otimes \id_{\mu_2})}) 
\end{align*}
by $\partial \mu_1 = \partial \mu_2$, $\partial \lambda_1(\partial \overline{\mu_1}) = \partial \lambda_2(\partial \overline{\mu_1})$ and Lemmata \ref{conjugatecomposition}, \ref{conjugatemonoidalproduct} and \ref{conjugatebraiding}. By the cyclicity of the trace and the Fourier expansion in the basis $\{ \tilde{\phi}_{l'}^{L(\partial \overline{\mu_1})} \}_{l'}$, we get
\begin{align*}
&\sum_{l'} \frac{1}{d_{\lambda_1}} \tr(\alpha^{G+}_A(\tilde{u}_{\partial \mu_1}^{\lambda_1(\partial \overline{\mu_1})}) \tilde{\phi}_l^{L} \alpha^{G-}_A(u_{\partial \mu_1}^{\lambda_2(\partial \overline{\mu_1})}) \tilde{\eta}^{\partial \mu_1-1}_{\lambda_2(\partial \overline{\mu_1})} {}^{\mu_1} \phi_{l'}^{L(\partial \overline{\mu_1})} \tilde{\eta}^{\partial \mu_1}_{\lambda_1(\partial \overline{\mu_1})}) \\
&\hspace{1cm} \times \tr(\alpha^{G+}_A(\tilde{e}_i^{\nu_1, \mu_1 (\lambda_1(\partial \overline{\mu_1}))}) (\tilde{\phi}_m^{M} \otimes_A \tilde{\phi}_{l'}^{L(\partial \overline{\mu_1})}) \alpha^{G-}_A(e_j^{\nu_2, \mu_2 (\lambda_2(\partial \overline{\mu_1}))}) \phi_n^N) \\
&= \tr(\alpha^{G+}_A(\tilde{e}_i^{\nu_1, \mu_1 (\lambda_1(\partial \overline{\mu_1}))}) (\tilde{\phi}_m^{M} \otimes_A {}^{\overline{\mu_1}}(\tilde{\eta}^{\partial \mu_1}_{\lambda_1(\partial \overline{\mu_1})} \alpha^{G+}_A(\tilde{u}_{\partial \mu_1}^{\lambda_1(\partial \overline{\mu_1})}) \tilde{\phi}_{l}^{L} \alpha_A^{G-}(u_{\partial \mu_1}^{\lambda_2(\partial \overline{\mu_1})})\tilde{\eta}^{\partial \mu_1-1}_{\lambda_2(\partial \overline{\mu_1})})) \\
&\hspace{1cm} \alpha^{G-}_A(e_j^{\nu_2, \mu_2 (\lambda_2(\partial \overline{\mu_1}))}) \phi_n^N).
\end{align*}
Then, since 
\begin{align*}
&(\tilde{\phi}_m^{M} \otimes_A {}^{\overline{\mu_1}}(\tilde{\eta}^{\partial \mu_1}_{\lambda_1(\partial \overline{\mu_1})} \alpha^{G+}_A(\tilde{u}_{\partial \mu_1}^{\lambda_1(\partial \overline{\mu_1})}) \tilde{\phi}_{l}^{L} \alpha_A^{G-}(u_{\partial \mu_1}^{\lambda_2(\partial \overline{\mu_1})})\tilde{\eta}^{\partial \mu_1-1}_{\lambda_2(\partial \overline{\mu_1})})) \alpha^{G-}_A (b^{\calc-}_{\lambda_2(\partial \overline{\mu_1}),\mu_2}) \\
&= \alpha^{G+}_A (b^{\calc-}_{\lambda_1(\partial \overline{\mu_1}),\mu_1})((\alpha^{G+}_A(\tilde{u}_{\partial \mu_1}^{\lambda_1(\partial \overline{\mu_1})}) \tilde{\phi}_{l}^{L} \alpha_A^{G-}(u_{\partial \mu_1}^{\lambda_2(\partial \overline{\mu_1})})) \otimes_A \tilde{\phi}_m^{M})
\end{align*}
by $\partial \mu_1 = \partial \mu_2$ and Lemma \ref{lem_alpha_braiding}, the left hand side of Figure \ref{graphicalgcocommutativity} is equal to
\begin{align*}
&\bigoplus_{L,M,N} \sum_{i,j} \frac{d_{\lambda_2} d_{\mu_2}}{d_{\nu_2} d_{\nu_1}} \tr( \alpha^{G+}_A(\tilde{e}_i^{\nu_1, \mu_1 \lambda_1(\partial \overline{\mu_1})} b^{\calc-}_{\lambda_1(\partial \overline{\mu_1}),\mu_1} (\tilde{u}_{\partial \mu_1}^{\lambda_1(\partial \overline{\mu_1})} \otimes \id_{\mu_1})) (\tilde{\phi}_l^{L}\otimes_A \tilde{\phi}_{m}^{M}) \\
&\hspace{3.5cm} \alpha^{G-}_A((u_{\partial \mu_1}^{\lambda_2(\partial \overline{\mu_1})} \otimes \id_{\mu_2}) b^\calc_{\mu_2, \lambda_2(\partial \overline{\mu_1})} e_j^{\nu_2, \mu_2 \lambda_2(\partial \overline{\mu_1})}) \phi_n^N) \\
&\hspace{1cm} \times ((u_{\partial \mu_1}^{\lambda_1(\partial \overline{\mu_1})} \otimes \id_{\mu_1}) b^\calc_{\mu_1, \lambda_1(\partial \overline{\mu_1})} e_i^{\nu_1, \mu_1 \lambda_1(\partial \overline{\mu_1})} \boxtimes \overline{\tilde{e}_j^{\nu_2, \mu_2 (\lambda_2(\partial \overline{\mu_1}))} b^{\calc-}_{\lambda_2(\partial \overline{\mu_1}),\mu_2} (\tilde{u}_{\partial \mu_1}^{\lambda_2(\partial \overline{\mu_1})} \otimes \id_{\mu_2})}),
\end{align*}
which is equal to $\Delta$ since $\{(u_{\partial \mu_1}^{\lambda_k(\partial \overline{\mu_1})} \otimes \id_{\mu_k}) b_{\mu_k, \lambda_k(\partial \overline{\mu_1})} e_i^{\nu_k, \mu_k \lambda_k(\partial \overline{\mu_1})}\}_i$ is a basis of $\hom(\nu_k, \lambda_k \mu_k)$ and $\{\tilde{e}_j^{\nu_k, \mu_k (\lambda_k(\partial \overline{\mu_1}))} b^{\calc-}_{\lambda_k(\partial \overline{\mu_1}),\mu_k} (\tilde{u}_{\partial \mu_1}^{\lambda_k(\partial \overline{\mu_1})} \otimes \id_{\mu_k})\}_j$ is its dual basis for $k=1,2$.
\end{proof}
\end{prop}

Then, we show that $\Theta_\alpha^G(A)$ is symmetric.

\begin{prop}
	\label{prop_alphafrob_sym}
	Let $A$ be a neutral symmetric special simple $G$-equivariant Frobenius algebra in a split spherical $G$-braided fusion category $\calc$. Then, $\Theta_\alpha^G(A)$ is symmetric.

	\begin{proof}
		Let us follow the notation in Theorem \ref{mainthm1}. For every $\lambda \in I$, let $\overline{\lambda}$ denote the element in $I$ with $\overline{\lambda} = \lambda^\vee$. Since $\calc$ is split by assumption, $\hom(\mathbf{1}_\calc, \mu \lambda)$ for $\lambda,\mu \in I$ is zero if $\mu \neq \overline{\lambda}$ and one-dimensional if $\mu = \overline{\lambda}$. Since $\Theta_\alpha^G(A)$ does not depend on a choice of a basis of $\hom(\mathbf{1}_\calc, \lambda \overline{\lambda})$, we may take it to be the right coevaluation $\mathrm{coev}_\lambda'$. Then, the corresponding dual basis is given by $\mathrm{ev}_\lambda$. 
		
		Now, let us calculate $\delta^{D(\calc)}_\theta$. By the proof of Proposition \ref{prop_semidir_def}, we can take $\theta^\vee$ to be 
		\begin{align*}
			\bigoplus_L \langle \alpha^{G+}_A(\lambda_1), \alpha^{G-}_A(\lambda_2) \rangle \overline{\lambda_1} \boxtimes {}^{\lambda_1} \overline{\lambda_2}^{\vee}.
		\end{align*}
		Since $\theta$ is a Frobenius algebra and therefore itself is a left dual, there exists a canonical isomorphism $\theta^{\vee \vee} \cong \theta$, which is given in Figure \ref{graphical_alphafrob_dual_canonical}.
		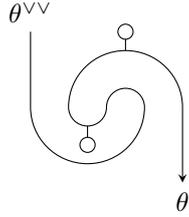
\begin{figure}[htb]
			\centering
			\begin{tikzpicture}
				\draw[->] (-1.5,1) --(-1.5,0) arc (180:360:0.75) arc (0:180:0.25) arc (360:180:0.25) arc (180:0:0.75) -- (0.5,-1);
				\draw (-0.75,-0.25) -- (-0.75,-0.5);
				\draw (-0.25,0.75) -- (-0.25,1);
				\draw[fill=white] (-0.75,-0.5) circle [radius=0.1];
				\draw[fill=white] (-0.25,1) circle [radius=0.1];
				\node at (-1.5,1.25){$\theta^{\vee \vee}$};
				\node at (0.5,-1.25){$\theta$};
			\end{tikzpicture}
			\caption{The canonical isomorphism $\theta^{\vee \vee} \cong \theta$}
			\label{graphical_alphafrob_dual_canonical}
		\end{figure}

		Then, the composition of $\delta^{D(\calc)}_\theta$ and the canonical isomorphism $\theta^{\vee \vee} \cong \theta$ are given by
		\begin{align*}
			\bigoplus_L \sum_{l'} \tr(\alpha^{G+}_A(\mathrm{ev}_{\lambda_1}) (\tilde{\phi}^{\overline{L}}_{l'} \otimes_A \tilde{\phi}^L_l) \alpha^{G-}_A(\mathrm{coev}'_{\lambda_2})) \tr(\alpha^{G-}_A(\mathrm{ev}_{\lambda_2}) (\phi^{\overline{L}}_{l'} \otimes_A \phi^L_l) \alpha^{G+}_A(\mathrm{coev}'_{\lambda_1})) f^L_1 \boxtimes f^L_2,
		\end{align*}
		where $\overline{L} = (\overline{\lambda_1}, \overline{\lambda_2}, l')$ and $f^L_1 \in \en(\lambda_1)$ and $f^L_2 \in \en(\lambda_2^\vee)$ is given in Figure \ref{graphical_proof_sym}, which are indeed identical. We used Figure \ref{graphical_reidemeister} in the proof of $f_2^L = \id_{\lambda_2^\vee}$. The coefficient cancels as in the proof of $\tilde{m} \circ \Delta = d_\theta$, see the proof of \cite[Theorem 1.4]{MR1754521}. Thus, we obtain $\delta^{D(\calc)}_\theta = \id_\theta$ under the canonical identification $\theta^{\vee \vee} \cong \theta$.
		\begin{figure}[htb]
			\centering
			\begin{tikzpicture}
				\draw[->] (-1.5,1) --(-1.5,0) arc (180:360:0.75) arc (0:180:0.25) arc (360:180:0.25) arc (180:0:0.75) -- (0.5,-1);
				\node at (-1.5,1.25){$\lambda_1$};
				\node at (0.5,-1.25){$\lambda_1$};
				\node at (-2.25,0){$f_1^L=$};
				\node at (1.25,0){$=\id_{\lambda_1}$};
				\node at (2.25,0){and};
				\begin{scope}[shift={(6,0)}]
					\draw (0,0) -- (-0.5,-0.5) arc (360:180:0.25) -- (-1,0) -- (-0.5,0.5) -- (-1,1);
					\draw[cross] (-2,1.5) -- (-2,0) -- (-0.5,-1.5) arc (180:360:0.25) -- (1,-0.5) -- (1,0) -- (0.5,0.5) arc (0:180:0.25) -- (-0.5,0) -- (0,-0.5);
					\draw[cross,->] (-1,0.5) -- (0,1.5) arc (180:0:0.25) -- (1.5,0.5) -- (1.5,-1.5);
					\draw (0,-0.5) arc (180:360:0.25) -- (0.5,0) arc (0:180:0.25);
					\draw (-1,0.5) arc (360:180:0.25) -- (-1.5,1) arc (180:0:0.25);
					\node at (-2,1.75){$\lambda_2^\vee$};
					\node at (1.5,-1.75){$\lambda_2^\vee$};
					\node at (-2.75,0){$f_2^L=$};
					\node at (2.25,0){$=\id_{\lambda_2^\vee}$};
				\end{scope}
			\end{tikzpicture}
			\caption{The proof of the symmetry of $\Theta_\alpha^G(A)$}
			\label{graphical_proof_sym}
		\end{figure}
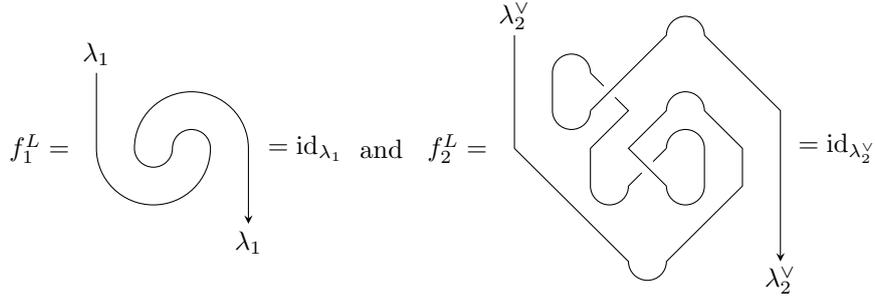
	\end{proof}
\end{prop}

Finally, we give an application of equivariant $\alpha$-induction Frobenius algebras to algebraic quantum field theory. Indeed, some of them give modular invariants of fixed point nets (Theorem \ref{theorem_modular}).

Let $\calc$ be a multitensor category with an action of a finite group $G$ and let $H$ be a subgroup of $G$. Then, the \emph{induction functor} $\ind_H^G : \calc^H \to \calc^G$, where $\calc^G$ denotes the $G$-equivariantization of $\calc$ (see e.g. \cite[Sections 2.7 and 4.15]{egno}), is defined by fixing a complete set of representatives of $G/H$ and putting $\ind_H^G(\lambda, z^\lambda) \coloneqq (\bigoplus_{g \in G/H} {}^g \lambda, z^{\ind_H^G(\lambda)})$. Here, $z^{\ind_H^G(\lambda)}_{g'}$ is defined as follows. For $g \in G$, there is a unique $\tilde{g} \in G$ with $g' g \in \tilde{g} H$. Let $h \in H$ be a unique element with $g' g = \tilde{g} h$. Then, the component of $z^{\ind_H^G(\lambda)}_{g'}$ on ${}^g \lambda$ is defined to be ${}^{\tilde{g}}z_h^\lambda : {}^{g' g} \lambda = {}^{\tilde{g} h} \lambda \cong {}^{\tilde{g}} \lambda$. 

By definition, for an algebra (resp. coalgebra, resp. Frobenius algebra, resp. Q-system) $A$ in $\calc$, $\ind_H^G(A)$ is naturally endowed with a structure of an algebra (resp. coalgebra, resp. Frobenius algebra, resp. Q-system).  

\begin{lem}
	\label{lemma_induction_commutative}
	Let $\calc$ be a $G$-braided multitensor category for a finite group $G$, and let $H$ be a subgroup of $G$. Let $\calc_H$ denote the $H$-braided multitensor category obtained by restricting $\calc$. If $A$ is an $H$-commutative algebra (resp. $H$-cocommutative coalgebra) in $\calc_H$, then $\ind_H^G(A)$ is $G$-commutative (resp. $G$-cocommutative). 

	\begin{proof}
		We only prove $G$-cocommutativity. Consider a component ${}^g A$ of $\ind_H^G(A)$. Then, the label $g'$ on the dashed line in Figure \ref{graphicalgcocommutativity} survives only when $g^{-1} g' g \in H$ since $A \in \calc_H$. By definition, for $h \coloneqq g^{-1} g' g \in H$, we have $z^{\ind_H^G(A)}_{g'} = {}^g z_h^A$. Thus, the $H$-cocommutativity of $A$ implies the $G$-cocommutativity of $\ind_H^G(A)$.
	\end{proof}
\end{lem}

\begin{lem}
	\label{lemma_equiv_reverse}
	Let $\calc$ be a $G$-braided multitensor category. Then, $(\calc^{\mathrm{rev}})^G \simeq (\calc^G)^{\mathrm{rev}}$ as braided multitensor categories.

	\begin{proof}
		$J_{(\lambda, z^\lambda), (\mu, z^\mu)}$ defined as the extension of $z_{g}^\lambda \otimes \id_\mu$ for $(\lambda, z^\lambda), (\mu, z^\mu) \in \obj((\calc^{\mathrm{rev}})^G)$, where a homogeneous decomposition of $\mu$ is labeled by $g \in G$, makes the identity into a tensor equivalence since
			$(z_{g}^\lambda \otimes \id_{\mu \nu}) (z_{h}^{{}^g \lambda \mu} \otimes \id_\nu) = (z_{g}^\lambda \otimes \id_{\mu \nu}) ({}^g z_{g^{-1} h g}^\lambda \otimes z^\mu_{h} \otimes \id_\nu) = z_{hg}^\lambda \otimes z^\mu_{h} \otimes \id_\nu$
		for $(\nu, z^\nu) \in \obj((\calc^{\mathrm{rev}})^G)$ and $g,h \in G$. This $J$ intertwines the braiding $(z_g^{\mu-1} \otimes \id_\lambda) b^{\calc-}_{\lambda, \mu}$ of $(\calc^{\mathrm{rev}})^G$ and $b^{\calc-}_{\lambda, \mu} (z_h^{\lambda-1} \otimes \id_\mu)$ of $(\calc^G)^{\mathrm{rev}}$.
	\end{proof}
\end{lem}

\begin{lem}
	\label{lemma_application_Frobenius}
	Let $G$ be a finite group. Let $A$ be a neutral symmetric special simple $G$-equivariant Frobenius algebra in a split spherical $G$-braided fusion category $\calc$. Then, $\ind_G^{G \times G} \Theta_\alpha^G(A)$ is a commutative Frobenius algebra in $\calc^G \boxtimes (\calc^G)^{\mathrm{rev}}$. If $\calc^G$ is modular, then $\ind_G^{G \times G} \Theta_\alpha^G(A)$ is maximal if and only if $\operatorname{FPdim} \Theta_\alpha^G(A) = \operatorname{FPdim} \calc$. 
	
	\begin{proof}
		$\ind_G^{G \times G} \Theta_\alpha^G(A)$ is a commutative Frobenius algebra in $\calc^G \boxtimes (\calc^{\mathrm{rev}})^G \simeq \calc^G \boxtimes (\calc^G)^{\mathrm{rev}}$ by Proposition \ref{prop_alphafrob_gcomm} and Lemmata \ref{lemma_induction_commutative} and \ref{lemma_equiv_reverse}. The last statement follows from \cite[Corollary 3.32]{MR3039775} since $\operatorname{FPdim} \ind_G^{G \times G} \Theta_\alpha^G(A) = \abs{G} \operatorname{FPdim} \Theta_\alpha^G(A)$ by definition.
	\end{proof}
\end{lem}

\begin{rem}
	If $\calc^G$ is modular and $\dim \Theta_\alpha^G(A) = \abs{G} \dim \rep \cala$, then there exists a simple symmetric special Frobenius algebra $\tilde{A}$ in $\calc^G$ such that $\ind_G^{G \times G} \Theta_\alpha^G(A)$ is isomorphic to the (ordinary) $\alpha$-induction Frobenius algebra $\Theta_\alpha(\tilde{A})$ by \cite[Theorem 3.22]{MR2551797} and \cite[Proposition 4.18]{MR3424476} (see also Theorem \ref{mainthm2} below). 
\end{rem}

Note that the $G$-equivariant Longo--Rehren Frobenius algebra $\Theta^G_{\mathrm{LR}}$ always satisfies the condition $\operatorname{FPdim} \Theta_{\mathrm{LR}}^G = \operatorname{FPdim} \calc$.

\begin{thm}
	\label{theorem_modular}
	Let $\cala$ be a completely rational irreducible local M\"{o}bius covariant net on $S^1$ with an action of a finite group $G$. Let $A$ be a simple $G$-equivariant Q-system in $\rep \cala$. Then, $\ind_G^{G \times G} \Theta_\alpha^G(A)$ is a commutative Q-system in $\rep \cala^G \boxtimes (\rep \cala^G)^{\mathrm{rev}}$. Therefore, it gives a finite index local standard extension of $\cala^G \otimes \cala^G$. It is modular invariant if and only if $\dim \Theta_\alpha^G(A) = \abs{G} \dim \rep \cala$. 

	\begin{proof}
		The first statement follows from \cite[Theorem 3.12]{MR2183964} and Lemma \ref{lemma_application_Frobenius}. The last statement follows from \cite[Proposition 6.6]{MR3424476}.
	\end{proof}
\end{thm}

\section{Equivariant full centers}
\label{section_fullcenter}

Let $A$ be a neutral symmetric special simple $G$-equivariant Frobenius algebra in a split spherical $G$-braided fusion category $\calc$. For $\lambda, \mu \in \obj(\calc)$, we have
\begin{align*}
	\hom (\alpha_A^{G+}(\lambda), \alpha_A^{G-}(\mu)) \subset \hom ({}_A A \lambda, {}_A A \mu) \cong \hom_\calc(\lambda, A \mu)
\end{align*}
by definition, where ${}_A A \lambda$ denotes $A \lambda$ as a left $A$-module. Then, it is natural to ask if $\Theta_\alpha^G(A)$ can be realized as a subalgebra of a $G$-equivariant Frobenius algebra with coefficients $\langle \lambda, A \mu \rangle$. In the case where $G$ is trivial, it is known \cite[Proposition 4.18]{MR3424476} that the answer is yes and the $\alpha$-induction Frobenius algebra $\Theta_\alpha(A)$ is realized as the full center $Z(A)$ \cite[Definition 4.9]{MR2443266} of $A$. In this section, we give the equivariant generalization of this theorem.

\subsection{The equivariant full center of a $G$-equivariant Frobenius algebra}

In this subsection, we define the equivariant generalization of a full center (Definition \ref{def_fullcenter}). Because it is defined to be a subalgebra of a product algebra as in the case where $G$ is trivial, we begin with products and subalgebras.

\begin{prop}
Let $A,B$ be $G$-equivariant Frobenius algebras in a $G$-braided multitensor category $\calc$. Then, $(AB, m_{AB}, \eta_A \otimes \eta_B, \Delta_{AB}, \varepsilon_A \otimes \varepsilon_B, z^A \otimes z^B)$ is a $G$-equivariant Frobenius algebra, where $m_{AB}$ and $\Delta_{AB}$ are given in Figure \ref{equivfrobalgproddef}. If $A$ and $B$ are special (resp. symmetric), then so is $A B$. If $A$ and $B$ are Q-systems in a $G$-braided ${}^\ast$-multitensor category, then $A B$ is again a $G$-equivariant Q-system.

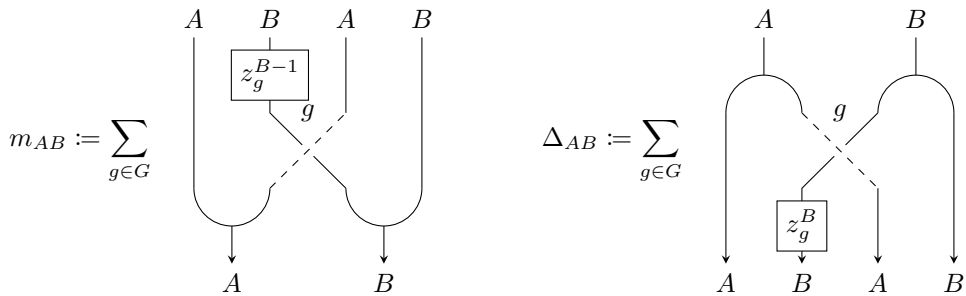
\begin{figure}[htb]
\centering
\begin{tikzpicture}
\node at (-2,-1){$\Delta_{A B} \coloneqq \displaystyle{\sum_{g \in G}}$};
\begin{scope}[shift={(2,0)},xscale=-1]
\draw (0,0.5) -- (0,0);
\draw (0.5,-0.5) arc (0:180:0.5);
\draw[->] (-0.5,-0.5) -- (-0.5,-2.5);
\draw[->] (0.5,-0.5) -- (1.5,-1.5) -- (1.5,-2.5);
\node at (0,0.75){$B$};
\node at (-0.5,-2.75){$B$};
\node at (1.5,-2.75){$B$};
\node[block] at (1.5, -2){$z^B_g$};
\end{scope}
\draw (0,0.5) -- (0,0);
\draw (0.5,-0.5) arc (0:180:0.5);
\draw[->] (-0.5,-0.5) -- (-0.5,-2.5);
\draw[dashed,cross] (0.5,-0.5) -- (1.5,-1.5);
\draw[->] (1.5,-1.5) -- (1.5,-2.5);
\node at (1,-0.5){$g$};
\node at (0,0.75){$A$};
\node at (-0.5,-2.75){$A$};
\node at (1.5,-2.75){$A$};
\begin{scope}[shift={(-7,-2)},yscale=-1]
\node at (-2,-1){$m_{A B} \coloneqq \displaystyle{\sum_{g \in G}}$};
\begin{scope}[shift={(2,0)},xscale=-1]
\draw[<-] (0,0.5) -- (0,0);
\draw (0.5,-0.5) arc (0:180:0.5);
\draw (-0.5,-0.5) -- (-0.5,-2.5);
\draw (0.5,-0.5) -- (1.5,-1.5) -- (1.5,-2.5);
\node at (0,0.75){$B$};
\node at (-0.5,-2.75){$B$};
\node at (1.5,-2.75){$B$};
\node[block] at (1.5, -2){$z^{B-1}_g$};
\end{scope}
\draw[<-] (0,0.5) -- (0,0);
\draw (0.5,-0.5) arc (0:180:0.5);
\draw (-0.5,-0.5) -- (-0.5,-2.5);
\draw[dashed,cross] (0.5,-0.5) -- (1.5,-1.5);
\draw (1.5,-1.5) -- (1.5,-2.5);
\node at (1,-1.5){$g$};
\node at (0,0.75){$A$};
\node at (-0.5,-2.75){$A$};
\node at (1.5,-2.75){$A$};
\end{scope}
\end{tikzpicture}
\caption{The product equivariant Frobenius algebra structure}
\label{equivfrobalgproddef}
\end{figure}

\begin{proof}
By definition, $(A B,z^{A B}) = (A, z^A) (B,z^B)$ in $\calc^G$ and the right-hand side induces the stated Frobenius algebra structure on the left-hand side by \cite[Proposition 3.22]{MR1940282}. Therefore, $A B$ is a Frobenius algebra. We can also give a direct proof with some graphical calculations in $\calc$. It follows from the naturality of dashed crossings that $z_g^{A B}$ is a Frobenius algebra isomorphism. The remaining statements follow easily.
\end{proof}
\end{prop}

\begin{defi}
Let $A$ be a $G$-equivariant Frobenius algebra in a $G$-crossed multitensor category, and let $e$ be a Frobenius idempotent for $A$ i.e. an idempotent in $\en_\calc(A)$ that satisfies \cite[Equations 2.54 and 2.55]{MR2187404}. We say $e$ is a \emph{$G$-equivariant Frobenius idempotent} for $A$ if $e z_g^A = z_g^A {}^g e$ for every $g \in G$. If $A$ is a $G$-equivariant Q-system in a $G$-crossed ${}^\ast$-multitensor category and $e$ is a projection, we say $e$ is a \emph{$G$-equivariant Frobenius projection} for $A$.
\end{defi}

\begin{prop}
\label{equivsubalgprop}
Let $A$ be a $G$-equivariant Frobenius algebra in a $G$-crossed multitensor category $\calc$ and let $e$ be a $G$-equivariant Frobenius idempotent for $A$ with a retract $(B,s,r)$. Then, for any nonzero scalar $\zeta$, the tuple $B_\zeta \coloneqq (B, rm_A(s \otimes s), r\eta_A, \zeta (r \otimes r)\Delta_A s, \zeta^{-1} \varepsilon_A s, \{r z_g^A {}^g s \}_g)$ is a $G$-equivariant Frobenius algebra. If moreover $\calc$ is a $G$-crossed ${}^\ast$-multitensor category, $A$ is a $G$-equivariant Q-system and $e$ is a $G$-equivariant Frobenius projection, $B_\zeta$ with $r =s^\ast$ is a $G$-equivariant Q-system. 

\begin{proof}
It is already known that $B_\zeta$ is a Frobenius algebra, see the proof of \cite[Proposition 2.37]{MR2187404}. Since $e$ is a $G$-equivariant Frobenius idempotent, we obtain
\begin{align*}
m_B (z_g^B \otimes z_g^B) &= r m_A (e z_g^A {}^g s \otimes e z_g^A {}^g s) = re z_g^A {}^g m_A ({}^g s \otimes {}^g s) = r z_g^A {}^g e {}^g(m_A (s \otimes s)) \\
&= r z_g^A {}^g s {}^g(r m_A (s \otimes s)) = z_g^B {}^g m_B \\
z_g^B {}^g \eta_B &= r z_g^A {}^g(e \eta_A) = rz_g^A {}^g \eta_A = r \eta_A = \eta_B
\end{align*}
for every $g \in G$. We can similarly check $\Delta_B z_g^B = (z_g^B \otimes z_g^B) {}^g \Delta_B$ and $\varepsilon^B z_g^B = {}^g \varepsilon^B$. Therefore $z_g^B$ is a Frobenius algebra endomorphism of $B$. Moreover,
\begin{align*}
z_g^B {}^g r z_g^{A-1} s = r z_g^A {}^g e z_g^{A-1} s = res = \id_{B} = {}^g r z_g^{A-1} s z_g^B
\end{align*}
and therefore $z_g^B$ is an isomorphism. Finally,
\begin{align*}
z_g^B {}^g z_h^B = r z_g^A {}^g(e z_h^A {}^h s) = r e z_g^A {}^g z_h^A {}^{g h} s = r z_{gh}^A {}^{g h} s = r z_{gh}^A {}^{gh}s = z_{gh}^B
\end{align*}
for any $g,h \in G$, and therefore $\{ z_g^B \}_g$ gives a $G$-equivariant structure on $B$. The last statement follows from \cite[Lemma 4.1]{bklr} and $z_g^{B-1} = {}^g r z_g^{A-1} s = z_g^{B \ast}$.
\end{proof}
\end{prop}

Note that the Frobenius algebra structure of $B_\zeta$ depends on $\zeta$ in general. Indeed, if $A$ is symmetric and special, then $\varepsilon_{B_\zeta} \eta_{B_\zeta} = \zeta^{-1} \dim A$. However, we do not write $\zeta$ when it is not important.

\begin{prop}
	\label{proposition_p_lambda}
Let $A$ be a symmetric special $G$-equivariant Frobenius algebra in a $G$-braided multitensor category $\calc$. Then, $P_A^{G}(\lambda) \in \en(A\lambda)$ defined in Figure \ref{equivleftcentralidempotent} for every $(\lambda, z^\lambda) \in \calc^G$ is an idempotent, where $\calc^G$ denotes the $G$-equivariantization of $\calc$ (see e.g. \cite[Sections 2.7 and 4.15]{egno}). Moreover, $P_A^{G} \coloneqq P_A^{G}(\mathbf{1}^\calc)$ is a $G$-equivariant Frobenius idempotent for $A$. When $\calc$ is a $G$-braided ${}^\ast$-multitensor category and $A$ is a Q-system, $P_A^G(\lambda)$ is a projection.

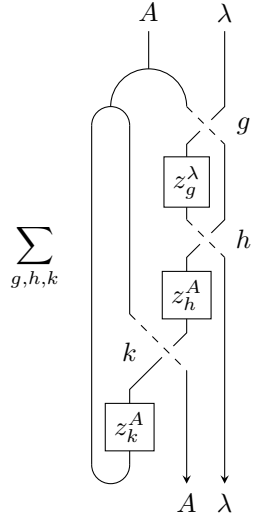
\begin{figure}[htb]
\centering
\begin{tikzpicture}
\draw (0,0) -- (0,0.5);
\draw (0.5,-0.5) arc (0:180:0.5);
\draw (-0.25,-0.75) arc (0:180:0.25);
\draw (1,0.5) -- (1,-0.5) -- (0.5,-1) -- (0.5,-2);
\node[block] at (0.5,-1.5){$z_g^\lambda$};
\draw[dashed,cross] (0.5, -0.5) -- (1,-1);
\draw (1,-1) -- (1,-2) -- (0.5,-2.5)--(0.5,-3.5)--(-0.25,-4.25) -- (-0.25,-5.25);
\node[block] at (-0.25,-4.75){$z_k^A$};
\draw (-0.25,-5.25) arc (360:180:0.25);
\draw (-0.75,-0.75) -- (-0.75,-5.25);
\draw[dashed,cross] (0.5, -2) -- (1,-2.5);
\node[block] at (0.5,-3){$z_h^A$};
\draw[->] (1,-2.5) -- (1,-5.5);
\draw[dashed,cross] (-0.25, -3.25) -- (0.5,-4);
\draw[->] (0.5,-4) -- (0.5,-5.5);
\draw (-0.25,-0.75) --(-0.25,-3.25);
\node at (0,0.75){$A$};
\node at (1,0.75){$\lambda$};
\node at (1.25,-0.75){$g$};
\node at (1.25,-2.25){$h$};
\node at (-0.25,-3.75){$k$};
\node at (0.5,-5.75){$A$};
\node at (1,-5.75){$\lambda$};
\node at (-1.5,-2.5){$\displaystyle \sum_{g, h, k}$};
\end{tikzpicture}
\caption{$P_A^{G}(\lambda)$}
\label{equivleftcentralidempotent}
\end{figure}

\begin{proof}
First, $P_A^{G}(\lambda) (z_g^A \otimes z_g^\lambda) = (z_g^A \otimes z_g^\lambda) {}^g P_A^{G}(\lambda)$ follows from Figure \ref{leftcentralprojisequiv}. Then, the remaining part of the statement follows from \cite[Lemma 5.2]{MR1940282}, \cite[Lemma 3.10]{MR2187404} and \cite[Lemma 4.6]{MR3424476} since $P_A^G(\lambda) = P_{(A,z^A)}((\lambda,z^\lambda))$ in $\calc^G$. Note that $\mathrm{ev}_A = \varepsilon_A m_A$ when $A^\vee$ is taken to be $A$ itself. We can also give a direct proof with some graphical calculations, see Figures \ref{graphical_ptilde_idempotent} and \ref{graphicallocalhomleftinvproof} below.
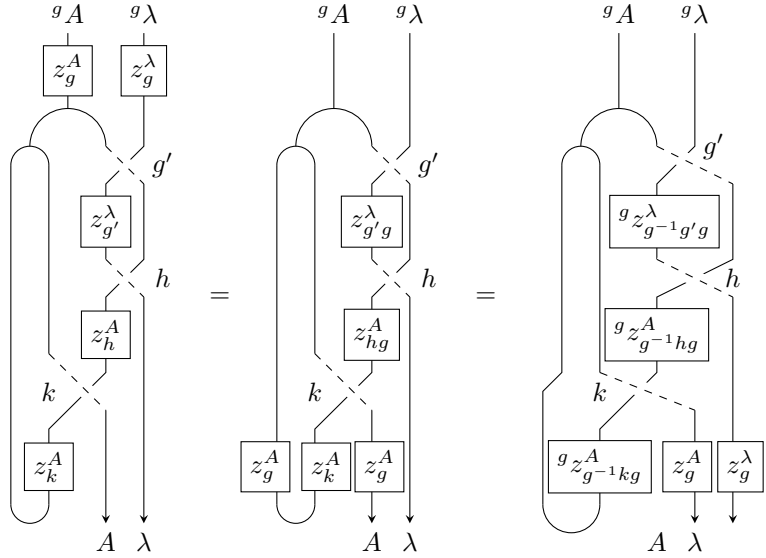
\begin{figure}[H]
\centering
\begin{tikzpicture}
\draw (0,0) -- (0,1);
\draw (0.5,-0.5) arc (0:180:0.5);
\draw (-0.25,-0.75) arc (0:180:0.25);
\draw (1,1) -- (1,-0.5) -- (0.5,-1) -- (0.5,-2);
\node[block] at (0.5,-1.5){$z_{g'}^\lambda$};
\draw[dashed,cross] (0.5, -0.5) -- (1,-1);
\draw (1,-1) -- (1,-2) -- (0.5,-2.5)--(0.5,-3.5)--(-0.25,-4.25) -- (-0.25,-5.25);
\node[block] at (-0.25,-4.75){$z_k^A$};
\draw (-0.25,-5.25) arc (360:180:0.25);
\draw (-0.75,-0.75) -- (-0.75,-5.25);
\draw[dashed,cross] (0.5, -2) -- (1,-2.5);
\node[block] at (0.5,-3){$z_h^A$};
\draw[->] (1,-2.5) -- (1,-5.5);
\draw[dashed,cross] (-0.25, -3.25) -- (0.5,-4);
\draw[->] (0.5,-4) -- (0.5,-5.5);
\draw (-0.25,-0.75) --(-0.25,-3.25);
\node at (0,1.25){${}^g A$};
\node at (1,1.25){${}^g \lambda$};
\node[block] at (0,0.5){$z_g^A$};
\node[block] at (1,0.5){$z_g^\lambda$};
\node at (1.25,-0.75){$g'$};
\node at (1.25,-2.25){$h$};
\node at (-0.25,-3.75){$k$};
\node at (0.5,-5.75){$A$};
\node at (1,-5.75){$\lambda$};
\node at (2,-2.5){$=$};
\begin{scope}[shift={(3.5,0)}]%equation 1
\draw (0,0) -- (0,1);
\draw (0.5,-0.5) arc (0:180:0.5);
\draw (-0.25,-0.75) arc (0:180:0.25);
\draw (1,1) -- (1,-0.5) -- (0.5,-1) -- (0.5,-2);
\node[block] at (0.5,-1.5){$z_{g' g}^\lambda$};
\draw[dashed,cross] (0.5, -0.5) -- (1,-1);
\draw (1,-1) -- (1,-2) -- (0.5,-2.5)--(0.5,-3.5)--(-0.25,-4.25) -- (-0.25,-5.25);
\node[block] at (-0.1,-4.75){$z_k^A$};
\draw (-0.25,-5.25) arc (360:180:0.25);
\draw (-0.75,-0.75) -- (-0.75,-5.25);
\draw[dashed,cross] (0.5, -2) -- (1,-2.5);
\node[block] at (0.5,-3){$z_{h g}^A$};
\draw[->] (1,-2.5) -- (1,-5.5);
\draw[dashed,cross] (-0.25, -3.25) -- (0.5,-4);
\draw[->] (0.5,-4) -- (0.5,-5.5);
\draw (-0.25,-0.75) --(-0.25,-3.25);
\node at (0,1.25){${}^g A$};
\node at (1,1.25){${}^g \lambda$};
\node[block] at (-0.9,-4.75){$z_g^A$};
\node[block] at (0.6,-4.75){$z_g^A$};
\node at (1.25,-0.75){$g'$};
\node at (1.25,-2.25){$h$};
\node at (-0.25,-3.75){$k$};
\node at (0.5,-5.75){$A$};
\node at (1,-5.75){$\lambda$};
\node at (2,-2.5){$=$};
\end{scope}
\begin{scope}[shift={(7.25,0)}]%equation 2
\draw (0,0) -- (0,1);
\draw (0.5,-0.5) arc (0:180:0.5);
\draw (-0.25,-0.75) arc (0:180:0.25);
\draw (1,1) -- (1,-0.5) -- (0.5,-1) -- (0.5,-2);
\node[block] at (0.6,-1.5){${}^g z_{g^{-1} g' g}^\lambda$};
\draw[dashed,cross] (0.5, -0.5) -- (1.5,-1);
\draw (1.5,-1) -- (1.5,-2) -- (0.5,-2.5)--(0.5,-3.5)--(-0.25,-4.25) -- (-0.25,-5.25);
\node[block] at (-0.25,-4.75){${}^g z_{g^{-1} kg}^A$};
\draw (-0.25,-5.25) arc (360:180:0.375);
\draw (-0.75,-0.75) -- (-0.75,-3.5) -- (-1,-3.75) -- (-1,-5.25);
\draw[dashed,cross] (0.5, -2) -- (1.5,-2.5);
\node[block] at (0.5,-3){${}^g z_{g^{-1}h g}^A$};
\draw[->] (1.5,-2.5) -- (1.5,-5.5);
\draw[dashed,cross] (-0.25, -3.5) -- (1,-4);
\draw[->] (1,-4) -- (1,-5.5);
\draw (-0.25,-0.75) --(-0.25,-3.5);
\node at (0,1.25){${}^g A$};
\node at (1,1.25){${}^g \lambda$};
\node[block] at (0.9,-4.75){$z_g^A$};
\node[block] at (1.6,-4.75){$z_g^\lambda$};
\node at (1.25,-0.5){$g'$};
\node at (1.5,-2.25){$h$};
\node at (-0.25,-3.75){$k$};
\node at (0.5,-5.75){$A$};
\node at (1,-5.75){$\lambda$};
\end{scope}
\end{tikzpicture}
\caption{$P_A^{G}(\lambda)$ is a $G$-equivariant morphism}
\label{leftcentralprojisequiv}
\end{figure}
\end{proof}
\end{prop}

\begin{defi}
Let $A$ be a symmetric special $G$-equivariant Frobenius algebra in a $G$-braided multitensor category $\calc$. We call $P_A^{G}$ the \emph{$G$-equivariant left central idempotent} of $A$. If $P_A^{G}$ is split, then we call the corresponding reduced Frobenius algebra, which is denoted by $C^{G}(A)$, the \emph{$G$-equivariant left center} of $A$.  
\end{defi}

\begin{prop}
	Let $A$ be a symmetric special $G$-equivariant Frobenius algebra in a $G$-braided multitensor category $\calc$. Then, $C^{G}(A)$ is symmetric, $G$-commutative and $G$-cocommutative. If $A$ is simple, then so is $C^{G}(A)$.

	\begin{proof}
		$(C^{G}(B),z^{C^{G}(B)})$ is equal to the left center $C((B,z^B))$ since $r$ and $s$ can be regarded as morphisms in $\calc^G$ for a retract $(C^{G}(B),r,s)$ by definition (see Proposition \ref{equivsubalgprop}), and $s r = P^{G}_B$ is $P_{(B,z^B)}$ in $\calc^G$. Then, the statement follows from \cite[Proposition 2.37 (i) and (ii)]{MR2187404}.
	\end{proof}
\end{prop}

\begin{rem}
	\begin{enumerate}
		\item We can also give the $G$-equivariant version of \cite[Proposition 2.25(iii)]{MR2187404}: a symmetric $G$-equivariant Frobenius algebra in a $G$-braided multitensor category is $G$-commutative if and only if it is $G$-cocommutative.
		\item The specialness of $C^{G}(A)$ is nontrivial in general, see \cite[Proposition 2.37(iii)]{MR2187404}.
	\end{enumerate}
\end{rem}

We need the following lemmata for the next subsection. 

\begin{lem}
	\label{prodidempotentlem}
	For symmetric special $G$-equivariant Frobenius algebras $A$ and $B$ in a $G$-braided multitensor category $\calc$, we have $P_{AB}^{G} = P_A^{G}(C^{G}(B))$.
	
	\begin{proof}
	$P_{A B}^{G}$ turns into the left central idempotent $P_{(A,z^A)(B,z^B)}$ in $\calc^G$. Then, the statement follows from \cite[Proposition 3.14(i)]{MR2187404}.
	\end{proof}
\end{lem}

\begin{lem}
	\label{lem_commutative_center}
	For a $G$-commutative (resp. $G$-cocommutative) symmetric special $G$-equivariant Frobenius algebra $A$ in a $G$-braided multitensor category, we have $P^G(A) = \id_A$.

	\begin{proof}
		As in the proof of Lemma \ref{prodidempotentlem}, it follows from the classical result \cite[Lemma 2.30]{MR2187404} since $(A,z^A)$ is commutative (resp. cocommutative) in $\calc^G$.
	\end{proof}
\end{lem}

\begin{defi}
	For a split $G$-braided fusion category $\calc$, we call $\Theta_{\mathrm{LR}}^G \coloneqq \Theta_\alpha^G(\mathbf{1}^\calc)$ the \emph{Longo--Rehren $G$-equivariant Frobenius algebra} in $D(\calc)$.
\end{defi}

Finally, note that for a neutral $G$-equivariant Frobenius algebra $A$ in a $G$-crossed multitensor category $\calc$, the tuple $(A \boxtimes \mathbf{1}^\calc, m_A \boxtimes \id, \eta_A \boxtimes \id, \Delta_A \boxtimes \id, \varepsilon_A \boxtimes \id, z^A \boxtimes \id)$ is a $G$-equivariant Frobenius algebra in $D(\calc)$. By definition, if $A$ is symmetric (resp. special), then so is $A \boxtimes \mathbf{1}^\calc$. Then, we can give the following definition, which is what we want and is the $G$-equivariant version of \cite[Definition 4.9]{MR2443266}.

\begin{defi}
	\label{def_fullcenter}
	Let $A$ be a neutral symmetric special $G$-equivariant Frobenius algebra in a split spherical $G$-braided fusion category $\calc$ with $\dim \calc \neq 0$. Then, we refer to the $G$-equivariant Frobenius algebra $Z^G(A) \coloneqq C^{G}((A \boxtimes \mathbf{1}^\calc) \Theta_{\mathrm{LR}}^G)$ in $D(\calc)$ as the \emph{$G$-equivariant full center} of $A$.
\end{defi}

Note that we need the sphericality of $\calc$ and $\dim \calc \neq 0$ for $\Theta_{\mathrm{LR}}^G$ to be symmetric special and therefore for $Z^G(A)$ to be well-defined, see Theorem \ref{mainthm1} and Proposition \ref{prop_alphafrob_sym}. In particular, when $\calc$ is a $G$-braided ${}^\ast$-fusion category, the assumption is always satisfied.

\subsection{Equivariant $\alpha$-induction Frobenius algebras as equivariant full centers}
\label{subsection_mainthm2}

In this subsection, we prove our second main theorem (Theorem \ref{mainthm2}). Indeed, the proof is given as the equivariant generalization of that of \cite[Proposition 4.18]{MR3424476}, and for this we need to define crossings that arise from $\alpha$-induction as in \cite[Proposition 3.1]{MR1729094}.

\begin{lem}
	Let $A$ be a neutral symmetric special $G$-equivariant Frobenius algebra in a $G$-braided multitensor category $\calc$. Then, for $\lambda \in \obj(\calc_g)$ and a left $A$-module $\mu$ in $\calc_h$, the morphisms in Figures \ref{graphicalibfeleftplus} and \ref{graphicalibfeleftminus} are left $A$-module isomorphisms and natural in $\lambda$ and $\mu$. Similarly, for a right $A$-module $\rho$, the morphisms in Figures \ref{graphicalibferightplus} and \ref{graphicalibferightminus} are right $A$-module isomorphisms and natural in $\lambda$ and $\rho$. If $\calc$ is a $G$-braided ${}^\ast$-multitensor category, $A$ is a $G$-equivariant Q-system and $\mu$ and $\rho$ are standard, then $B^\pm_{\lambda,\mu}$ and $B^{\pm}_{\rho,\lambda}$ are unitary when we take $s_{\alpha^{G \pm}_A(\lambda),\mu}$ and $s_{\rho, \alpha^{G \pm}_A(\lambda)}$ to be isometric. 

	\begin{figure}[H]
		\begin{tabular}{cc}
			\begin{minipage}[b]{0.45\hsize}
				\centering
				\begin{tikzpicture}
					\draw[->] (1,1) -- (0,0);
					\draw[cross,->] (0,1) -- (1,0);
					\draw[ultra thick] (0,1) -- (0.5,0.5);
					\node at (0,1.25){$\alpha^{G+}_A(\lambda)$};
					\node at (1,1.25){$\mu$};
					\node at (0,-0.25){${}^g \mu$};
					\node at (1,-0.25){$\lambda$};
					\node at (1.5,0.5){$\coloneqq$};
					\begin{scope}[shift={(3,-0.5)}]
						\draw[<-] (0,0) -- (0,0.5) -- (0.5,1) -- (0.5,1.5);
						\draw[cross,<-] (0.5,0) -- (0.5,0.5) -- (0,1) -- (0,2);
						\draw[rounded corners] (0,0.25) -- (-0.5,0.25) -- (-0.5,1.5);
						\draw (0,0.35) arc (90:270:0.1);
						\node[block] at (0,1.5){$s_{\alpha^{G+}_A(\lambda),\mu}$}; 
						\node at (0,2.25){$\alpha^{G+}_A(\lambda) \otimes_A \mu$};
						\node at (0,-0.25){${}^g \mu$};
						\node at (0.5,-0.25){$\lambda$};
						\node at (-0.5,0){$A$};
					\end{scope}
				\end{tikzpicture}
				\caption{$B^+_{\lambda, \mu}$}
				\label{graphicalibfeleftplus}
			\end{minipage}
			\begin{minipage}[b]{0.45\hsize}
				\centering
				\begin{tikzpicture}
					\draw[->] (0,1) -- (1,0);
					\draw[ultra thick] (0,1) -- (0.5,0.5);
					\draw[cross,->] (1,1) -- (0,0);
					\node at (0,1.25){$\alpha^{G-}_A({}^h \lambda)$};
					\node at (1,1.25){$\mu$};
					\node at (0,-0.25){$\mu$};
					\node at (1,-0.25){$\lambda$};
					\node at (1.5,0.5){$\coloneqq$};
					\begin{scope}[shift={(3,-0.5)}]
						\draw[<-] (0.5,0) -- (0.5,0.5) -- (0,1) -- (0,2);
						\draw[cross,<-] (0,0) -- (0,0.5) -- (0.5,1) -- (0.5,1.5);
						\draw[rounded corners] (0,0.25) -- (-0.5,0.25) -- (-0.5,1.5);
						\draw (0,0.35) arc (90:270:0.1);
						\node[block] at (0,1.5){$s_{\alpha^{G-}_A({}^h \lambda),\mu}$}; 
						\node at (0,2.25){$\alpha^{G-}_A({}^h \lambda) \otimes_A \mu$};
						\node at (0,-0.25){$\mu$};
						\node at (0.5,-0.25){$\lambda$};
						\node at (-0.5,0){$A$};
					\end{scope}
				\end{tikzpicture}
				\caption{$B^-_{\lambda, \mu}$}
				\label{graphicalibfeleftminus}	
			\end{minipage}
		\end{tabular}
	\end{figure}

	\begin{figure}[H]
		\begin{tabular}{cc}
			\begin{minipage}[b]{0.45\hsize}
				\centering
				\begin{tikzpicture}
					\draw[->] (0,1) -- (1,0);
					\draw[cross,->] (1,1) -- (0,0);
					\draw[ultra thick] (1,1) -- (0.5,0.5);
					\node at (1,1.25){$\alpha^{G+}_A(\lambda)$};
					\node at (0,1.25){${}^g \rho$};
					\node at (1,-0.25){$\rho$};
					\node at (0,-0.25){$\lambda$};
					\node at (1.5,0.5){$\coloneqq$};
					\begin{scope}[shift={(3,-0.5)}]
						\draw (0,1) -- (0,2);
						\draw[->] (-0.5,1.5) -- (-0.5,0.5) -- (0.5,-0.5);
						\draw[cross,<-] (-0.5,-0.5) -- (0.5,0.5) -- (0.5,1.5);
						\draw[rounded corners] (-0.5,0.75) -- (0,0.75) -- (0,1);
						\draw (-0.5,0.85) arc (90:-90:0.1);
						\node[block] at (0,1.5){$s_{{}^g \rho,\alpha^{G+}_A(\lambda)}$}; 
						\node at (0,2.25){${}^g \rho \otimes_A \alpha^{G+}_A(\lambda)$};
						\node at (0.5,-0.75){$\rho$};
						\node at (-0.5,-0.75){$\lambda$};
						\node at (0,0.5){$A$};
					\end{scope}
				\end{tikzpicture}
				\caption{$B^+_{\rho, \lambda}$}
				\label{graphicalibferightplus}
			\end{minipage}
			\begin{minipage}[b]{0.45\hsize}
				\centering
				\begin{tikzpicture}
					\draw[->] (1,1) -- (0,0);
					\draw[ultra thick] (1,1) -- (0.5,0.5);
					\draw[cross,->] (0,1) -- (1,0);
					\node at (1,1.25){$\alpha^{G-}_A(\lambda)$};
					\node at (0,1.25){$\rho$};
					\node at (0,-0.25){${}^h \lambda$};
					\node at (1,-0.25){$\rho$};
					\node at (1.5,0.5){$\coloneqq$};
					\begin{scope}[shift={(3,-0.5)}]
						\draw (0,1) -- (0,2);
						\draw[<-] (-0.5,-0.5) -- (0.5,0.5) -- (0.5,1.5);
						\draw[cross,->] (-0.5,1.5) -- (-0.5,0.5) -- (0.5,-0.5);
						\draw[rounded corners] (-0.5,0.75) -- (0,0.75) -- (0,1);
						\draw (-0.5,0.85) arc (90:-90:0.1);
						\node[block] at (0,1.5){$s_{\rho,\alpha^{G-}_A(\lambda)}$}; 
						\node at (0,2.25){$\rho \otimes_A \alpha^{G-}_A(\lambda)$};
						\node at (0.5,-0.75){$\rho$};
						\node at (-0.5,-0.75){${}^h \lambda$};
						\node at (0,0.5){$A$};
					\end{scope}
				\end{tikzpicture}
				\caption{$B^-_{\rho,\lambda}$}
				\label{graphicalibferightminus}	
			\end{minipage}
		\end{tabular}
	\end{figure}
	
	\begin{proof}
		The naturality of the morphisms follows by definition. We only show the remaining statement for $B^+_{\lambda, \mu}$ and $B^+_{\rho,\lambda}$ because the proof for $B_{\lambda, \mu}^-$ and $B^-_{\rho, \lambda}$ can be obtained just by replacing a crossing by its reverse. The left $A$-modularity of $B^+_{\lambda, \mu}$ follows from that of $s_{\alpha^{G+}_A(\lambda),\mu}$ and the left $A$-modularity of ${}^g \mu$. The morphism in Figure \ref{graphicalibfeinvleft} is the left inverse of $B^+_{\lambda,\mu}$ since $(B^+_{\lambda,\mu})^{-1} B^+_{\lambda,\mu}$ is equal to $r_{\alpha_A^{G+}(\lambda),\mu} f s_{\alpha_A^{G+}(\lambda),\mu}$ with $f$ given in Figure \ref{graphicalibfeinvleftproofleft}, where $z \coloneqq z^A$, and is equal to $(r_{\alpha_A^{G+}(\lambda),\mu} \circ s_{\alpha_A^{G+}(\lambda),\mu})^2 = \id_{\alpha_A^{G+}(\lambda) \otimes_A \mu}$ by the graphical calculation there. It is also the right inverse by Figure \ref{graphicalibfeinvleftproofright}. Next, the right $A$-modularity of $B_{\rho,\lambda}^+$ follows from Figure \ref{graphicalibferightproof}. Then, similarly to the argument for $B_{\lambda, \mu}^{+}$, we can show by some graphical calculations that the morphism in Figure \ref{graphicalibfeinvright} is the inverse of $B^+_{\rho, \lambda}$. The final statement follows from \cite[Lemma 3.23]{bklr} and the definitions of $(B^\pm_{\lambda, \mu})^{-1}$ and $(B^\pm_{\rho, \lambda})^{-1}$.
		\begin{figure}[htb]
			\begin{tabular}{cc}
				\begin{minipage}[b]{0.45\hsize}
					\centering
					\begin{tikzpicture}
						\draw[->] (0,1) -- (1,0);
						\draw[cross,->] (1,1) -- (0,0);
						\draw[ultra thick] (0,0) -- (0.5,0.5);
						\node at (0,-0.25){$\alpha^{G+}_A(\lambda)$};
						\node at (1,-0.25){$\mu$};
						\node at (0,1.25){${}^g \mu$};
						\node at (1,1.25){$\lambda$};
						\node at (1.5,0.5){$\coloneqq$};
						\begin{scope}[shift={(3,-1)}]
							\draw[<-] (0,0) -- (0,0.5);
							\draw (0.5,0.5) -- (0.5,1) -- (0,1.5) -- (0,2.75);
							\draw[cross,<-] (0,0.5) -- (0,1) -- (0.5,1.5) -- (0.5,2.75);
							\draw[rounded corners] (0,1.75) -- (-0.25,1.75) -- (-0.25,2);
							\draw (0,1.85) arc (90:270:0.1);
							\draw (-0.25,2) arc (0:180:0.125) -- (-0.5,0.5);
							\draw (-0.375,2.125) -- (-0.375,2.5);
							\draw[fill=white] (-0.375,2.5) circle [radius=0.1];
							\node[block] at (0,0.5){$r_{\alpha^{G+}_A(\lambda),\mu}$}; 
							\node at (0,-0.25){$\alpha^{G+}_A(\lambda) \otimes_A \mu$};
							\node at (0,3){${}^g \mu$};
							\node at (0.5,3){$\lambda$};
							\node at (-0.75,1.5){$A$};
						\end{scope}
					\end{tikzpicture}
					\caption{The inverse of $B_{\lambda,\mu}^+$}
					\label{graphicalibfeinvleft}
				\end{minipage}
				\begin{minipage}[b]{0.45\hsize}
					\centering
					\begin{tikzpicture}
						\draw[->] (1,1) -- (0,0);
						\draw[cross,->] (0,1) -- (1,0);
						\draw[ultra thick] (1,0) -- (0.5,0.5);
						\node at (1,-0.25){$\alpha^{G+}_A(\lambda)$};
						\node at (0,-0.25){${}^g \rho$};
						\node at (1,1.25){$\rho$};
						\node at (0,1.25){$\lambda$};
						\node at (1.5,0.5){$\coloneqq$};
						\begin{scope}[shift={(3,-1)}]
							\draw (0,0.5) -- (0,1.25) arc (0:180:0.125);
							\draw[rounded corners] (-0.25,1.25) -- (-0.25,1) -- (-0.5,1);
							\draw (-0.125,1.375) -- (-0.125,1.75);
							\draw[fill=white] (-0.125,1.75) circle [radius=0.1];
							\draw[->] (0,0.5) -- (0,0);
							\draw (-0.5,0.5) -- (-0.5,2) -- (0.5,3);
							\draw[cross] (-0.5,3) -- (0.5,2) -- (0.5,0.5);
							\node[block] at (0,0.5){$r_{{}^g \rho,\alpha^{G+}_A(\lambda)}$}; 
							\node at (0,-0.25){${}^g \rho \otimes_A \alpha^{G+}_A(\lambda)$};
							\node at (0.5,3.25){$\rho$};
							\node at (-0.5,3.25){$\lambda$};
							\node at (-0.75,2){${}^g \rho$};
							\node at (0.25,1.5){$A$};
							\draw (-0.5,1.1) arc (90:-90:0.1);
						\end{scope}
					\end{tikzpicture}
					\caption{The inverse of $B_{\rho, \lambda}^+$}
					\label{graphicalibfeinvright}
				\end{minipage}
			\end{tabular}
		\end{figure}
		\begin{figure}[htb]
			\centering
			\begin{tikzpicture}
				\draw (0,0) -- (0,0.5) -- (0.5,1) -- (0.5,1.5);
				\draw[cross] (0.5,0) -- (0.5,0.5) -- (0,1) -- (0,1.5);
				\draw[rounded corners] (0,0.25) -- (-0.5,0.25) -- (-0.5,1.5);
				\draw (0,0.35) arc (90:270:0.1);
				\node at (-0.5,1.75){$A$};
				\node at (1.25,-0.25){$=$};
				\node at (0.5,1.75){$\mu$};
				\node at (0,1.75){$\lambda$};
				\begin{scope}[shift={(0,-2.75)}]
					\draw[<-] (0.5,0.5) -- (0.5,1) -- (0,1.5) -- (0,2.75);
					\draw[cross,<-] (0,0.5) -- (0,1) -- (0.5,1.5) -- (0.5,2.75);
					\draw[rounded corners] (0,1.75) -- (-0.25,1.75) -- (-0.25,2);
					\draw (0,1.85) arc (90:270:0.1);
					\draw[->] (-0.25,2) arc (0:180:0.125) -- (-0.5,0.5);
					\draw (-0.375,2.125) -- (-0.375,2.5);
					\draw[fill=white] (-0.375,2.5) circle [radius=0.1];
					\node at (-0.5,0.25){$A$};
					\node at (0,0.25){$\lambda$};
					\node at (0.5,0.25){$\mu$};
				\end{scope}	
				\begin{scope}[shift={(3,0)}]
					\draw[<-] (0.75,-2) -- (0.75,1.5);
					\draw[<-] (-0.75,-2) -- (-0.75,0.5) arc (180:0:0.25) -- (-0.25,-0.5) -- (0.25,-1) -- (0.75,-1);
					\draw (-0.5,0.75) -- (-0.5,1.5);
					\draw (0.75,-0.9) arc (90:270:0.1);
					\draw[cross,->] (0.25,1.5) -- (0.25,-0.5) -- (-0.25,-1) -- (-0.25,-2);
					\node at (1.5,-0.25){$=$};
					\node at (-0.5,1.75){$A$};
					\node at (0.25,1.75){$\lambda$};
					\node at (0.75,-2.25){$\mu$};
					\node at (0.5,-1.25){$A$};
					\node at (-0.75,-2.25){$A$};
					\node at (0.75,1.75){$\mu$};
					\node at (-0.25,-2.25){$\lambda$};
					\node[block] at (-0.25,0){$z_g^{-1}$};
				\end{scope}
				\begin{scope}[shift={(6,0)}]
					\draw[<-] (0.75,-2) -- (0.75,1.5);
					\draw (-0.75,1.5) -- (-0.75,-1.5) arc (180:360:0.25) -- (-0.25,-0.5) -- (0.25,0) -- (0.25,0.25) arc (180:0:0.125);
					\draw[rounded corners] (0.5,0.25) -- (0.5,0) -- (0.75,0);
					\draw (0.375,0.375) -- (0.375,0.75);
					\draw[fill=white] (0.375,0.75) circle [radius=0.1]; 
					\draw[->] (-0.5,-1.75) -- (-0.5,-2);
					\draw (0.75,0.1) arc (90:270:0.1);
					\draw[cross,->] (-0.25,1.5) -- (-0.25,0) -- (0.25,-0.5) -- (0.25,-2);
					\node at (-0.5,-2.25){$A$};
					\node at (0.25,-2.25){$\lambda$};
					\node at (0.75,-2.25){$\mu$};
					\node at (-0.75,1.75){$A$};
					\node at (0.5,-0.25){$A$};
					\node at (0.75,1.75){$\mu$};
					\node at (-0.25,1.75){$\lambda$};
					\node[block] at (-0.25,-1){$z_g$};
				\end{scope}
			\end{tikzpicture}
			\caption{A proof of the left invertibility of $B_{\lambda,\mu}^+$}
			\label{graphicalibfeinvleftproofleft}
		\end{figure}
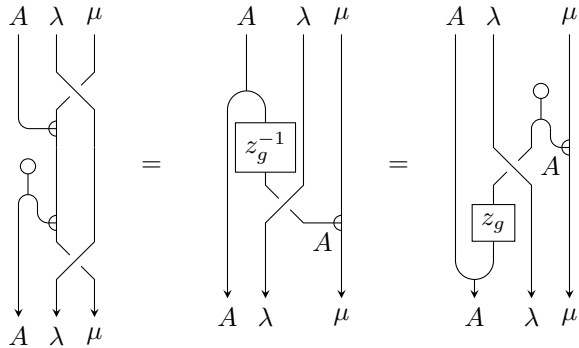
		\begin{figure}[htb]
			\centering
			\begin{tikzpicture}
				\draw (0,0) -- (0,0.5);
				\draw (0.5,0.5) -- (0.5,1) -- (0,1.5) -- (0,2.75);
				\draw[cross,<-] (0,0.5) -- (0,1) -- (0.5,1.5) -- (0.5,2.75);
				\draw[rounded corners] (0,1.75) -- (-0.25,1.75) -- (-0.25,2);
				\draw (0,1.85) arc (90:270:0.1);
				\draw (-0.25,2) arc (0:180:0.125) -- (-0.5,0.5);
				\draw (-0.375,2.125) -- (-0.375,2.5);
				\draw[fill=white] (-0.375,2.5) circle [radius=0.1];
				\node[block] at (0,0.5){$r_{\alpha^{G+}_A(\lambda),\mu}$}; 
				\node at (0,3){${}^g \mu$};
				\node at (0.5,3){$\lambda$};
				\node at (-0.75,1.5){$A$};
				\begin{scope}[shift={(0,-2)}]
					\draw[<-] (0,0) -- (0,0.5) -- (0.5,1) -- (0.5,1.5);
					\draw[cross,<-] (0.5,0) -- (0.5,0.5) -- (0,1) -- (0,2);
					\draw[rounded corners] (0,0.25) -- (-0.5,0.25) -- (-0.5,1.5);
					\draw (0,0.35) arc (90:270:0.1);
					\node[block] at (0,1.5){$s_{\alpha^{G+}_A(\lambda),\mu}$};
					\node at (0,-0.25){${}^g \mu$};
					\node at (0.5,-0.25){$\lambda$};
					\node at (-0.5,0){$A$};
				\end{scope}
				\node at (1.5,0){$=$};
				\begin{scope}[shift={(3,0)}]
					\begin{scope}[shift={(0,0.25)}]
						\draw (0.75,0.75) -- (0,1.5) -- (0,2.75);
						\draw[rounded corners] (0,1.75) -- (-0.25,1.75) -- (-0.25,2);
						\draw (0,1.85) arc (90:270:0.1);
						\draw (-0.25,2) arc (0:180:0.25) -- (-0.75,0.75);
						\draw (-0.5,2.25) -- (-0.5,2.5);
						\draw[fill=white] (-0.5,2.5) circle [radius=0.1];
						\node at (0,3){${}^g \mu$};
						\node at (0.5,3){$\lambda$};
					\end{scope}
					\draw (0.75,-2) -- (0.75,1);
					\draw (-0.75,1) -- (-0.75,-1.5) arc (180:360:0.25) -- (-0.25,-0.5) -- (0.25,0) -- (0.25,0.25) arc (180:0:0.125);
					\draw[rounded corners] (0.5,0.25) -- (0.5,0) -- (0.75,0);
					\draw (0.375,0.375) -- (0.375,0.75);
					\draw[fill=white] (0.375,0.75) circle [radius=0.1]; 
					\draw (-0.5,-1.75) -- (-0.5,-2);
					\draw (0.75,0.1) arc (90:270:0.1);
					\node at (-1,-1.75){$A$};
					\node at (0.5,-0.25){$A$};
					\node at (1,1){$\mu$};
					\node[block] at (-0.25,-1){$z_g$};
					\node at (1.5,0){$=$};
					\begin{scope}[shift={(0,-3.25)}]
						\draw[<-] (0,0) -- (0,0.5) -- (0.75,1.25);
						\draw[rounded corners] (0,0.25) -- (-0.5,0.25) -- (-0.5,1.5);
						\draw (0,0.35) arc (90:270:0.1);
						\node at (0,-0.25){${}^g \mu$};
						\node at (0.75,-0.25){$\lambda$};
					\end{scope}
					\draw[cross,->] (0.5,3) --(0.5,1.75) -- (-0.25,1) -- (-0.25,0) -- (0.25,-0.5) -- (0.25,-2) -- (0.75,-2.5) -- (0.75,-3.25);
				\end{scope}
				\begin{scope}[shift={(6,0)}]
					\begin{scope}[shift={(0,0.25)}]
						\draw (0.75,0.75) -- (0,1.5) -- (0,2.75);
						\draw[rounded corners] (0,1.75) -- (-0.25,1.75) -- (-0.25,2);
						\draw (0,1.85) arc (90:270:0.1);
						\draw (-0.25,2) arc (0:180:0.25) -- (-0.75,0.75);
						\draw (-0.5,2.25) -- (-0.5,2.5);
						\draw[fill=white] (-0.5,2.5) circle [radius=0.1];
						\node at (0,3){${}^g \mu$};
						\node at (0.5,3){$\lambda$};
					\end{scope}
					\draw (0.75,-2) -- (0.75,1);
					\draw (-0.75,1) -- (-0.75,-1.5) arc (180:360:0.25) -- (-0.25,-0.5) -- (0.25,0) -- (0.25,0.25) arc (180:0:0.125);
					\draw[rounded corners] (0.5,0.25) -- (0.5,0) -- (0.75,0);
					\draw (0.375,0.375) -- (0.375,0.75);
					\draw[fill=white] (0.375,0.75) circle [radius=0.1]; 
					\draw (-0.5,-1.75) -- (-0.5,-2);
					\draw (0.75,0.1) arc (90:270:0.1);
					\node at (-1,-1.75){$A$};
					\node at (0.25,-0.5){${}^{g^{-1}} A$};
					\node at (1,1){$\mu$};
					\node at (1.5,0){$=$};
					\begin{scope}[shift={(0,-3.25)}]
						\draw[<-] (0,0) -- (0,0.5) -- (0.75,1.25);
						\draw[rounded corners] (0,0.25) -- (-0.5,0.25) -- (-0.5,1.5);
						\draw (0,0.35) arc (90:270:0.1);
						\node at (0,-0.25){${}^g \mu$};
						\node at (0.75,-0.25){$\lambda$};
					\end{scope}
					\draw[cross,->] (0.5,3) --(0.5,1.75) -- (-0.5,0.75) -- (-0.5,-0.75) -- (0.25,-1.5) -- (0.25,-2) -- (0.75,-2.5) -- (0.75,-3.25);
				\end{scope}
				\begin{scope}[shift={(9,-1.75)}]
					\draw[rounded corners] (0.5,0.25) -- (-0.5,0.25) -- (-0.5,0.75);
					\draw (-0.25,3) arc (0:180:0.25) -- (-0.75,1) arc (180:360:0.25) -- (-0.25,1.75) arc (180:0:0.125);
					\draw[rounded corners] (0,1.75) -- (0,1.5) -- (0.5,1.5);
					\draw (-0.125,1.875) -- (-0.125,2.25);
					\draw[fill=white] (-0.125,2.25) circle [radius=0.1];
					\draw[rounded corners] (-0.25,3) -- (-0.25,2.5) -- (0.5,2.5);
					\draw (-0.5,3.25) -- (-0.5,3.5);
					\draw[fill=white] (-0.5,3.5) circle [radius=0.1];
					\draw[->] (0.5,3.75) -- (0.5,0);
					\draw[->] (1,3.75) -- (1,0);
					\draw (0.5,2.6) arc (90:270:0.1);
					\draw (0.5,1.6) arc (90:270:0.1);
					\draw (0.5,0.35) arc (90:270:0.1);
					\node at (-1,3){$A$};
					\node at (0.5,4){${}^g \mu$};
					\node at (1,4){$\lambda$};
					\node at (0.5,-0.25){${}^g \mu$};
					\node at (1,-0.25){$\lambda$};
					\node at (1.75,1.75){$= \id$};
				\end{scope}
			\end{tikzpicture}
			\caption{A proof of the right invertibility of $B_{\lambda,\mu}^+$}
			\label{graphicalibfeinvleftproofright}
		\end{figure}
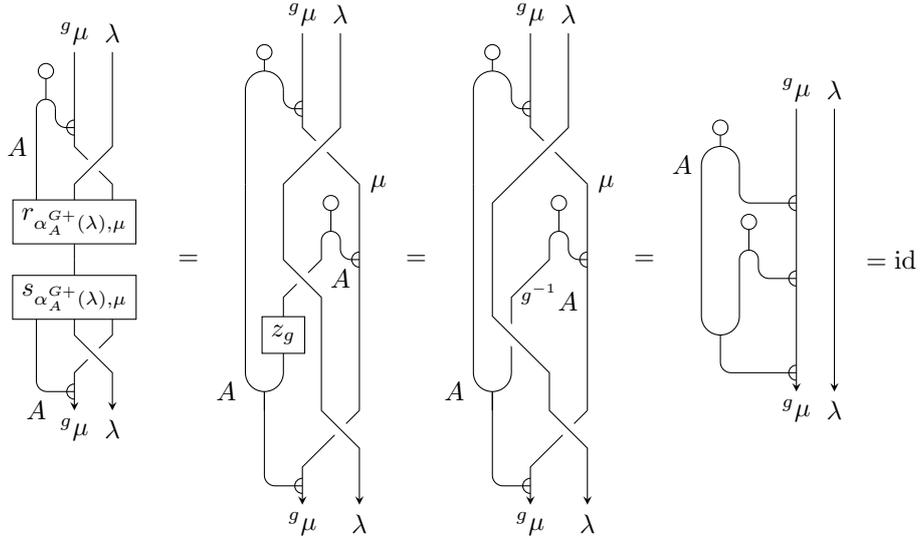
		\begin{figure}[htb]
			\centering
			\begin{tikzpicture}
				\draw (0,1) -- (0,2.25);
				\draw[rounded corners] (0,2) -- (1,2) -- (1,2.25);
				\draw[->] (-0.5,1.5) -- (-0.5,0.5) -- (0.5,-0.5);
				\draw[cross,<-] (-0.5,-0.5) -- (0.5,0.5) -- (0.5,1.5);
				\draw[rounded corners] (-0.5,0.75) -- (0,0.75) -- (0,1);
				\draw (-0.5,0.85) arc (90:-90:0.1);
				\draw (0,2.1) arc (90:-90:0.1);
				\node[block] at (0,1.5){$s_{{}^g \rho,\alpha^{G+}_A(\lambda)}$}; 
				\node at (-0.5,2.5){${}^g \rho \otimes_A \alpha^{G+}_A(\lambda)$};
				\node at (0.5,-0.75){$\rho$};
				\node at (-0.5,-0.75){$\lambda$};
				\node at (0,0.5){$A$};
				\node at (1,2.5){$A$};
				\node at (1.5,0.75){$=$};
				\node at (-0.75,1){${}^g \rho$};
				\begin{scope}[shift={(3,0)}]
					\draw (0,0.5) -- (0,2.25);
					\draw (0,0.5) arc (180:360:0.25) -- (0.5,1) -- (1,1.5) -- (1,2.25);
					\draw[->] (-0.5,1.5) -- (-0.5,-0.25) -- (0.5,-1.25);
					\draw[cross,<-] (-0.25,-1.25) -- (1,0) -- (1,1) -- (0.5,1.5);
					\draw[rounded corners] (-0.5,0) -- (0.25,0) -- (0.25,0.25);
					\draw (-0.5,0.1) arc (90:-90:0.1);
					\node[block] at (0,1.75){$s_{{}^g \rho,\alpha^{G+}_A(\lambda)}$}; 
					\node[block] at (0.5,0.75){$z_g$};
					\node at (-0.25,2.5){${}^g \rho \otimes_A \alpha^{G+}_A(\lambda)$};
					\node at (0.5,-1.5){$\rho$};
					\node at (-0.25,-1.5){$\lambda$};
					\node at (0,-0.25){$A$};
					\node at (1,2.5){$A$};
					\node at (1.5,0.75){$=$};
					\node at (-0.75,1.25){${}^g \rho$};
				\end{scope}
				\begin{scope}[shift={(6,0)}]
					\draw[rounded corners] (-0.5,1) -- (0,1) -- (0,2.25);
					\draw (0.5,1) -- (1,1.5) -- (1,2.25);
					\draw[->] (-0.5,1.5) -- (-0.5,-0.25) -- (0.5,-1.25);
					\draw[cross,<-] (-0.25,-1.25) -- (1,0) -- (1,1) -- (0.5,1.5);
					\draw[rounded corners] (-0.5,0) -- (0.5,0) -- (0.5,1);
					\draw (-0.5,0.1) arc (90:-90:0.1);
					\draw (-0.5,1.1) arc (90:-90:0.1);
					\node[block] at (0,1.75){$s_{{}^g \rho,\alpha^{G+}_A(\lambda)}$}; 
					\node[block] at (0.5,0.75){$z_g$};
					\node at (-0.25,2.5){${}^g \rho \otimes_A \alpha^{G+}_A(\lambda)$};
					\node at (0.5,-1.5){$\rho$};
					\node at (-0.25,-1.5){$\lambda$};
					\node at (0,-0.25){$A$};
					\node at (-0.25,0.75){$A$};
					\node at (1,2.5){$A$};
					\node at (1.5,0.75){$=$};
					\node at (-0.75,1.25){${}^g \rho$};
				\end{scope}
				\begin{scope}[shift={(9,0)}]
					\draw[rounded corners] (-0.5,1.25) -- (0,1.25) -- (0,2.25);
					\draw[->] (-0.5,1.5) -- (-0.5,1) -- (0.5,0) -- (0.5,-0.5);
					\draw[cross,<-] (-0.5,-0.5) --(-0.5,0) -- (0.5,1) -- (0.5,1.5);
					\draw[rounded corners] (0.5,-0.25) -- (1,-0.25) -- (1,2.25);
					\draw (0.5,-0.15) arc (90:-90:0.1);
					\draw (-0.5,1.35) arc (90:-90:0.1);
					\node[block] at (0,1.75){$s_{{}^g \rho,\alpha^{G+}_A(\lambda)}$}; 
					\node at (-0.25,2.5){${}^g \rho \otimes_A \alpha^{G+}_A(\lambda)$};
					\node at (0.5,-0.75){$\rho$};
					\node at (-0.5,-0.75){$\lambda$};
					\node at (1,2.5){$A$};
					\node at (-0.75,1.25){${}^g \rho$};
				\end{scope}
			\end{tikzpicture}
			\caption{The right $A$-modularity of $B_{\rho,\lambda}^+$}
			\label{graphicalibferightproof}
		\end{figure}
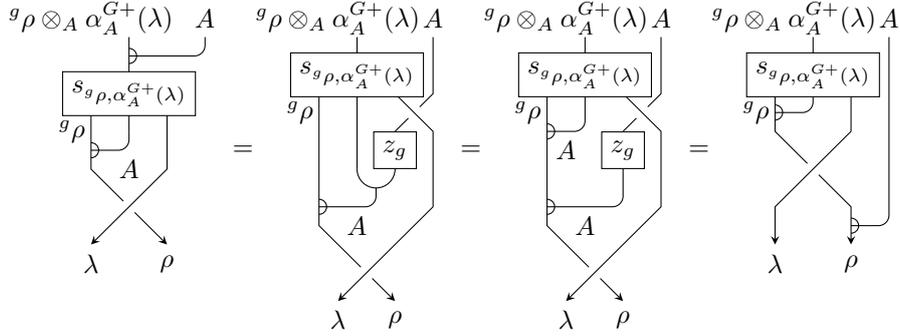
	\end{proof}
\end{lem}

We give the following lemmata for later use.

\begin{lem}
	\label{lem_ibfe_action}
	Let $A$ be a neutral symmetric special $G$-equivariant Frobenius algebra in a $G$-braided multitensor category $\calc$. Then, $B^+_{{}^g \lambda, {}^g \mu} = {}^g B_{\lambda, \mu}^+ (\tilde{\eta}^{g}_\lambda \otimes_A \id_{{}^g \mu})$ and $B^-_{{}^g \lambda, {}^g \mu} = {}^g B_{\lambda, \mu}^- ({}^{g \partial \mu g^{-1}} \tilde{\eta}^{g}_\lambda \otimes_A \id_{{}^g \mu})$, where $\tilde{\eta}^g_\lambda \coloneqq z^{A-1}_g \otimes_\calc \id_{{}^g \lambda}$, for $\lambda \in \homog(\calc)$, a homogeneous left $A$-module $\mu$ in $\calc$ and $g \in G$. Similarly, $B^+_{{}^g \rho, {}^g \lambda} = {}^g B^+_{\rho, \lambda} (\id_{{}^g \rho} \otimes_A \tilde{\eta}^{g}_\lambda)$ and $B^-_{{}^g \rho, {}^g \lambda}={}^g B^-_{\rho, \lambda} (\id_{{}^g \rho} \otimes_A \tilde{\eta}^{g}_\lambda)$ for a homogeneous right $A$-module $\rho$ in $\calc$.

	\begin{proof}
		Put $z \coloneqq z^A$. Since ${}^g z_\lambda^{-1} z_g^{-1} = z_{g \partial \lambda}^{-1}$, we obtain ${}^g B_{\lambda, \mu}^+ = B^+_{{}^g \lambda, {}^g \mu} (\tilde{\eta}^g_\lambda \otimes_A \id_\mu)$ from the naturality of $s$ and the definition of $m^{\mathrm{L}}_{{}^\lambda \mu}$. The proof for $B^\pm_{\rho, \lambda}$ is similar. The statement for $B^-_{\lambda, \mu}$ follows from $z_g^{-1} z_{g \partial \mu g^{-1}} {}^{g \partial \mu g^{-1}} z_g = z_g^{-1} z_{g \partial \mu} = {}^g z_{\partial \mu}$. 
	\end{proof}
\end{lem}

\begin{lem}
	\label{lem_alpha_ordinary_braiding}
	Let $A$ be a neutral special symmetric $G$-equivariant Frobenius algebra in a $G$-braided multitensor category $\calc$. Then, the equations in Figures \ref{graphical_alpha_ordinary_braiding} hold for $\lambda \in \obj(\calc_g)$, a homogeneous left $A$-module $\mu$ and a homogeneous right $A$-module $\rho$.
	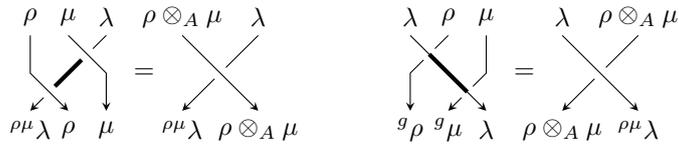
\begin{figure}[H]
		\centering
		\begin{tikzpicture}
			\draw[->] (1,0) -- (0,-1);
			\draw[ultra thick] (0.75,-0.25) -- (0.25,-0.75);
			\draw[cross,->] (0.5,0) -- (1,-0.5) -- (1,-1);
			\draw[cross,->] (0,0) -- (0,-0.5) -- (0.5,-1);
			\node at (0,0.25){$\rho$};
			\node at (0.5,0.25){$\mu$};
			\node at (1,0.25){$\lambda$};
			\node at (0,-1.25){${}^{\rho \mu} \lambda$};
			\node at (0.5,-1.25){$\rho$};
			\node at (1,-1.25){$\mu$};
			\node at (1.5,-0.5){$=$};
			\begin{scope}[shift={(2,0)}]
				\draw[->] (1,0) -- (0,-1);
				\draw[cross,->](0,0) -- (1,-1);
				\node at (0,0.25){$\rho \otimes_A \mu$};
				\node at (1,0.25){$\lambda$};
				\node at (0,-1.25){${}^{\rho \mu} \lambda$};
				\node at (1,-1.25){$\rho \otimes_A \mu$};
			\end{scope}
			\begin{scope}[shift={(5,0)}]
				\draw[->] (0.5,0) -- (0,-0.5) -- (0,-1);
				\draw[->] (1,0) -- (1,-0.5) -- (0.5,-1);
				\draw[cross,->] (0,0) -- (1,-1);
			\draw[ultra thick] (0.25,-0.25) -- (0.75,-0.75);
			\node at (0,0.25){$\lambda$};
			\node at (0.5,0.25){$\rho$};
			\node at (1,0.25){$\mu$};
			\node at (0,-1.25){${}^{g} \rho$};
			\node at (0.5,-1.25){${}^g \mu$};
			\node at (1,-1.25){$\lambda$};
			\node at (1.5,-0.5){$=$};
			\begin{scope}[shift={(2,0)}]
				\draw[->] (1,0) -- (0,-1);
				\draw[cross,->](0,0) -- (1,-1);
				\node at (1,0.25){$\rho \otimes_A \mu$};
				\node at (0,0.25){$\lambda$};
				\node at (1,-1.25){${}^{\rho \mu} \lambda$};
				\node at (0,-1.25){$\rho \otimes_A \mu$};
			\end{scope}
			\end{scope}
		\end{tikzpicture}
		\caption{Thick crossings and ordinary crossings}
		\label{graphical_alpha_ordinary_braiding}
	\end{figure} 
	\begin{proof}
		We only show the first equality because the proof of the other is similar. By explicitly writing down canonical morphisms in $\frob^G(\calc)$ as morphisms in $\calc$, we find that $r_{\alpha^{G-}_A({}^\mu \lambda),\mu}$ and $s_{\rho,\alpha^{G-}_A({}^\mu \lambda)}$ in the definition of $(B^+_{\lambda,\mu})^{-1}$ and $B^-_{\rho, \lambda}$ cancel. Then, the left-hand side of the first equation in Figure \ref{graphical_alpha_ordinary_braiding} is equal to $(\id_{{}^{\rho \mu} \lambda} \otimes r_{\rho, \mu}) b_{\rho \mu, \lambda}^\calc (e_{\rho, \mu} s \otimes \id_\lambda) = b^\calc_{\rho \otimes_A \mu, \lambda}$. 
	\end{proof}
\end{lem}

Thanks to this lemma, we can move an arc along a crossing as in Figure \ref{graphicalmovearc} when it contains a thick line segment inside. On the other hand, when thick line segments are outside an arc, we need more arguments to move the arc. 

\begin{lem}
	\label{lem_arc_alpha_braiding}
	Let $A$ be a neutral symmetric special $G$-equivariant Frobenius algebra in a $G$-braided multitensor category $\calc$. Then, the equations in Figure \ref{graphical_arc_alpha_braiding} hold for $\lambda \in \obj(\calc_g)$, a homogeneous left $A$-module $\mu$ in $\calc$ and a homogeneous right $A$-module $\rho$ in $\calc$. We also have similar equations for $\alpha_A^{G-}$.
	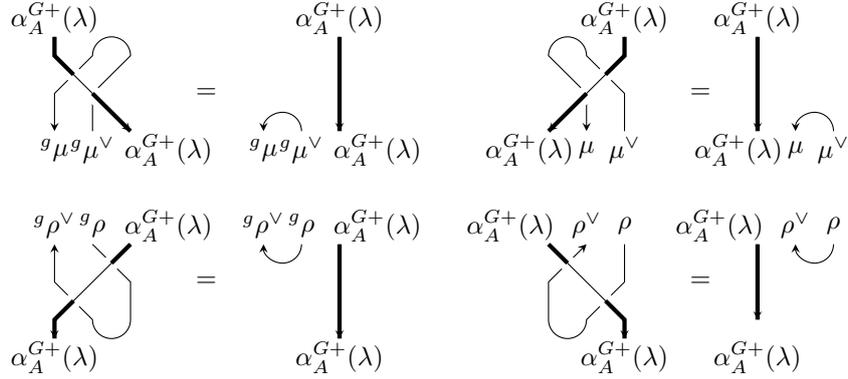
\begin{figure}[htb]
				\centering
				\begin{tikzpicture}
					\draw[->] (0.5,-1) -- (0.5,-0.5) -- (1,0) arc (0:180:0.25) -- (0,-0.5) -- (0,-1);
					\draw[cross,->] (0,0.25) -- (0,0) -- (1,-1);
					\draw[ultra thick] (0,0.25) -- (0,0) -- (0.25,-0.25);
					\draw[ultra thick] (0.5,-0.5) -- (1,-1);
					\node at (0,0.5){$\alpha_A^{G+}(\lambda)$};
					\node at (1.5,-1.25){$\alpha_A^{G+}(\lambda)$};
					\node at (0,-1.25){${}^g \mu$};
					\node at (0.5,-1.25){${}^g \mu^\vee$};
					\node at (2,-0.5){$=$};
					\begin{scope}[shift={(2.75,0)}]
						\draw[->] (0.5,-1) arc (0:180:0.25);
						\draw[cross,->] (1,0.25) -- (1,-1);
						\draw[ultra thick] (1,0.25) -- (1,-1);
						\node at (1,0.5){$\alpha_A^{G+}(\lambda)$};
						\node at (1.5,-1.25){$\alpha_A^{G+}(\lambda)$};
						\node at (0,-1.25){${}^g \mu$};
						\node at (0.5,-1.25){${}^g \mu^\vee$};
					\end{scope}
					\begin{scope}[shift={(6.5,0)}]
						\draw[->] (1,-1) -- (1,-0.5) -- (0.5,0) arc (0:180:0.25) -- (0.5,-0.5) -- (0.5,-1);
					\draw[cross,->] (1,0.25) -- (1,0) -- (0,-1);
					\draw[ultra thick] (1,0.25) -- (1,0) -- (0.75,-0.25);
					\draw[ultra thick] (0.5,-0.5) -- (0,-1);
					\node at (1,0.5){$\alpha_A^{G+}(\lambda)$};
					\node at (-0.25,-1.25){$\alpha_A^{G+}(\lambda)$};
					\node at (0.5,-1.25){$\mu$};
					\node at (1,-1.25){$\mu^\vee$};
					\node at (2,-0.5){$=$};
					\begin{scope}[shift={(2.75,0)}]
						\draw[->] (0,0.25) -- (0,-1);
						\draw[ultra thick] (0,0.25) -- (0,-1);
						\draw[cross,->] (1,-1) arc (0:180:0.25);
						\node at (0,0.5){$\alpha_A^{G+}(\lambda)$};
						\node at (-0.25,-1.25){$\alpha_A^{G+}(\lambda)$};
						\node at (0.5,-1.25){$\mu$};
						\node at (1,-1.25){$\mu^\vee$};
					\end{scope}	
					\end{scope}
					\begin{scope}[shift={(0,-2.5)}]
						\draw[<-] (0,0) -- (0,-0.5) -- (0.5,-1) arc (180:360:0.25) -- (1,-0.5) -- (0.5,0);
						\draw[cross,->] (1,0) -- (0,-1) -- (0,-1.25);
						\draw[ultra thick] (1,0) -- (0.75,-0.25);
						\draw[ultra thick] (0.25,-0.75) -- (0,-1) -- (0,-1.25);
						\node at (0,0.25){${}^g \rho^\vee$};
						\node at (0.5,0.25){${}^g \rho$};
						\node at (1.5,0.25){$\alpha_A^{G+}(\lambda)$};
						\node at (0,-1.5){$\alpha_A^{G+}(\lambda)$};
						\node at (2,-0.5){$=$};
						\begin{scope}[shift={(2.75,0)}]
							\draw[<-] (0,0) arc (180:360:0.25);
							\draw[cross,->] (1,0) -- (1,-1.25);
							\draw[ultra thick] (1,0) -- (1,-1.25);
							\node at (0,0.25){${}^g \rho^\vee$};
							\node at (0.5,0.25){${}^g \rho$};
							\node at (1.5,0.25){$\alpha_A^{G+}(\lambda)$};
							\node at (1,-1.5){$\alpha_A^{G+}(\lambda)$};	
						\end{scope}
						\begin{scope}[shift={(6.5,0)}]
							\draw[<-] (0.5,0) -- (0,-0.5) -- (0,-1) arc (180:360:0.25) -- (1,-0.5) -- (1,0);
							\draw[cross,->] (0,0) -- (1,-1) -- (1,-1.25);
							\draw[ultra thick] (0,0) -- (0.25,-0.25);
							\draw[ultra thick] (0.75,-0.75) -- (1,-1) -- (1,-1.25);
							\node at (-0.5,0.25){$\alpha_A^{G+}(\lambda)$};
							\node at (0.5,0.25){$\rho^\vee$};
							\node at (1,0.25){$\rho$};
							\node at (1,-1.5){$\alpha_A^{G+}(\lambda)$};
							\node at (2,-0.5){$=$};
							\begin{scope}[shift={(2.75,0)}]
								\draw[<-] (0.5,0) arc (180:360:0.25);
								\draw[cross,->] (0,0) -- (0,-1);
								\draw[ultra thick] (0,0) -- (0,-1);
								\node at (-0.5,0.25){$\alpha_A^{G+}(\lambda)$};
							\node at (0.5,0.25){$\rho^\vee$};
							\node at (1,0.25){$\rho$};
							\node at (0,-1.5){$\alpha_A^{G+}(\lambda)$};
							\end{scope}
						\end{scope}
					\end{scope}
				\end{tikzpicture}
				\caption{Arcs and thick crossings}
				\label{graphical_arc_alpha_braiding}
	\end{figure}
	\begin{proof}
		We only show the first equation because the proofs of other equations are similar. By explicitly writing down canonical morphisms in $\frob^G(\calc)$ as morphisms in $\calc$, we find that the left-hand and right-hand sides of the first equation in Figure \ref{graphical_arc_alpha_braiding} are equal to $r_{{}^g \mu {}^g \mu^\vee, \alpha_A^{G+}(\lambda)} f s_{\alpha^{G+}_A(\lambda),A}$ with $f$ given respectively by the leftmost and rightmost diagrams in Figure \ref{graphical_arc_alpha_braiding_proof}, which proves the statement. Note that we used at the second equality that $\mathrm{coev}_\mu$ is given by Figure \ref{graphicalbimodcoev} as in the case of bimodules.
		\begin{figure}[htb]
			\centering
			\begin{tikzpicture}
				\draw (0.25,0) arc (0:180:0.125) -- (-0.5,-0.5) -- (-0.5,-1) arc (0:-180:0.25) -- (-1,0) -- (-1,1);
				\node[block] at (-0.5,-0.75){$z_g$};
				\draw[rounded corners] (-0.75,-1.25) -- (-0.75,-1.5) -- (-0.25,-1.5);
				\draw[rounded corners](0.25,0) -- (0.25,-0.25) -- (0.5,-0.25);
				\draw (0.125,0.125) -- (0.125,0.35);
				\draw[rounded corners] (0.5,0.5) -- (0.25,0.5) -- (0.25,1);
				\draw[fill=white] (0.125,0.35) circle [radius=0.1];
				\draw[->] (0,-3.5) -- (0,-2.75) --(0.25,-2.5) -- (0.25,-1.75) -- (1,-1) -- (1,0.75) arc (0:180:0.25) -- (0.5,-0.5) -- (-0.25,-1.25) -- (-0.25,-2.5) -- (-0.5,-2.75) -- (-0.5,-3.5);
				\draw[rounded corners] (0.25,-2.25) -- (0.5,-2.25) -- (0.5,-2);
				\draw (0.5,-2) arc (180:0:0.125) -- (0.75,-3) arc (360:180:0.125) arc (0:180:0.125);
				\draw[rounded corners] (0.25,-3) -- (0.25,-3.25) -- (0,-3.25);
				\draw (0.375,-2.875) -- (0.375,-2.65);
				\draw[fill=white] (0.375,-2.65) circle [radius=0.1];
				\draw (0.625,-1.875) -- (0.625,-1.65);
				\draw[fill=white] (0.625,-1.65) circle [radius=0.1];
				\draw[->] (0.625,-3.125) -- (0.625,-3.5);
				\draw[cross,->] (-0.5,1) -- (-0.5,0) -- (1,-1.5) -- (1,-3.5);
				\draw (0.5,0.6) arc (90:270:0.1);
				\draw (0.5,-0.15) arc (90:270:0.1);
				\draw (-0.25,-1.4) arc (90:270:0.1);
				\draw (0.25,-2.15) arc (90:-90:0.1);
				\draw (0,-3.15) arc (90:-90:0.1);
				\node at (-1,1.25){$A$};
				\node at (-0.5,1.25){$\lambda$};
				\node at (0.25,1.25){$A$};
				\node at (-0.5,-3.75){${}^g \mu$};
				\node at (0,-3.75){${}^g \mu^\vee$};
				\node at (0.625,-3.75){$A$};
				\node at (1,-3.75){$\lambda$};
				\node at (1.5,-1.25){$=$};
				\begin{scope}[shift={(3,0)}]
					\draw (0,-1) arc (0:180:0.25) -- (-0.5,-1.5) arc (0:-180:0.25) -- (-1,0.5) -- (-1,1);
					\draw (-0.25,-0.75) -- (-0.25,0.25) -- (0.5,1);
					\node[block] at (-0.5,-1.25){$z_g$};
					\draw[rounded corners] (-0.75,-1.75) -- (-0.75,-2) -- (-0.25,-2);
					\draw[rounded corners](0,-1) -- (0,-1.25) -- (0.25,-1.25);
					\draw[->] (0,-3.5) -- (0,-2.75) --(0.25,-2.5) -- (0.25,-2.25) -- (0.75,-1.5) -- (0.75,-1) arc (0:180:0.25) -- (0.25,-1.5) -- (-0.25,-1.75) -- (-0.25,-2.5) -- (-0.5,-2.75) -- (-0.5,-3.5);
					\draw[<-] (0.5,-3.5) -- (0.5,-3) arc (0:180:0.125);
					\draw[rounded corners] (0.25,-3) -- (0.25,-3.25) -- (0,-3.25);
					\draw (0.375,-2.875) -- (0.375,-2.65);
					\draw[fill=white] (0.375,-2.65) circle [radius=0.1];
					\draw[cross,->] (-0.5,1) -- (1,-0.5) -- (1,-3.5);
					\draw (-0.25,-1.9) arc (90:270:0.1);
					\draw (0.25,-1.15) arc (90:270:0.1);
					\draw (0,-3.15) arc (90:-90:0.1);
					\node at (-1,1.25){$A$};
					\node at (-0.5,1.25){$\lambda$};
					\node at (0.5,1.25){$A$};
					\node at (-0.5,-3.75){${}^g \mu$};
					\node at (0,-3.75){${}^g \mu^\vee$};
					\node at (0.5,-3.75){$A$};
					\node at (1,-3.75){$\lambda$};
					\node at (1.5,-1.25){$=$};
				\end{scope}
				\begin{scope}[shift={(6,0)}]
					\draw[<-] (0.5,-3.5) -- (0.5,-1.5) arc (0:180:0.5) arc (0:-180:0.25) -- (-1,0.5) -- (-1,1);
					\draw (0,-1) -- (0,0) -- (1,1);
					\node[block] at (0,-0.5){$z_g$};
					\draw[rounded corners] (-0.75,-1.75) -- (-0.75,-3.25) -- (-0.5,-3.25);
					\draw[->] (0,-3.5) -- (0,-3) arc (0:180:0.25) -- (-0.5,-3.5);
					\draw[cross,->] (-0.5,1) -- (1,-0.5) -- (1,-3.5);
					\draw (-0.5,-3.15) arc (90:270:0.1);
					\node at (-1,1.25){$A$};
					\node at (-0.5,1.25){$\lambda$};
					\node at (1,1.25){$A$};
					\node at (-0.5,-3.75){${}^g \mu$};
					\node at (0,-3.75){${}^g \mu^\vee$};
					\node at (0.5,-3.75){$A$};
					\node at (1,-3.75){$\lambda$};
					\node at (1.5,-1.25){$=$};	
				\end{scope}
				\begin{scope}[shift={(9,0)}]
					\draw (-1,1) -- (-1,0) arc (180:360:0.25);
					\draw (-0.5,0) -- (-0.5,0.5) -- (0,1);
					\node[block] at (-0.5,0.25){$z_g$};
					\draw[rounded corners] (-0.75,-2.75) -- (-0.75,-3.25) -- (-0.5,-3.25);
					\draw[->] (0,-3.5) -- (0,-3) arc (0:180:0.25) -- (-0.5,-3.5);
					\draw[cross,->] (-0.5,1) -- (1,-0.5) -- (1,-3.5);
					\draw (-0.5,-3.15) arc (90:270:0.1);
					\draw (-1,-1.5) -- (-1,-2.5) arc (180:360:0.25) arc (180:0:0.25) arc (180:360:0.25);
					\draw (-1,-1.5) arc (180:0:0.25) -- (0.25,-2.25);
					\draw (0.75,-2.25) arc (360:180:0.25);
					\draw (-0.75,-0.25) -- (-0.75,-0.5) -- (0.75,-2) -- (0.75,-2.25);
					\draw (-0.75,-1.25) -- (-0.75,-1);
					\draw[fill=white] (-0.75,-1) circle [radius=0.1];
					\draw (-0.25,-2.25) -- (-0.25,-2);
					\draw[fill=white] (-0.25,-2) circle [radius=0.1];
					\node at (-1,1.25){$A$};
					\node at (-0.5,1.25){$\lambda$};
					\node at (0,1.25){$A$};
					\node at (-0.5,-3.75){${}^g \mu$};
					\node at (0,-3.75){${}^g \mu^\vee$};
					\node at (0.5,-3.75){$A$};
					\node at (1,-3.75){$\lambda$};	
				\end{scope}
			\end{tikzpicture}
			\caption{The proof of Figure \ref{graphical_arc_alpha_braiding}}
			\label{graphical_arc_alpha_braiding_proof}
		\end{figure}
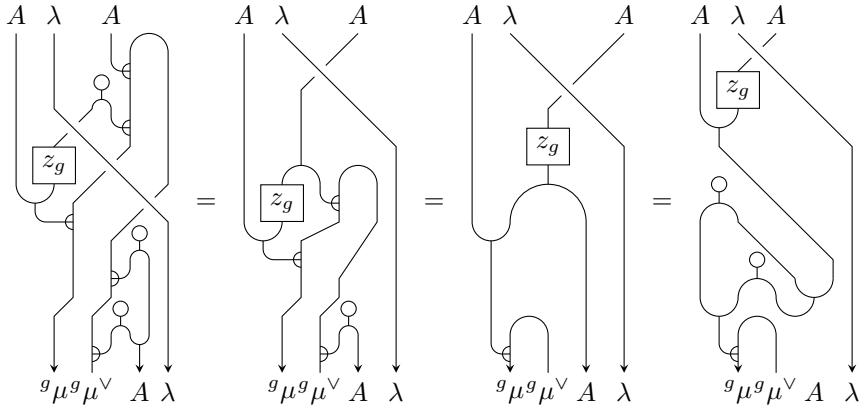
	\end{proof}
\end{lem}

We can also move thick arcs by the following lemma since $\alpha^{G \pm}_A$ are tensor functors and therefore preserve arcs.

\begin{lem}
	\label{lem_thick_crossing_product}
	Let $A$ be a neutral symmetric special $G$-equivariant Frobenius algebra in a $G$-braided multitensor category $\calc$. Then, the equations in Figure \ref{graphical_thick_crossing_product} hold for $\lambda_1 \in \obj(\calc_{g_1})$, $\lambda_2 \in \obj (\calc_{g_2})$, a homogeneous left $A$-module $\mu$ in $\calc$ and a homogeneous right $A$-module $\rho$. Similar statements hold for $\alpha_A^{G-}$.
	\begin{figure}[htb]
		\centering
		\begin{tikzpicture}
			\draw[->] (1,0) -- (0,-1);
			\draw[cross,->] (0,0) -- (0,-0.5) -- (0.5,-1);
			\draw[cross,->] (0.5,0) -- (1,-0.5) -- (1,-1);
			\draw[ultra thick] (0,0) -- (0,-0.5) -- (0.25,-0.75);
			\draw[ultra thick] (0.5,0) -- (0.75,-0.25);
			\node at (0.5,0.25){$\alpha_A^{G+}(\lambda_1)\alpha_A^{G+}(\lambda_2)\mu$};
			\node at (-0.25,-1.25){${}^{g_1 g_2} \mu$};
			\node at (0.5,-1.25){$\lambda_1$};
			\node at (1,-1.25){$\lambda_2$};
			\node at (2,-0.5){$=$};
			\begin{scope}[shift={(3,0)}]
				\draw[->] (1,0) -- (0,-1);
				\draw[cross,->] (0,0) -- (1,-1);
				\draw[ultra thick] (0,0) -- (0.5,-0.5);
				\node at (0,0.25){$\alpha_A^{G+}(\lambda_1 \lambda_2)$};
				\node at (1,0.25){$\mu$};
				\node at (0,-1.25){${}^{g_1 g_2} \mu$};
				\node at (1,-1.25){$\lambda_1 \lambda_2$};
			\end{scope}
			\begin{scope}[shift={(6.5,0)}]
				\draw[->] (0,0) -- (1,-1);
				\draw[cross,->] (0.5,0) -- (0,-0.5) -- (0,-1);
				\draw[cross,->] (1,0) -- (1,-0.5) -- (0.5,-1);
				\draw[ultra thick] (0.5,0) -- (0.25,-0.25);
				\draw[ultra thick] (1,0) -- (1,-0.5) -- (0.75,-0.75);
				\node at (0.5,0.25){${}^{g_1 g_2} \rho \alpha_A^{G+}(\lambda_1) \alpha_A^{G+}(\lambda_2)$};
				\node at (0,-1.25){$\lambda_1$};
				\node at (0.5,-1.25){$\lambda_2$};
				\node at (1,-1.25){$\rho$};
				\node at (2,-0.5){$=$};
				\begin{scope}[shift={(3,0)}]
					\draw[->] (0,0) -- (1,-1);
					\draw[cross,->] (1,0) -- (0,-1);
					\draw[ultra thick] (1,0) -- (0.5,-0.5);
					\node at (-0.25,0.25){${}^{g_1 g_2} \rho$};
					\node at (1.25,0.25){$\alpha_A^{G+}(\lambda_1 \lambda_2)$};
					\node at (0,-1.25){$\lambda_1 \lambda_2$};
					\node at (1,-1.25){$\rho$};
				\end{scope}
			\end{scope}
		\end{tikzpicture}
		\caption{Thick crossings and monoidal products}
		\label{graphical_thick_crossing_product}
	\end{figure}
	\begin{proof}
		We only show the first equation because the proofs of other equations are similar. By a similar argument to that in the proof of Lemma \ref{lem_alpha_ordinary_braiding}, we see that the left-hand side of the first equation in Figure \ref{graphical_thick_crossing_product} is equal to $f (\id_{A \lambda} \otimes s_{\alpha_A^{G+}(\lambda_2), \mu}) s_{\alpha_A^{G+}(\lambda_1), \alpha_A^{G+}(\lambda_2) \otimes_A \mu}$ with $f$ given by the leftmost diagram in Figure \ref{graphical_thick_crossing_product_proof}, where $z \coloneqq z^A$. On the other hand, since $r_{\alpha^{G+}_A(\lambda_1), \alpha^{G+}_A(\lambda_2)}$ is given in Figure \ref{graphicalalphaindsplit}, the right-hand side of the first equation in Figure \ref{graphical_thick_crossing_product} is equal to $f (\id_{A \lambda} \otimes s_{\alpha_A^{G+}(\lambda_2), \mu}) s_{\alpha_A^{G+}(\lambda_1), \alpha_A^{G+}(\lambda_2) \otimes_A \mu}$ with $f$ given by the leftmost diagram in Figure \ref{graphical_thick_crossing_product_proof}, which proves the statement. 
		\begin{figure}[htb]
			\centering
			\begin{tikzpicture}
				\draw[->] (1,0) -- (0.5,-0.5) -- (0.5,-1) -- (0,-1.5) -- (0,-2);
				\draw[cross,->] (0.5,0) -- (1,-0.5) -- (1,-2);
				\draw[rounded corners] (0,0) -- (0,-0.75) -- (0.5,-0.75);
				\draw[cross,->] (-0.5,0) -- (-0.5,-1) -- (0.5,-1.5) -- (0.5,-2);
				\draw[rounded corners] (-1,0) -- (-1,-1.75) -- (0,-1.75);
				\draw (0.5,-0.65) arc (90:270:0.1);
				\draw (0,-1.65) arc (90:270:0.1);
				\node at (-1,0.25){$A$};
				\node at (-0.5,0.25){$\lambda_1$};
				\node at (0,0.25){$A$};
				\node at (0.5,0.25){$\lambda_2$};
				\node at (1,0.25){$\mu$};
				\node at (-0.25,-2.25){${}^{g_1 g_2} \mu$};
				\node at (0.5,-2.25){$\lambda_1$};
				\node at (1,-2.25){$\lambda_2$};
				\node at (1.5,-1){$=$};
				\begin{scope}[shift={(3,0)}]
					\draw[->] (1,0) -- (0,-1) -- (0,-1.5) -- (0,-2);
					\draw[cross,->] (0.5,0) -- (1,-0.5) -- (1,-2);
					\draw[rounded corners] (0,0) -- (-0.5,-0.5) -- (-0.5,-1.25) -- (0,-1.25);
					\draw[cross,->] (-0.5,0) -- (0.5,-1) -- (0.5,-1.5) -- (0.5,-2);
					\draw[rounded corners] (-1,0) -- (-1,-1.75) -- (0,-1.75);
					\draw (0,-1.15) arc (90:270:0.1);
					\draw (0,-1.65) arc (90:270:0.1);
					\node at (-1,0.25){$A$};
					\node at (-0.5,0.25){$\lambda_1$};
					\node at (0,0.25){$A$};
					\node at (0.5,0.25){$\lambda_2$};
					\node at (1,0.25){$\mu$};
					\node at (-0.25,-2.25){${}^{g_1 g_2} \mu$};
					\node at (0.5,-2.25){$\lambda_1$};
					\node at (1,-2.25){$\lambda_2$};
					\node at (-0.5,-1.5){${}^{g_1} A$};
					\node at (1.5,-1){$=$};
				\end{scope}
				\begin{scope}[shift={(6,0)}]
					\draw[->] (1,0) -- (0,-1) -- (0,-1.5) -- (0,-2);
					\draw[cross,->] (0.5,0) -- (1,-0.5) -- (1,-2);
					\draw[rounded corners] (0,0) -- (-0.5,-0.5) -- (-0.5,-1.25) -- (0,-1.25);
					\draw[cross,->] (-0.5,0) -- (0.5,-1) -- (0.5,-1.5) -- (0.5,-2);
					\draw[rounded corners] (-1,0) -- (-1,-1.75) -- (0,-1.75);
					\draw (0,-1.15) arc (90:270:0.1);
					\draw (0,-1.65) arc (90:270:0.1);
					\node at (-1,0.25){$A$};
					\node at (-0.5,0.25){$\lambda_1$};
					\node at (0,0.25){$A$};
					\node at (0.5,0.25){$\lambda_2$};
					\node at (1,0.25){$\mu$};
					\node at (-0.25,-2.25){${}^{g_1 g_2} \mu$};
					\node at (0.5,-2.25){$\lambda_1$};
					\node at (1,-2.25){$\lambda_2$};
					\node at (-0.5,-1.5){$A$};
					\node at (1.5,-1){$=$};
					\node[block] at (-0.5,-0.75){$z_{g_1}$};
				\end{scope}
				\begin{scope}[shift={(9,0)}]
					\draw[->] (1,0) -- (0,-1) -- (0,-1.5) -- (0,-2);
					\draw[cross,->] (0.5,0) -- (1,-0.5) -- (1,-2);
					\draw[rounded corners] (0,0) -- (-0.5,-0.5) -- (-0.5,-1);
					\draw[cross,->] (-0.5,0) -- (0.5,-1) -- (0.5,-1.5) -- (0.5,-2);
					\draw (-1,0) -- (-1,-1) arc (180:360:0.25);
					\draw[rounded corners] (-0.75,-1.25) -- (-0.75,-1.5) -- (0,-1.5);
					\draw (0,-1.4) arc (90:270:0.1);
					\node at (-1,0.25){$A$};
					\node at (-0.5,0.25){$\lambda_1$};
					\node at (0,0.25){$A$};
					\node at (0.5,0.25){$\lambda_2$};
					\node at (1,0.25){$\mu$};
					\node at (-0.25,-2.25){${}^{g_1 g_2} \mu$};
					\node at (0.5,-2.25){$\lambda_1$};
					\node at (1,-2.25){$\lambda_2$};
					\node[block] at (-0.5,-0.75){$z_{g_1}$};
				\end{scope}
			\end{tikzpicture}
			\caption{The proof of Figure \ref{graphical_thick_crossing_product}}
			\label{graphical_thick_crossing_product_proof}
		\end{figure}
	\end{proof}
\end{lem}

\begin{cor}
	\label{cor_thick_rotation}
	Let $A$ be a neutral symmetric special $G$-equivariant Frobenius algebra in a $G$-braided multitensor category $\calc$. Then, the equations in Figure \ref{graphical_rotation_thick} hold for $\lambda \in \obj(\calc_g)$, a homogeneous left $A$-module $\mu$ in $\calc$ and a homogeneous right $A$-module $\rho$. Similar statements hold for $\alpha_A^{G-}$.
	\begin{figure}[htb]
		\centering
		\begin{tikzpicture}
			\draw[->] (0,0) -- (-1,-1);
			\draw[->, cross] (-1,0) -- (0,-1);
			\draw[ultra thick] (-1,0) -- (-0.5,-0.5);
			\node at (-1, 0.25){$\alpha_A^{G+} (\lambda)$};
			\node at (0, 0.25){$\mu$};
			\node at (-1, -1.25){${}^g \mu$};
			\node at (0, -1.25){$\lambda$};
			\node at (-1.75,-0.5) {$=$};
			\node at (0.75,-0.5) {$=$};
			\begin{scope}[shift={(-3,0)}]
			\draw[->] (-1,0) -- (0,-1);
			\draw[->, cross] (-1.5,0) -- (-1.5,-1) arc (180:360:0.25) -- (0,0) arc (180:0:0.25) -- (0.5, -1);
			\draw[ultra thick] (-1.5,0) -- (-1.5,-1) arc (180:360:0.25) -- (-0.5,-0.5); 
			\node at (0.5, -1.25){$\lambda$};
			\node at (-1, 0.25){$\mu$};
			\node at (0, -1.25){${}^g \mu$};
			\node at (-1.75, 0.25){$\alpha_A^{G+}( \lambda)$};
			\end{scope}
			\begin{scope}[shift={(3,0)}]
			\draw[->] (0.5,0) -- (0.5,-1) arc (360:180:0.25) -- (-1,0) arc (0:180:0.25) -- (-1.5,-1);
			\draw[->, cross] (0,0) -- (-1,-1);
			\draw[ultra thick] (0,0) -- (-0.5,-0.5);
			\node at (-0.25, 0.25){$\alpha_A^{G+} (\lambda)$};
			\node at (0.5, 0.25){$\mu$};
			\node at (-1.5, -1.25){${}^g \mu$};
			\node at (-1, -1.25){$\lambda$};
			\end{scope}
			\begin{scope}[shift={(0,-2.5)}]
				\draw[->] (-1,0) -- (0,-1);
				\draw[->,cross] (0,0) -- (-1,-1);
				\draw[ultra thick] (0,0) -- (-0.5,-0.5);
				\node at (0, 0.25){$\alpha_A^{G+} (\lambda)$};
				\node at (-1, 0.25){${}^g \rho$};
				\node at (0, -1.25){$\rho$};
				\node at (-1, -1.25){$\lambda$};
				\node at (-1.75,-0.5) {$=$};
				\node at (0.75,-0.5) {$=$};
				\begin{scope}[shift={(-3,0)}]
				\draw[->] (-1.5,0) -- (-1.5,-1) arc (180:360:0.25) -- (0,0) arc (180:0:0.25) -- (0.5,-1);
				\draw[->,cross] (-1,0) -- (0,-1);
				\draw[ultra thick] (-1,0) -- (-0.5,-0.5); 
				\node at (0.5, -1.25){$\rho$};
				\node at (-1.75, 0.25){${}^g \rho$};
				\node at (0, -1.25){$\lambda$};
				\node at (-0.75, 0.25){$\alpha_A^{G+}( \lambda)$};
				\end{scope}
				\begin{scope}[shift={(3,0)}]
				\draw[->] (0,0) -- (-1,-1);
				\draw[->,cross] (0.5,0) -- (0.5,-1) arc (360:180:0.25) -- (-1,0) arc (0:180:0.25) -- (-1.5,-1);
				\draw[ultra thick] (0.5,0) -- (0.5,-1) arc (360:180:0.25) -- (-0.5,-0.5);
				\node at (0.75, 0.25){$\alpha_A^{G+} (\lambda)$};
				\node at (-0.25, 0.25){${}^g \rho$};
				\node at (-1.5, -1.25){$\lambda$};
				\node at (-1, -1.25){$\rho$};
				\end{scope}
			\end{scope}
		\end{tikzpicture}
		\caption{Rotation of thick crossings}
		\label{graphical_rotation_thick}
	\end{figure}
	\begin{proof}
		The same argument as that in the proof of Lemma \ref{lem_crossing_rotation} works by Lemmata \ref{lem_alpha_ordinary_braiding}, \ref{lem_arc_alpha_braiding} and \ref{lem_thick_crossing_product}.
	\end{proof}
\end{cor}

Next, we describe $\hom(\alpha_A^{G+}(\lambda),\alpha_A^{G-}(\mu))$ as the $G$-equivariant version of the space of local morphisms \cite[Definition 4.14]{MR3424476}.

\begin{defi}
	Let $A$ be a neutral symmetric special $G$-equivariant Frobenius algebra in a $G$-braided multitensor category $\calc$. For $\lambda \in \obj(\calc)$, we define a morphism $\tilde{P}^{G}_A(\lambda) \in \en(A\lambda)$ to be the morphism in Figure \ref{ptildeplusdef}, where $z \coloneqq z^A$.
	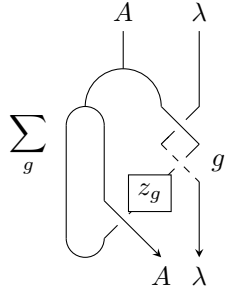
\begin{figure}[htb]
		\centering
		\begin{tikzpicture}
		\draw (0,0) -- (0,0.5);
		\draw (0.5,-0.5) arc (0:180:0.5);
		\draw (-0.25,-0.75) arc (0:180:0.25);
		\draw (1,0.5) -- (1,-0.5) -- (0.5,-1);
		\draw[cross] (0.5, -0.5) -- (1,-1);
		\draw (1,-1) -- (-0.25,-2.25);
		\draw (-0.25,-2.25) arc (360:180:0.25);
		\draw (-0.75,-0.75) -- (-0.75,-2.25);
		\draw[cross,dashed] (0.5, -1) -- (1,-1.5);
		\draw (0.5,-1) -- (0.6,-0.9);
		\draw[->] (1,-1.5) -- (1,-2.5);
		\draw[cross,->] (-0.25, -1.75) -- (0.5,-2.5);
		\draw (-0.25,-0.75) --(-0.25,-1.75);
		\node at (0,0.75){$A$};
		\node at (1,0.75){$\lambda$};
		\node at (0.5,-2.75){$A$};
		\node at (1,-2.75){$\lambda$};
		\node at (1.25,-1.25){$g$};
		\node at (-1.25,-1){$\displaystyle \sum_g$};
		\node[block] at (0.35,-1.65){$z_g$};
		\end{tikzpicture}
		\caption{$\tilde{P}_{A}^{G}(\lambda)$}
		\label{ptildeplusdef}
	\end{figure}
\end{defi}

\begin{rem}
Note that the definition of $P_A^G(\lambda)$ in Proposition \ref{proposition_p_lambda} requires $\lambda \in \obj(\calc^G)$ but does not require that $A$ is neutral, while the definition of $\tilde{P}_A^G(\lambda)$ requires that $A$ is neutral but does not require $\lambda \in \obj(\calc^G)$. If $\lambda \in \obj(\calc^G)$ and $A$ is neutral, then $P_A^G(\lambda) = \tilde{P}_A^{G}(\lambda)$ by definition.
\end{rem}

\begin{lem}
	\label{lem_tildep_idempotent}
	Let $A$ be a neutral symmetric special $G$-equivariant Frobenius algebra in a $G$-braided multitensor category $\calc$. Then, $\tilde{P}_A^{G}(\lambda)$ is an idempotent for any $\lambda \in \obj(\calc)$. 

	\begin{proof}
		The proof is similar to that of \cite[Lemma 4.26]{bklr}. Indeed, the statement follows from the graphical calculation in Figure \ref{graphical_ptilde_idempotent}, where we only gave calculations on equivariant structures. 
		\begin{figure}[htb]
			\centering
			\begin{tikzpicture}
				\draw (0,0) -- (0,0.5);
				\draw (0.5,-0.5) arc (0:180:0.5);
				\draw (-0.25,-0.75) arc (0:180:0.25);
				\draw (1,0.5) -- (1,-0.5) -- (0.5,-1);
				\draw[cross] (0.5, -0.5) -- (1,-1);
				\draw (1,-1) -- (-0.25,-2.25);
				\draw (-0.25,-2.25) arc (360:180:0.25);
				\draw (-0.75,-0.75) -- (-0.75,-2.25);
				\draw[cross,dashed] (0.5, -1) -- (1,-1.5);
				\draw (0.5,-1) -- (0.6,-0.9);
				\draw (1,-1.5) -- (1,-2.5);
				\draw[cross] (-0.25, -1.75) -- (0,-2) -- (0,-2.5);
				\draw (-0.25,-0.75) --(-0.25,-1.75);
				\node at (0,0.75){$A$};
				\node at (1,0.75){$\lambda$};
				\node at (1.25,-1.25){$g$};
				\node at (-1.25,-2.5){$\displaystyle \sum_g$};
				\node at (1.75,-2.5){$=$};
				\node[block] at (0.35,-1.65){$z_g$};
				\begin{scope}[shift={(0,-3)}]
					\draw (0,0) -- (0,0.5);
					\draw (0.5,-0.5) arc (0:180:0.5);
					\draw (-0.25,-0.75) arc (0:180:0.25);
					\draw (1,0.5) -- (1,-0.5) -- (0.5,-1);
					\draw[cross] (0.5, -0.5) -- (1,-1);
					\draw (1,-1) -- (-0.25,-2.25);
					\draw (-0.25,-2.25) arc (360:180:0.25);
					\draw (-0.75,-0.75) -- (-0.75,-2.25);
					\draw[cross,dashed] (0.5, -1) -- (1,-1.5);
					\draw (0.5,-1) -- (0.6,-0.9);
					\draw[->] (1,-1.5) -- (1,-2.5);
					\draw[cross,->] (-0.25, -1.75) -- (0.5,-2.5);
					\draw (-0.25,-0.75) --(-0.25,-1.75);
					\node at (0.5,-2.75){$A$};
					\node at (1,-2.75){$\lambda$};
					\node at (1.25,-1.25){$g$};
					\node[block] at (0.35,-1.65){$z_g$};
				\end{scope}
				\begin{scope}[shift={(3.75,-0.75)}]
					\draw (0,0) -- (0,0.5);
					\draw (0.5,-0.5) arc (0:180:0.5);
					\draw (-0.25,-0.75) arc (0:180:0.25);
					\draw (1,0.5) -- (1,-0.5) -- (0.5,-1);
					\draw[cross] (0.5, -0.5) -- (1,-1);
					\draw (1,-1) -- (0.5,-1.5) -- (0.5,-1.75);
					\draw (0.25,-2) arc (180:0:0.25) -- (0.75,-2.5) arc (360:180:0.25) -- (0.25,-2);
					\draw (0.5,-2.75) -- (0.5,-3) -- (-0.25,-3.75);
					\draw (-0.25,-3.75) arc (360:180:0.25);
					\draw (-0.75,-0.75) -- (-0.75,-3.75);
					\draw[cross,dashed] (0.5, -1) -- (1.25,-1.75);
					\draw (0.5,-1) -- (0.6,-0.9);
					\draw[->] (1.25,-1.75) -- (1.25,-4);
					\draw[cross,->] (-0.25, -3) -- (0.5,-3.75) -- (0.5,-4);
					\draw (-0.25,-0.75) --(-0.25,-3);
					\node at (0,0.75){$A$};
					\node at (1,0.75){$\lambda$};
					\node at (1.25,-1.25){$g$};
					\node at (-1.25,-1.75){$\displaystyle \sum_g$};
					\node[block] at (0.15,-2.25){$z_g$};
					\node[block] at (0.85,-2.25){$z_g$};
					\node at (0.5,-4.25){$A$};
					\node at (1.25,-4.25){$\lambda$};
					\node at (1.75,-1.75){$=$};
				\end{scope}
				\begin{scope}[shift={(7.5,-1.5)}]
					\draw (0,0) -- (0,0.5);
					\draw (0.5,-0.5) arc (0:180:0.5);
					\draw (-0.25,-0.75) arc (0:180:0.25);
					\draw (1,0.5) -- (1,-0.5) -- (0.5,-1);
					\draw[cross] (0.5, -0.5) -- (1,-1);
					\draw (1,-1) -- (-0.25,-2.25);
					\draw (-0.25,-2.25) arc (360:180:0.25);
					\draw (-0.75,-0.75) -- (-0.75,-2.25);
					\draw[cross,dashed] (0.5, -1) -- (1,-1.5);
					\draw (0.5,-1) -- (0.6,-0.9);
					\draw[->] (1,-1.5) -- (1,-2.5);
					\draw[cross,->] (-0.25, -1.75) -- (0.5,-2.5);
					\draw (-0.25,-0.75) --(-0.25,-1.75);
					\node at (0,0.75){$A$};
					\node at (1,0.75){$\lambda$};
					\node at (0.5,-2.75){$A$};
					\node at (1,-2.75){$\lambda$};
					\node at (1.25,-1.25){$g$};
					\node at (-1.25,-1){$\displaystyle \sum_g$};
					\node[block] at (0.35,-1.65){$z_g$};
				\end{scope}
			\end{tikzpicture}
			\caption{$\tilde{P}_A^{G}(\lambda)$ is an idempotent}
			\label{graphical_ptilde_idempotent}
		\end{figure}
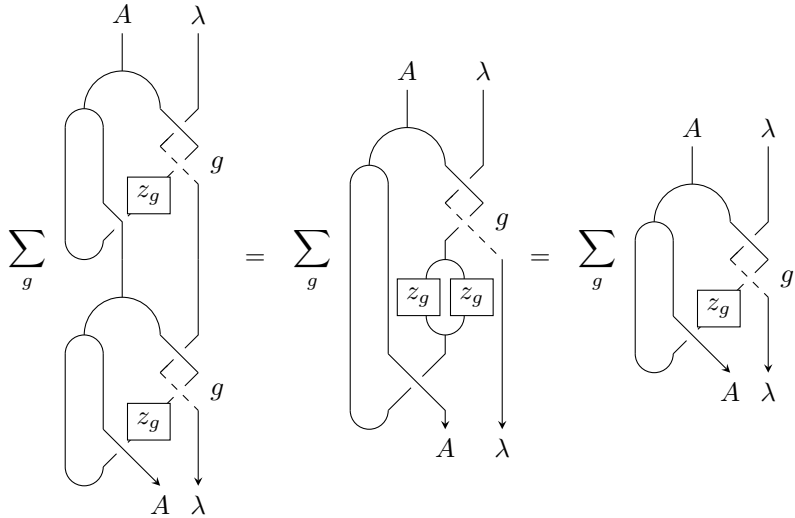
	\end{proof}
\end{lem}

\begin{lem}
	\label{localhomlem}
	Let $A$ be a neutral symmetric special $G$-equivariant Frobenius algebra in a $G$-ribbon multitensor category $\calc$. For $\lambda,\mu \in \obj(\calc_g)$, define finite dimensional vector spaces
	\begin{align*}
		\hom^G_{\mathrm{loc}} (A\lambda, \mu) &\coloneqq \{ f \in \hom(A \lambda, \mu) \mid f \circ \tilde{P}^{G}_A(\lambda) = f \} \\
		\hom^G_{\mathrm{loc}} (\mu, A\lambda) &\coloneqq \{ f \in \hom(\mu, A\lambda) \mid \tilde{P}^{G}_A(\lambda) \circ f = f \}.
	\end{align*}
	Then, $\hom_{\mathrm{loc}}^G(\lambda, A\mu) \cong \hom(\alpha_A^{G+}(\lambda),\alpha^{G-}_A(\mu))$ by the map
	\begin{align*}
		I_{\lambda,\mu}^A: \hom_{\mathrm{loc}}^G(\lambda,A\mu) \to \hom(\alpha_A^{G+}(\lambda),\alpha^{G-}_A(\mu))
	\end{align*}
	defined in Figure \ref{graphicallocalhom} for $f \in \hom_{\mathrm{loc}}^G(\lambda,A\mu)$, where ${}_A A$ (resp. ${}_A {}^g A$, resp. $A_A$) denotes $A$ as a left $A$-module (resp. ${}^g A$ as a left $A$-module, resp. $A$ as a right $A$-module). Note that ${}_A A^\vee = {}^\vee_A A = A_A$ by \cite[Section 5]{MR2075605}. We similarly have the isomorphism
	\begin{align*}
		\tilde{I}^A_{\lambda, \mu}: \hom_{\mathrm{loc}}^G(A\mu, \lambda) \to \hom(\alpha_A^{G-}(\mu),\alpha^{G+}_A(\lambda))
	\end{align*}
	defined in Figure \ref{graphical_localhom_dual} for $\tilde{f} \in \hom_{\mathrm{loc}}^G(A \mu, \lambda)$.
	\begin{figure}[htb]
		\begin{tabular}{cc}
			\begin{minipage}[b]{0.45\hsize}
				\centering
				\begin{tikzpicture}
					\draw[<-] (0,-0.5)--(0,0) -- (0.5,0.5);
					\draw[ultra thick] (0,-0.5)--(0,0) -- (0.25,0.25);
					\draw[cross] (-0.5,0.5) arc (360:180:0.25) -- (-1,1) -- (-0.25,1.75) arc (180:0:0.25) -- (1,1) -- (1,0) arc (360:180:0.25) -- (0,0.5);
					\draw[cross] (-1,2) -- (0,1); 
					\draw[ultra thick] (-1,2) -- (-0.5,1.5);
					\draw (-0.65,1) rectangle (0.75,0.5) node[midway]{$f$};
					\node[block] at (-1,0.75){$z_g$};
					\node at (-1,2.25){$\alpha^{G+}_A(\lambda)$};
					\node at (0,-0.75){$\alpha^{G-}_A(\mu)$};
					\node at (0,1.35){$\lambda$};
					\node at (0.75,0.25){$\mu$};
					\node at (0.75,1.75){$A_A$};
					\node at (-1.25,1.25){${}_A {}^g A$};
					\node at (-1.25,0.25){${}_A A$};
				\end{tikzpicture}
				\caption{$I^A_{\lambda,\mu}(f)$}
				\label{graphicallocalhom}
			\end{minipage}
			\begin{minipage}[b]{0.45\hsize}
				\centering
				\begin{tikzpicture}
					\draw[<-] (0.5,0) -- (0,0.5) -- (0,1.5);
					\draw[ultra thick] (0,1) -- (0,0.5) -- (0.25,0.25);
					\draw[cross,<->] (0,0) -- (0.5,0.5) -- (0.5,2) arc (0:180:0.25) -- (-0.75,1.25) -- (-0.75,0);
					\draw[cross] (-0.5,2.25) -- (-0.5,2) -- (0,1.5) -- (0,1);
					\draw[ultra thick] (-0.25,1.75) -- (0,1.5) -- (0,1);
					\node[block] at (0,1){$f'$}; 
					\node[block] at (-0.75,1){$z_g$};
					\node at (0,-0.25){${}_A A$};
					\node at (0.5,-0.25){$\mu$};
					\node at (-0.75,-0.25){$A_A$};
					\node at (-0.5,2.5){$\lambda$};
				\end{tikzpicture}
				\caption{$I^{A-1}_{\lambda,\mu}(f')$}
				\label{graphicallocalhominv}
			\end{minipage}
		\end{tabular}
	\end{figure}
	\begin{figure}[htb]
		\begin{tabular}{cc}
			\begin{minipage}[b]{0.45\hsize}
				\centering
				\begin{tikzpicture}
					\draw (0.5,1) -- (0,1.5); 
					\draw[ultra thick] (0,2)--(0,1.5) -- (0.25,1.25);
					\draw[cross] (0,1) -- (0.5,1.5) arc (180:0:0.25) -- (1,0.5) -- (0,-0.5) -- (0,-1) arc (360:180:0.25) -- (-1,-0.5) -- (-1,1) arc (180:0:0.25);
					\draw[cross,->] (0,0.5) -- (1,-0.5) -- (1,-1.5);
					\draw[ultra thick] (0.5,0) -- (1,-0.5) -- (1,-1.5);
					\draw (-0.65,1) rectangle (0.75,0.5) node[midway]{$\tilde{f}$};
					\node[block] at (0,-0.75){$z_g$};
					\node at (0,2.25){$\alpha^{G-}_A(\mu)$};
					\node at (1,-1.75){$\alpha^{G+}_A(\lambda)$};
					\node at (-0.25,0.25){$\lambda$};
					\node at (1.35,0.75){$A_A$};
					\node at (-1.35,0.75){${}_A A$};
				\end{tikzpicture}
				\caption{$\tilde{I}^A_{\lambda,\mu}(f)$}
				\label{graphical_localhom_dual}
			\end{minipage}
			\begin{minipage}[b]{0.45\hsize}
				\centering
				\begin{tikzpicture}
					\draw (0,0) -- (0,1) -- (0.5,1.5);
					\draw[ultra thick] (0.25,1.25) -- (0,1) -- (0,0);
					\draw[cross] (0,1.5) -- (0.5,1) -- (0.5,0) -- (0,-0.5) -- (0,-1) arc (360:180:0.25) -- (-0.5,1.5);
					\draw[cross,->] (0,0) -- (0.5,-0.5) -- (0.5,-1.5);
					\draw[ultra thick] (0,0) -- (0.25,-0.25);
					\node[block] at (0,0.5){$\tilde{f}'$}; 
					\node[block] at (0,-0.75){$z_g$};
					\node at (0,1.75){${}_A A$};
					\node at (0.5,-1.75){$\lambda$};
					\node at (-0.75,1.75){$A_A$};
					\node at (0.5,1.75){$\mu$};
				\end{tikzpicture}
				\caption{$\tilde{I}^{A-1}_{\lambda,\mu}(\tilde{f}')$}
				\label{graphical_localhom_dual_inv}
			\end{minipage}
		\end{tabular}
	\end{figure}

	\begin{proof}
		We only show the statement for $\hom_{\mathrm{loc}}^G(\lambda,A \mu)$ because the proof for $\hom_{\mathrm{loc}}^G(A\lambda,\mu)$ is easier. We show that $I^{A-1}_{\lambda, \mu}$ is given by the morphism in Figure \ref{graphicallocalhominv} for $f' \in \hom(\alpha^{G+}_A(\lambda),\alpha^{G-}_A(\mu))$ up to the canonical isomorphism $A_A \otimes_A {}_A A \cong A$. By the graphical calculation in Figure \ref{graphicalequivstrtensor}, we can do a similar graphical calculation to that in the proof of \cite[Lemma 4.16]{MR3424476} using Lemma \ref{lem_alpha_ordinary_braiding} and Corollary \ref{cor_thick_rotation} and see that $I_{\lambda,\mu}^{A-1} \circ I_{\lambda,\mu}^A (f)$ is equal to the morphism in Figure \ref{graphicallocalhomleftinvproof} by the graphical calculation in Figure \ref{graphicalequivstrtensor}. Then, from the graphical calculation in Figure \ref{graphicallocalhomleftinvproof}, where the third equality follows from a similar argument to the proof of \cite[Lemma 4.26]{bklr}, we obtain $I_{\lambda,\mu}^{A-1} \circ I_{\lambda,\mu}^A (f) = \tilde{P}^{G}_A(\lambda) f = f$. We can also check that $I_{\lambda,\mu}^{A-1}$ is indeed the right inverse as in the proof of \cite[Lemma 4.16]{MR3424476} using Lemma \ref{lem_arc_alpha_braiding} and the graphical calculation in Figure \ref{graphicallocalhomrightinvproof}. 
		\begin{figure}[htb]
			\centering
			\begin{tikzpicture}
				\draw[<-] (0,0) -- (0,1);
				\draw[<-] (0.75,0) -- (0.75,1);
				\node[block] at (0,0.5){$z_g$};
				\node[block] at (0.75,0.5){$z_g$};
				\node at (0,-0.25){$A_A$};
				\node at (0,1.25){${}^g A_A$};
				\node at (0.75,-0.25){${}_A A$};
				\node at (0.75,1.25){${}_A^g A$};
				\node at (1.5,0.5){$=$};
				\begin{scope}[shift={(2.25,0)}]
					\draw[->] (0.5,1.5) -- (0.5,-0.5);
					\draw[fill=white] (0.5,0.5)circle [radius=0.5];
					\node[block] at (0,0.5){$z_g$}; 
					\node[block] at (1,0.5){$z_g$};
					\node at (0.5,-0.75){$A$};
					\node at (0.5,1.75){${}^g A$};
					\node at (1.75,0.5){$=$};
				\end{scope}
				\begin{scope}[shift={(4.25,0)}]
					\draw[->] (0.5,1.5) -- (0.5,-0.5);
					\draw[fill=white] (0.5,0)circle [radius=0.25];
					\node[block] at (0.5,1){$z_g$}; 
					\node at (0.5,-0.75){$A$};
					\node at (0.5,1.75){${}^g A$};
					\node at (1.25,0.5){$=$};
					\begin{scope}[shift={(2,0)}]
						\draw[<-] (0,0) -- (0,1);
						\node[block] at (0,0.5){$z_g$};
						\node at (0,-0.25){$A$};
						\node at (0,1.25){${}^g A$};	
					\end{scope}
				\end{scope}
			\end{tikzpicture}
			\caption{$z_g \otimes_A z_g = z_g$}
			\label{graphicalequivstrtensor}
		\end{figure}
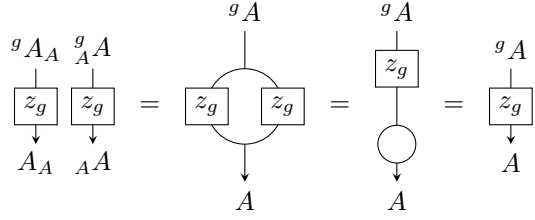
		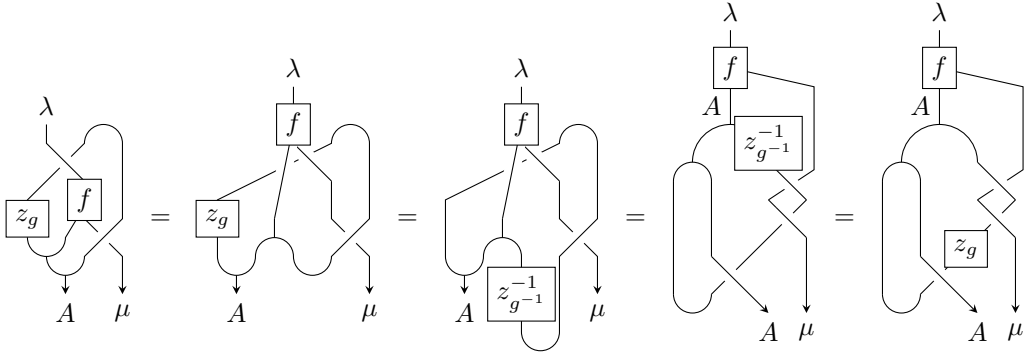
\begin{figure}[htb]
			\centering
			\begin{tikzpicture}
				\draw[->] (0.25,1) -- (0.75,0.5) -- (0.75,0);
				\draw[cross] (-0.5,0.75) --(-0.5,1.25)--(0.25,2) arc (180:0:0.25) --(0.75,1) -- (0.25,0.5) arc (360:180:0.25);
				\draw[->] (0,0.25) -- (0,0);
				\draw (-0.5,0.75) arc (180:360:0.25) -- (0.15,1);
				\draw[cross] (-0.25,2.25) --(-0.25,2) -- (0.25,1.5); 
				\node[block] at (0.25,1.25){$f$};
				\node[block] at (-0.5,1){$z_g$};
				\node at (0,-0.25){$A$};
				\node at (0.75,-0.25){$\mu$};
				\node at (-0.25,2.5){$\lambda$};
				\node at (1.25,1){$=$};
				\begin{scope}[shift={(3.25,0)}]
					\draw[->] (0.25,1.5) --(0.25,1) -- (0.75,0.5) -- (0.75,0);
					\draw[cross] (-1.25,0.5) --(-1.25,1.25)--(0.25,2) arc (180:0:0.25) --(0.75,1) -- (0.25,0.5) arc (360:180:0.25) arc (0:180:0.25) arc (360:180:0.25);
					\draw[->] (-1,0.25) -- (-1,0);
					\draw[cross] (-0.25,2.75) --(-0.25,2) -- (0.25,1.5);
					\draw[cross] (-0.25,2.25) --(-0.5,1);
					\draw (-0.5,1) -- (-0.5,0.75);
					\node[block] at (-0.25,2.25){$f$};
					\node[block] at (-1.25,1){$z_g$};
					\node at (-1,-0.25){$A$};
					\node at (0.75,-0.25){$\mu$};
					\node at (-0.25,3){$\lambda$};
					\node at (1.25,1){$=$};
				\end{scope}
				\begin{scope}[shift={(6.25,0)}]
					\draw[->] (0.25,1.5) --(0.25,1) -- (0.75,0.5) -- (0.75,0);
					\draw[cross] (-1.25,0.5) --(-1.25,1.25)--(0.25,2) arc (180:0:0.25) --(0.75,1) -- (0.25,0.5) -- (0.25,-0.5) arc (360:180:0.25) -- (-0.25,0.5) arc (0:180:0.25) arc (360:180:0.25);
					\draw[->] (-1,0.25) -- (-1,0);
					\draw[cross] (-0.25,2.75) --(-0.25,2) -- (0.25,1.5);
					\draw[cross] (-0.25,2.25) --(-0.5,1);
					\draw (-0.5,1) -- (-0.5,0.75);
					\node[block] at (-0.25,2.25){$f$};
					\node[block] at (-0.25,0){$z_{g^{-1}}^{-1}$};
					\node at (-1,-0.25){$A$};
					\node at (0.75,-0.25){$\mu$};
					\node at (-0.25,3){$\lambda$};
					\node at (1.25,1){$=$};
				\end{scope}
				\begin{scope}[shift={(8.75,2.25)}]
					\draw (0,0) -- (0,1.25);
					\draw (0.5,-0.5) arc (0:180:0.5);
					\draw (-0.25,-0.75) arc (0:180:0.25);
					\draw (0,0.75)--(1.1,0.5) -- (1.1,-0.6) -- (0.5,-1);
					\draw[cross] (0.5, -0.5) -- (1,-1);
					\draw (1,-1) -- (-0.25,-2.25);
					\draw (-0.25,-2.25) arc (360:180:0.25);
					\draw (-0.75,-0.75) -- (-0.75,-2.25);
					\draw[cross] (0.6, -1.1) -- (1,-1.5);
					\draw (0.5,-1) -- (0.6,-1.1);
					\draw[->] (1,-1.5) -- (1,-2.5);
					\draw[cross,->] (-0.25, -1.75) -- (0.5,-2.5);
					\draw (-0.25,-0.75) --(-0.25,-1.75);
					\node at (-0.25,0.25){$A$};
					\node at (0,1.5){$\lambda$};
					\node at (0.5,-2.75){$A$};
					\node at (1,-2.75){$\mu$};
					\node[block] at (0.5,-0.25){$z_{g^{-1}}^{-1}$};
					\node[block] at (0,0.75){$f$};
				\end{scope}
				\begin{scope}[shift={(11.5,2.25)}]
					\draw (0,0) -- (0,1.25);
					\draw (0.5,-0.5) arc (0:180:0.5);
					\draw (-0.25,-0.75) arc (0:180:0.25);
					\draw (0,0.75)--(1.1,0.5) -- (1.1,-0.6) -- (0.5,-1);
					\draw[cross] (0.5, -0.5) -- (1,-1);
					\draw (1,-1) -- (-0.25,-2.25);
					\draw (-0.25,-2.25) arc (360:180:0.25);
					\draw (-0.75,-0.75) -- (-0.75,-2.25);
					\draw[cross] (0.6, -1.1) -- (1,-1.5);
					\draw (0.5,-1) -- (0.6,-1.1);
					\draw[->] (1,-1.5) -- (1,-2.5);
					\draw[cross,->] (-0.25, -1.75) -- (0.5,-2.5);
					\draw (-0.25,-0.75) --(-0.25,-1.75);
					\node at (-0.25,0.25){$A$};
					\node at (0,1.5){$\lambda$};
					\node at (0.5,-2.75){$A$};
					\node at (1,-2.75){$\mu$};
					\node[block] at (0,0.75){$f$};
					\node[block] at (0.35,-1.65){$z_g$};
					\node at (-1.25,-1.25){$=$};
				\end{scope}
			\end{tikzpicture}
			\caption{$I_{\lambda,\mu}^{A-1} \circ I_{\lambda,\mu}^A (f)$}
			\label{graphicallocalhomleftinvproof}
		\end{figure}
		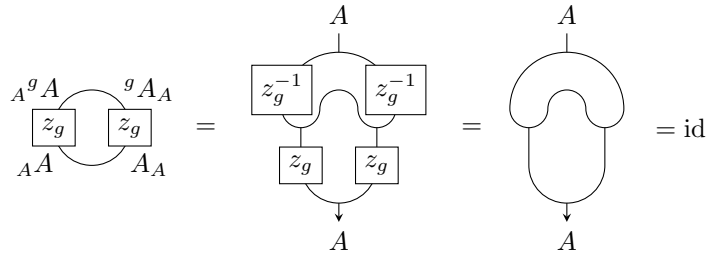
\begin{figure}[htb]
			\centering
			\begin{tikzpicture}
				\draw (0,0) circle [radius=0.5];
				\node[block] at (-0.5,0){$z_g$};
				\node[block] at (0.5,0){$z_g$};
				\node at (-0.75,0.5){${}_A {}^g A$};
				\node at (0.75,0.5){${}^g A_A$};
				\node at (-0.75,-0.5){${}_A A$};
				\node at (0.75,-0.5){$A_A$};
				\node at (1.5,0){$=$};
				\begin{scope}[shift={(2.5,0.25)}]
					\draw (0,0) arc (180:0:0.75) arc (360:180:0.25) arc (0:180:0.25) arc (360:180:0.25); 
					\draw (0.25,-0.25) -- (0.25,-0.75) arc (180:360:0.5) -- (1.25,-0.25);
					\draw (0.75,0.75) -- (0.75,1);
					\draw[->] (0.75,-1.25) -- (0.75,-1.5);
					\node[block] at (0.25,-0.75){$z_g$};
					\node[block] at (1.25,-0.75){$z_g$};
					\node[block] at (0,0.25){$z_g^{-1}$};
					\node[block] at (1.5,0.25){$z_g^{-1}$};
					\node at (0.75,1.25){$A$};
					\node at (0.75,-1.75){$A$};
					\node at (2.5,-0.25){$=$};
				\end{scope}
				\begin{scope}[shift={(5.5,0.25)}]
					\draw (0,0) arc (180:0:0.75) arc (360:180:0.25) arc (0:180:0.25) arc (360:180:0.25); 
					\draw (0.25,-0.25) -- (0.25,-0.75) arc (180:360:0.5) -- (1.25,-0.25);
					\draw (0.75,0.75) -- (0.75,1);
					\draw[->] (0.75,-1.25) -- (0.75,-1.5);
					\node at (0.75,1.25){$A$};
					\node at (0.75,-1.75){$A$};
					\node at (2.25,-0.25){$= \id$};
				\end{scope}
			\end{tikzpicture}
			\caption{$z_g \otimes z_g = \id_A$}
			\label{graphicallocalhomrightinvproof}
		\end{figure}
	\end{proof}
\end{lem}

\begin{rem}
	When $\calc$ is a $G$-braided ${}^\ast$-multitensor category and $A$ is a $G$-equivariant Q-system, Figure \ref{graphicallocalhomleftinvproof} with $\lambda = A \mu$ and $f = \id_{A \mu}$ shows that $\tilde{P}^{G+}_A(\mu)$ is a projection as in \cite[Lemma 4.26]{bklr}.
\end{rem}

\begin{lem}
	\label{lem_i_action}
	Let $A$ be a neutral symmetric special $G$-equivariant Frobenius algebra in a $G$-ribbon multitensor category $\calc$ and let $\lambda, \mu \in \homog(\calc)$. Then, ${}^g I_{\lambda, \mu}^{A-1}(f') = \tilde{\eta}^g_\mu I_{{}^g \lambda, {}^g \mu}^{A-1}(\tilde{\eta}^{g-1}_\mu {}^g f' \tilde{\eta}^g_\lambda)$ for $f' \in \hom(\alpha_A^{G+}(\lambda), \alpha^{G-}_A(\mu))$ and $g \in G$.

	\begin{proof}
		Put $z \coloneqq z^A$ and $\partial \lambda = h$. By Lemma \ref{lem_ibfe_action}, $\tilde{\eta}^{g-1}_\mu {}^g I_{\lambda, \mu}^{A-1}(f')$ is given by the morphism in Figure \ref{graphical_lem_i_action_proof}. Then, the graphical calculation there shows the statement. We used that $\mathrm{coev}_{A_A} = \eta_A$ by \cite[Section 5]{MR2075605} at the first equality and $z_g {}^g z_h {}^{ghg^{-1}} z_g^{-1} = z_g {}^g z_h {}^{gh} z_{g^{-1}} = z_{ghg^{-1}}$ at the second equality. 
		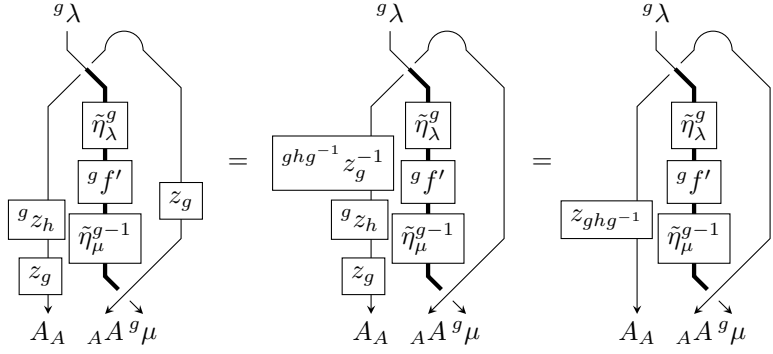
\begin{figure}[htb]
			\centering
			\begin{tikzpicture}
				\draw[<-] (0.5,-1.5) -- (0,-1) -- (0,1.5);
				\draw[ultra thick] (0,1) -- (0,-1) -- (0.25,-1.25);
				\draw[cross,<->] (0,-1.5) -- (1,-0.5) -- (1,1.5) -- (0.5,2) arc (0:180:0.25) -- (-0.75,1.25) -- (-0.75,-1.5);
				\draw[cross] (-0.5,2.25) -- (-0.5,2) -- (0,1.5) -- (0,1);
				\draw[ultra thick] (-0.25,1.75) -- (0,1.5) -- (0,1);
				\node[block] at (0,1){$\tilde{\eta}_{\lambda}^g$};
				\node[block] at (0,0.25){${}^g f'$}; 
				\node[block] at (0,-0.5){$\tilde{\eta}_{\mu}^{g-1}$};
				\node[block] at (-0.9,-0.25){${}^g z_h$};
				\node[block] at (-0.85,-1){$z_g$};
				\node[block] at (1,0){$z_g$};
				\node at (0,-1.75){${}_A A$};
				\node at (0.5,-1.75){${}^g \mu$};
				\node at (-0.75,-1.75){$A_A$};
				\node at (-0.5,2.5){${}^g \lambda$};
				\node at (1.75,0.5){$=$};
				\begin{scope}[shift={(4.25,0)}]
					\draw[<-] (0.5,-1.5) -- (0,-1) -- (0,1.5);
					\draw[ultra thick] (0,1) -- (0,-1) -- (0.25,-1.25);
					\draw[cross,<->] (0,-1.5) -- (1,-0.5) -- (1,1.5) -- (0.5,2) arc (0:180:0.25) -- (-0.75,1.25) -- (-0.75,-1.5);
					\draw[cross] (-0.5,2.25) -- (-0.5,2) -- (0,1.5) -- (0,1);
					\draw[ultra thick] (-0.25,1.75) -- (0,1.5) -- (0,1);
					\node[block] at (0,1){$\tilde{\eta}_{\lambda}^g$};
					\node[block] at (0,0.25){${}^g f'$}; 
					\node[block] at (0,-0.5){$\tilde{\eta}_{\mu}^{g-1}$};
					\node[block] at (-0.9,-0.25){${}^g z_h$};
					\node[block] at (-0.85,-1){$z_g$};
					\node[block] at (-1.25,0.5){${}^{ghg^{-1}} z_g^{-1}$};
					\node at (0,-1.75){${}_A A$};
					\node at (0.5,-1.75){${}^g \mu$};
					\node at (-0.75,-1.75){$A_A$};
					\node at (-0.5,2.5){${}^g \lambda$};
					\node at (1.5,0.5){$=$};
				\end{scope}
				\begin{scope}[shift={(7.75,0)}]
					\draw[<-] (0.5,-1.5) -- (0,-1) -- (0,1.5);
					\draw[ultra thick] (0,1) -- (0,-1) -- (0.25,-1.25);
					\draw[cross,<->] (0,-1.5) -- (1,-0.5) -- (1,1.5) -- (0.5,2) arc (0:180:0.25) -- (-0.75,1.25) -- (-0.75,-1.5);
					\draw[cross] (-0.5,2.25) -- (-0.5,2) -- (0,1.5) -- (0,1);
					\draw[ultra thick] (-0.25,1.75) -- (0,1.5) -- (0,1);
					\node[block] at (0,1){$\tilde{\eta}_{\lambda}^g$};
					\node[block] at (0,0.25){${}^g f'$}; 
					\node[block] at (0,-0.5){$\tilde{\eta}_{\mu}^{g-1}$};
					\node[block] at (-1.15,-0.25){$z_{ghg^{-1}}$};
					\node at (0,-1.75){${}_A A$};
					\node at (0.5,-1.75){${}^g \mu$};
					\node at (-0.75,-1.75){$A_A$};
					\node at (-0.5,2.5){${}^g \lambda$};
				\end{scope}
			\end{tikzpicture}
			\caption{The proof of Lemma \ref{lem_i_action}}
			\label{graphical_lem_i_action_proof}
		\end{figure}
	\end{proof}
\end{lem}

\begin{lem}
	\label{localhomnondeglem}
	Let $A$ be a neutral symmetric special $G$-equivariant Frobenius algebra in a $G$-ribbon multitensor category $\calc$ and let $\lambda, \mu \in \homog(\calc)$. Then, the nondegenerate pairing $\hom(A \mu, \lambda) \times \hom(\lambda, A \mu)$ restricts to $\hom_{\mathrm{loc}}^G(A \mu, \lambda) \times \hom_{\mathrm{loc}}^G(\lambda, A\mu)$. Moreover, the maps $I^A_{\lambda, \mu}$ and $\tilde{I}^A_{\lambda, \mu}$ defined in Lemma \ref{localhomlem} are isometric with respect to this pairing. If $\calc$ is a $G$-braided ${}^\ast$-multitensor category and $A$ is a $G$-equivariant Q-system, then they are unitary.

	\begin{proof}
		The first statement follows from Lemma \ref{lem_tildep_idempotent}. The second statement follows from a similar graphical calculation to that in the proof of \cite[Lemma 4.16]{MR3424476} by Lemma \ref{lem_alpha_ordinary_braiding}. Note that $\alpha_A^{G+}$ is pivotal and therefore preserves traces. The final statement follows from the unitarity of thick crossings.
	\end{proof}
\end{lem}

Finally, we can prove our second main theorem in this article, which is the equivariant generalization of \cite[Proposition 4.18]{MR3424476}.

\begin{thm}
	\label{mainthm2}
Let $A$ be a neutral symmetric special simple $G$-equivariant Frobenius algebra in a split spherical $G$-braided fusion category $\calc$ with $\dim \calc \neq 0$. Then, the $G$-equivariant full center $Z^G(A)_{\dim A}$ of $A$ (Definition \ref{def_fullcenter}), where the subscript $\dim A$ means that we take $\zeta$ in Proposition \ref{equivsubalgprop} to be $\dim A$, is isomorphic to the $G$-equivariant $\alpha$-induction Frobenius algebra $\Theta_{\alpha}^G(A)$ associated with $A$ (Theorem \ref{mainthm1}). If $\calc$ is a $G$-braided ${}^\ast$-multitensor category and $A$ is a $G$-equivariant Q-system, then we have an isomorphism between Q-systems.

\begin{proof}
By Lemmata \ref{prodidempotentlem} and \ref{lem_commutative_center}, we have $P^{G}_{(A \boxtimes \mathbf{1})\Theta_{\mathrm{LR}}^G} = P^{G}_{A \boxtimes \mathbf{1}}(C^{G}(\Theta_{\mathrm{LR}}^G)) = P^{G}_{A \boxtimes \mathbf{1}}(\Theta_{\mathrm{LR}}^G)$ since $\Theta^G_{\mathrm{LR}}$ is $G$-commutative by Proposition \ref{prop_alphafrob_gcomm}. Hence, by \cite[Lemma 2.4(ii)]{MR2187404} and semisimplicity, we have
\begin{align*}
Z^G(A) &= \bigoplus_{\lambda_1,\lambda_2 \in \Delta} \langle \lambda_1, A\lambda_2 \rangle^G_{\mathrm{loc}} \lambda_1 \boxtimes \lambda_2^\vee
\end{align*}
as objects for a complete system $\Delta$ of representatives of the simple objects of $\calc$. Hence, by Lemma \ref{localhomlem}, $Z^G(A) \cong \Theta_\alpha^G(A)$ as objects. 

By Lemma \ref{localhomnondeglem}, we can see that for a basis $\{ \varphi_l^L \}_{l}$ of $\hom_{\mathrm{loc}}^G(\lambda_1, A \lambda_2)$ and its dual basis $\{ \tilde{\varphi}_l^{L}\}_l$ of $\hom_{\mathrm{loc}}^G(A \lambda_2,\lambda_1)$, the morphisms $s_{Z^G(A)} \coloneqq \bigoplus_{L} \varphi_l^L \boxtimes \id_{\lambda_2}$ and $ r_{Z^G(A)} \coloneqq \bigoplus_{L} \tilde{\varphi}_l^{L} \boxtimes \id_{\lambda_2}$, where $L \coloneqq (\lambda_1, \lambda_2, l)$, split the idempotent $P_{A \boxtimes \mathbf{1}}^{G+}(\Theta_{\mathrm{LR}}^G)$ as in \cite[Lemma 4.15]{MR3424476}. Then, by a similar argument to that in the proof of \cite[Proposition 4.18]{MR3424476}, we can see that $Z^{G}(A)_{\dim A} \cong \Theta_{\alpha}^G(A)$ as Frobenius algebras. Namely, their units and counits coincide up to an isomorphism thanks to the normalization $\zeta = \dim A$. In order to show the coincidence of coproducts, it suffices to calculate the quantity in Figure \ref{graphical_mainthm2_proof_coproduct}, which is the coefficient of $e_i^{\nu_1, \lambda_1 \mu_1} \boxtimes \overline{\tilde{e}_j^{\nu_2,\lambda_2 \mu_2}}$ up to $d_A d_{\lambda_2} d_{\mu_2}/d_{\nu_2} d_{\nu_1}$, as in the proof of \cite[Proposition 4.18]{MR3424476}. Then, by the graphical calculation there, the coproducts coincide. We put $\phi^{L}_l \coloneqq I^A_{\lambda_1, \lambda_2}(\varphi^L_l)$ and used Lemma \ref{localhomnondeglem} at the first equality. Note that the second equality follows since the quantity is nonzero only if $\partial \nu_1 = \partial \lambda_1 \partial \mu_1$. The third equality follows as in the proof of \cite[Proposition 4.18]{MR3424476}. The proof for products is similar.
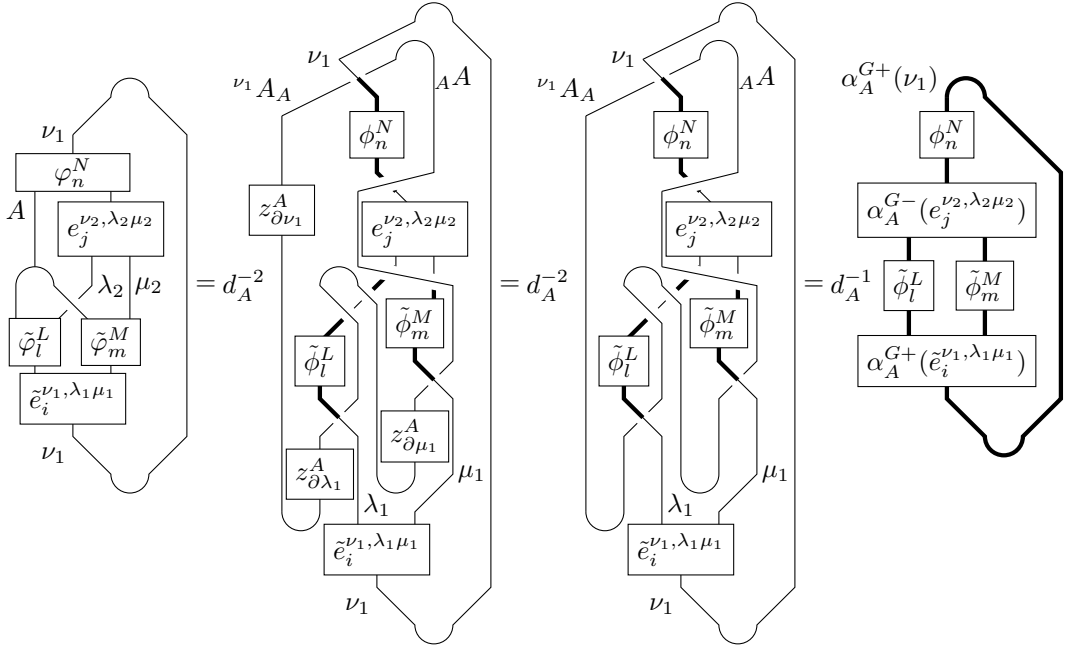
\begin{figure}[H]
	\centering
	\begin{tikzpicture}
		\draw (0,1) -- (0,0.5);
		\draw (0.5,-0.25) -- (0.5,0.5);
		\draw (0.25,-0.5) -- (0.25,-1) -- (-0.25,-1.5);
		\draw (0.75,-0.5) -- (0.75,-1.5);
		\draw[cross] (-0.75,-1.5) -- (-0.75,-1) arc (180:0:0.25) -- (0.25,-1.5);
		\draw (-0.5,0.25) -- (-0.5,-0.75);
		\draw (0.5,-1.75) -- (0.5,-2.5);
		\draw (-0.5,-1.75) -- (-0.5,-2.5);
		\draw (0,-2.5) -- (0,-3) -- (0.5,-3.5) arc (180:360:0.25) -- (1.5,-3) -- (1.5,1) -- (1,1.5) arc (0:180:0.25) -- (0,1);
		\draw[fill=white] (-0.75,0.75) rectangle (0.75,0.25) node[midway]{$\varphi^N_n$};
		\node[block] at (0.5,-0.25){$e_j^{\nu_2,\lambda_2 \mu_2}$};
		\node[block] at (-0.5,-1.75){$\tilde{\varphi}^L_l$};
		\node[block] at (0.5,-1.75){$\tilde{\varphi}^M_m$};
		\node[block] at (0,-2.5){$\tilde{e}_i^{\nu_1,\lambda_1\mu_1}$};
		\node at (1.75,-1){$=$};
		\node at (-0.25,1){$\nu_1$};
		\node at (-0.25,-3.25){$\nu_1$};
		\node at (-0.75,0){$A$};
		\node at (1,-1){$\mu_2$};
		\node at (0.5,-1){$\lambda_2$};
		\begin{scope}[shift={(4,0)}]
			\draw (0.5,-0.25) -- (0.5,0) -- (0,0.5)-- (0,1);
			\draw[ultra thick] (0,1) -- (0,0.5) -- (0.25,0.25);
			\draw[cross] (-0.25,0.25) -- (0.75,0.5) -- (0.75,2) arc (0:180:0.25) -- (-0.75,1.5) -- (-1.25,1.25); 
			\draw[cross] (0,1) -- (0,1.5) -- (-0.5,2);
			\draw[ultra thick] (-0.25,1.75) -- (0,1.5) -- (0,1);
			\draw (0.25,-0.5) -- (0.25,-0.75) -- (-0.75,-1.75) -- (-0.75,-2.25);
			\draw[ultra thick] (0.25,-0.75) -- (-0.1,-1.1);
			\draw[ultra thick] (-0.35,-1.35) -- (-0.75,-1.75);
			\draw (0.75,-0.5) -- (0.75,-1.5);
			\draw[ultra thick] (0.75,-1) -- (0.75,-1.5);
			\draw[cross] (-0.25,-0.75) -- (1,-1) -- (1,-2) -- (0.5,-2.5) -- (0.5,-3.5) arc (360:180:0.25) -- (0,-1.25) -- (-0.25,-1) arc (0:180:0.25) -- (-0.25,-1.5) -- (-0.25,-2.5) -- (-0.75,-3) -- (-0.75,-4) arc (360:180:0.25);
			\draw[cross] (-0.75,-2.25) -- (-0.75,-2.5) -- (-0.25,-3) -- (-0.25,-4.5);
			\draw[ultra thick] (-0.75,-2.25) -- (-0.75,-2.5) -- (-0.5,-2.75);
			\draw (-1.25,-4) -- (-1.25,1.25);
			\draw (-0.25,0.25) -- (-0.25,-0.75);
			\draw[cross] (0.5,-2) -- (1,-2.5) -- (1,-3.5) -- (0.5,-4) -- (0.5,-4.5);
			\draw[ultra thick] (0.5,-1.75) -- (0.5,-2) -- (0.75,-2.25);
			\draw (0.5,-1.75) -- (0.5,-2);
			\draw (0,-4.5) -- (0,-5) -- (0.5,-5.5) arc (180:360:0.25) -- (1.5,-5) -- (1.5,2) -- (1,2.5) arc (0:180:0.25) -- (-0.5,2);
			\node[block] at (0,1){$\phi^N_n$};
			\node[block] at (0.5,-0.25){$e_j^{\nu_2,\lambda_2 \mu_2}$};
			\node[block] at (-0.75,-2){$\tilde{\phi}^L_l$};
			\node[block] at (0.5,-1.5){$\tilde{\phi}^M_m$};
			\node[block] at (0,-4.5){$\tilde{e}_i^{\nu_1,\lambda_1\mu_1}$};
			\node[block] at (0.5,-3){$z^A_{\partial \mu_1}$};
			\node[block] at (-0.75,-3.5){$z^A_{\partial \lambda_1}$};
			\node[block] at (-1.25,0){$z^A_{\partial \nu_1}$};
			\node at (-0.75,2){$\nu_1$};
			\node at (-0.25,-5.25){$\nu_1$};
			\node at (0,-3.9){$\lambda_1$};
			\node at (1.25,-3.5){$\mu_1$};
			\node at (-1.5,1.6){${}^{\nu_1}A_A$};
			\node at (1,1.75){${}_A A$};
			\node at (1.75,-1){$=$};
			\node at (-1.75,-1){$d_A^{-2}$};
		\end{scope}
		\begin{scope}[shift={(8,0)}]
			\draw (0.5,-0.25) -- (0.5,0) -- (0,0.5)-- (0,1);
			\draw[ultra thick] (0,1) -- (0,0.5) -- (0.25,0.25);
			\draw[cross] (-0.25,0.25) -- (0.75,0.5) -- (0.75,2) arc (0:180:0.25) -- (-0.75,1.5) -- (-1.25,1.25); 
			\draw[cross] (0,1) -- (0,1.5) -- (-0.5,2);
			\draw[ultra thick] (-0.25,1.75) -- (0,1.5) -- (0,1);
			\draw (0.25,-0.5) -- (0.25,-0.75) -- (-0.75,-1.75) -- (-0.75,-2.25);
			\draw[ultra thick] (0.25,-0.75) -- (-0.1,-1.1);
			\draw[ultra thick] (-0.35,-1.35) -- (-0.75,-1.75);
			\draw (0.75,-0.5) -- (0.75,-1.5);
			\draw[ultra thick] (0.75,-1) -- (0.75,-1.5);
			\draw[cross] (-0.25,-0.75) -- (1,-1) -- (1,-2) -- (0.5,-2.5) -- (0.5,-3.5) arc (360:180:0.25) -- (0,-1.25) -- (-0.25,-1) arc (0:180:0.25) -- (-0.25,-1.5) -- (-0.25,-2.5) -- (-0.75,-3) -- (-0.75,-4) arc (360:180:0.25);
			\draw[cross] (-0.75,-2.25) -- (-0.75,-2.5) -- (-0.25,-3) -- (-0.25,-4.5);
			\draw[ultra thick] (-0.75,-2.25) -- (-0.75,-2.5) -- (-0.5,-2.75);
			\draw (-1.25,-4) -- (-1.25,1.25);
			\draw (-0.25,0.25) -- (-0.25,-0.75);
			\draw[cross] (0.5,-2) -- (1,-2.5) -- (1,-3.5) -- (0.5,-4) -- (0.5,-4.5);
			\draw[ultra thick] (0.5,-1.75) -- (0.5,-2) -- (0.75,-2.25);
			\draw (0.5,-1.75) -- (0.5,-2);
			\draw (0,-4.5) -- (0,-5) -- (0.5,-5.5) arc (180:360:0.25) -- (1.5,-5) -- (1.5,2) -- (1,2.5) arc (0:180:0.25) -- (-0.5,2);
			\node[block] at (0,1){$\phi^N_n$};
			\node[block] at (0.5,-0.25){$e_j^{\nu_2,\lambda_2 \mu_2}$};
			\node[block] at (-0.75,-2){$\tilde{\phi}^L_l$};
			\node[block] at (0.5,-1.5){$\tilde{\phi}^M_m$};
			\node[block] at (0,-4.5){$\tilde{e}_i^{\nu_1,\lambda_1\mu_1}$};
			\node at (-0.75,2){$\nu_1$};
			\node at (-0.25,-5.25){$\nu_1$};
			\node at (0,-3.9){$\lambda_1$};
			\node at (1.25,-3.5){$\mu_1$};
			\node at (-1.5,1.6){${}^{\nu_1}A_A$};
			\node at (1,1.75){${}_A A$};
			\node at (1.75,-1){$=$};
			\node at (-1.75,-1){$d_A^{-2}$};
		\end{scope}
		\begin{scope}[shift={(11.5,0)}]
			\draw[ultra thick] (0,1.5) -- (0,0);
			\draw[ultra thick] (-0.5,0) -- (-0.5,-2);
			\draw[ultra thick] (0.5,0) -- (0.5,-2);
			\draw[ultra thick] (0,-2) -- (0,-2.5)  -- (0.5,-3) arc (180:360:0.25) -- (1.5,-2.5) -- (1.5,0.5) -- (0.5,1.5) arc (0:180:0.25);
			\node[block] at (0,1){$\phi^N_n$};
			\node[block] at (0,0){$\alpha^{G-}_A(e_j^{\nu_2,\lambda_2 \mu_2})$};
			\node[block] at (-0.5,-1){$\tilde{\phi}^L_l$};
			\node[block] at (0.5,-1){$\tilde{\phi}^M_m$};
			\node at (-1.25,-1){$d_A^{-1}$};
			\node at (-0.75,1.75){$\alpha_A^{G+}(\nu_1)$};
			\node[block] at (0,-2){$\alpha_A^{G+} (\tilde{e}_i^{\nu_1,\lambda_1\mu_1})$};
		\end{scope}
	\end{tikzpicture}
	\caption{The coincidence of coproducts}
	\label{graphical_mainthm2_proof_coproduct}
\end{figure}

We show that they are indeed isomorphic as $G$-equivariant Frobenius algebras. Since
\begin{align*}
	z^{Z^G(A)}_g = \bigoplus_{L,m} \tilde{\varphi}^{L(g)}_m (z_g \otimes u_g^{\lambda_2}) {}^g \varphi^L_l \boxtimes \overline{\tilde{u}_g^{\lambda_2}} = \bigoplus_{L,m} \frac{1}{d_{\lambda_1}} \tr (\tilde{u}_g^{\lambda_1} \tilde{\varphi}^{L(g)}_m (z_g \otimes u_g^{\lambda_2}) {}^g \varphi^L_l) u_g^{\lambda_1} \boxtimes \overline{\tilde{u}_g^{\lambda_2}}
\end{align*}
for $g \in G$. Then, the traces in coefficients are represented in Figure \ref{graphicalmainthmproof}, which gives the coincidence of the equivariant structures. We used Lemma \ref{lem_i_action} at the first equality and $\partial {}^g \lambda_1 = \partial \lambda_1(g)$ at the second equality. The final statement follows since then $\{ \phi^L_l \}_l$ is an orthonormal basis by Lemma \ref{localhomnondeglem}.
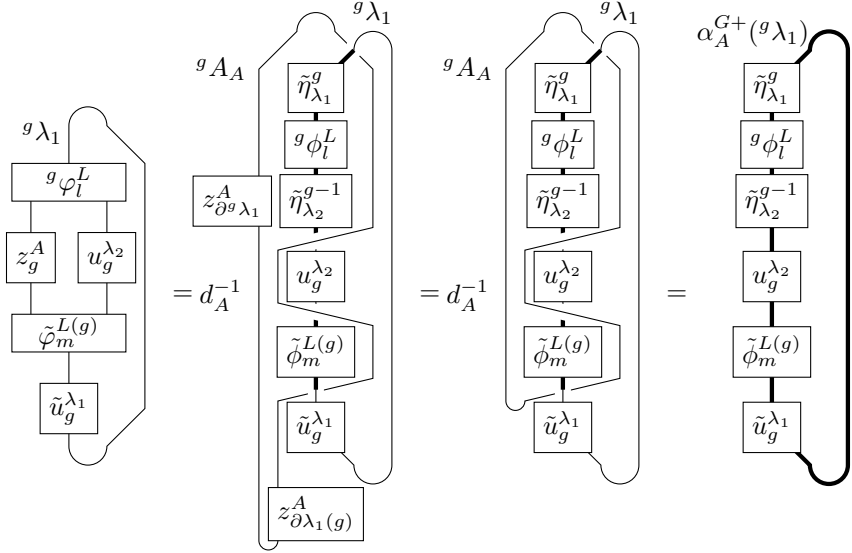
\begin{figure}[htb]
	\centering
	\begin{tikzpicture}
		\draw (-0.5,1) -- (-0.5,-1);
		\draw (0.5,1) -- (0.5,-1);
		\draw (0,-1) -- (0,-2.5) arc (180:360:0.25) -- (1,-2) -- (1,1.25) -- (0.5,1.75) arc (0:180:0.25) -- (0,1);
		\draw[fill=white] (-0.75,1.25) rectangle (0.75,0.75) node[midway]{${}^g \varphi^L_l$};
		\node [block] at (0.5,0){$u_g^{\lambda_2}$};
		\node [block] at (-0.5,0){$z_g^{A}$};
		\draw[fill=white] (-0.75,-0.75) rectangle (0.75,-1.25) node[midway]{$\tilde{\varphi}^{L(g)}_m$};
		\node [block] at (0,-2){$\tilde{u}_g^{\lambda_1}$};
		\node at (1.5,-0.5){$=$};
		\node at (-0.35,1.65){${}^g \lambda_1$};
		\begin{scope}[shift={(2.75,-0.25)}]
			\draw (0.5,1) -- (0.5,-1);
			\draw[ultra thick] (0.5,2.5) -- (0.5,0.5);
			\draw[ultra thick] (0.5,-1) -- (0.5,-0.5);
			\draw[cross] (0,0) -- (0,0.35) -- (1.25,0.65) -- (1.25,2.75) -- (0.75,3.25) arc (0:180:0.25) -- (-0.25,2.75) -- (-0.25,-3.5) arc (180:360:0.125) -- (0,-1.65) -- (1.25,-1.35) -- (1.25,-0.65) -- (0,-0.35) -- (0,0);
			\draw[cross] (0.5,-1) -- (0.5,-2) -- (1,-2.5) arc (180:360:0.25) -- (1.5,-1.75) -- (1.5,3) arc (0:180:0.25) -- (0.5,2.5);
			\draw[ultra thick] (0.5,2.5) -- (1,3);
			\draw[ultra thick] (0.5,-1) -- (0.5,-1.5);
			\node[block] at (0.5,1.75){${}^g \phi^L_l$};
			\node[block] at (0.5,1){$\tilde{\eta}_{\lambda_2}^{g-1}$};
			\node[block] at (0.5,2.5){$\tilde{\eta}_{\lambda_1}^g$};
			\node [block] at (0.5,0){$u_g^{\lambda_2}$};
			\node[block] at (0.5,-1){$\tilde{\phi}^{L(g)}_m$};
			\node [block] at (0.5,-2){$\tilde{u}_g^{\lambda_1}$};
			\node[block] at (0.5,-3.15){$z_{\partial \lambda_1(g)}^A$};
			\node[block] at (-0.6,1){$z_{\partial {}^g \lambda_1}^A$};
			\node at (2,-0.25){$=$};
			\node at (-0.75,-0.25){$d_A^{-1}$};
			\node at (-0.75,2.75){${}^g A_A$};
			\node at (1.25,3.5){${}^g \lambda_1$};
		\end{scope}
		\begin{scope}[shift={(6,-0.25)}]
			\draw (0.5,1) -- (0.5,-1);
			\draw[ultra thick] (0.5,2.5) -- (0.5,0.5);
			\draw[ultra thick] (0.5,-1) -- (0.5,-0.5);
			\draw[cross] (0,0) -- (0,0.35) -- (1.25,0.65) -- (1.25,2.75) -- (0.75,3.25) arc (0:180:0.25) -- (-0.25,2.75) -- (-0.25,-1.65) arc (180:360:0.125) -- (1.25,-1.35) -- (1.25,-0.65) -- (0,-0.35) -- (0,0);
			\draw[cross] (0.5,-1) -- (0.5,-2) -- (1,-2.5) arc (180:360:0.25) -- (1.5,-1.75) -- (1.5,3) arc (0:180:0.25) -- (0.5,2.5);
			\draw[ultra thick] (0.5,2.5) -- (1,3);
			\draw[ultra thick] (0.5,-1) -- (0.5,-1.5);
			\node[block] at (0.5,1.75){${}^g \phi^L_l$};
			\node[block] at (0.5,1){$\tilde{\eta}_{\lambda_2}^{g-1}$};
			\node[block] at (0.5,2.5){$\tilde{\eta}_{\lambda_1}^g$};
			\node [block] at (0.5,0){$u_g^{\lambda_2}$};
			\node[block] at (0.5,-1){$\tilde{\phi}^{L(g)}_m$};
			\node [block] at (0.5,-2){$\tilde{u}_g^{\lambda_1}$};
			\node at (2,-0.25){$=$};
			\node at (-0.75,-0.25){$d_A^{-1}$};
			\node at (-0.75,2.75){${}^g A_A$};
			\node at (1.25,3.5){${}^g \lambda_1$};
		\end{scope}
		\begin{scope}[shift={(8.75,-0.25)}]
			\draw[ultra thick] (0.5,2.5) -- (0.5,-1) -- (0.5,-2) -- (1,-2.5) arc (180:360:0.25) -- (1.5,-1.75) -- (1.5,3) arc (0:180:0.25) -- (0.5,2.5);
			\node[block] at (0.5,1.75){${}^g \phi^L_l$};
			\node[block] at (0.5,1){$\tilde{\eta}_{\lambda_2}^{g-1}$};
			\node[block] at (0.5,2.5){$\tilde{\eta}_{\lambda_1}^g$};
			\node [block] at (0.5,0){$u_g^{\lambda_2}$};
			\node[block] at (0.5,-1){$\tilde{\phi}^{L(g)}_m$};
			\node [block] at (0.5,-2){$\tilde{u}_g^{\lambda_1}$};
			\node at (0.25,3.25){$\alpha_A^{G+} ({}^g \lambda_1)$};
		\end{scope}
	\end{tikzpicture}
	\caption{The coincidence of equivariant structures}
	\label{graphicalmainthmproof}
\end{figure}
\end{proof}
\end{thm}

\section*{Acknowledgements}

The author thanks his supervisor Yasuyuki Kawahigashi for his constant support and advice. In particular, he thanks him for suggesting studying equivariant $\alpha$-induction. He thanks Roberto Longo for answering his question. He thanks Sebastiano Carpi for comments on Remark \ref{nonfaithrem}. He thanks Kan Kitamura for comments on an early draft of this article. He thanks the reviewers for their detailed comments and suggestions. This work was supported by RIKEN Junior Research Associate Program and JSPS KAKENHI Grant Number JP23KJ0540. He thanks their financial support.

\section*{Declarations}

\textbf{Funding:} This work was supported by RIKEN Junior Research Associate Program and JSPS KAKENHI Grant Number JP23KJ0540.

\noindent
\textbf{Conflicts of interest/Competing interests:} The author has no conflicts of interest to declare that are relevant to the content of this article.

\bibliography{oikawa_induction}
\bibliographystyle{alpha}

\end{document}